\documentclass[12pt]{article}
\usepackage[dvips]{graphicx}
\usepackage{subfigure}
\usepackage{tikz}
\usepackage{ pstricks, amssymb, xr, multind, multicol,textcomp}
\def\hfl#1{\smash{\mathop{\hbox to 10mm{\rightarrowfill}}\limits^{\textstyle
#1}}}
\newtheorem{proposition}[equation]{Proposition} 
 
\newtheorem{theorem}[equation]{Theorem} 
\newtheorem{exa}[equation]{Example} 
\newtheorem{ex}[equation]{Exercise} 
\newtheorem{s-ex}[equation]{Side-exercise} 
\newtheorem{exas}[equation]{Examples} 
\newtheorem{lemma}[equation]{Lemma} 
\newtheorem{sublemma}[equation]{Sublemma} 
\newtheorem{remar}[equation]{Remark} 
\newtheorem{remars}[equation]{Remarks} 
\newtheorem{nota}[equation]{Notation} 
\newtheorem{sremar}[equation]{Side-remark} 
\newtheorem{definitio}[equation]{Definition}

\newenvironment{remark}{\begin{remar} \rm }{\end{remar}}

\newenvironment{example}{\begin{exa} \rm }{\end{exa}} 
\newenvironment{definition}{\begin{definitio} \rm }{\end{definitio}}

\newcommand{\KK}{\mathbb{K}}

\newcommand{\CA}{{\cal A}}

\newcommand{\TCM}{\tilde{C}_2(M)}
\newcommand{\TCMD}{\widetilde{M^2}} 
\newcommand{\CB}{{\cal B}}

\newcommand{\CI}{{\cal I}} 
 
\newcommand{\CL}{{\cal L}}

\newcommand{\CT}{{\cal T}} 
 
\newcommand{\CQ}{{\cal Q}}
\newcommand{\CS}{{\cal S}}
\newcommand{\projconf}{Q}
\newcommand{\fvar}{x}
\newcommand{\svar}{y}
\newcommand{\tvar}{z}
\newcommand{\fvarM}{X}
\newcommand{\svarM}{Y}
\newcommand{\tvarM}{Z}
\newcommand{\qvarM}{W}
\newcommand{\cvarM}{V}
\newcommand{\Bor}{B}
\newcommand{\Borp}{B^{\prime}}
\newcommand{\degr }{d }
\newcommand{\ID}{I_{\Delta}}
\newcommand{\ZZ}{\mathbb{Z}} 
\newcommand{\RR}{\mathbb{R}} 
\newcommand{\QQ}{\mathbb{Q}} 
\newcommand{\CC}{\mathbb{C}} 
\newcommand{\NN}{\mathbb{N}}

\newcommand{\bp}{\noindent {\sc Proof: }} 
\newcommand{\eop}{\nopagebreak \hspace*{\fill}{$\diamond$} \medskip}

\begin{document} 
\title{On the cube of the equivariant linking pairing\\ for knots and $3$--manifolds of rank one}
\author{Christine Lescop \thanks{Institut Fourier, UJF Grenoble, CNRS}}
\maketitle
\begin{abstract} 
Let $M$ be a closed oriented $3$-manifold with first Betti number one.
Its equivariant linking pairing may be seen as a two-dimensional cohomology class in an appropriate infinite cyclic covering of the space of ordered pairs of distinct points of $M$. We show how to define
the equivariant cube $\CQ(\KK)$ of this Blanchfield pairing with respect to a framed knot $\KK$ that generates $H_1(M;\ZZ)/\mbox{Torsion}$.

This article is devoted to the study of the invariant $\CQ$. We prove many properties for this
invariant including two surgery formulae.

Via surgery, the invariant $\CQ$ is equivalent to an invariant $\hat{\CQ}$ of null-homologous knots in rational homology spheres,
that coincides with the two-loop part of the Kricker
rational lift of the Kontsevich integral, at least for knots with trivial Alexander polynomial in integral homology spheres.

The invariant $\CQ$
takes its values in a polynomial ring $R$ over $\QQ$.
We determine the rational sub-vector space of $R$ generated by the variations $(\CQ(\KK^{\prime})-\CQ(\KK))$ for two framed knots that generate $H_1(M;\ZZ)/\mbox{Torsion}$, and we show that the invariant $\overline{\CQ}(M)$ of $M$, that is the class of $\CQ(\KK)$ in the quotient
of $R$ by this subspace, detects the connected sums with rational homology spheres with non trivial Casson-Walker invariant, when the Alexander polynomial of $M$ has no multiple roots.
Conjecturally, the invariant $\overline{\CQ}(M)$ refines the two-loop part of an invariant
recently defined by Ohtsuki.

\vskip.5cm

\end{abstract}

\newpage
\section{Introduction}
\setcounter{equation}{0}
\label{secintro}

\subsection{Short introduction}
The study of $3$--manifold invariants built from
integrals over configuration spaces started
after the work of Witten on Chern-Simons theory in 1989 \cite{wit},
with work of Axelrod, Singer~\cite{as1,as2}, Kontsevich~\cite{ko}, Bott, Cattaneo~\cite{BC,bc2}, Taubes~\cite{taubes}.

The simplest non-trivial of these invariants is an invariant of homology 3-spheres $N$ (equipped with an appropriate trivialisation) that may be written as the integral of the cube of a closed $2$-form $\omega$ over a configuration space $C_2(N)$, where the cohomology class of $\omega$ represents the linking form of $N$. This invariant is associated to the $\theta$--graph,
G. Kuperberg and D. Thurston proved that it is the Casson invariant \cite{kt}.

We shall present a similar invariant for closed oriented $3$-manifolds with first Betti number one, in an equivariant setting, already investigated by Julien March\'e in \cite{Ju} in the knot case.

Let $M$ be such a closed oriented $3$-manifold with first Betti number one.
Its equivariant linking pairing may be seen as a two-dimensional cohomology class in an appropriate infinite cyclic covering of the space of ordered pairs of distinct points of $M$. We shall show how to define
the equivariant cube $\CQ(\KK)$ of this Blanchfield pairing with respect to a framed knot $\KK$ that generates $H_1(M;\ZZ)/\mbox{Torsion}$. We show many properties
for the invariant $\CQ$ including two surgery formulae.

Via surgery, the invariant $\CQ$ is equivalent to an invariant $\hat{\CQ}$ of null-homologous knots in rational homology spheres.
For a knot $\hat{K}$ with trivial Alexander polynomial in an integral homology sphere, one of our surgery formulae combined with results of Garoufalidis and Rozansky \cite{GR} shows that $\hat{\CQ}(\hat{K})$
coincides with the two--loop part of the Kricker
rational lift of the Kontsevich integral of $\hat{K}$.
This lift was defined in \cite{Kr,GK}, from the LMO invariant of Le, Murakami and Ohtsuki \cite{lmo}, following conjectures of Rozansky \cite{Ro}, who may have had in mind the construction that is presented in this article.  Its two-loop part is also called the two-loop polynomial (up to normalization). It was introduced by Rozansky in \cite[Section 6]{Ro} and has been extensively studied by Ohtsuki in \cite{oht2}.

I think that $\hat{\CQ}$
is always equivalent to the two--loop polynomial in the sense that if two knots with equivalent equivariant linking pairings are distinguished by one of these two invariants, then they are distinguished by the other one.

The invariant $\CQ$
takes its values in a polynomial ring $R$ over $\QQ$.
We determine the rational subvector space of $R$ generated by the variations $(\CQ(\KK^{\prime})-\CQ(\KK))$ for pairs of framed knots that generate $H_1(M;\ZZ)/\mbox{Torsion}$, and we define an invariant $\overline{\CQ}(M)$ of $M$ that is the class of $\CQ(\KK)$ in the quotient
of $R$ by this subspace.

Conjecturally, the invariant $\overline{\CQ}(M)$ refines the two-loop part of the invariant
recently defined by Ohtsuki using the LMO invariant in \cite{ohtb}. It detects connected sums with rational homology spheres with non trivial Casson-Walker invariant, when the Alexander polynomial of $M$ has no multiple roots.

The constructions contained in this article apply to the higher loop degree case, and allow us to define an invariant conjecturally equivalent to the whole rational lift of the Kontsevich integral for null-homologous knots in rational homology spheres. This is discussed
in \cite{lesbonn}.

I started to work on this project after a talk of Tomotada Ohtsuki at
a workshop at the CTQM in {\AA}rhus in Spring 2008. I wish to thank Joergen Andersen and Bob Penner for organizing this very stimulating meeting.
Theorem~\ref{thmaugcas} answers a question that George Thompson asked me at the conference {\em Chern-Simons Gauge theory : 20 years after, Hausdorff center for Mathematics\/} in Bonn in August 2009. I thank him for asking and I thank the organizers Joergen Andersen, Hans Boden, Atle Hahn, Benjamin Himpel of this great conference.

\subsection{On the equivariant configuration space $\TCM$}
\label{subtcm}
Consider a closed oriented $3$-manifold $M$ with first Betti number one, and its standard infinite cyclic covering
$\tilde{M}$.
Let $\theta_M$ denote one of the two homeomorphisms that generate the covering group of $\tilde{M}$, and let
$\TCMD$ be the quotient of $\tilde{M}^2$ by the equivalence relation that identifies $(p,q)$ to $(\theta_M(p),\theta_M(q))$.
It is an infinite cyclic covering of $M^2$ whose covering group is generated by 
$\theta$ that acts as follows on the equivalence class $\overline{(p,q)}$ of $(p,q)$
$$\theta\overline{(p,q)}=\overline{(\theta_M(p),q)}=\overline{(p,\theta_M^{-1}(q))}.$$ 

The diagonal of $\tilde{M}^2$ projects to a preferred lift of the diagonal of $M^2$ in $\TCMD$.

The configuration space $C_2(M)$ is the compactification of $\left( M^2 \setminus \mbox{diagonal}\right)$ that is obtained from $M^2$ by {\em blowing-up\/} the diagonal of $M^2$ in the following sense:
The diagonal  of $M^2$ is replaced by the total space of its unit normal bundle. Thus, $C_2(M)$ is diffeomorphic to the complement of an open tubular neighborhood of the diagonal of $M^2$.

The normal bundle $\left(TM^2/\mbox{diag}\right)$ of the diagonal of $M^2$ in $M^2$ is identified to the tangent bundle $TM$ of $M$ via $((x,y) \mapsto y-x)$.
Therefore, $\partial C_2(M)$ is canonically identified to the unit tangent bundle $ST(M)$ of $M$.

The {\em equivariant configuration space\/} $\TCM$ is the $\ZZ$-covering of $C_2(M)$ that is obtained from $\TCMD$ by blowing-up all the lifts of the diagonal of $M^2$ as above.

$$\partial \TCM =\ZZ \times ST(M)$$
where $n \times ST(M)$ stands for $\theta^n(ST(M))$ and $ST(M)$ is the preimage of the preferred lift of the diagonal under the blow-up map.

The configuration space $\TCM$ is a smooth $6$--dimensional manifold in which we can define the equivariant algebraic intersections of two chains $C$ and $A$, with complementary dimensions (i.e. the sum of dimensions is $6$), whose projections in $C_2(M)$ are transverse,
as follows:

$$\langle C,A \rangle_e=\sum_{i\in \ZZ} \langle \theta^{-i}(C),A\rangle_{\TCM} t^i \in \QQ[t^{\pm1}]$$
where $\langle \cdot ,\cdot\rangle$ stands for the algebraic intersection and a $c$-dimensional {\em chain\/} is a linear combination of oriented compact smooth $c$-submanifolds with boundaries and corners. Definitions and conventions about algebraic intersections are recalled in Subsection~\ref{subinteq}.

We can similarly define the equivariant algebraic triple intersection of $3$ codimension $2$ rational chains $C_{\fvarM}$, $C_{\svarM}$, $C_{\tvarM}$ of $\TCM$,
 whose projections in $C_2(M)$ are transverse, as the following polynomial in $\QQ[\svar^{\pm 1},\tvar^{\pm 1}]$:

$$\langle C_{\fvarM},C_{\svarM},C_{\tvarM} \rangle_e=\sum_{(i,j)\in \ZZ^2} \langle C_{\fvarM},\theta^{-i}(C_{\svarM}),\theta^{-j}(C_{\tvarM})\rangle_{\TCM} \svar^i\tvar^j \in \QQ[\svar^{\pm1},\tvar^{\pm1}].$$

\subsection{Some other notation}

In this article, all manifolds are oriented.
\begin{itemize}
\item Boundaries are oriented by the {\em outward normal first\/} convention.
\item $K$ is a knot in $M$ whose homology class generates $H_1(M;\ZZ)/\mbox{Torsion}$ and acts on $\tilde{M}$ by $\theta_M$.
\item $S$ is a closed surface of $M$ that intersects $K$ transversally at one point with a positive sign.
\item Unless otherwise mentioned, homology coefficients are in $\QQ$.
\item $\Delta=\Delta(M)$ is the Alexander polynomial of $M$ and $\delta=\delta(M)$ is the annihilator of $H_{1}(\tilde{M})$. (The polynomial $\delta(M)$ divides $\Delta(M)$ and it has the same roots.) These polynomials are normalized so that $\Delta(1)=\delta(1)=1$, $\Delta(t_M)=\Delta(t_M^{-1})$, and $\delta(t_M)=\delta(t_M^{-1})$. The definitions of $\Delta$ and $\delta$ are recalled in Lemma~\ref{lemhomtilM}.
\item $$\ID=\ID(t)=\frac{1+t}{1-t} + \frac{t\Delta^{\prime}(M)(t)}{\Delta(M)(t)}.$$
\item $$R_{\delta}= \frac{\QQ[\fvar^{\pm 1},\svar^{\pm 1}, \tvar^{\pm 1},\frac{1}{\delta(\fvar)},\frac{1}{\delta(\svar)},\frac{1}{\delta(\tvar)}]}{(\fvar \svar \tvar=1)}.$$
\item A {\em rational homology sphere\/} or {\em $\QQ$--homology sphere\/} is a closed $3$-manifold with the same rational homology as the standard sphere $S^3$.
\item The Casson-Walker invariant will be denoted by
$\lambda$ and normalized like the Casson invariant in \cite{akmc,gm,mar}. If $\lambda_W$ denotes the Walker invariant normalized as in \cite{wal}, then $\lambda=\frac{\lambda_W}{2}$.
\end{itemize}

\subsection{Definition of the invariant $\CQ(\KK)$}
\label{subannmain}

The invariant that we shall study in this article is defined by the following theorem that is a direct consequence of Propositions~\ref{propuninvtripl}, \ref{propvartau}, \ref{propdenwithoutz} and \ref{propGrat}.

\begin{theorem}
\label{thminvtripl}
Let $\tau: TM \rightarrow M \times \RR^3$ be a trivialisation of $TM$
and let $p_1(\tau)$ be its first Pontrjagin class ($p_1(\tau)$ is an integer whose definition is recalled in Subsection~\ref{submanpar}).
Assume that $\tau$ maps oriented unit tangent vectors of $K$ to some fixed
$\qvarM \in S^2$. Then $\tau$ induces a parallelisation of $K$. Let $K_{\parallel}$ be a parallel of $K$ with respect to this parallelisation. Let $K_{\fvarM}$, $K_{\svarM}$, $K_{\tvarM}$ be three disjoint parallels of $K$, on the boundary $\partial N(K)$ of a tubular neighborhood of $K$, that induce the same parallelisation of $K$ as $K_{\parallel}$.
Consider the continuous map
$$\begin{array}{llll}\check{A}(K) \colon &(S^1=[0,1]/(0\sim 1)) \times [0,1] &\rightarrow & C_2(M)\\
&(t,u \in]0,1[)& \mapsto &(K(t),K(t+u)),
\end{array}$$
and its lift ${A}(K) \colon S^1 \times [0,1] \rightarrow  \TCM$
such that the lift of $(K(t),K(t+\varepsilon))$ is in a small neighborhood of the blown-up canonical lift of the diagonal, for a small positive $\varepsilon$. Let $A(K)$ also denote the $2$--chain $A(K)(S^1 \times [0,1])$.

For $\cvarM \in S^2$, let $s_{\tau}(M;\cvarM)=\tau^{-1}(M \times \cvarM) \subset ST(M) \subset \partial \TCM$.

Let $\fvarM$, $\svarM$, $\tvarM$ be three distinct points in $S^2 \setminus\{\qvarM, -\qvarM\}$.

There exist three transverse rational\/ $4$--dimensional chains $G_{\fvarM}$, $G_{\svarM}$ and $G_{\tvarM}$ of $\TCM$ whose boundaries
are 
$$\partial G_{\fvarM} =(\theta -1)\delta(\theta)\left(s_{\tau}(M;\fvarM)-\ID(\theta) ST(M)_{|K_{\fvarM}}\right) ,$$ 
$$\partial G_{\svarM} =(\theta -1)\delta(\theta)\left(s_{\tau}(M;\svarM)-\ID(\theta) ST(M)_{|K_{\svarM}}\right) \mbox{and}$$ $$\partial G_{\tvarM} =(\theta -1)\delta(\theta)\left(s_{\tau}(M;\tvarM)-\ID(\theta) ST(M)_{|K_{\tvarM}}\right)$$
and such that the following equivariant algebraic intersections vanish
$$\langle G_{\fvarM}, A(K)\rangle_e=\langle G_{\svarM}, A(K)\rangle_e=\langle G_{\tvarM}, A(K)\rangle_e=0.$$
Then $$\CQ(K,K_{\parallel})=\frac{\langle G_{\fvarM},G_{\svarM},G_{\tvarM} \rangle_e}{(\fvar -1)(\svar -1)(\tvar -1)\delta(\fvar)\delta(\svar)\delta(\tvar)} - \frac{p_1(\tau)}{4} \in R_{\delta}$$ only depends on the isotopy class of the knot $K$ and on its parallelisation.
\end{theorem}

Fix $$\KK=(K,K_{\parallel}).$$
We shall furthermore show
that
 $$\delta(M)(\fvar)\delta(M)(\svar)\delta(M)(\tvar)\CQ(\KK) \in \frac{\QQ[\fvar^{\pm 1},\svar^{\pm 1}, \tvar^{\pm 1}]}{(\fvar \svar \tvar=1)}$$
in Proposition~\ref{propdenwithoutz}, and that
$$\CQ(\KK)(\fvar,\svar,\tvar)=\CQ(\KK)(\svar,\fvar,\tvar)=\CQ(\KK)(\tvar,\svar,\fvar)=\CQ(\KK)(\fvar^{-1},\svar^{-1},\tvar^{-1})$$
in Subsection~\ref{subsecsym}.
We shall also see that $\CQ((S^1 \times \qvarM,S^1 \times \qvarM^{\prime})\subset S^1 \times S^2) =0$ in Remark~\ref{remvarpar}.

\begin{remark}
In the definition of $\CQ$ by the statement of Theorem~\ref{thminvtripl}, replace the assumption that the knots $K_{\fvarM}$, $K_{\svarM}$, $K_{\tvarM}$ are parallel to $K$ by the assumption that these three knots are rationally homologous to $K$ and that $K_{\fvarM}$, $K_{\svarM}$, $K_{\tvarM}$ and $K$ are disjoint. This defines an invariant of
the $4$-uples $(K_{\fvarM}, K_{\svarM}, K_{\tvarM},\KK)$ that could be interesting to study. See Subsection~\ref{subinvtripl}.
\end{remark}

\subsection{On the equivariant linking number and its cube}

By definition, the  {\em equivariant linking number\/} $lk_e(\alpha,\beta)$ of two cycles $\alpha$ and $\beta$ in $\tilde{M}$ of respective dimensions $\mbox{dim}(\alpha)$ and $\mbox{dim}(\beta)$, such that
\begin{itemize}
\item $\mbox{dim}(\alpha) + \mbox{dim}(\beta)=2$ and 
\item the projections of $\alpha$ and $\beta$ in $M$ do not intersect
\end{itemize} is equal to the equivariant intersection
$$lk_e(\alpha,\beta)=\langle \alpha,B \rangle_{e,\tilde{M}}=\sum_{i\in \ZZ} \langle \theta_M^{-i}(\alpha),B\rangle_{\TCM} t_M^i \in \QQ[t_M^{\pm1}]$$ if $\beta =\partial B$.
In general, $(\theta_M-1)\delta(\theta_M)(\beta)$ bounds a chain $(t_M-1)\delta(t_M) B$ and $$lk_e(\alpha,\beta)=\frac{\langle \alpha,(t_M-1)\delta(t_M)B \rangle_{e,\tilde{M}}}{(t_M^{-1}-1)\delta(t_M^{-1})}$$
and, for a polynomial $P \in \QQ[t_M,t_M^{-1}]$,
$$lk_e(P(\theta_M)(\alpha),\beta)=P(t_M)lk_e(\alpha,\beta)=lk_e(\alpha,P(\theta_M^{-1})(\beta)).$$

See Section~\ref{secblanch} for more definitions and properties of $lk_e$.

Let the multiplication by $t$ stand for morphisms induced by $\theta$ and allow division by polynomials in $t$.
We shall see (Proposition~\ref{propdeflkeq} and Theorem~\ref{thmstauM}) that the assumptions on 
$$F_{\fvarM}=\frac{G_{\fvarM}}{(t -1)\delta(t)},$$ $F_{\svarM}=\frac{G_{\svarM}}{(t -1)\delta(t)}$ and $F_{\tvarM}=\frac{G_{\tvarM}}{(t -1)\delta(t)}$ in Theorem~\ref{thminvtripl} imply that
$$lk_e(\alpha,\beta)=\langle \alpha \times \beta,F_{\fvarM}\rangle_e=\langle \alpha \times \beta,F_{\svarM}\rangle_e=\langle \alpha \times \beta,F_{\tvarM}\rangle_e$$
so that the $4$-dimensional chains $F_{\fvarM}$, $F_{\svarM}$, $F_{\tvarM}$ may be seen as three parallel representatives 
of the equivariant linking number, and $\CQ(\KK)$ may be thought of as the cube of the equivariant linking number with respect to $\KK$.

\subsection{Introduction to the properties of $\CQ$}
\label{subintroprop}

According to Theorem~\ref{thminvtripl}, $\CQ(\KK \subset M)$ is an invariant of framed knots $\KK$ generating $\frac{H_1(M;\ZZ)}{\mbox{\tiny Torsion}}$.

Let $M_{\KK}$ denote the rational homology sphere obtained from $M$ by surgery along $\KK$. This manifold is obtained from $M$ by replacing a tubular neighborhood $N(K)$ of $K$ by another solid torus $N(\hat{K})$ whose meridian is the given parallel of $K$. The core $\hat{K}$ of the new torus $N(\hat{K})$ is a null-homologous knot in $M_{\KK}$.
Since $M$ is obtained from $M_{\KK}$ by $0$-surgery on $\hat{K}$, the data $(M,\KK)$ are equivalent to $(M_{\KK},\hat{K})$.
In particular, the formula
$$\hat{\CQ}(\hat{K} \subset M_{\KK})=\CQ(\KK \subset M)$$
defines an invariant of null-homologous knots $\hat{K}$ in rational homology spheres that will be shown to share a lot of properties with the two--loop polynomial. For example, it will be shown to have the same variation under a surgery on the Garoufalidis and Rozansky degree $2$ null-claspers of \cite{GR} in Lemma~\ref{lemnullclasp}.

We shall see in Propositions~\ref{propsym3} and \ref{propsym4} that $\CQ(\KK)$ is independent of the orientation of $K$, and that
$$\CQ(\KK \subset -M)=-\CQ(\KK \subset M).$$ The following behaviour of $\CQ$ under connected sum with a rational homology sphere will be
established in Subsection~\ref{proofconcas}.

\begin{proposition}
\label{propconcas}
 Let $N$ be a rational homology sphere, then 
$$\CQ(\KK \subset M \sharp N)=\CQ(\KK \subset M) + 6 \lambda(N)$$
where $\sharp$ stands for the connected sum and $\lambda$ is the Casson-Walker invariant.
\end{proposition}

The following property will be proved in Section~\ref{secaug}.
 $$\CQ(\KK \subset M)(1,1,1)=6\lambda(M_{\KK}).$$

The invariant $\CQ$ satisfies the following surgery formula for surgeries on knots. The Dehn surgery formula of Theorem~\ref{thmDehn} will be proved in Subsection~\ref{subDehnproof} from an LP surgery formula (Theorem~\ref{thmLP}) for Lagrangian-Preserving replacements of rational homology handlebodies, that generalize the null borromean surgeries. The LP surgery formula is precisely stated in Subsection~\ref{substateLP}. According to results of Garoufalidis and Rozansky \cite{GR}, it implies that $\hat{Q}$ lifts the (primitive) two-loop part of the Kontsevich integral for knots with trivial Alexander polynomials in integral homology spheres. See Theorem~\ref{thmKont}.

\begin{theorem}
\label{thmDehn}
Let $J$ be a knot of $M$ that bounds
a Seifert surface $\Sigma$ disjoint from $K$
such that $H_1(\Sigma)$ is mapped to $0$ in $H_1(M)$.\\
Let $p/q$ be a nonzero rational number.
Let $(c_i,d_i)_{i=1,\dots,g}$ be a symplectic basis of $H_1(\Sigma)$.
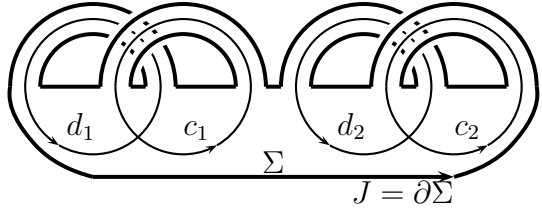
\begin{center}
\begin{pspicture}[shift=-0.1](0,.2)(5.4,2.9)
\psarc[linewidth=1.5pt](1.1,1.6){1.1}{0}{180}
\psarc[linewidth=1.5pt](1.1,1.6){.7}{0}{180}
\psarc{->}(1.1,1.6){.9}{95}{240}
\psarc(1.1,1.6){.9}{240}{95}
\psarc[linewidth=1.5pt,border=2pt](2.3,1.6){1.1}{0}{180}
\psarc[linewidth=1.5pt,border=2pt](2.3,1.6){.7}{0}{180}
\psarc[border=1pt](2.3,1.6){.9}{85}{180}
\psarc{->}(2.3,1.6){.9}{180}{300}
\psarc(2.3,1.6){.9}{-60}{85}
\psecurve[linewidth=1.5pt](1.1,2.7)(0,1.6)(.33,.83)(1.1,.4)(1.3,.4)
\psline[linewidth=1.5pt]{->}(1.1,.4)(5.9,.4)
\rput[b](3.5,.45){$\Sigma$}
\rput[rt](5.9,.35){$J=\partial \Sigma$}
\rput[rb](2.65,.9){$c_1$} 
\rput[lb](.75,.9){$d_1$}
\psline[linewidth=1.5pt](.4,1.6)(1.2,1.6)
\psline[linewidth=1.5pt](1.6,1.6)(1.8,1.6)
\psline[linewidth=1.5pt](2.2,1.6)(3,1.6) 
\psarc[linewidth=1.5pt](4.7,1.6){1.1}{0}{180}
\psarc[linewidth=1.5pt](4.7,1.6){.7}{0}{180}
\psarc{->}(4.7,1.6){.9}{95}{240}
\psarc(4.7,1.6){.9}{240}{95}
\psarc[linewidth=1.5pt,border=2pt](5.9,1.6){1.1}{0}{180}
\psarc[linewidth=1.5pt,border=2pt](5.9,1.6){.7}{0}{180}
\psarc[border=1pt](5.9,1.6){.9}{85}{180}
\psarc{->}(5.9,1.6){.9}{180}{300}
\psarc(5.9,1.6){.9}{-60}{40}
\psarc[border=1pt](5.9,1.6){.9}{40}{85}
\psecurve[linewidth=1.5pt](5.9,2.7)(7,1.6)(6.67,.83)(5.9,.4)(5.7,.4)
\rput[rb](6.25,.9){$c_2$} 
\rput[lb](4.35,.9){$d_2$}
\psline[linewidth=1.5pt](3.4,1.6)(3.6,1.6)
\psline[linewidth=1.5pt](4,1.6)(4.8,1.6)
\psline[linewidth=1.5pt](5.2,1.6)(5.4,1.6)
\psline[linewidth=1.5pt](5.8,1.6)(6.6,1.6)
\end{pspicture}
\end{center}
Let $$\lambda_e^{\prime}(J)=\frac{1}{12}\sum_{(i,j) \in \{1,\dots,g\}^2}\sum_{\mathfrak{S}_3(\fvar,\svar,\tvar)} \left(\alpha_{ij}(\fvar,\svar)+\alpha_{ij}(\fvar^{-1},\svar^{-1})+\beta_{ij}(\fvar,\svar)\right)$$ where 
$$\alpha_{ij}(\fvar,\svar)=lk_e(c_i ,c_j^+)(\fvar)lk_e(d_i ,d_j^+)(\svar) -lk_e(c_i ,d_j^+)(\fvar)lk_e(d_i ,c_j^+)(\svar),$$
$$\beta_{ij}(\fvar,\svar)=\left(lk_e(c_i,d^+_i )(\fvar)-lk_e(d^+_i,c_i)(\fvar)\right)\left(lk_e(c_j,d^+_j )(\svar)-lk_e(d^+_j,c_j)(\svar)\right)$$
and \/ $\sum_{\mathfrak{S}_3(\fvar,\svar,\tvar)}$ stands for the sum over the $6$ terms
obtained by replacing $\fvar$ and $\svar$ by their images under the $6$ permutations of $\{\fvar,\svar,\tvar\}$,
then $$\CQ(\KK \subset M(J;p/q))-\CQ(\KK \subset M)= 6\frac{q}{p}\lambda_e^{\prime}(J) + 6 \lambda(S^3(U ;p/q))$$
where $S^3(U ;p/q)$ is the lens space $L(p,-q)$ obtained from $S^3$ by $p/q$--surgery on the unknot $U$.
\end{theorem}

Since $H_1(\Sigma)$ goes to $0$ in $H_1(M)$
in the above statement, $\Sigma$ lifts as homeomorphic copies of $\Sigma$ in $\tilde{M}$ and $lk_e(c_i ,c_j^+)$ denotes the equivariant linking number of a lift of $c_i$ in $\tilde{M}$ in some lift of $\Sigma$ and a lift of $c^+_j$ near the same lift of $\Sigma$. The superscript $+$ means that 
$c_j$ is pushed in the direction of the positive normal to $\Sigma$.

When the above knot $J$ is inside a rational homology ball, $\lambda_e^{\prime}(J)$ coincides with $\frac{1}{2}\Delta^{\prime\prime}(J)$, where $\Delta(J)$ is the Alexander polynomial of $J$, and the right-hand side is nothing but $6$ times the variation of the Casson-Walker invariant 
under a $p/q$--surgery on $J$. This is consistent with Proposition~\ref{propconcas}.

Recall $I_{\Delta}(t)=\frac{1+t}{1-t} +  \frac{t\Delta^{\prime}(t)}{\Delta(t)}$.
The following theorem will be proved in Section~\ref{secvarcob}.
\begin{theorem}
\label{thmfrakcha}
Let $\KK^{\prime}$ be another framed knot of $M$ such that $H_1(M)/\mbox{Torsion} = \ZZ[K^{\prime}]$.
Then there exists an antisymmetric polynomial 
${\cal V}(\KK,\KK^{\prime})$ in $\QQ[t,t^{-1}]$ such that 
$$\CQ(\KK^{\prime}) - \CQ(\KK)=\sum_{\mathfrak{S}_3(\fvar,\svar,\tvar)}\frac{{\cal V}(\KK,\KK^{\prime})(x)}{\delta(x)}I_{\Delta}(y).$$

Furthermore, for any $k \in \ZZ$, there exists a pair of framed knots $(\KK,\KK^{\prime})$ such that ${\cal V}(\KK,\KK^{\prime})=q(t^k-t^{-k})$ for some nonzero rational number $q$.
\end{theorem}

The following two propositions use the notation of Theorem~\ref{thmfrakcha}. The first one is a direct consequence of Proposition~\ref{propvarpar},
and the second one is proved right after Lemma~\ref{lemkeyvarsur}.
\begin{proposition} If $\KK=(K,K_{\parallel})$ and if $\KK^{\prime}=(K,K^{\prime}_{\parallel})$, where $K^{\prime}_{\parallel}$ is another parallel of $K$ such that the difference $(K^{\prime}_{\parallel} -K_{\parallel})$ is homologous to a positive meridian of $K$ in $\partial N(K)$, then
$${\cal V}(\KK,\KK^{\prime})(t)=-\frac{\delta (t)}{2}\frac{t\Delta^{\prime}(t)}{\Delta(t)}.$$
\end{proposition}

\begin{proposition}
\label{propcorcalvarb} 
If $K$ and $K^{\prime}$ coincide along an interval, if $(K^{\prime}-K)$ bounds a surface $\Bor$ that lifts in $\tilde{M}$ such that $(K^{\prime}_{\parallel} - K_{\parallel})$ is homologous to a curve of $\Bor$ in the complement of $(\partial \Bor \cup K)$ in a regular neighborhood of $\Bor$, and
if $(a_i,b_i)_{i \in \{1,\dots,g\}}$ is a symplectic basis of $H_1(\Bor;\ZZ)$, then
$$\frac{{\cal V}(\KK,\KK^{\prime})(t)}{\delta (t)}=\sum_{i=1}^g\left(lk_e(a_i,b_i^+) -\overline{lk_e(a_i,b_i^+)}\right).$$
\end{proposition}

\subsection{The derived $3$-manifold invariant}

\begin{definition}{\em Definition of an invariant for $3$-manifolds of rank one:\/}\\
Let $Q_k(\delta,\Delta)=\sum_{\mathfrak{S}_3(\fvar,\svar,\tvar)}\frac{x^k-x^{-k}}{\delta(x)}I_{\Delta}(y)$ for $k \in (\NN \setminus \{0\})$.
For a fixed $(\delta,\Delta)$, define ${\overline{\CQ}}(M)$ in the 
quotient of $R_{\delta}$ by the vector space generated by the $Q_k(\delta,\Delta)$ for $k \in (\NN \setminus \{0\})$ as the class of $\CQ(\KK \subset M)$.
\end{definition}

For $3$-manifolds $M$ such that $H_1(M;\ZZ)=\ZZ$ and $\Delta(M)=1$, the invariant ${\overline{\CQ}}$ coincides with the two-loop part of the invariant defined by Ohtsuki in 2008, combinatorially in \cite{ohtb}
for $3$-manifolds of rank one, up to normalization. 
I think that the two invariants are equivalent, in the sense that they distinguish the same pairs of manifolds with isomorphic equivariant pairings, when $\Delta=\delta$, and that ${\overline{\CQ}}$ refines the two-loop part of the Ohtsuki invariant when $\Delta\neq\delta$.

The following proposition is a direct consequence of Proposition~\ref{propinj} that is proved in Section~\ref{sectarg}.

\begin{proposition}
If $\Delta$ has only simple roots and
if $N$ is a rational homology sphere such that $\lambda(N)\neq 0$,
then $\overline{\CQ}(M) \neq \overline{\CQ}(M\sharp N)$.
\end{proposition}

\subsection{On the derived knot invariant $\hat{Q}$}

Here, we list a few properties of the invariant $\hat{Q}$ of null-homologous knots in rational homology spheres that was defined in the beginning of Subsection~\ref{subintroprop}. These properties obviously follow from the corresponding properties of $\CQ$. Let us first recall the definition of $\hat{Q}$.
Let $\hat{K}$ be a null-homologous knot in a rational homology sphere $N$. Let $N(\hat{K};0)$ be the $3$-manifold obtained by $0$-surgery on $\hat{K}$.
Let $K$ be the core of the solid torus in $N(\hat{K};0)$ that is glued during the surgery, and let $K_{\parallel}$ be the meridian of $\hat{K}$ on the boundary of this torus.
Then $$\hat{\CQ}(\hat{K} \subset N)=\CQ((K,K_{\parallel}) \subset N(\hat{K};0)) \in R_{\delta(N(\hat{K};0))}.$$
The {\em Alexander polynomial\/} of $\hat{K}$ is defined like the Alexander polynomial of $M$ in Lemma~\ref{lemhomtilM}, it is normalized so that
$\Delta(\hat{K})=\Delta(N(\hat{K};0))$.

\begin{proposition}
The invariant $\hat{\CQ}$ satisfies the following properties for any null-homologous knot $\hat{K}$ in a rational homology sphere $N$.
\begin{enumerate}
 \item $\delta(M)(\fvar)\delta(M)(\svar)\delta(M)(\tvar)\hat{\CQ}(\hat{K}) \in \frac{\QQ[\fvar^{\pm 1},\svar^{\pm 1}, \tvar^{\pm 1}]}{(\fvar \svar \tvar=1)}$.
\item $\hat{\CQ}(\hat{K})(\fvar,\svar,\tvar)=\hat{\CQ}(\hat{K})(\svar,\fvar,\tvar)=\hat{\CQ}(\hat{K})(\tvar,\svar,\fvar)=\hat{\CQ}(\hat{K})(\fvar^{-1},\svar^{-1},\tvar^{-1})$.
\item $\hat{\CQ}(\hat{K})(1,1,1)=6\lambda(N)$.
\item $\hat{\CQ}(-\hat{K})=\hat{\CQ}(\hat{K})$, $\hat{\CQ}(\hat{K} \subset (-N) )=-\hat{\CQ}(\hat{K} \subset N)$.
\item For any other rational homology sphere $N_2$, $\hat{\CQ}(\hat{K}\subset N\sharp N_2)=\hat{\CQ}(\hat{K}\subset N)+6\lambda(N_2)$.
\item For any knot $J$ of $N$ that bounds a Seifert surface $\Sigma$ disjoint from $\hat{K}$ such that $H_1(\Sigma)$ maps to $0$ in $H_1(N \setminus \hat{K})$, for any $r$ in $\QQ \setminus \{0\}$,
$$\hat{\CQ}(\hat{K} \subset N(J;r))-\hat{\CQ}(\hat{K} \subset N)= \frac{6}{r}\lambda_e^{\prime}(J) + 6 \lambda(S^3(U ;r))$$
where $\lambda_e^{\prime}(J)$ is defined like in Theorem~\ref{thmDehn}, and the other used notation can be found there, too.
\item If $J$ is the parallel of $\hat{K}$ that bounds in $N \setminus \hat{K}$, then 
$$\hat{\CQ}(\hat{K} \subset N(J;1))-\hat{\CQ}(\hat{K} \subset N)=\sum_{\mathfrak{S}_3(\fvar,\svar,\tvar)}\frac{\fvar\Delta^{\prime}(\fvar)}{2\Delta(\fvar)}I_{\Delta}(\svar).$$
\item The invariant $\hat{\CQ}$ satisfies the $LP$-surgery formula of Theorem~\ref{thmLP} when $A$ and $B$ are disjoint rational homology handlebodies of $N \setminus \hat{K}$ whose $H_1$ map to $0$ in $H_1(N \setminus \hat{K})$.
\item If $N$ is an integral homology sphere, and if $\Delta(\hat{K})=1$, then
$\hat{\CQ}(\hat{K})$ lifts the (primitive) two-loop part of the Kontsevich integral. See Theorem~\ref{thmKont}.
\end{enumerate}
\end{proposition}
Note that the fourth property implies that $\hat{\CQ}(U)=0$ for the trivial knot $U$ of $S^3$.

\subsection{Organization of the article}

The claims contained in the statement of Theorem~\ref{thminvtripl} are
mostly homological:

The existence of the chain $G_{\fvarM}$ will be proved as soon as we have proved that the $3$--cycle 
$$(\theta -1)\delta(\theta)\left(s_{\tau}(M;\fvarM)-\ID(\theta) ST(M)_{|K_{\fvarM}}\right)$$
vanishes in $H_3(\TCM)$.

The fact that the equivariant triple algebraic intersection 
$\langle G_{\fvarM},G_{\svarM},G_{\tvarM} \rangle_e$ is well-defined by the statement comes from the fact that $G_{\fvarM}$ is well-defined up to the boundary of a $5$-chain in the interior of $\TCM$, and from the similar facts for $G_{\svarM}$ and $G_{\tvarM}$.

Therefore, we shall begin with a study of the homology of 
$\TCM$ in Section~\ref{sechomtcm}. We shall recall known facts about equivariant intersections and equivariant linking numbers in Section~\ref{secblanch}. Then we shall compare various homology classes of $\TCM$ and we shall prove the existence of $G_{\fvarM}$, $G_{\svarM}$ and $G_{\tvarM}$ (but not yet their rationality) in Section~\ref{seccomparhom}.
In Section~\ref{sectriple}, we shall define an invariant $\CQ(K_{\fvarM},K_{\svarM},K_{\tvarM},\KK,\tau)$, determine the dependence on $\tau$ to complete the definition of $\CQ(\KK)$, and study an example with $M=S^1\times S^2$.

In Section~\ref{secconcas}, we give a configuration space definition for the Casson-Walker invariant $\lambda$, and we prove that
$\CQ(\KK \subset M \sharp N)=\CQ(\KK \subset M)+6 \lambda(N)$, for a rational homology sphere $N$.
The dependence of $\CQ$ on the parallelisation of $K$ is discussed in Section~\ref{secvar}.

In Section~\ref{secvarcob}, we compute $(\CQ(\KK^{\prime})-\CQ(\KK))$ when $ \KK^{\prime}$ is another framed knot whose homology class generates $H_1(M)/\mbox{Torsion}$. From Section~\ref{secsur} to \ref{secproofLP}, we prove our two surgery formulae that generalize properties of the two-loop polynomial described in \cite{oht2}. In order to prove the LP surgery formula in full generality in Section~\ref{secproofLP}, we give a more flexible definition of $\CQ$ in Section~\ref{secpseudotriv} by introducing {\em pseudo-trivialisations\/} that generalize the trivialisations used in Theorem~\ref{thminvtripl}.

In Section~\ref{secspecK}, we specialize our chains $G_{\fvarM}$,
$G_{\svarM}$ and $G_{\tvarM}$ in the preimage of $K \times M \cup M \times K$. We say that we show how $K$ {\em interacts\/} via the chains $G_{\fvarM}$,
$G_{\svarM}$ and $G_{\tvarM}$. This allows us to prove that these chains are rational, in Proposition~\ref{propGrat}, and that $\delta(x)\delta(y)\delta(z)\CQ(\KK)$ is in $\QQ[\fvar^{\pm1},\svar^{\pm1}]$.
To prove this result, we give an alternative definition of $\CQ(\KK)$ as an equivariant triple algebraic intersection in an equivariant configuration space associated with the complement of $K$ in Subsection~\ref{subdefbord}. This alternative definition is very close to work of Julien March\'e \cite{Ju}.
This definition also allows us to compute the evaluation at $(1,1,1)$ of $\CQ(\KK)$ and to find that it is $6\lambda(M_{\KK})$, in Section~\ref{secaug}.
Finally, in Section~\ref{sectarg}, we prove that the induced invariant of $3$-manifolds detects the connected sums with rational homology spheres when $\Delta$ has no multiple roots. The various sections are as independent as possible.

\bigskip

More precisely, we could partition this big article into $5$ separate blocks.
The first block is the {\em definition block\/} from Section 1 to Section~\ref{sectriple}, together with Section~\ref{secconcas} that provides useful links with the construction of the Casson invariant.
The other blocks rely on this first block but they do not rely on other blocks
except for a few specific statements, so that they can be read independently of each other.
The {\em dependence on the framed knot $K$ in $M$ block\/} goes from Section~\ref{secvar} to Section~\ref{secvarcob}, 
the {\em surgery formula block\/} goes from Section~\ref{secsur} to Section~\ref{secproofLP}, 
then
Section~\ref{secspecK} and Section~\ref{secaug} also form a block even if Section~\ref{secaug} also depends on the properties of pseudo-trivialisations studied in Section~\ref{secpseudotriv}, 
and Section~\ref{sectarg} is the final small block.

\newpage 
\section{On the homology of $\TCM$}
\setcounter{equation}{0}
\label{sechomtcm}

\subsection{On the homology of $M\setminus S$}
Recall that $S$ is a closed surface of $M$ that intersects $K$ transversally at one point with a positive sign.
Identify $H_{\ast}(S^-=S \times \{-1\})$, $H_{\ast}(S^+=S \times \{+1\})$ and $H_{\ast}(S \times [-1,1])$, naturally, and, write the element of $H_{\ast}(S^{+})$ (resp. $H_{\ast}(S^{-})$) corresponding to an element $c$ of $H_{\ast}(S \times [-1,1])$ as $c^+$ (resp. $c^-$). Let $\ast$ be a basepoint of $S$.

\begin{lemma}
\label{lemhomMsetminusS}
$H_i(M\setminus S;\ZZ)=0$ for any $i \geq 3$,
$H_2(M\setminus S;\ZZ)=\ZZ[S^+]$ and $H_0(M\setminus S;\ZZ)=\ZZ[\ast^+]$.\\
Let $(z_i)_{i=1, \dots 2g}$ and $(z^{\ast}_i)_{i=1, \dots, 2g}$ be two dual bases
of $H_1(S;\ZZ)$ such that 
$$\langle z_i, z_j^{\ast}\rangle=\delta_{ij}=\left\{ \begin{array}{ll}0 \; & \mbox{if}\; i \neq j\\
1 & \mbox{if}\; i = j.\end{array}\right.$$
Then $$H_1(M \setminus S;\QQ) =\bigoplus_{i=1}^{2g} \QQ[z_i^{+} -z_i^{-}]$$
and for any $v\in H_1(M \setminus S;\QQ)$,
$$v=\sum_{i=1}^{2g}lk(v,z_i^{\ast})(z_i^{+} -z_i^{-}).$$
\end{lemma}
\bp
Let $\CB$ be a basis of $H_{i-1}(S)$ such that any $[b] \in \CB$ is represented by an $(i-1)$-submanifold $b$ of $S$. Then
$$H_{i}(M, M\setminus S;\ZZ) \cong H_{i}(S \times [-1,1],S^+ \sqcup S^-;\ZZ)$$
is freely generated by the classes $[b\times[-1,1]]$, for $[b] \in \CB$.
The natural maps $$H_{i}(M;\ZZ) \rightarrow H_{i}(M, M\setminus S;\ZZ)$$ are therefore
surjective for $i=1$ and $i=3$, and the kernel is the torsion of $H_{1}(M;\ZZ)$ when $i=1$.
In particular, the natural map $H_{2}(M\setminus S;\ZZ) \rightarrow (H_{2}(M;\ZZ)=\ZZ[S])$ is injective, and since it is obviously surjective, it is an isomorphism.
Then the boundary map
$$\partial_{2} \colon H_{2}(M, M\setminus S;\ZZ) \rightarrow H_{1}( M\setminus S;\ZZ)$$
is injective, and it is an isomorphism if we take coefficients in $\QQ$.
This gives the announced expression for $H_1(M \setminus S;\QQ)$.

Recall that the linking number of two disjoint null-homologous links makes sense in $M$, and if $u$ is a curve of $S$, $lk(z_i^{+} -z_i^{-},u)=\langle z_i,u\rangle$. This allows us to express the coordinates of an element $v\in H_1(M \setminus S;\QQ)$ as in the statement.
\eop

\subsection{On the homology of $\tilde{M}$, the Alexander polynomial and the annihilator}
\label{subsechomtilM}

Let $p_M \colon \tilde{M} \rightarrow M$ be the infinite cyclic covering of $M$. $$p_M^{-1}(S)=\sqcup_{n\in \ZZ} \theta_M^n(\hat{S})$$
where $\hat{S}$ is a homeomorphic lift of $S$ in $\tilde{M}$. 
Let $\ast \in S$ and let $\hat{\ast}$ be its lift in $\hat{S}$.
Here, we compute the homology of $\tilde{M}$ with coefficients in $\QQ$ and we see it as
endowed with a structure of $\QQ[t_M,t_M^{-1}]$-module where the multiplication by $t_M$ is induced by the action of $\theta_M$ on  $\tilde{M}$.
Set 
$$\Lambda_M=\QQ[t_M,t_M^{-1}].$$

For a polynomial $P$ of $t_M^{1/2}\Lambda_M$, $\frac{\Lambda_M}{(P)}$ will denote the $\Lambda_M$-module $\frac{\Lambda_M}{(t_M^{1/2}P)}$ quotient of $\Lambda_M$ by the $\Lambda_M$-ideal generated by $t_M^{1/2}P$. (We may need polynomials of $t_M^{1/2}\Lambda_M$ when $H_1(M;\ZZ)$ has $2$--torsion. See Remark~\ref{remtwotor}.)
The homology of $\tilde{M}$ is given by the following lemma.

\begin{lemma}
\label{lemhomtilM}
Let $(z_i)_{i=1, \dots 2g}$ be a basis
of $H_1(S;\ZZ)$.
Set
$$\Delta(M)=\mbox{det}\left([t_M^{1/2}lk(z^{+}_j,z_i)-t_M^{-1/2}lk(z^{-}_j,z_i)]_{(i,j)\in \{1,2,\dots,2g\}^2}\right).$$
$$H_{0}(\tilde{M}) = \frac{\Lambda_M}{(t_M-1)}[\hat{\ast}],\;\;\; H_{2}(\tilde{M}) = \frac{\Lambda_M}{(t_M-1)}[\hat{S}],\;\;\; H_{3}(\tilde{M}) = 0$$
and
$$H_{1}(\tilde{M}) \cong \oplus_{i=1}^k\frac{\Lambda_M}{(\delta_i(M))}$$
for an integer $k\geq 0$ and polynomials $\delta_i(M)$ of $\Lambda_M \cup t_M^{1/2}\Lambda_M$ such that
\begin{itemize}
 \item  $\delta_i(M)(1)=1$,
\item $\delta_i(M)(t_M)=\delta_i(M)(t_M^{-1})$,
\item $\delta_i(M)$ divides $\delta_{i+1}(M)$,
\item  $\delta_1(M)$ is not a unit,
\item $\Delta(M)=\prod_{i=1}^k\delta_i(M)$.
\end{itemize}
Then the {\em Alexander polynomial of $M$\/} is $\Delta(M)$, it satisfies $$\Delta(M)(t_M)=\Delta(M)(t_M^{-1})$$
and the {\em annihilator\/} of $H_{1}(\tilde{M})$ is $\delta(M)=\delta_k(M)$ (or $1$ if $\Delta(M)=1$, when $k=0$).
\end{lemma}

Like the previous one, the following proof is classical but it is useful to recall it to introduce notation.

\noindent{\sc Proof of Lemma~\ref{lemhomtilM}:}
In $\tilde{M}$, both $S$ and $M\setminus S$ lift homeomorphically.
Consider a lift $(M\setminus S)_0$ of $M\setminus S$ and a lift $(S \times [-1,1])_0$
of $S \times [-1,1]$ that intersect along a lift of $S \times [-1,0[$ so that
$\theta_M((M\setminus S)_0)$ intersects $(S \times [-1,1])_0$ along a lift of $S \times ]0,1]$.

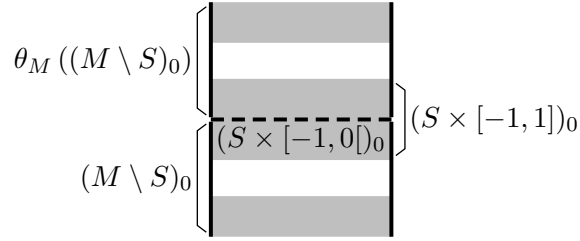
\begin{center}
\begin{pspicture}[shift=-0.1](0,-.2)(7.2,3.2)
\psset{xunit=.6cm,yunit=.6cm}
\pspolygon*[linecolor=lightgray](4,1.7)(8,1.7)(8,3.5)(4,3.5)
\pspolygon*[linecolor=lightgray](4,4.3)(8,4.3)(8,5.2)(4,5.2)
\pspolygon*[linecolor=lightgray](4,0)(8,0)(8,.9)(4,.9)
\psline[linewidth=1.5pt](4,0)(4,2.55)
\psline[linewidth=1.5pt](4,2.65)(4,5.2)
\psline[linewidth=1.5pt](8,0)(8,2.55)
\psline[linewidth=1.5pt](8,2.65)(8,5.2)
\psline[linewidth=1.5pt, linestyle=dashed](4,2.6)(8,2.6)
\psline[linewidth=.5pt](3.9,.1)(3.7,.25)(3.7,2.35)(3.9,2.5)
\rput[r](3.6,1.3){\small $(M\setminus S)_0$}
\psline[linewidth=.5pt](3.9,2.7)(3.7,2.85)(3.7,4.95)(3.9,5.1)
\rput[r](3.6,3.9){\small $\theta_M\left((M\setminus S)_0\right)$}
\psline[linewidth=.5pt](8.1,3.4)(8.3,3.25)(8.3,1.95)(8.1,1.8)
\rput[l](8.4,2.6){\small $(S \times [-1,1])_0$}
\rput(6,2.15){\small $(S \times [-1,0[)_0$}
\end{pspicture}
\end{center}

$\tilde{M}=p_M^{-1}(M\setminus S) \cup p_M^{-1}(S \times [-1,1])$
where
$$p_M^{-1}(M\setminus S)=\sqcup_{i\in \ZZ} \theta_M^i((M\setminus S)_0)\;\;\;
\mbox{and}\;\;\;
p_M^{-1}(S \times [-1,1])= \sqcup_{i\in \ZZ} \theta_M^i((S \times [-1,1])_0) .$$
$$H_{\ast}(p_M^{-1}(M\setminus S))= H_{\ast}(M \setminus S)\otimes_{\QQ} \Lambda_M, \;\;\; H_{\ast}(p_M^{-1}(S \times [-1,1]))= H_{\ast}(S \times [-1,1])\otimes_{\QQ} \Lambda_M, $$
$$H_{\ast}(\tilde{M},p_M^{-1}(M\setminus S))= H_{\ast}(p_M^{-1}(S \times [-1,1]),p_M^{-1}(S^+ \sqcup S^-))=H_{\ast}((S \times [-1,1])_0,S^+ \sqcup S^-)\otimes_{\QQ} \Lambda_M .$$
Consider the boundary maps
$$\partial_{i} \colon H_{i}(\tilde{M},p_M^{-1}(M\setminus S)) \rightarrow H_{i-1}( M\setminus S)\otimes_{\QQ} \Lambda_M.$$
Recall that $\ast$ is a basepoint of $S$, let $\ast^+$ and $\ast^-$ denote the corresponding basepoints of $S^+$ and $S^-$. Let $\hat{\ast}^+$ and $\hat{\ast}^-$ be their respective lifts in $(M\setminus S)_0$.
The map $\partial_1$ maps the preferred generator of $H_{1}(\tilde{M},p_M^{-1}(M\setminus S))$ to
$t_M[\hat{\ast}^+] - [\hat{\ast}^-]=(t_M-1)[\hat{\ast}^+]$. It is therefore injective and
$$H_{0}(\tilde{M}) = \frac{\Lambda_M}{(t_M-1)}[\hat{\ast}^+] (\cong \QQ).$$
Similarly, $\partial_3$ maps the preferred generator of $H_{3}(\tilde{M},p_M^{-1}(M\setminus S))$ to
$t_M[\hat{S}^+] - [\hat{S}^-]=(t_M-1)[\hat{S}^+]$ and is injective.
The map $\partial_2$ maps a basis of $H_{\ast}(\tilde{M},p_M^{-1}(M\setminus S))$
to the system $$(t_M\hat{z}_i^{+} -\hat{z}_i^{-})_{i \in \{1,2,\dots, 2g\}}$$
of $H_{i-1}( M\setminus S)\otimes_{\QQ} \Lambda_M$. 
Its determinant reads
$$\mbox{det}(\partial_2)=\frac{\bigwedge_{j=1}^{2g}(t_M\hat{z}^{+}_j-\hat{z}^{-}_j)}{\bigwedge_{j=1}^{2g}(\hat{z}^{+}_j-\hat{z}^{-}_j)}.$$
In particular, this determinant maps $t_M=1$ to $1$ and it does not vanish.
Therefore, $\partial_2$ is injective, too.

Hence, for any $i$, $H_{i}(\tilde{M})$ is the cokernel of $\partial_{i+1}$.
This gives the result for $i=2$, and we are left with the computation of $H_{1}(\tilde{M})$.

According to the theory of modules over principal domains (see \cite[Chapter 7]{bou}), $H_{1}(\tilde{M})$ reads $\oplus_{i=1}^k\frac{\Lambda_M}{\delta_i}$
for polynomials $\delta_i$ of $\Lambda_M$ such that $\delta_i$ divides $\delta_{i+1}$ and  $\delta_1$ is not a unit, and these polynomials are well-defined up to multiplications by units $(qt_M^{\pm k})$ of $\Lambda_M$. Furthermore, 
$\prod_{i=1}^k\delta_i=\mbox{det}(\partial_2)$ up to units of $\Lambda_M$.
Since $\mbox{det}(\partial_2)$ maps $t_M=1$ to $1$, we can assume that
$\delta_i(1)=1$ for any $i$. The $\delta_i$ are now well-defined up to multiplications by units $(t^{\pm k})$ of $\QQ[t^{\pm 1}]$.

The matrix of $\partial_2$ is a presentation matrix for $H_{1}(\tilde{M})$, it is equivalent to the matrix $$[t_Mlk(z_j^+,z_i)-lk(z_j^-,z_i)]_{(i,j)\in \{1,2,\dots,2g\}^2},$$
it is also equivalent to its transposed matrix $[t_Mlk(z_i^+,z_j)-lk(z_i^-,z_j)]$ and by multiplication by the unit $(-t_M^{-1})$ to
$[t_M^{-1}lk(z_j,z_i^-)-lk(z_j,z_i^+)]$ which is the initial matrix where
$t_M$ has been changed into $t_M^{-1}$. Therefore, $\delta_i(t_M)=\delta_i(t_M^{-1})$ up to units and it suffices to multiply $\delta_i$ by some $t_M^{r/2}$ to normalize it as wanted.
\eop 

\begin{remark}
\label{remtwotor}
We could be more precise and notice that 
$$\delta_i(M) \in \Lambda_M \cup (t_M^{1/2}+t_M^{-1/2})\Lambda_M $$
and that if $H_1(M;\ZZ)$ has no $2$--torsion, then $\delta_i(M) \in \Lambda_M$, for any $i$.

Indeed, if the degree (difference between highest degree and lowest degree) of $\delta_i(M)$ is even, then $\delta_i(M) \in \Lambda_M$. Since
 $\delta_i(M)(t_M)=\delta_i(M)(t_M^{-1})$, if $x$ is a root of $\delta_i(M)$, then
$x^{-1}$ is a root with the same order. In particular, if $x\neq x^{-1}$ for all the roots $x$ of $\delta_i(M)$, then the degree of $\delta_i(M)$ is even. Since $1$ cannot be a root of 
$\delta_i(M)$, the only annoying root could be $(-1)$ and this leads to the first part of the statement.
To prove the second part,
note that $(-1)$ cannot be a root of $\Delta$ if there is no two-torsion in $H_1(M)$. Indeed, in this case, there exists an odd integer $r$ such that $r\Delta(M)$ has integral coefficients. Then
$r\Delta(M)(1)-r\Delta(M)(-1)$ is even and since $r\Delta(M)(1)=r$,
$\Delta(M)(-1)\neq 0$. See \cite{moussard} for more precise results about the occuring polynomials.
\end{remark}

\subsection{First remarks on the homology of $\TCMD$}
\label{subhomTCMD}

The homology of $\TCMD$ with coefficients in $\QQ$ is
endowed with a structure of $\QQ[t,t^{-1}]$-module where the multiplication by $t$ is induced by the action of $\theta$ on  $\TCMD$.
Set 
$$\Lambda=\QQ[t,t^{-1}].$$

Recall that $\delta=\delta(M)$ is the annihilator of $H_1(\tilde{M})$. 
Since $\delta$ has the same roots as $\Delta$, $\delta \frac{\Delta^{\prime}}{\Delta} \in \Lambda \cup t^{1/2}\Lambda$.

\begin{proposition}
\label{proptilMtwotor}
For any element $x \in H_{\ast}(\TCMD;\QQ)$,
$(t-1)^2 \delta(M)(t)^2 x = 0$.
\end{proposition}
To prove this proposition -which would be sufficient to define the invariant $\CQ$ but not sufficient to get the correct denominators-, we shall use the following lemma.

\begin{lemma}
\label{lemdesctcmd}
Let $\tilde{M}_1$ and $\tilde{M}_2$ be two identified copies of $\tilde{M}$.  Let $S^{I-}=(S \times [-1,0[)$ and $S^{I+}=(S \times ]0,1])$.
Then $$\TCMD=\left((S \times [-1,1])\times \tilde{M}_1\right) \sqcup \left((M \setminus S) \times \tilde{M}_2\right)/\sim$$
where $\sim$ identifies 
$$(\sigma^- \in (S^{I-} \subset S \times [-1,1]) , \mu \in \tilde{M}_1) \sim (\sigma^- \in M \setminus S , \mu \in \tilde{M}_2)$$ and
$$(\sigma^+ \in (S^{I+} \subset S \times [-1,1]) , \mu \in \tilde{M}_1) \sim (\sigma^+ \in M \setminus S , \theta_M^{-1}(\mu) \in \tilde{M}_2)=\theta(\sigma^+,\mu)=(\theta_M(\sigma^+),\mu).$$
\end{lemma}
\bp
$\TCMD$ can be constructed as the union of 
$$\TCMD=\left((S \times [-1,1])\times \tilde{M}_1\right) \cup \left((M \setminus S) \times \tilde{M}_2\right)$$
where $\tilde{M}_1$ and $\tilde{M}_2$ are two copies of $\tilde{M}$ that contain a preferred copy $(S \times [-1,1])_1$ of $S \times [-1,1]$, and, a preferred copy $(M \setminus S)_2$ of $M \setminus S$, respectively. These preferred copies are such that the canonical lift of the diagonal of $(S \times [-1,1])^2$ is inside
$(S \times [-1,1]) \times (S \times [-1,1])_1$ and the canonical lift of the diagonal of $(M\setminus S)^2$ is inside
$(M \setminus S) \times (M \setminus S)_2$.

In $\tilde{M}_2$, $(M \setminus S)_2$ intersects a lift $(S \times [-1,1])_2$ of $S \times [-1,1]$ along a lift of $S \times [-1,0[$ and the lift $\theta_M^{-1}((S \times [-1,1])_2)$ along a lift of $S \times ]0,1]$. Identify $(\tilde{M}_2,(S \times [-1,1])_2)$ with $(\tilde{M}_1,(S \times [-1,1])_1)$.

Then  $(S \times [-1,1])\times \tilde{M}_1$ and $(M \setminus S) \times \tilde{M}_2$ intersect along $ (S^{I-} \times \tilde{M}_1) \sqcup  (S^{I+} \times \tilde{M}_1) $
that are seen as natural parts of $(S \times [-1,1])\times \tilde{M}_1$ and that 
map to $(M \setminus S) \times \tilde{M}_2$ by sending the diagonals of $S^{I+}$
and $S^{I-}$ to $(M \setminus S) \times (M \setminus S)_2$ and therefore to
$(M \setminus S) \times \theta_M^{-1}((S \times ]0,1])_2)$, and to $(M \setminus S) \times ((S \times [-1,0[)_2)$, respectively.
This gives the result.
\eop

\noindent{\sc Proof of Proposition~\ref{proptilMtwotor} :}
Let $p \colon \TCMD \rightarrow M^2$ be the covering map.
$$\TCMD=p^{-1}((M \setminus S) \times M) \cup p^{-1}((S \times [-1,1]) \times M)$$
where
$$p^{-1}((M \setminus S) \times M)=(M \setminus S)\times \tilde{M}\;\;\;\;
\mbox{and}
\;\;\;\;p^{-1}((S \times [-1,1]) \times M)= (S \times [-1,1])\times \tilde{M}.$$
Thanks to the K\"unneth formula,
$$\begin{array}{ll}H_{r}((M \setminus S)\times \tilde{M})&= \left(H_{\ast}(M \setminus S)\otimes_{\QQ} H_{\ast}(\tilde{M})\right)_r \\&= \bigoplus_{p,q;p+q=r}\left(H_{p}(M \setminus S)\otimes_{\QQ} H_{q}(\tilde{M})\right)\end{array}$$
is a torsion module over $\QQ[t,t^{-1}]$. More precisely, all its homogeneous elements are 
either $(t-1)$-torsion elements or $\delta(t)$-torsion elements.
$$\begin{array}{ll}H_{i}(\TCMD,(M \setminus S)\times \tilde{M};\QQ)&= H_{i}((S \times [-1,1])\times \tilde{M},(S^+ \sqcup S^-)\times \tilde{M};\QQ)\\&\cong \left(H_{\ast}(S \times [-1,1], S^+ \sqcup S^-)\otimes_{\QQ} H_{\ast} (\tilde{M})\right)_i\end{array}$$
is also a torsion module over $\QQ[t,t^{-1}]$ whose homogeneous elements are annihilated by $(t-1)$ or $\delta(t)$.
Then the long exact sequence associated to the pair $(\TCMD,(M \setminus S)\times \tilde{M})$ allows us to conclude that
$(t-1)^2 \delta(M)^2 x = 0$ for any $x \in H_{\ast}(\TCMD;\QQ)$.
\eop

\subsection{More on the homology of $\TCMD$}
\label{subhomTCMDmore}

The involution that exchanges the two factors in $\tilde{M}^2$ induces an involution $\iota$ of $\TCMD$.
See $S^1$ as the unit circle of $\CC$.
Define
$$\begin{array}{llll}
f\colon & M & \rightarrow &S^1\\
&(M\setminus (S\times ]-1,1[))& \mapsto & -1\\
&(\sigma,u) \in S \times [-1,1]& \mapsto & \exp(i\pi u)
\end{array}$$
Assume that $K$ meets $S \times [-1,1]$ as $\ast \times [-1,1]$, and define a homeomorphism $f_K \colon K \rightarrow S^1$ that coincides with $f$ on $\ast \times [-1/2,1/2]$.

Define $\mbox{\rm diag}_u(K \times M)$  as the lift in $\TCMD$
of $(f_K \times f)^{-1}(\mbox{\rm diag}((S^1)^2)$ that meets the preferred lift of the diagonal of $M^2$. Define $\mbox{\rm diag}_u(M \times K)$, similarly. $$\mbox{\rm diag}_u(M \times K)=\pm \iota(\mbox{\rm diag}_u(K \times M)).$$

Note that if $\gamma$ is a curve homotopic to the diagonal of $(S^1)^2$, then the preimage of the homotopy yields a cobordism between $\mbox{\rm diag}_u(K \times M)$ and a lift of the preimage of $\gamma$ under $(f_K \times f)$. In particular the homotopy $$(u \in [0,1],z\in S^1) \mapsto (z\exp(2i\pi u),z)$$ 
shows that $\theta\left(\mbox{\rm diag}_u(K \times M)\right)$
is homologous to $\mbox{\rm diag}_u(K \times M)$, that is homologous to
all the lifts of $(f_K \times f)^{-1}(\gamma)$.

Similarly, all the lifts in $\TCMD$ of the preimages $(f \times f)^{-1}(\gamma)$ for the curves $\gamma$ of $(S^1)^2$ homotopic to the diagonal of $(S^1)^2$, such that $(f \times f)^{-1}(\gamma)$ is a $5$--chain, are homologous to each other. 
Denote their homology class by $[\mbox{\rm diag}_u(M \times M)]$ and note that $$t[\mbox{\rm diag}_u(M \times M)]=[\mbox{\rm diag}_u(M \times M)].$$ Recall that boundaries are oriented by the {\em outward normal first\/} convention.

\begin{proposition}
\label{proptilMtwocomp}
Recall that
$$H_{1}(\tilde{M}) = \oplus_{i=1}^k\frac{\Lambda_M}{(\delta_i)}[c_i]$$
where $\delta_1$ is not a unit, $\delta_i(t_M)=\delta_i(t_M^{-1})$, $\delta_i$ divides $\delta_{i+1}$ for any $i$ such that $i<k$ and $\delta_k=\delta$.
Then the rational homology of $\TCMD$ reads as follows\\
$\begin{array}{lll}
H_0(\TCMD)&=&\frac{\Lambda}{(t-1)}[\ast \times \ast]\\
H_1(\TCMD)&=&\frac{\Lambda}{(t-1)}[\mbox{\rm diag}(K^2)]\\
H_2(\TCMD)&=&\frac{\Lambda}{(t-1)}[S \times \ast] \oplus \frac{\Lambda}{(t-1)}[ \ast \times S] \oplus \oplus_{(i,j) \in \{1,2,\dots,k\}^2}\frac{\Lambda}{\left(\delta_{\mbox{\tiny \rm min}(i,j)}(t)\right)}[c_i \times c_j]\\
H_3(\TCMD)&=&\frac{\Lambda}{(t-1)}[\mbox{\rm diag}_u(K \times M)] \oplus \frac{\Lambda}{(t-1)}[\mbox{\rm diag}_u(M \times K)] \oplus \oplus_{(i,j) \in \{1,2,\dots,k\}^2}\frac{\Lambda}{\left(\delta_{\mbox{\tiny \rm min}(i,j)}\right)}[C(\Sigma_i \times \Sigma_j)]\\
H_4(\TCMD)&=&\frac{\Lambda}{(t-1)}[S \times S]\\
H_5(\TCMD)&=&\frac{\Lambda}{(t-1)}[\mbox{\rm diag}_u(M \times M)]\\
\end{array}$

where $\ast \times \ast$, $S \times \ast$, $\ast \times S$, $c_i \times c_j$ and $S \times S$ abusively denote the projections in $\TCMD$ of $\hat{\ast} \times \hat{\ast}$, $\hat{S} \times \hat{\ast}$, $\hat{\ast} \times \hat{S}$, $c_i \times c_j$ and $\hat{S} \times \hat{S}$, respectively, 
$\Sigma_i$ is a rational $2$--chain of $\tilde{M}$ whose boundary is $\delta_i c_i$ and $C(\Sigma_i \times \Sigma_j)$ denotes
$$C(\Sigma_i \times \Sigma_j)=\frac{\delta_i(t)}{\delta_{\mbox{\tiny \rm min}(i,j)}(t)} c_i \times \Sigma_j \cup \frac{\delta_j(t)}{\delta_{\mbox{\tiny \rm min}(i,j)}(t)} \Sigma_i \times c_j$$
or its projection in $\TCMD$.
\end{proposition}
\bp First notice that all the mentioned chains are cycles that are annihilated by the corresponding polynomials. For example if $i \leq j$, 
$\partial (\Sigma_i \times \Sigma_j)=\delta_i C(\Sigma_i \times \Sigma_j)$.
Now, continue the proof of Proposition~\ref{proptilMtwotor} to determine
$H_i(\TCMD)$ from the short exact sequences associated with the long exact sequence of the pair $(\TCMD,(M \setminus S)\times \tilde{M})$
$$ 0 \rightarrow \mbox{Coker}(\partial_{i+1}) \rightarrow H_i(\TCMD) \rightarrow \mbox{Ker}(\partial_i) \rightarrow 0$$
where 
$$\partial_{i+1} \colon H_{i+1}(\TCMD,(M \setminus S)\times \tilde{M}) \rightarrow H_{i}((M \setminus S)\times \tilde{M})$$
will be rewritten as
$$F_i =\oplus_{(r,s);r+s=i} F_{r,s}\colon H_{i+1}((S \times [-1,1])\times \tilde{M},(S^+ \sqcup S^-)\times \tilde{M})\rightarrow H_{i}((M \setminus S)\times \tilde{M})$$ 
thanks to the excision isomorphism with
$$\begin{array}{llll}
   F_{r,s}\colon & H_{r}(S)\otimes_{\QQ} H_{s} (\tilde{M})&\rightarrow &H_{r}(M \setminus S)\otimes_{\QQ} H_{s}(\tilde{M})\\
&\sigma \otimes x &\mapsto & \sigma^+ \otimes t_M^{-1}x -\sigma^- \otimes x
  \end{array}$$
thanks to Lemma~\ref{lemdesctcmd}.
Thus, the above short exact sequence may be rewritten as 
$$ 0 \rightarrow \oplus_{(r,s);r+s=i} \mbox{Coker}(F_{r,s}) \rightarrow H_i(\TCMD) \rightarrow \oplus_{(r,s);r+s=i-1} \mbox{Ker}(F_{r,s}) \rightarrow 0$$
and we are left with the computation of the $F_{r,s}$.

\begin{lemma}
Let $\sigma_0 =[\ast]$ and $\sigma_2= [S]$.\\
When $r$ and $s$ belong to $\{0,2\}$, 
$F_{r,s} \colon \frac{\Lambda}{(t-1)}[\sigma_r] \otimes [\sigma_s] \rightarrow \frac{\Lambda}{(t-1)}[\sigma_r] \otimes [\sigma_s]$ is the multiplication by $t-1$ that is zero.
Therefore, $\mbox{Coker}(F_{r,s})=\frac{\Lambda}{(t-1)}[\sigma_r] \otimes [\sigma_s]$.
Furthermore, there is a section $$s_{r,s} \colon \mbox{Ker}(F_{r,s})=\frac{\Lambda}{(t-1)}[\sigma_r] \otimes [\sigma_s]
\rightarrow H_{r+s+1}(\TCMD)$$
and
$$\begin{array}{lll}s_{0,0}(\mbox{Ker}(F_{0,0}))&=&\frac{\Lambda}{(t-1)}[\mbox{\rm diag}(K^2)]\\
s_{0,2}(\mbox{Ker}(F_{0,2}))&=&\frac{\Lambda}{(t-1)}[\mbox{\rm diag}_u(K\times M)]\\
s_{2,0}(\mbox{Ker}(F_{2,0}))&=&\frac{\Lambda}{(t-1)}[\mbox{\rm diag}_u(M\times K)]\\
s_{2,2}(\mbox{Ker}(F_{2,2}))&=&\frac{\Lambda}{(t-1)}[\mbox{\rm diag}_u(M^2)].
  \end{array}$$
\end{lemma}
\bp The first assertions are obvious. Let us prove the existence of a section $s_{0,2}$
such that
 $$s_{0,2}(\mbox{Ker}(F_{0,2}))=\frac{\Lambda}{(t-1)}[\mbox{\rm diag}_u(K\times M)].$$
The class $[\mbox{\rm diag}_u(K\times M)]$ is the class of a lift of the preimage under
$f_K \times f$ of any curve $\gamma$ of $(S^1)^2$ homologous to the diagonal. There exists such a curve $\gamma$ that intersects $\exp(i\pi[-1/2,1/2]) \times S^1$ as $\exp(i\pi[-1/2,1/2])\times 1$. This makes clear that the
image of $[\mbox{\rm diag}_u(K\times M)]$ in $H_3(\TCMD,(M \setminus S)\times \tilde{M})$ comes from $\pm[(\sigma_0 \times[-1/2,1/2]) \times \sigma_2] \in H_{3}((S \times [-1/2,1/2])\times \tilde{M},(S^+ \sqcup S^-)\times \tilde{M})$, and this shows how to define the wanted section. The other cases can be treated similarly.
\eop

\begin{lemma}
If $r\in\{0,2\}$, then $F_{r,1}$ is an isomorphism so that $\mbox{Ker}(F_{r,1}) =0$ and $\mbox{Coker}(F_{r,1}) =0$.
\end{lemma}
\bp
$$F_{r,1} \colon \QQ \otimes H_1(\tilde{M}) \rightarrow H_1(\tilde{M})$$
is the multiplication by $(t-1)$ that is an isomorphism on each $\frac{\Lambda}{(\delta_i)}$ since $(t-1)$ is coprime with $\delta_i$.
\eop

To treat the symmetric case $r=1$, we shall use the following standard diagonalisation lemma over principal domains, see \cite[Chapter 7]{bou}.
\begin{lemma}
Recall
$$H_{1}(\tilde{M}) = \oplus_{i=1}^k\frac{\Lambda_M}{(\delta_i)}[c_i].$$
 There exist a basis $(d_i)_{i \in \{1,2,\dots, 2g\}}$ of $H_1(S) \otimes_{\QQ} \Lambda_M$ and a basis $(c_i)_{i \in \{1,2,\dots, 2g\}}$ of $H_1(M\setminus S) \otimes_{\QQ} \Lambda_M$ such that
$$\begin{array}{lll}
t_Md_i^+ -d_i^- &=\delta_i c_i \;\;\;\;&\mbox{\rm if}\; i \leq k\\
&= c_i \;\;\;\;&\mbox{\rm if}\; i > k
\end{array}$$
in $H_1(M\setminus S) \otimes_{\QQ} \Lambda_M$.
\end{lemma}
\eop

Using these bases, if $s \in \{0,2\}$, $F_{1,s}$ reads
$$\begin{array}{llll}F_{1,s} \colon &\bigoplus_{i=1}^{2g}\frac{\Lambda}{(t-1)}[d_i] \otimes  [\sigma_s] &\rightarrow &\bigoplus_{i=1}^{2g}\frac{\Lambda}{(t-1)}[c_i]  \otimes [\sigma_s]\\
   & [d_i] \otimes  [\sigma_s] &\mapsto & \delta_i[c_i] \otimes  [\sigma_s]
\end{array}$$
where $\delta_i=1$ if $i>k$. Therefore we have the following lemma.

\begin{lemma}
 If $s \in \{0,2\}$, then $F_{1,s}$ is an isomorphism, $\mbox{Ker}(F_{1,s}) =0$ and $\mbox{Coker}(F_{1,s}) =0$.
\end{lemma}
\eop

\begin{lemma}
$F_{1,1}$ reads $\bigoplus_{i=1}^{2g}\bigoplus_{j=1}^{k} f_{i,j}$ where
$$\begin{array}{llll}
f_{i,j} \colon&\frac{\Lambda}{(\delta_j)}[d_i \times c_j]&\rightarrow& 
\frac{\Lambda}{(\delta_j)}[c_i \times c_j]\\
&[d_i \times c_j]&\mapsto&\delta_i[c_i \times c_j]
\end{array}$$
so that $\mbox{Ker}(F_{1,1})=\bigoplus_{i=1}^{k}\bigoplus_{j=1}^{k}\mbox{Ker}(f_{i,j})$, $\mbox{Coker}(F_{1,1})=\bigoplus_{i=1}^{k}\bigoplus_{j=1}^{k}\mbox{Coker}(f_{i,j})$,
$\mbox{Coker}(f_{i,j})=\frac{\Lambda}{(\delta_{\mbox{\tiny \rm min}(i,j)})}[c_i \times c_j]$ and there is a section 
 $$\begin{array}{llll}u_{i,j} \colon & \mbox{Ker}(f_{i,j})=\frac{\Lambda}{(\delta_{\mbox{\tiny \rm min}(i,j)})}\frac{\delta_j}{\delta_{\mbox{\tiny \rm min}(i,j)}}[d_i \times c_j]
&\rightarrow &H_{3}(\TCMD)\\
&\frac{\delta_j}{\delta_{\mbox{\tiny \rm min}(i,j)}}[d_i \times c_j]&\mapsto& C(\Sigma_i \times \Sigma_j).\end{array}$$

\end{lemma}
\bp
Assume that the $c_i$ are outside $p_M^{-1}(S\times [-1,1])$.
Since $S \times c_i$ is rationally null-homologous (because $(t-1)S \times c_i$ and $\delta_iS \times c_i$ bound), the homology class of $C(\Sigma_i \times \Sigma_j)$ is independent of the chains $\Sigma_i$
and $\Sigma_j$ that have the given boundaries.
Since $\delta_i c_i$ cobounds with $(\theta_M(d_i^+) -d_i^-)$ in $p_M^{-1}(M\setminus S)$, there is a $2$-chain $\pm\Sigma_i$ in $\tilde{M}$ with boundary $\pm \delta_i c_i$
that intersects $p_M^{-1}(S\times [-1,1])$ as $(d_i \times [-1,1])$.

Then $C(\Sigma_i \times \Sigma_j)$ intersects $(S\times [-1,1]) \times \tilde{M}$ as $\pm\frac{\delta_j}{\delta_{\mbox{\tiny \rm min}(i,j)}}[(d_i \times [-1,1]) \times c_j]$ and its class in $H_3(\TCMD,(M\setminus S)\times \tilde{M})$ is therefore the image of this element
under the excision isomorphism. Therefore the class $[C(\Sigma_i \times \Sigma_j)]$ in $H_3(\TCMD)$ is mapped to the given generator of
$\mbox{Ker}(f_{i,j})$. Since 
$\delta_{\mbox{\tiny \rm min}(i,j)} C(\Sigma_i\times \Sigma_j)=\partial(\Sigma_i\times \Sigma_j)$, the wanted section exists.
\eop

Now, Proposition~\ref{proptilMtwocomp} is proved.
\eop

\subsection{On the homology of $\TCM$}
 \label{subhomTCM}

Recall from Subsection~\ref{subtcm} that $\TCM$ is obtained from $\TCMD$ by blowing-up the lifts of the diagonal, and that $\TCM$ is diffeomorphic to the complement of an open tubular neighborhood of  the lifts of the diagonal.
Again 
we consider the homology of $\TCM$ with coefficients in $\QQ$
endowed with the structure of $\QQ[t,t^{-1}]$-module where the multiplication by $t$ is induced by the action of $\theta$ on $\TCMD$ or $\TCM$.
Recall
$\Lambda=\QQ[t,t^{-1}]$ and
let $\QQ(t)$ be the field of fractions of $\Lambda$.

Then $H_{\ast}(\TCM)$ is a graded $\Lambda$-module, set 
$$H_{\ast}(C_2(M);\QQ(t)) = H_{\ast}(\TCM;\QQ) \otimes_\Lambda \QQ(t).$$

Also recall from Subsection~\ref{subtcm} that $\partial \TCM =\ZZ \times ST(M)$ where $ST(M)$ is the unit tangent bundle of $M$ that
 is diffeomorphic to $M\times S^2$. When $N$ is a submanifold of $M$, $ST(N)$ will denote the restriction of $ST(M)$ to $N$
viewed either as a submanifold of $\partial C_2(M)$ or as a submanifold of $\partial \TCM$ that sits inside the preimage of the preferred lift of the diagonal under the blow-up map.

\begin{proposition}
\label{prophomTCM}
$$H_{i}(C_2(M);\QQ(t))=H_{i-2}(M;\QQ) \otimes_{\QQ} \QQ(t)$$
for any $i \in \ZZ$.
$$H_{2}(C_2(M);\QQ(t))=\QQ(t)[ST(\ast) (\cong \ast \times S^2)]$$
$$H_{3}(C_2(M);\QQ(t))=\QQ(t)[ST(K) (\cong K \times S^2)]$$
$$H_{4}(C_2(M);\QQ(t))=\QQ(t)[ST(S) (\cong S \times S^2)]$$
$$H_{5}(C_2(M);\QQ(t))=\QQ(t)[ST(M) (\cong M \times S^2)].$$
\end{proposition}
\bp
We know that $$H_{\ast}(M^2;\QQ(t)) = H_{\ast}(\TCMD) \otimes_\Lambda \QQ(t)=0,$$
according to Proposition~\ref{proptilMtwotor},
and we compute 
$H_{\ast}(C_2(M);\QQ(t)) = H_{\ast}(\TCM) \otimes_\Lambda \QQ(t).$
In this proof, think of $C_2(M)$ as $(M^2\setminus \mbox{diag})$ that has the same homotopy type.
$$H_{i}(\TCMD,\TCM;\QQ)= H_{i}(M \times B^3,M \times S^2;\QQ)\otimes_{\QQ} \Lambda
=H_{i-3}(M;\QQ)\otimes_{\QQ} \Lambda$$
This module is always free and therefore the natural maps from $H_i(\TCMD;\QQ)$
to  $H_{i}(\TCMD,\TCM;\QQ)$ vanish.
Thus, we have the exact sequence of $\Lambda$-modules
$$0 \rightarrow H_{i+1}(\TCMD,\TCM;\QQ) \rightarrow H_{i}(\TCM;\QQ)\rightarrow H_{i}(\TCMD;\QQ) \rightarrow 0.$$
that gives rise to an isomorphism 
$$H_{i+1}(M^2,C_2(M);\QQ(t)) \cong H_{i}(C_2(M);\QQ(t)).$$
that proves the result.
\eop

\newpage
\section{On the equivariant linking number}
\setcounter{equation}{0}
\label{secblanch}

\subsection{On the equivariant intersection}
\label{subinteq}
In this text, unless otherwise mentioned, all the diffeomorphisms preserve the orientation, and the order of appearance of coordinates induces the orientation. Recall that boundaries are oriented with the {\em outward normal first\/} convention.

The fiber $N_u(A)$ of the normal bundle $N(A)$ of an oriented submanifold $A$ in an oriented manifold $C$ at $u \in A$ is oriented so that $T_uC=N_u(A)\oplus T_uA$ as oriented vector spaces.
For two oriented transverse submanifolds $A$ and $B$ of $C$, $A \cap B$ is oriented so that $N_u(A \cap B)=N_u(A)\oplus N_u(B)$. 
In particular, when the sum of the dimensions of $A$ and $B$ is the dimension of $C$, the sign of an intersection point $u$ of $A$ and $B$ is positive if and only if
$T_uC=N_u(A)\oplus N_u(B)$, that is if and only if $T_uC=T_uA\oplus T_uB$ as oriented vector spaces.

When $C$ is equipped with a free action of an abelian group $G$, the equivariant intersection has its coefficients in the group ring
$$\ZZ[G]=\oplus_{g \in G} \ZZ\exp(g)$$
and reads
$$A \cap_e B=\cup_{g \in G} \exp(g)A \cap g.B,$$
when $A$ is transverse to $g.B$ for any $g$,
so that $$g_1.A \cap_e g_2.B=\exp(g_1-g_2) g_1.(A \cap_e B).$$
When the sum of the dimensions of $A$ and $B$ is the dimension of $C$,
the {\em equivariant intersection number\/} $\langle A,B \rangle_{e,C}$ is the sum of the coefficients of the points.
Using the linear involution $$\exp(g) \mapsto \overline{\exp(g)}=\exp(-g),$$
$$\langle A,B \rangle_{e,C}=(-1)^{\mbox{\tiny dim}(A)\mbox{\tiny dim}(B)} \overline{\langle B,A \rangle_{e,C}}. $$

\subsection{On the linking number}

\begin{proposition}
\label{proplk}
Let $C$ be a compact connected manifold of dimension $c$.
Let $A$ and $B$ be two transverse submanifolds of $C$ with boundaries, and with respective dimensions $\mbox{dim}(A)=a$ and $\mbox{dim}(B)=b$, such that $a+b=c+1$, and $\partial A \cap \partial B =\emptyset$.
Then the {\em linking number\/} of $\partial A$ and $\partial B$ can be defined in the following equivalent ways:
$$lk(\partial A,\partial B)=\langle \partial A,B \rangle_C=(-1)^{a}\langle A ,\partial B \rangle_C=(-1)^{ab}lk(\partial B,\partial A).$$
Let $\ast \in C$, and let $S(\ast)$ denote the oriented boundary of a ball that contains $\ast$. The homology class of $\ast \times S(\ast)$ in $H_{c-1}(C^2 \setminus \mbox{diag}(C^2))$ is denoted by $[S^{c-1}(C)]$.
Then $\partial A \times \partial B$ is homologous to
$lk(\partial A,\partial B)[S^{c-1}(C)]$ in $H_{c-1}(C^2 \setminus \mbox{diag}(C^2))$.
\end{proposition}
Note
that $lk(\partial A,\partial B)=(-1)^alk(\partial B,\partial A)$ when $c$ is odd. Also note that the above formulae show that $\langle \partial A,B\rangle_C$ is independent of $B$ and that $\langle A,\partial B \rangle_C$ is independent of $A$.

\noindent{\sc Proof of Proposition~\ref{proplk}:}
Since $N(\partial A)$ is oriented as $N(A) \oplus N_{ext, A}(\partial A)$,
$$\partial (A\cap B)= (-1)^{c-b}\partial A \cap B + A \cap \partial B.$$
Therefore, 
$$\langle \partial A,B \rangle_C =(-1)^{c-b+1}\langle A,\partial B \rangle_C=(-1)^{a}\langle A ,\partial B \rangle_C.$$

Let us now prove that with this definition, $\partial A \times \partial B$ is homologous to $lk(\partial A,\partial B)[S^{c-1}(C)]$ in $H_{c-1}(C^2 \setminus \mbox{diag}(C^2))$. The manifold $A$ induces a cobordism between $\partial A$ and boundaries
 of neighborhoods in $A$ of points of $A \cap \partial B$ in $C \setminus \partial B$. This cobordism allows us to reduce the proof to the case where $A$ is a disk of dimension $a$ that intersects $\partial B$ once. Using a similar cobordism induced by $B$ allows us to reduce the proof to the case when $B$ is also a small disk of dimension $b$, and $A \cap B$ is an interval. We compute the homology class of $\partial A \times \partial B$ in this case of generalized Hopf links where
$C=\RR^c$, $A=[-1,1]^{a-1}\times [-2,0]\times(0)^{c-a}$, $B=(0)^{a-1} \times[-1,1]^{b}$.
$$\partial A= \left( \partial [-1,1]^{a-1} \times [-2,0] \times (0)^{c-a} \right)\bigcup \left( (-1)^{a-1}[-1,1]^{a-1}\times(\partial [-2,0]=0-(-2)) \times (0)^{c-a}\right)$$
$$\langle \partial A,B \rangle_C=(-1)^{a-1}.$$
Split $\partial A$ into two topological disks
 $$D^{a-1}_1= (-1)^{a-1}[-1,1]^{a-1}\times (0)^{b}\;\;\mbox{and}\;\;
D^{a-1}_2=\partial A \setminus D^{a-1}_1 ,$$
so that $D^{a-1}_2$ does not intersect $B$ and 
$\partial (D^{a-1}_2 \times B) =\partial (D^{a-1}_2) \times B + (-1)^{a-1} (D^{a-1}_2 \times \partial B)$.
 
Then $\partial A \times \partial B$ is homologous to 
$$D^{a-1}_1 \times \partial B  +(-1)^{a} \partial (D^{a-1}_2) \times B=D^{a-1}_1 \times \partial B  +(-1)^{a-1} \partial (D^{a-1}_1) \times B.$$
Now, change 
$\partial (D^{a-1}_1) \times B = \{((x,0),(0,\beta \in B)) \in C^2 ;x \in \partial (D^{a-1}_1), \beta \in B\}$ to
$$E=\{((x,0),(-x d(\beta,\partial B),\beta)\in C^2 \setminus \mbox{diag}(C^2);x \in \partial (D^{a-1}_1), \beta \in B \}$$ by the obvious homotopy where $d(\beta,\partial B)$ denotes the distance between $\beta$ and $\partial B$ (for an arbitrary continuous distance of $B$). Now, the first factor in
$$D^{a-1}_1 \times \partial B +(-1)^{a-1} E$$
can be contracted to $0$ in $D^{a-1}_1$ without meeting the diagonal.
Therefore $\partial A \times \partial B$ is homologous to
$$(-1)^{a-1}\{0\} \times\{(-x d(\beta,\partial B),\beta);(x,\beta) \in \partial (D^{a-1}_1 \times B)\}$$
that is homologous to $\{0\} \times \partial (D^{a-1}_1 \times B)$, that is in turn homologous to $(-1)^{a-1}[S^{c-1}(C)]$.
\eop

\begin{exa}
In $\RR^3$, let $B$ be the unit ball, let $x\in \mbox{Int}(B)$ and let $y \notin B$, 
then $$lk(x-y,\partial B)=1=-lk(\partial B,x-y).$$
\end{exa}

\subsection{On the considered equivariant linking number}
\label{subeqlk}

\begin{proposition}
\label{propdeflkeq}
Let $\alpha$ and $\beta$ be two submanifolds of $\tilde{M}$ of respective dimensions $\mbox{dim}(\alpha)$ and $\mbox{dim}(\beta)$, such that 
$\mbox{dim}(\alpha) + \mbox{dim}(\beta)=2$, whose projections in $M$ do not intersect.
There exist rational chains $A$ and $B$ such that
$\Delta(\theta_M)(\theta_M-1)\alpha=\partial A$ and $\Delta(\theta_M)(\theta_M-1)\beta=\partial B$.

The two following definitions for
the equivariant linking number $lk_e(\alpha,\beta)$ of $\alpha$ and $\beta$ are equivalent.
\begin{enumerate}
 \item $$lk_e(\alpha,\beta)(t_M)=\frac{\langle \alpha,B \rangle_{e,\tilde{M}}}{\Delta(t_M^{-1})(t_M^{-1}-1)}=(-1)^{\mbox{\tiny dim}(\alpha)+1}\frac{\langle A,\beta \rangle_{e,\tilde{M}}}{\Delta(t_M)(t_M-1)}$$
where $\exp(\theta_M^n)$ is denoted by $t_M^n$ so that $\langle \alpha,B \rangle_{e,\tilde{M}}=\sum_{n \in \ZZ}\langle \alpha,\theta_M^n(B) \rangle_{\tilde{M}}t_M^n$.
\item The class of $\alpha \times \beta$ in $H_2(\TCM;\QQ) \otimes_\Lambda \QQ(t)$ is equal to $lk_e(\alpha,\beta)(t)[ST(\ast)]$.
\end{enumerate}

The equivariant linking number has the following properties
$$lk_e(\alpha,\beta)=(-1)^{\mbox{\tiny dim}(\alpha)+1}\overline{lk_e(\beta,\alpha)}$$
and, when $\lambda \in \QQ(t_M)$,
$$lk_e(\lambda \alpha,\beta)=\lambda lk_e(\alpha,\beta) = lk_e(\alpha,\overline{\lambda}\beta).$$
\end{proposition}
\bp Applying Proposition~\ref{proplk} to $A$ and $\theta_M^n(B)$, for any $n$, shows that 
$$\langle \partial A ,B \rangle_{e,\tilde{M}}=(-1)^{\mbox{\tiny dim}(\alpha)+1}\langle A , \partial B \rangle_{e,\tilde{M}}$$
and allows us to define $lk_e(\partial A, \partial B)$ consistently by this expression so that for any $\lambda \in \Lambda_M$
$$lk_e(\lambda \partial A,\partial B)=\lambda lk_e(\partial A,\partial B) = lk_e(\partial A,\overline{\lambda}\partial B).$$
Extend the definition of $lk_e$ so that this proposition holds for $\lambda \in \QQ(t_M)$.
The proof of Proposition~\ref{proplk} easily adapts to show that the class of $\alpha \times \beta$ in $H_2(C_2(M);\QQ(t))$ is equal to $lk_e(\alpha,\beta)(t)[ST(\ast)]$. According to Proposition~\ref{prophomTCM}, this yields a definition of $lk_e(\alpha,\beta)$.
\eop

\begin{example}
\label{exeqlk}
The lifts $\hat{S}^+$, $\hat{S}^-$, $\hat{\ast}^+$ and $\hat{\ast}^+$ of submanifolds of $S^+$, $S^-$, $\ast^+$ and $\ast^+$ are in $(M\setminus S)_0$.
The surface $\hat{S}$ is between $\hat{S}^-$ and $\theta_M(\hat{S}^+)$. Following the preimage of $K$ under the covering map according to the orientation of $K$, we successively meet $\theta_M^{-1}(\hat{\ast})$, $\hat{\ast}^+$, $\hat{\ast}^{-}$, $\hat{\ast}$, $\theta_M(\hat{\ast}^+)$ as in the following figure.
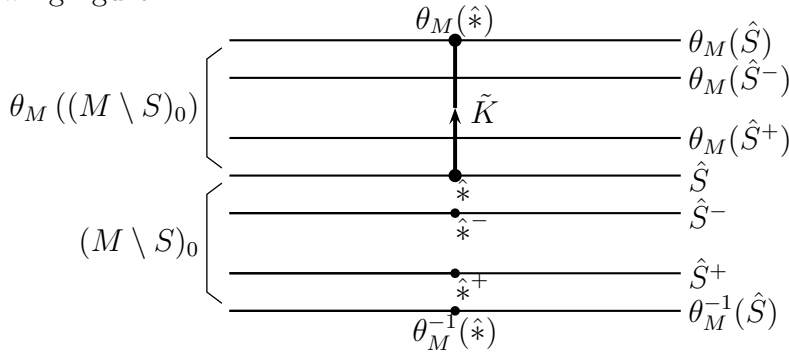
\begin{center}
\begin{pspicture}[shift=-0.1](0,-.2)(12,3.8)
\psline(3,0)(9,0)
\rput[l](9.1,0){$\theta_M^{-1}(\hat{S})$}
\psline(3,.5)(9,.5)
\rput[l](9.1,.5){$\hat{S}^+$}
\psline(3,1.3)(9,1.3)
\rput[l](9.1,1.3){$\hat{S}^-$}
\psline(3,1.8)(9,1.8)
\rput[l](9.1,1.8){$\hat{S}$}
\psline(3,2.3)(9,2.3)
\rput[l](9.1,2.3){$\theta_M(\hat{S}^+)$}
\psline(3,3.1)(9,3.1)
\rput[l](9.1,3.1){$\theta_M(\hat{S}^-)$}
\psline(3,3.6)(9,3.6)
\rput[l](9.1,3.6){$\theta_M(\hat{S})$}
\psline[linewidth=1.5pt]{*->}(6,1.8)(6,2.7)
\rput[lt](6,1.75){$\hat{\ast}$}
\psline[linewidth=1.5pt]{-*}(6,2.7)(6,3.6)
\rput[l](6.2,2.7){$\tilde{K}$}
\rput[b](6,3.65){$\theta_M(\hat{\ast})$}
\psline{-*}(3,.5)(6,.5)
\rput[lt](6,.5){$\hat{\ast}^+$}
\psline{-*}(3,0)(6,0)
\rput[t](6,-.05){$\theta_M^{-1}(\hat{\ast})$}
\psline{-*}(3,1.3)(6,1.3)
\rput[lt](6,1.3){$\hat{\ast}^-$}
\psline[linewidth=.5pt](2.9,.1)(2.7,.25)(2.7,1.55)(2.9,1.7)
\rput[r](2.6,.9){$(M\setminus S)_0$}
\psline[linewidth=.5pt](2.9,1.9)(2.7,2.05)(2.7,3.35)(2.9,3.5)
\rput[r](2.6,2.7){$\theta_M\left((M\setminus S)_0\right)$}
\end{pspicture}
\end{center}

In $\tilde{M}$, there is a lift $\tilde{K}$ of $K$ whose boundary is 
$\partial \tilde{K} =(\theta_M-1)\hat{\ast}$. The boundary of the closure of some lift of $(M \setminus S)$ is $(\theta_M-1)\hat{S}$. Then
$$lk_e(\hat{\ast},t_M\hat{S}^+)=\frac{1}{1-t_M}$$
and $$lk_e(\hat{S},t_M\hat{\ast}^+)=\frac{1}{1-t_M}.$$

\end{example}

The following lemma follows from Blanchfield duality \cite{Blanch}.

\begin{lemma}
 \label{lemblanchnondeg}
The equivariant linking number is non-degenerate in the following sense.
For any class $\beta$ of $H_1(\tilde{M})$ of order $\delta(\beta)\in \Lambda_M$, there exists a two-component link $(J,J^{\ast})$ of $\tilde{M}$ such that $J$ represents $\beta$, $p_M(J) \cap p_M(J^{\ast}) = \emptyset$ and
$$lk_e(J,J^{\ast})=\frac{q}{\delta(\beta)}$$
for some nonzero rational number $q$.
\end{lemma}
\bp
The polynomial $\delta(\beta)$ is defined up to units $q t_M^{\pm n}$ 
of $\Lambda_M$. Realize $\beta$ by a link $L$ of $p_M^{-1}\left(M \setminus (S\times ]-1,1[)\right)$.
By assumption, $\delta(\beta)L$ is rationally homologous to 
$$\sum_{i=1}^{2g}\alpha_i \left(z_i\times \partial [-1,1] \subset p_M^{-1}(S\times [-1,1])\right)$$
in $p_M^{-1}\left(M \setminus (S\times ]-1,1[)\right)$, for some polynomials $\alpha_i$ of $\Lambda_M$, where $\delta(\beta)$ is coprime with the greatest common divisor $d$ of the $\alpha_i$.
Then the B\'ezout Identity allows us to find a link $L^{\ast}$ of 
$p_M^{-1}(S\times ]-1,1[)$ such that 
$$lk_e(L,L^{\ast})=\frac{q_1 d}{\delta(\beta)}$$ for some $q_1 \in \QQ\setminus\{0\}.$
The B\'ezout Identity also implies the existences of $u$ and $v$ in $\ZZ[t_M,t_M^{-1}]$ such that
$$u\delta(\beta) + vd=q_2$$ for some $q_2 \in \QQ\setminus\{0\}.$
Then there exists a link $L^{\prime\ast}$ homologous to 
$\overline{v}L^{\ast}$ such that
$$lk_e(L,L^{\prime\ast})=\frac{q_1 q_2}{\delta(\beta)}-uq_1.$$
Adding to $L^{\prime\ast}$ meridians of $\theta_M^n(L)$
transforms $L^{\prime\ast}$ into $L^{\prime\prime\ast}$ such that $lk_e(L,L^{\prime\prime\ast})=\frac{q_1 q_2}{\delta(\beta)}$.
Finally make $L^{\prime\prime\ast}$ connected by connecting components
by bands that don't intersect $\cup_{n\in \ZZ}\theta_M^n(L)$ and make $L$ connected in a similar way.
\eop

\newpage
\section{Comparing homology classes in $\TCM$}
\setcounter{equation}{0}
\label{seccomparhom}

\subsection{On the logarithmic derivative of the Alexander polynomial}
\label{subsecderal}

Recall that $(z_i)_{i=1, \dots 2g}$ and $(z^{\ast}_i)_{i=1, \dots, 2g}$ are two dual bases of $H_1(S;\ZZ)$ such that 
$\langle z_i, z^{\ast}_j\rangle=\delta_{ij}$.
We shall need the following proposition.

\begin{proposition}
\label{proplogdereq}
Let $\tilde{z}^{\ast}_i$, $\tilde{z}^+_i$ and $\tilde{z}^-_i$ denote the lifts of $z^{\ast}_i$, ${z}^+_i$
and ${z}^-_i$ in the lift $(S\times[-1,1])_0$ of $S\times[-1,1]$ in $\tilde{M}$, respectively.
Then 
$$\frac{t \Delta^{\prime}(t)}{\Delta(t)}=\sum_{i=1}^{2g}lk_e(\frac{\tilde{z}^+_i+\tilde{z}^-_i}{2},\tilde{z}^{\ast}_i)=\sum_{i=1}^{2g}lk_e(\tilde{z}^+_i,\tilde{z}^{\ast}_i)-g=\sum_{i=1}^{2g}lk_e(\tilde{z}^-_i,\tilde{z}^{\ast}_i)+g.$$
\end{proposition}

\begin{remark}
\label{remlog}
The above proposition implies that if $(a_i,b_i)_{i \in \{1,\dots,g\}}$ is a symplectic basis of $H_1(S;\ZZ)$, ($\langle a_i,a_j \rangle=\langle b_i,b_j \rangle=0$, $\langle a_i,b_j \rangle=\delta_{ij}$), then
$$\frac{t \Delta^{\prime}(t)}{\Delta(t)}=\sum_{i=1}^g\left(lk_e(\tilde{a}_i,\tilde{b}_i^+) -\overline{lk_e(\tilde{a}_i,\tilde{b}_i^+)}\right).$$
Indeed, according to the proposition,\\
$\begin{array}{ll}\frac{t \Delta^{\prime}(t)}{\Delta(t)}&=\sum_{i=1}^{g} lk_e(\tilde{a}_i,\tilde{b}_i^+)-\sum_{i=1}^{g} lk_e(\tilde{b}_i, \tilde{a}^{+}_i) +g\\
&=\sum_{i=1}^{g} lk_e(\tilde{a}_i,\tilde{b}_i^+)-\sum_{i=1}^{g} lk_e(\tilde{b}_i, \tilde{a}^{+}_i) + \sum_{i=1}^{g}lk_e(\tilde{b}_i, \tilde{a}^{+}_i - \tilde{a}^{-}_i)\\
&=\sum_{i=1}^{g} lk_e(\tilde{a}_i,\tilde{b}_i^+) -\sum_{i=1}^{g}lk_e(\tilde{b}^{+}_i, \tilde{a}_i).
\end{array}$ 
\end{remark}

To prove Proposition~\ref{proplogdereq}, we shall prove the following proposition that gives another topological expression of the logarithmic derivative of the Alexander polynomial and the following lemma that may be of independent interest.
\begin{proposition}
\label{proplogderlk}
Let $\hat{z}^+_i$ and $\hat{z}^-_i$ denote the lifts of ${z}^+_i$
and ${z}^-_i$ in the lift $(M \setminus S)_0$ of $M \setminus S$ in $\tilde{M}$, so that $\hat{z}^-_i=\tilde{z}^-_i$ and $\tilde{z}^+_i=\theta_M(\hat{z}^+_i)$ with the notation of the previous statement.
Then $(\hat{z}^+_i-\hat{z}^-_i)_{i=1, \dots 2g}$ and $t_M^{-1/2}(\tilde{z}^+_i-\tilde{z}^-_i)_{i=1, \dots 2g}$ are two bases of 
$H_1(p_M^{-1}(M\setminus S);\QQ)\otimes_{\Lambda_M} \QQ(t_M)[t_M^{1/2}]$ and
$$t_M^{-1/2}(\tilde{z}^+_j-\tilde{z}^-_j)=\sum_{i=1}^{2g}a_{ij}(\hat{z}_i^+-\hat{z}_i^-)$$
where $$a_{ij}=t_M^{1/2} lk(z^{+}_j,z^{\ast}_i) -t_M^{-1/2}lk(z^{-}_j,z^{\ast}_i).$$
Let $A$ denote the matrix
$A=[a_{ij}]_{(i,j) \in \{1,\dots,2g\}^2}$. Then
$\det(A)=\Delta(M)$, and the trace $\mbox{tr}(A^{-1})$ of $A^{-1}$ is related to $\frac{\Delta^{\prime}}{\Delta}$ as follows:
$$\frac{t_M \Delta^{\prime}(t_M)}{\Delta(t_M)}= g \frac{t_M+1}{t_M-1} 
-\frac{\mbox{tr}(A^{-1})}{t^{1/2}_M-t^{-1/2}_M} .$$
\end{proposition}

\begin{lemma}
\label{leminutile} Set $A^{-1}=[b_{ij}]_{(i,j) \in \{1,\dots,2g\}^2}$,
then
 $$lk_e(\tilde{z}^+_j,\tilde{z}^{\ast}_i)=\frac{\delta_{ij}}{(1-t_M^{-1})} - \frac{b_{ij}}{(t_M^{1/2}-t_M^{-1/2})}, \;\; lk_e(\tilde{z}^-_j,\tilde{z}^{\ast}_i)=\frac{\delta_{ij}}{(t_M-1)} - \frac{b_{ij}}{(t_M^{1/2}-t_M^{-1/2})}\;\;\mbox{and}$$
$$lk_e(\frac{\tilde{z}^+_j+\tilde{z}^-_j}{2},\tilde{z}^{\ast}_i)=\frac{\delta_{ij}(t_M+1)}{2(t_M-1)} - \frac{b_{ij}}{(t_M^{1/2}-t_M^{-1/2})}.$$
\end{lemma}

\noindent{\sc Proof of Proposition~\ref{proplogderlk}:}
In $H_1(M\setminus S)$, $z^{\varepsilon}_j=\sum_{i=1}^{2g}lk(z_j^{\varepsilon},z^{\ast}_i)(z_i^+-z_i^-).$
Thus, $$(\theta_M(\hat{z}_j^+)-\hat{z}_j^-)
=\sum_{i=1}^{2g}\left(t_M lk(z_j^{+},z^{\ast}_i)-lk(z_j^{-},z^{\ast}_i)\right)
(\hat{z}_i^+-\hat{z}_i^-)$$
and 
$t_M^{-1/2}(\tilde{z}^+_j-\tilde{z}^-_j)=\sum_{i=1}^{2g}a_{ij}(\hat{z}_i^+-\hat{z}_i^-).$
Therefore, $A$ is a transition matrix from one basis of the statement to the other one. Since it coincides with the matrix of the statement of Lemma~\ref{lemhomtilM}, up to a linear transformation with determinant $1$ on its rows, $\det(A)=\Delta(M)$.

Set $z=t_M^{1/2}-t_M^{-1/2}$.
$$\begin{array}{ll}a_{ij}(t_M)&=t_M^{1/2} lk(z^{+}_j,z^{\ast}_i) -t_M^{-1/2}lk(z^{-}_j,z^{\ast}_i)\\
&=zlk(z^{+}_j,z^{\ast}_i)+t_M^{-1/2}lk(z^{+}_j-z^{-}_j,z^{\ast}_i)\\
&=zlk(z^{+}_j,z^{\ast}_i)+t_M^{-1/2}\delta_{ij}.\end{array}$$

$$\begin{array}{ll}t_M a^{\prime}_{ij}(t_M)&=\frac{1}{2}(t_M^{1/2} +t_M^{-1/2})lk(z^{+}_j,z^{\ast}_i)
- \frac{1}{2}t_M^{-1/2}\delta_{ij}\\
&=\frac{t_M^{1/2} +t_M^{-1/2}}{2z}a_{ij}(t_M) 
- \frac{t_M^{1/2} +t_M^{-1/2}}{2z}t_M^{-1/2}\delta_{ij}
- \frac{1}{2}t_M^{-1/2}\delta_{ij}\\
&=\frac{t_M^{1/2} +t_M^{-1/2}}{2z}a_{ij}(t_M) - \frac{ 1}{z}\delta_{ij}.\end{array}$$
Let $\Delta_{ii}$ denote the cofactor of $(i,i)$ in the matrix $A$.
$$\begin{array}{ll}t_M\Delta^{\prime}(M)=2g\frac{t_M^{1/2} +t_M^{-1/2}}{2z} \Delta(M)
-\frac{1}{z}\sum_{i=1}^{2g}\Delta_{ii}.\end{array}$$
Recall $A^{-1}=[b_{ij}]_{(i,j) \in \{1,\dots,2g\}^2}$,
$b_{ii}=\frac{\Delta_{ii}}{\Delta(M)}$, then
$$t_M\frac{\Delta^{\prime}(M)}{\Delta(M)}=g\frac{t_M^{1/2} +t_M^{-1/2}}{z} -\frac{1}{z}\mbox{tr}(A^{-1}).$$
Proposition~\ref{proplogderlk} is proved.
\eop

\noindent{\sc Proof of Lemma~\ref{leminutile}:}
We compute $lk_e(\tilde{z}^+_j,\tilde{z}^{\ast}_i)$, using that $lk_e(\tilde{z}^+_j-\tilde{z}^-_j,\tilde{z}^{\ast}_i)=\delta_{ij}$.
$$\tilde{z}^+_j=t_M\hat{z}^+_j= \frac{\tilde{z}^+_j-\tilde{z}^-_j - (\hat{z}^+_j-\hat{z}^-_j)}{1-t_M^{-1}}.$$
According to Proposition~\ref{proplogderlk},
$$\hat{z}^+_j-\hat{z}^-_j
=\sum_{i=1}^{2g}b_{ij}t_M^{-1/2}(\tilde{z}^+_i-\tilde{z}^-_i).$$
Thus, $(1-t_M^{-1})lk_e(\tilde{z}^+_j,\tilde{z}^{\ast}_i)=\delta_{ij} - b_{ij}t_M^{-1/2}$. The other expressions follow easily.
\eop

\noindent{\sc Proof of Proposition~\ref{proplogdereq}:}
Since $lk_e(\tilde{z}^+_i-\tilde{z}^-_i,\tilde{z}^{\ast}_i)=1$,
$$\sum_{i=1}^{2g}lk_e(\tilde{z}^+_i,\tilde{z}^{\ast}_i)-g=\sum_{i=1}^{2g}lk_e(\tilde{z}^-_i,\tilde{z}^{\ast}_i)+g=\sum_{i=1}^{2g}lk_e(\frac{\tilde{z}^+_i+\tilde{z}^-_i}{2},\tilde{z}^{\ast}_i).$$
According to Lemma~\ref{leminutile},
$$\sum_{i=1}^{2g}lk_e(\frac{\tilde{z}^+_i+\tilde{z}^-_i}{2},\tilde{z}^{\ast}_i)=\frac{g(t_M+1)}{(t_M-1)} - \frac{\mbox{tr}(A^{-1})}{(t_M^{1/2}-t_M^{-1/2})}$$
and Proposition~\ref{proplogderlk} allows us to conclude.
\eop

\subsection{Homology class of $[s_{\tau}(S)]$}

\begin{proposition}
\label{propsplusstau}
Let $\tau: TM \rightarrow M \times \RR^3$ be a trivialisation of $TM$, let $S$ be a closed (oriented) surface of $M$, let $\cvarM \in S^2$ and let $s_{\tau}(S;\cvarM)=\tau^{-1}(S \times \cvarM) \subset ST(M)$ be a section induced by $\tau$.
Let $s_+(S)$ and $s_-(S)$ be the sections of $ST(M)_{|S}$ induced by the positive normal of $S$ and the negative normal of $S$, respectively.
$$[s_+(S)]=[s_{\tau}(S;\cvarM)] + \frac{\chi(S)}{2}[ST(\ast)]$$
$$[s_-(S)]=[s_{\tau}(S;\cvarM)] - \frac{\chi(S)}{2}[ST(\ast)]$$
\end{proposition}
\bp
First notice that the homology class of $s_{\tau}(S;\cvarM)$ in $ST(M)_{\tau}$ does not depend on $\cvarM$. It will be denoted by $[s_{\tau}(S)]$. Next, since the generator of $\pi_1(SO(3))$ can be realized by rotations around the $\cvarM$-axis and since $\pi_2(SO(3))=0$, the homology class of $s_{\tau}(S)$ in $ST(M)_{|S}$ does not depend on $\tau$ either.

Embed $S\times [-1,1]$ in $\RR^3$, then the tangent bundles of $S\times [-1,1]$ 
 in $\RR^3$ and in $M$ are isomorphic (they are both isomorphic to the direct sum of the tangent bundle of $S$ and the trivial normal bundle).
 
 Using the trivialisation $\tau$ of $ST(M)_{|S}$ induced by the standard trivialisation of $\RR^3$, the positive normal section of $ST(S)$ is a map from $S$ to $S^2$ that can be homotoped to a constant outside an open disk. Then  $[s_+(S)]-[s_{\tau}(S)] \in H_2(D^2 \times S^2)$ and $[s_+(S)]-[s_{\tau}(S)] =c [ST(\ast)]$, where $c$ is the degree of the Gauss map from $S$ to $S^2$ that maps a point to the direction of the positive normal of $S$.
 \eop

\begin{proposition}
\label{prophomdiagS}
Let $S$ be a closed (oriented) surface.
Let $S$ and $S^+$ be two copies of $S$, let $(z_i)_{i=1, \dots 2g}$ and $(z^{\ast}_i)_{i=1, \dots, 2g}$ be two dual bases
of $H_1(S;\ZZ)$ such that 
$\langle z_i, z^{\ast}_j\rangle=\delta_{ij}.$ Let $\ast \in S$.
Let $\mbox{diag}(S \times S^+)=\{(x,x^+); x \in S\}$.
We have the following equality in $H_2(S \times S^+)$
$$[\mbox{diag}(S \times S^+)]=[\ast \times S^+] + [S \times \ast^+]+ \sum_{i=1}^{2g} [z_i \times z^{\ast+}_i].$$
\end{proposition}
\bp
$$H_2(S \times S^+)=\ZZ[\ast \times S^+] \oplus \ZZ[S \times \ast^+] \oplus \bigoplus_{(i,j) \in \{1,2,\dots,2g\}^2}\ZZ[z_i \times z^{\ast +}_j].$$
The dual basis of the above basis with respect to the intersection form is $$\left([S \times \ast^+],[\ast \times S^+], ([z^{\ast}_i \times z^{+}_j])_{(i,j) \in \{1,2,\dots,2g\}^2}\right).$$
To get the coordinates of $[\mbox{diag}(S \times S^+)]$ in the first decomposition
we compute the intersection numbers with the second one.
$\langle[\mbox{diag}(S \times S^+)],[z^{\ast}_i \times z^{+}_i]\rangle=\pm 1$
where the tangent space to $\mbox{diag}(S \times S^+)$ is naturally parametrized by
$(u_i, v^{\ast}_i, u_i, v^{\ast}_i)$ and the tangent space to $[z^{\ast}_i \times z^{+}_i]$ is naturally parametrized by $(0,w^{\ast}_i,x_i,0)$, so the intersection sign is the sign of the permutation $(u,v,w,x) \mapsto (u,w,x,v)$ which is $+1$.
\eop

\begin{theorem}
\label{thmstauS}
Let $$\ID=\frac{1+t}{1-t} + \frac{t\Delta^{\prime}(M)}{\Delta(M)}.$$
Let $S$ be a surface generating $H_2(M)$, then
$$[s_{\tau}(S)]=\ID [ST(\ast)].$$
\end{theorem}
\bp According to Propositions~\ref{propdeflkeq}, \ref{propsplusstau} and \ref{prophomdiagS}, and to Example~\ref{exeqlk},
$$[s_{\tau}(S)]= \left(\frac{2}{1-t} -\frac{\chi(S)}{2}+ \sum_{i=1}^{2g}lk_e(\tilde{z}_i,\tilde{z}^{\ast +}_i)\right) [ST(\ast)],$$
with the notation of Proposition~\ref{proplogdereq}, according to which,
$$\sum_{i=1}^{2g}lk_e(\tilde{z}^-_i,\tilde{z}^{\ast}_i)=\frac{t \Delta^{\prime}(t)}{\Delta(t)}-g.$$
\eop

\subsection{Homology classes of $H_{\ast}(C_2(M),\partial C_2(M);\QQ(t))$}

We use equivariant intersections with homology classes of $H_{(6-\ast)}(C_2(M),\partial C_2(M);\QQ(t))$ to evaluate homology classes of $H_{\ast}(C_2(M);\QQ(t))$.

The configuration space $C_2(M)$ is oriented like $M^2$, $\partial C_2(M)$ is oriented like $ST(M)$.
Note that when $A$ is in $H_{6-i}(C_2(M),\partial C_2(M);\QQ(t))$ and when $B$ is in $H_i(\partial C_2(M);\QQ(t))$,
$$\langle A,B \rangle_{e,\TCM}=\langle \partial A,B \rangle_{e,\partial \TCM}.$$

\begin{theorem}
\label{thmstauM}
Let $\tau: TM \rightarrow M \times \RR^3$ be a trivialisation of $TM$ and let $s_{\tau}(M)=\tau^{-1}(M \times \cvarM) \subset ST(M)$ be a section induced by $\tau$.
Then $$[s_{\tau}(M) ]=\ID [ST(K)].$$
In particular, there exists a $4$--dimensional chain $F$ with coefficients in $\QQ(t)$ whose boundary is
$$\partial F =s_{\tau}(M)-\ID ST(K).$$
Such a $4$--chain satisfies
$$\langle ST(\ast),F \rangle_e=1.$$
\end{theorem}
\bp
As a consequence of Theorem~\ref{thmstauS}, $([s_{\tau}(S)]-\ID [ST(\ast)])$ bounds a $3$-dimensional chain
$A_{\tau}(S)$ and 
$$\begin{array}{ll}\langle ST(K), A_{\tau}(S)\rangle_{e,\TCM}&=
-\overline{\langle A_{\tau}(S) , [ST(K) ]\rangle_{e,\TCM}}\\
&=-\overline{\langle [s_{\tau}(S)]-\ID [ST(\ast)] , [ST(K) ]\rangle_{e,\partial \TCM}}=-1.
\end{array}$$
On the other hand, since
$$\langle A_{\tau}(S) , [s_{\tau}(M) ]\rangle_{e,\TCM}
=\langle [s_{\tau}(S)]-\ID [ST(\ast)] , [s_{\tau}(M) ]\rangle_{e,\partial \TCM}=-\ID,$$
$$\langle[s_{\tau}(M) ], A_{\tau}(S) \rangle_{e,\TCM}=\overline{\ID}=-\ID.$$
This shows that $[s_{\tau}(M) ]=\ID [ST(K)]$ thanks to Proposition~\ref{prophomTCM}. The other assertions follow.
\eop

\begin{lemma}
\label{lemak}
See the knot $K$ as a map from $\RR/\ZZ$ to $M$.
Consider the continuous map
$$\begin{array}{llll}\check{A}(K) \colon &(S^1=[0,1]/(0\sim 1)) \times ]0,1[ &\rightarrow & C_2(M)\\
&(t,u \in]0,1[)& \mapsto &(K(t),K(t+u)),
\end{array}$$
and its lift ${A}(K) \colon S^1 \times ]0,1[ \rightarrow  \TCM$
such that the lift of $(K(t),K(t+\varepsilon))$ is in a small neighborhood of the canonical lift of the diagonal, for a small positive $\varepsilon$.
This map $A(K)$ extends to the closed annulus $S^1 \times [0,1]$. Its extension, still denoted by $A(K)$, maps $S^1 \times \{0\}$ to the section $s_{TK}(K)$ of 
$ST(K)$ given by the direction of the tangent vector to $K$ and $A(K)(S^1 \times \{1\})=\theta^{-1}(s_{(-TK)}(K))$. The image $A(K)$ of $A(K)$ is supported in $p^{-1}(K^2)$.
Then the homology class $[A(K)]$ of $A(K)$ in $H_{2}(C_2(M),\partial C_2(M);\QQ(t))$
satisfies:
$$\partial [A(K)] =[s_{TK}(K)]-t^{-1}[(s_{(-TK)}(K))]$$
and
$$\langle ST(S) , A(K)\rangle_{e,\TCM}=1-t.$$
\end{lemma}
\bp
$$\langle A(K),ST(S) \rangle_{e,\TCM}=\langle \partial A(K),ST(S) \rangle_{e,\partial \TCM}=1-t^{-1}.$$
\eop

\begin{proposition}
\label{propexistF}
Let $\tau: TM \rightarrow M \times \RR^3$ be a trivialisation of $TM$ that maps the vectors of $T^+K$ (tangent to $K$ and directed by $K$) to $K \times \RR^+\{\qvarM\}$, for some $\qvarM \in S^2$.
Let $K_{\fvarM}$ be a knot disjoint from $K$ that is rationally homologous to $K$.
Let $\fvarM \in (S^2 \setminus \{\qvarM,-\qvarM\})$ and
let $s_{\tau}(M;\fvarM)=\tau^{-1}(M \times \fvarM) \subset\partial ST(M)$.
Then there exists a $4$--dimensional chain $F_{\fvarM}$ with coefficients in $\QQ(t)$ whose boundary is
$$\partial F_{\fvarM} =s_{\tau}(M;\fvarM)-\ID ST(K_{\fvarM})$$
such that
$$\langle F_{\fvarM} , A(K)\rangle_{e,\TCM}=0.$$
Furthermore, if $F_{\fvarM}^{\prime}$ is another such, then the class of $(F_{\fvarM}-F_{\fvarM}^{\prime})$ vanishes in $H_4(C_2(M);\QQ(t))$.
\end{proposition}
\bp According to Theorem~\ref{thmstauM}, there exists a chain $F$ such that $\partial F =s_{\tau}(M;\fvarM)-\ID ST(K_{\fvarM})$. Now, $$F_{\fvarM}=F -\frac{\langle F ,A(K) \rangle_{e,\TCM}}{(1-t)}ST(S)$$ is a $4$-- chain with the announced properties.
The fact that $\langle F_{\fvarM}-F_{\fvarM}^{\prime} ,A(K) \rangle=0$ guarantees that the second assertion is true, thanks to Lemma~\ref{lemak} and to Proposition~\ref{prophomTCM}.
\eop 

It will be proved in Section~\ref{secspecK} that $\delta(M)(t-1)F_{\fvarM}$ can be furthermore assumed to be rational. See Proposition~\ref{propGrat}.

\newpage
\section{Definition of the invariant $\CQ$}
\setcounter{equation}{0}
\label{sectriple}

\subsection{On the equivariant triple intersection in $\TCM$}
\label{subdeftriple}

Let $C_{\fvarM}$, $C_{\svarM}$ and $C_{\tvarM}$ be three $4$-dimensional rational chains in $\TCM$ whose boundaries are in $\partial \TCM$, whose projections $p(C_{\fvarM})$, $p(C_{\svarM})$ and $p(C_{\tvarM})$ in $C_2(M)$ are transverse, and such that (therefore) $p(C_{\fvarM}) \cap p(C_{\svarM}) \cap p(C_{\tvarM}) \cap \partial C_2(M) =\emptyset$, define
$$\langle C_{\fvarM},C_{\svarM},C_{\tvarM} \rangle_e=\sum_{(i,j)\in \ZZ^2} \langle C_{\fvarM},\theta^{-i}(C_{\svarM}),\theta^{-j}(C_{\tvarM})\rangle_{\TCM} \svar^i\tvar^j \in \QQ[\svar^{\pm1},\tvar^{\pm1}]$$
This number does not change if $C_{\fvarM}$ (resp. $C_{\svarM}$ or $C_{\tvarM}$) is replaced by a chain $C^{\prime}_{\fvarM}$ (resp. $C^{\prime}_{\svarM}$ or $C^{\prime}_{\tvarM}$) with the same boundary such that $(C^{\prime}_{\fvarM}-C_{\fvarM})$ (resp. $(C^{\prime}_{\svarM}-C_{\svarM})$ or $(C^{\prime}_{\tvarM}-C_{\tvarM})$) is a rationally null-homologous cycle.

Indeed, in this case $(C^{\prime}_{\fvarM}-C_{\fvarM})$ is the boundary of a $5$-chain $A$ whose projection can be assumed to be transverse to $p(C_{\svarM})$ and $p(C_{\tvarM})$ and disjoint from $\partial C_2(M)$. Then $$\left(C^{\prime}_{\fvarM} \cap \theta^{-i}(C_{\svarM})\cap\theta^{-j}(C_{\tvarM})\right) -\left(C_{\fvarM} \cap \theta^{-i}(C_{\svarM})\cap\theta^{-j}(C_{\tvarM})\right) =\pm \partial (A \cap \theta^{-i}(C_{\svarM})\cap\theta^{-j}(C_{\tvarM})).$$
Set $$\fvar=\svar^{-1}\tvar^{-1}.$$

For any three one-variable polynomials $f$, $g$ and $h$ with rational coefficients,
$$\langle f(t)C_{\fvarM},g(t)C_{\svarM},h(t)C_{\tvarM} \rangle_e=f(\fvar)g(\svar)h(\tvar)\langle C_{\fvarM},C_{\svarM},C_{\tvarM} \rangle_e.$$

Therefore, the definition of the above {\em equivariant triple intersection\/} extends to
chains $Q_{\fvarM}$, $Q_{\svarM}$ and $Q_{\tvarM}$ of $C_2(M)$ with coefficients in $\QQ(t)$,
(whose boundaries are in $\partial \TCM$, and whose projections $p(Q_{\fvarM})$, $p(Q_{\svarM})$ and $p(Q_{\tvarM})$ in $C_2(M)$ are transverse) as follows.
There exist $f(t)$, $g(t)$ and $h(t)$
in $\Lambda=\QQ[t,t^{-1}]$ such that $f(t)Q_{\fvarM}$, $g(t)Q_{\svarM}$ and $h(t)Q_{\tvarM}$ are rational chains and
$$\langle Q_{\fvarM},Q_{\svarM},Q_{\tvarM} \rangle_e=\frac{1}{f(\fvar)g(\svar)h(\tvar)}\langle f(t)Q_{\fvarM},g(t)Q_{\svarM},h(t)Q_{\tvarM} \rangle_e.$$
This equivariant triple intersection takes values in the field of fractions $$\QQ(\svar,\tvar)=K(\QQ[\svar^{\pm1},\tvar^{\pm1}])=K\left(\frac{\QQ[\fvar^{\pm1},\svar^{\pm1},\tvar^{\pm1}]}{\fvar\svar\tvar=1}\right)$$
of $\QQ[\svar^{\pm1},\tvar^{\pm1}]$,
and it is such that for any $f$, $g$ and $h$ in $\QQ(t)$, 
$$\langle f(t)Q_{\fvarM},g(t)Q_{\svarM},h(t)Q_{\tvarM} \rangle_e=f(\fvar)g(\svar)h(\tvar)\langle Q_{\fvarM},Q_{\svarM},Q_{\tvarM} \rangle_e.$$

\begin{lemma}
\label{lemsym}
Set $P(\svar,\tvar)=\langle Q_{1},Q_{2},Q_{3} \rangle_e$.
Let $\sigma$ be a permutation of $\{1,2,3\}$ and let $\beta$ be the
natural bijection from $\{1,2,3\}$ to $\{\fvar,\svar,\tvar\}$, $\fvar=\beta(1)$, $\svar=\beta(2)$, $\tvar=\beta(3)$.
Then
$$\langle Q_{\sigma(1)},Q_{\sigma(2)},Q_{\sigma(3)} \rangle_e=
P(\beta \circ \sigma^{-1}(2 ),\beta \circ \sigma^{-1}(3 )).$$
\end{lemma}
\bp
View the ring $\QQ[\svar^{\pm1},\tvar^{\pm1}]$ as the quotient of the ring $\QQ[\fvar^{\pm1},\svar^{\pm1},\tvar^{\pm1}]$ by the relation $\fvar\svar\tvar=1$.
Then if an equivariant intersection point $\xi$ that contributes with a sign $\varepsilon$ to $\langle Q_1,Q_2,Q_3 \rangle_e$ lifts as $\xi_0 \in \TCM$, and if
$\xi_1=\theta^{\alpha_1}(\xi_0) \in Q_1$, $\xi_2=\theta^{\alpha_2}(\xi_0) \in Q_2$ and $\xi_3=\theta^{\alpha_3}(\xi_0) \in Q_3$, then
this point contributes as $c(\xi)=\varepsilon \fvar^{\alpha_1}\svar^{\alpha_2}\tvar^{\alpha_3}$.
Note that $\varepsilon$ does not depend on the order of $Q_1$, $Q_2$ and $Q_3$
and that this point will contribute to $\langle Q_{\sigma(1)},Q_{\sigma(2)},Q_{\sigma(3)} \rangle_e$ as 
$c_{\sigma}(\xi)=\varepsilon \fvar^{\alpha_{\sigma(1)}}\svar^{\alpha_{\sigma(2)}}\tvar^{\alpha_{\sigma(3)}}$. Then
$$c_{\sigma}(\xi)=\varepsilon \prod_{i=1}^3\beta(i)^{\alpha_{\sigma(i)}}=\varepsilon \prod_{i=1}^3\beta(\sigma^{-1}(i))^{\alpha_{i}}.$$
Assume without loss that 
$ \alpha_1=0$. 
Then $$c_{\sigma}(\xi)=\varepsilon \beta(\sigma^{-1}(2))^{\alpha_{2}}\beta(\sigma^{-1}(3))^{\alpha_{3}}.$$
\eop

\subsection{Equivariant triple intersection of representatives of the linking number}
\label{subinvtripl}
\begin{proposition}
\label{propuninvtripl}
Assume that $\tau$ maps the vectors of $T^+K$ 
to $K \times \RR^+\{\qvarM\}$ for some $\qvarM \in S^2$.
Let $\fvarM$, $\svarM$, $\tvarM$ be three distinct points in $(S^2 \setminus \{\qvarM,-\qvarM\})$, and let $K_{\fvarM}$, $K_{\svarM}$, $K_{\tvarM}$ and $K$ be four disjoint knots in $M$ that are rationally homologous to $K$.
For $\cvarM \in S^2$, $s_{\tau}(M;\cvarM)=\tau^{-1}(M \times \cvarM) \subset\partial ST(M)$.\\
There exist three $4$--dimensional chains $F_{\fvarM}$, $F_{\svarM}$ and $F_{\tvarM}$ of $C_2(M)$ with coefficients in $\QQ(t)$, whose boundaries
are 
$\partial F_{\fvarM} =s_{\tau}(M;\fvarM)-\ID ST(K_{\fvarM}) $, $\partial F_{\svarM} =s_{\tau}(M;\svarM)-\ID K_{\svarM} \times S^2$ and $\partial F_{\tvarM} =s_{\tau}(M;\tvarM)-\ID K_{\tvarM} \times S^2$, such that
$$\langle F_{\fvarM}, A(K)\rangle_e=\langle F_{\svarM}, A(K)\rangle_e=\langle F_{\tvarM}, A(K)\rangle_e=0.$$
Then $\langle F_{\fvarM},F_{\svarM},F_{\tvarM} \rangle_e$ only depends on the isotopy class of the link
$(K_{\fvarM},K_{\svarM},K_{\tvarM},K)$, on the homotopy class of $\tau$ and on the trivialisation of the tubular neighborhood of $K$ induced by $\tau_{|K}$.
It will be denoted by $\CQ(K_{\fvarM},K_{\svarM},K_{\tvarM},K,K_{ \parallel},\tau)$ where $K_{\parallel}$ is the parallel of $K$ induced by $\tau_{|K}$.
It will be denoted by $\CQ(K,K_{\parallel},\tau)$ when $K_{\fvarM}$, $K_{\svarM}$ and $K_{\tvarM}$ are parallels of $K$ with respect to the trivialisation induced by $K_{\parallel}$, living on the boundary of a tubular neighborhood of $K$.
\end{proposition}
\bp
The existence of the chain $F_{\fvarM}$ and the fact that if a chain $F^{\prime}_{\fvarM}$ satisfies the same conditions,
 then $(F^{\prime}_{\fvarM}-F_{\fvarM})$ is null-homologous come from
Proposition~\ref{propexistF}. Then according to the previous subsection, $\langle F_{\fvarM},F_{\svarM},F_{\tvarM} \rangle_e$ only depends on $\partial F_{\fvarM}$, $\partial F_{\svarM}$, and $\partial F_{\tvarM}$, and therefore on $(K_{\fvarM},K_{\svarM},K_{\tvarM},K,K_{ \parallel},\tau)$, $\qvarM$, $\fvarM$, $\svarM$ and $\tvarM$.

Hence, it is enough to see that $\langle F_{\fvarM},F_{\svarM},F_{\tvarM} \rangle_e$ does not vary 
\begin{itemize}
\item under a homotopy $\tau_t$ of $\tau$ such that $\tau_0(T^+K)=\tau_1(T^+K)=\RR^+\qvarM$, and $\tau_0$ and $\tau_1$ induce the same parallelisation of $K$ (where there is no loss in assuming that $\tau_t(T^+K)=\RR^+\qvarM$ for any $t$, as we show at the end of this proof),
\item under an isotopy of the link
$(K_{\fvarM},K_{\svarM},K_{\tvarM},K)$, extended to an isotopy $H \colon [0,1] \times M \rightarrow M$,
\item under a rotation of $SO(3)$ that acts on $\qvarM$, $\fvarM$, $\svarM$, $\tvarM$ and $\tau$, simultaneously,
\item when $(\fvarM,\svarM,\tvarM)$ varies among triples of distinct points of $S^2 \setminus\{\qvarM,-\qvarM\}$.
\end{itemize}
Without loss, assume that $F_{\fvarM}$, $F_{\svarM}$ and $F_{\tvarM}$ are products of their boundaries by the interval in a neighborhood $\partial \TCM \times [0,1]$ of $\partial \TCM$. Then each of these moves can be realized by global homeomorphisms of the whole picture in this neighborhood that change none of the intersection numbers $\langle F_{\fvarM}, A(K)\rangle_e$, $\langle F_{\svarM}, A(K)\rangle_e$, $\langle F_{\tvarM}, A(K)\rangle_e$ and
$\langle F_{\fvarM},F_{\svarM},F_{\tvarM} \rangle_e$.

To conclude, let us prove that there is no loss in assuming that $\tau_t(T^+K)=\RR^+\qvarM$ for any $t$, for the first item. It suffices to prove that --after
a possible composition by a homotopy $\tau^{\prime \prime}$ such that $\tau^{\prime \prime}_t(T^+K)=\RR^+\qvarM$ for any $t$-- $\tau_1$ is obtained from $\tau_0$ by a homotopy $\tau^{\prime}$ such that $\tau^{\prime}_t(T^+K)=\RR^+\qvarM$ for any $t$.

The indermediate trivialisation $\tau_t$ is the composition of $\tau_0$ with a map
$$\begin{array}{lll}
M \times \RR^3 &\rightarrow & M \times \RR^3\\
(m,v)& \mapsto & (m,\rho(t,m)(v))
\end{array}$$
for a map $\rho \colon [0,1]\times M \rightarrow SO(3)$ that
maps $\{0\} \times M$ to $1$, and that maps $\{1\} \times K(S^1)$ to a loop of the group $SO(2)(\qvarM)$ of rotations with axis $\RR \qvarM$. This loop is null-homotopic since $\tau_0$ and $\tau_1$ induce the same trivialisation of $K$.
Then there is a homotopy $\tau^{\prime \prime}$, such that $\tau^{\prime \prime}_t(T^+K)=\RR^+\qvarM$ between $\tau_1$ and a trivialisation $\tau^{\prime}_1$ such that $\tau^{\prime}_1=\tau_0$ on
$TM_{|K}$, and we assume from now on that $\tau_1=\tau_0$ on
$TM_{|K}$ without loss.
Change $\tau_t$ into $\rho(t,K(1))^{-1} \circ \tau_t$ again without loss so that $\tau_t=\tau_0$ at $K(1)$, for any $t$.
Now, the restriction of $\rho$ to $[0,1]\times K(S^1)$ maps
$$\left(\{0,1\} \times K(S^1) \right)\cup \left([0,1] \times K(1)\right)$$ to $1 \in SO(3)$, so that it can be seen as a map from $(S^2,\ast)$ to $(SO(3),1)$, that is homotopic to the constant map with value $1$ since $\pi_2(SO(3))$ is trivial.

See this homotopy as associated with a map 
$$\begin{array}{llll}r \colon& [0,1] \times [0,1]\times K(S^1)&\rightarrow & SO(3) \\
&(0,t,K(z)) &\mapsto &\rho(t,K(z)) \\
&(1,t,K(z)) &\mapsto &1 \\
&(u,t\in\{0,1\},K(z)) &\mapsto &1
\end{array}$$ as above.
Let $N(K)=K \times D^2$ be a tubular neighborhood of $K$ and see 
$D^2=\{dz_D; d \in [0,1], z_D \in S^1\}$.

Then replace our homotopy induced by $\rho$ by the homotopy induced by $r^{\prime}$ with
$$\begin{array}{llll}r^{\prime} \colon & [0,1] \times M & \rightarrow & SO(3)\\
&(t,m \notin N(K)) &\mapsto &\rho(t,m)\\
&(t,(K(z),dz_D)) &\mapsto &r(d,t,K(z))^{-1} \circ \rho(t,(K(z),dz_D)).
\end{array}$$
\eop

\subsection{An example in $S^1 \times S^2$}

In this subsection, we compute $\CQ(K,K_{\parallel},\tau)$ when $M=S^1 \times S^2$, $K =S^1 \times \qvarM$, $K_{\parallel} =S^1 \times \qvarM^{\prime}$, and 
$T(S^1 \times S^2)=S^1 \times \RR \times TS^2= S^1 \times (TS^2 \oplus \RR)$ is seen as the product by $S^1$ of the stabilisation $(TS^2 \oplus \RR)$ of the tangent bundle of $S^2$. The bundle $(TS^2 \oplus \RR)$ is the tangent bundle of $\RR^3$ restricted to its unit sphere that is trivialised as such, using the obvious trivialisation of $T\RR^3$.

\begin{proposition}
\label{propsonestwo}
Under the above assumptions,
$\CQ(K,K_{\parallel},\tau) =0.$
\end{proposition}
\bp
In this case, $$\TCMD=S^1 \times (S^2 \times \RR \times S^2)$$ where the projection map sends
$(z,v\in S^2,\alpha\in \RR,w \in S^2)$ to $((z,v),(z\exp(2i\pi \alpha),w)) \in M^2$.
The structure of a product by $S^1$ extends to the configuration space $\TCM$ that reads $S^1 \times \TCM_0$.
Note that the chain $A(K)$ of Lemma~\ref{lemak} also reads $S^1 \times A(K)_0$, for a $1$-chain  $A(K)_0$ in $\TCM_0$ (that is the closure of $\{\qvarM\} \times ]0,1[\times\{\qvarM\}$ in $\TCM_0$).
For a point $\fvarM \in S^2 \setminus \{\qvarM,-\qvarM\}$, we are going to contruct a rational $3$--chain $f_{\fvarM}$ in $\TCM_0$ with boundary
$$\partial f_{\fvarM}=(t^{-1}-1)s_{\tau}(\{1\} \times S^2;\fvarM) - (1+t^{-1})ST(\{1\} \times \fvarM)$$
explicitly, so that the boundary of
$F_{\fvarM}= \frac{1}{1-t^{-1}} (S^1 \times f_{\fvarM})$ will be $s_{\tau}(S^1 \times S^2;\fvarM) - \frac{1+t}{1-t}ST(S^1 \times \{\fvarM\})$.
Furthermore, since the dimension of $f_{\fvarM}$ is $3$ and the dimension of $A(K)_0$ is $1$ in the $5$ dimensional-manifold $\TCM_0$,
and since $\partial f_{\fvarM}$ does not meet $A(K)_0$,
we can assume that $f_{\fvarM}$ does not meet $A(K)_0$. 
Therefore, we can use chains $F_{\fvarM}$, $F_{\svarM}$ and $F_{\tvarM}$ that factor through $S^1$ to compute $\CQ(K,K_{\parallel},\tau)$. Since the projections of these chains in
$\TCM_0$ are of codimension $2$, again, we can assume that they have no equivariant triple intersection. Then $F_{\fvarM}$, $F_{\svarM}$ and $F_{\tvarM}$ have no equivariant triple intersection. Thus, we are left with the construction of $f_{\fvarM}$ to finish proving the proposition.

Consider the closure $G_0$ in $\TCM_0$
of (the natural lift of) $$\{((1,v),(\exp(2i\pi \alpha),v)); \alpha \in ]0,1[, v \in S^2\}$$
$$\partial(G_0) =t^{-1}s_-(S^2)-s_+(S^2)$$
where $s_+(S^2)$ is the section of $ST( S^2=\{1\} \times S^2)$ given by the direction of the $S^1$ factor and $s_-(S^2)$ is the opposite section.

We wish to change $s^+(S^2)$ to $s_{\tau}(S^2;\fvarM)$.

Let $v\in S^2$, note that $\tau(s_+(v))=(v,v)$.
If $v \neq - \fvarM$, there is a unique shortest arc $[\fvarM,v]$ of great circle from $\fvarM$ to $v$.
Let $s_{[\fvarM,+]}(S^2)$ denote the closure in $ST(S^2)$ of $$\mbox{Int}(s_{[\fvarM,+]}(S^2))=\cup_{v\in S^2\setminus\{-\fvarM\}}\tau^{-1}(\{v\}\times[\fvarM,v]).$$
Then 
$$\partial s_{[\fvarM,+]}(S^2)=s_+(S^2)-s_{\tau}(S^2;\fvarM) -ST(-\fvarM).$$
About $ST(-\fvarM)$, 
let $D(-\fvarM)$ be a disk centered at $(-\fvarM)$, $D(-\fvarM) \setminus \{-\fvarM\} =]0,1] \times \partial D(-\fvarM)$, $ST(-\fvarM)$ meets $s_{[\fvarM,+]}(S^2)$ on its boundary as $$-s_{\tau}(-\fvarM;lim_{\varepsilon \rightarrow 0 }\cup_{v\in\partial(\varepsilon D(-\fvarM))}[\fvarM,v])$$
that is as $-ST(-\fvarM)$. Note that the sign in front of $ST(-\fvarM)$ can also be deduced from the fact that $s_+(S^2)-s_{\tau}(S^2;\fvarM)$ is homologous to $ST(-\fvarM)$.

Similarly, let $s_{[\fvarM,-]}(S^2)$ denote the closure in $ST(S^2)$ of $\cup_{v\in S^2\setminus\{\fvarM\}}s_{\tau}(v;[\fvarM,-v])$.
$$\partial s_{[\fvarM,-]}(S^2)=s_-(S^2)-s_{\tau}(S^2;\fvarM) +ST(\fvarM).$$
Let $\alpha$ an arc of $S^2$ such that $\partial \alpha = \{-\fvarM\} -\{\fvarM\}$, $\partial ST(\alpha)=ST(-\fvarM)-ST(\fvarM)$.

Let $$G=G_0 + s_{[\fvarM,+]}(S^2) -t^{-1}s_{[\fvarM,-]}(S^2) + ST(\alpha).$$
$$\partial G=(t^{-1}-1)s_{\tau}(S^2;\fvarM) - ST(\fvarM) -t^{-1}ST(\fvarM)$$
and $G$ can be transformed into a $3$--chain $f_{\fvarM}$
transverse to the boundary.
\eop

\subsection{Symmetries}
\label{subsecsym}

\begin{proposition}
\label{propsym1}
Under the hypotheses of Proposition~\ref{propuninvtripl},
let $\CQ=\CQ(K_{\fvarM},K_{\svarM},K_{\tvarM},K,K_{ \parallel},\tau)$, then
$$\CQ(\fvar,\svar,\tvar)=\CQ(\fvar^{-1},\svar^{-1},\tvar^{-1}).$$
\end{proposition}
\bp
Consider the involution $\iota$ of $\TCM$ induced by
the involution that exchanges the two factors of $\tilde{M}^2$. This involution
reverses the orientation of $\TCM$ and is such that $\iota\theta=\theta^{-1}\iota$.
Let $\iota(F_{\fvarM})$ inherit its orientation from $F_{\fvarM}$.
Then $$\partial \iota(F_{\fvarM})=s_{\tau}(M;-\fvarM)-\ID(t^{-1})(-ST(K_{\fvarM}))=s_{\tau}(M;-\fvarM)-\ID(t)ST(K_{\fvarM})$$
and $$\iota(A(K))= -tA(K).$$
Therefore,
$$\CQ=\langle \iota(F_{\fvarM}),\iota(F_{\svarM}),\iota(F_{\tvarM}) \rangle_{e,\TCM}.$$
If an intersection point $u$ contributes to $\langle F_{\fvarM},F_{\svarM},F_{\tvarM} \rangle_{e,\TCM}$ as $\varepsilon \fvar^a \svar^b \tvar^c$ , the point $\iota(u)$ contributes to $\langle \iota(F_{\fvarM}),\iota(F_{\svarM}),\iota(F_{\tvarM}) \rangle_{e,\TCM}$ as  $\varepsilon \fvar^{-a} \svar^{-b} \tvar^{-c}$.
Let us check that the signs are indeed the same.
The oriented normal of $\iota(F_{\fvarM})$ in $\TCM$ is $-\iota(N(F_{\fvarM}))$. Therefore the sign for $\iota(u)$ will be positive if and only if $\iota(N(F_{\fvarM})) \oplus \iota(N(F_{\svarM})) \oplus \iota(N(F_{\tvarM}))$ induces the orientation of $(-\TCM)$ that is the orientation of $\iota(\TCM)$.
Hence, the sign for $\iota(u)$ will be positive if and only if the sign for $u$ is positive.
\eop

\begin{proposition}
\label{propsym2}
Under the hypotheses of Proposition~\ref{propuninvtripl},
let $\CQ=\CQ(K,K_{ \parallel},\tau)$, then
$$\CQ(\fvar,\svar,\tvar)=\CQ(\svar,\fvar,\tvar)=\CQ(\fvar,\tvar,\svar).$$
\end{proposition}
\bp
This is a direct consequence of Lemma~\ref{lemsym} and Proposition~\ref{propuninvtripl}.
\eop

\begin{proposition}
\label{propsym3}
Under the hypotheses of Proposition~\ref{propuninvtripl},
$\CQ(K,K_{\parallel},\tau)$ does not depend on the orientation of $K$.
\end{proposition}
\bp
When the orientation of $K$ is reversed, $t$ is changed into $t^{-1}$,
and $A(K)$ is changed by a multiplication by a unit. Therefore, we may use the same chains $F$ to compute $\CQ$ and Proposition~\ref{propsym1}
allows us to conclude.
\eop

\begin{proposition}
\label{propsym4}
Let $(-M)$ denote the manifold obtained from $M$ by reversing the orientation of $M$.
Let $\tau_{-M}$ be a trivialisation of $T(-M)$ obtained from $\tau$ by a composition by a fixed
orientation-reversing isomorphism of $\RR^3$.
Under the hypotheses of Proposition~\ref{propuninvtripl},
$$\CQ(K \subset (-M),K_{\parallel},\tau_{-M})=-\CQ(K,K_{\parallel},\tau).$$
\end{proposition}
\bp
$$\CQ(K \subset (-M),K_{\parallel},\tau_{-M})=\langle -F_X,-F_Y,-F_Z\rangle_e$$ 
in $\tilde{C}_2(-M)=\TCM$.
\eop

\subsection{Changing the manifold parallelisation}
\label{submanpar}

Any closed oriented $3$-manifold $M$ bounds an oriented compact $4$-dimensional manifold $W$ with signature $0$. Then $TW_{|M}=\RR \oplus TM$.
A trivialisation $\tau$ of $TM$ induces a trivialisation of $TW \otimes \CC$ on $M$. 
The {\em first Pontrjagin class\/} $p_1(\tau)$ of such a trivialisation $\tau$ of the tangent bundle of $M$ is the obstruction $p_1(W;\tau)$ to extend this trivialisation to $W$. It belongs to $H^4(W,M;\pi_3(SU(4)))\cong \ZZ$.  We use the notation and conventions of \cite{milnorsta}, see also \cite[Section 1.5]{lesconst}.

Since $p_1(X)=3\,\mbox{signature}(X)$ for a closed (oriented) $4$-manifold $X$, and since the signature of a closed $4$-manifold 
obtained by gluing two $4$-manifolds $W_1$ and $W_2$ along their common boundary is the sum of the signatures of $W_1$ and $W_2$ by Novikov additivity, $p_1(\tau)$ is well-defined.

\begin{proposition}
\label{propvartau}
 Let $\tau$ and $\tau^{\prime}$ be two trivialisations that induce the same parallelisation on $K$. Then
$$\CQ(K,K_{ \parallel},\tau^{\prime})-\CQ(K,K_{ \parallel},\tau)=\frac{p_1(\tau^{\prime})-p_1(\tau)}{4}.$$
In particular $$\CQ(K,K_{\parallel})=\CQ(K,K_{\parallel},\tau) -\frac{p_1(\tau)}{4}$$
does not depend on $\tau$.
\end{proposition}

\begin{remark}
\label{remvarpar}
Since $\pi_2(SU(4))=\{0\}$, the stabilised complexification of the trivialisation $\tau$ of $S^1 \times S^2$ involved in Proposition~\ref{propsonestwo} extends to $S^1 \times B^3$ whose signature is zero. Therefore $p_1(\tau)=0$, and $$\CQ((S^1 \times \qvarM,S^1 \times \qvarM^{\prime})\subset S^1 \times S^2) =0.$$
\end{remark}

\noindent{\sc Proof of Proposition~\ref{propvartau}:}
Let $\tau\colon TM \rightarrow M \times\RR^3$ be a parallelisation of $M$.
Any other parallelisation $\tau^{\prime}$ of $M$ reads $\tau^{\prime}=G\circ \tau$ for some $$\begin{array}{llll}G\colon &M \times\RR^3 &\rightarrow& M \times\RR^3\\
&(m,v)& \mapsto &(m,(g(m))^{-1}(v))\end{array}$$
associated with a map $g \colon M \rightarrow SO(3)$.

Since we work up to homotopy,
we shall assume that  
$G$ maps a tubular neighborhood $K \times D^2$ of $K$ that contains $K_{\fvarM}$, $K_{\svarM}$ and $K_{\tvarM}$ to $1$. 

Set $E=M \setminus \mbox{Int}\left(K \times D^2\right)$. In
$[0,3] \times ST(M),$
that is thought of as an external collar glued to $ST(M)=\{0\} \times ST(M)$,
consider the $3$-cycle 
$$c_{\fvarM}=1 \times G^{-1}(E \times \fvarM)-0 \times (E \times \fvarM) -[0,1] \times (K \times S^1) \times \fvarM.$$

\begin{lemma}
\label{lemhomhthree}
 The homology class of $c_{\fvarM}$ vanishes in $H_3([0,1] \times ST(E))$.
\end{lemma}
\bp
The algebraic intersection with $[0,1] \times (S \setminus \mbox{Int}(D^2))\times_{\tau} \svarM$ determines the elements of $H_3\left([0,1] \times ST(E)\right)=\QQ[\{1/2\} \times ST(K_{\parallel})].$
In particular, to prove that the homology class of $c_{\fvarM}$ is zero, it suffices to prove that the degree of $g(\cdot)(\fvarM)\colon S\rightarrow S^2$ at $\svarM$ is zero.
Of course, this degree only depends on the homotopy class of the restriction of $g$ to $S$. We may assume that $\svarM$ is a regular point of $g(\cdot)(\fvarM)$. The changes of homotopy class of $g\colon S\rightarrow SO(3)$ may be realised by compositions of rotations with axis $\svarM$ supported in neighborhoods of a geometric basis of $H_1(S)$ that avoid $g(\cdot)(\fvarM)^{-1}(\svarM)$ (since $\pi_2(SO(3))=\{0\}$ and 
$\pi_1(SO(3))$ is generated by a path of rotations with axis $\svarM$).
Therefore, the degree is independent of $g$ and equal to the degree for a constant map to $1\in SO(3)$ that is zero.
This ends the proof of the lemma. \eop

Back to the proof of Proposition~\ref{propvartau}, $c_{\fvarM}$ is the boundary of a chain $C_{0,\fvarM}$ in $[0,1] \times ST\left(E\right)$, and 
$$C_{\fvarM} =C_{0,\fvarM} + \left([0,1] \times (K \times D^2) \times \fvarM \right)$$ is a $4$--chain in $[0,1] \times M \times_{\tau} S^2$ such that
$$\partial C_{\fvarM} =\{1\} \times s_{\tau^{\prime}}(M;\fvarM)-\{0\} \times s_{\tau}(M;\fvarM)$$ and $C_{\fvarM}$ meets $[0,1] \times (K \times D^2) \times_{\tau} S^2$ as $[0,1] \times (K \times D^2) \times_{\tau} \fvarM$.

Set $$D_{\fvarM}=C_{\fvarM} + \left([1,3] \times s_{\tau^{\prime}}(M;\fvarM)\right) -\ID(t) [0,3] \times ST(K_{\fvarM}),$$
$$D_{\svarM}= \left([0,1] \times s_{\tau}(M;\svarM)\right) + C_{\svarM} + \left([2,3] \times s_{\tau^{\prime}}(M;\svarM)\right) -\ID(t) [0,3] \times ST(K_{\svarM}),$$
$$D_{\tvarM}= \left([0,2] \times s_{\tau}(M;\tvarM)\right) + C_{\tvarM} -\ID(t) [0,3] \times ST(K_{\tvarM}),$$
where \begin{itemize}
       \item $C_{\svarM} \subset [1,2] \times ST(M)$, $C_{\tvarM} \subset [2,3] \times ST(M)$
\item $\partial C_{\svarM} =\{2\} \times s_{\tau^{\prime}}(M;\svarM)-\{1\} \times s_{\tau}(M;\svarM)$, $\partial C_{\tvarM} =\{3\} \times s_{\tau^{\prime}}(M;\tvarM)-\{2\} \times s_{\tau}(M;\tvarM)$
\item $C_{\svarM} \cap \left([1,2] \times ST(K \times D^2)\right)=[1,2] \times (K \times D^2) \times_{\tau} \{\svarM\}$
and $C_{\tvarM} \cap [2,3] \times ST(K \times D^2)=[2,3] \times (K \times D^2) \times_{\tau} \{\tvarM\}$.
      \end{itemize}

Now $$
\delta=\CQ(K,K_{\parallel},\tau^{\prime})-\CQ(K,K_{\parallel},\tau)=
\langle D_{\fvarM},D_{\svarM},D_{\tvarM} \rangle_e.$$
First note that no triple intersection can occur in $[0,3] \times ST(K \times D^2)$. Then no triple intersection can occur in $([0,1] \cup [2,3]) \times ST(M)$ either and 
$$\delta=\langle [1,2] \times s_{G \circ \tau}(M;\fvarM),C_{\svarM},[1,2] \times s_{\tau}(M;\tvarM)\rangle_e.$$
It does not depend on $\tau$. Therefore $\delta$ will be denoted by $\delta(g)$.
Then $\delta(g^2)=2\delta(g)$. Since $\pi_1(SO(3))=\ZZ/2\ZZ$ and $\pi_2(SO(3))=0$, $g^2$ is homotopic to a map that sends the complement of a $3$-ball $B^3$ in $M$ to the identity, that is homotopic to $[p^k_{SO(3)}]$ where $[p_{SO(3)}]$ is defined as follows.
See $$B^3=[0,2\pi]\times S^2/ (0 \times S^2 \sim 0).$$
Then $p_{SO(3)}(\beta,\cvarM)$ is the rotation $R_{\beta,\cvarM}$ with angle $\beta$ and whose axis is directed by $\cvarM$. (Indeed, $S^3$ is the quotient of $B^3$ by its boundary $\{2\pi\} \times S^2$, $p_{SO(3)}$
factors through $S^3$ to define the covering map $\overline{p_{SO(3)}}$ whose homotopy class generates $\pi_3(SO(3))=\ZZ[\overline{p_{SO(3)}}]$.)
Then $\delta(g)=\frac{k}{2}\delta(p_{SO(3)})$.

\begin{lemma}
$\delta(p_{SO(3)})=1$.
\end{lemma}
\bp
When $g=p_{SO(3)}$, we can assume that 
$$C_{\svarM} \cap \left([1,2] \times ST(M \setminus B^3)\right)=[1,2] \times (M \setminus B^3)\times_{\tau} \{\svarM\} .$$
Then $$\delta(p_{SO(3)})=\langle [1,2] \times  B^3\times_{\tau^{\prime}} \{\fvarM\}, C_{\svarM},[1,2] \times  B^3 \times_{\tau}\{\tvarM\}\rangle_{e,\ZZ \times [1,2] \times B^3 \times S^2 }.$$
Assume that $\tvarM$ is the North Pole and that $\fvarM=-\tvarM$ is the South Pole. Then $R_{\beta,\cvarM}(-\tvarM)=\tvarM$ if and only if $\beta = \pi$ and $\cvarM \in \partial D(\tvarM )$, where $\partial D(\tvarM )$ is the equator that is oriented as the boundary of the northern hemisphere $D(\tvarM )$.  Consider a point $\cvarM \in \partial D(\tvarM )$, the great circle $(\tvarM,\cvarM,-\tvarM)$ containing $\cvarM$ and the poles -oriented by $(\tvarM,\cvarM,-\tvarM)$- and the orthogonal great circle $\partial D(\cvarM)$ containing the poles oriented as the boundary of the hemisphere $D(\cvarM)$ containing $\cvarM$. 
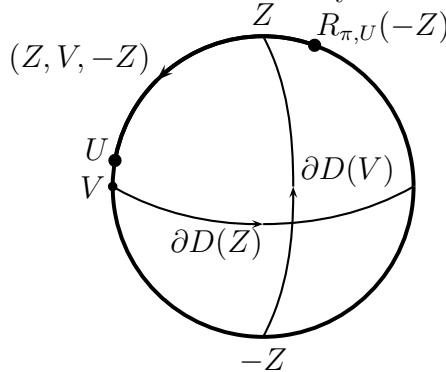
\begin{center}
\begin{pspicture}[shift=-0.1](-2,-.2)(6,4.2)
\rput[r](.5,3.7){$(\tvarM,\cvarM,-\tvarM)$}
\psarc[linewidth=1.5pt]{*-*}(2,2){2}{70}{170}
\rput[bl](2.7,3.9){$R_{\pi,U}(-\tvarM)$}
\rput[r](-0.05,2.5){$U$}
\rput[r](-.1,2){$\cvarM$}
\rput[b](2,4.15){$\tvarM$}
\rput[t](2,-.15){$-\tvarM$}
\psecurve{*->}(-1,3)(0,2)(2,1.5)(4,2)
\rput[rt](2,1.45){$\partial D(\tvarM)$}
\psecurve(0,2)(2,1.5)(4,2)(5,3)
\psecurve{->}(1,-1)(2,0)(2.4,2)(2,4)
\rput[lb](2.5,2){$\partial D(\cvarM)$}
\psecurve(2,0)(2.4,2)(2,4)(1,5)
\psarc[linewidth=1.5pt]{->}(2,2){2}{0}{135}
\psarc[linewidth=1.5pt](2,2){2}{135}{360}
\end{pspicture}
\end{center}

The normal bundle of $\pi \times \partial D(\tvarM )$ in $S^3$ at $\cvarM$ is $[0,2 \pi] \times$ $(\tvarM,\cvarM,-\tvarM)$. Along $[0,2 \pi]$, near $\pi$, $R_{\beta,\cvarM}(-\tvarM)$ moves along $\partial D(\cvarM)$.
Along $(\tvarM,\cvarM,-\tvarM)$, near $\cvarM$, $R_{\pi,U}(-\tvarM)$ moves along $(\tvarM,\cvarM,-\tvarM)$ near $\tvarM$. Therefore the normal bundle of $\pi \times \partial D(\tvarM )$ maps to $S^2$ in an orientation-preserving way and
$$\delta(p_{SO(3)})=\langle [1,2] \times \pi \times \partial D(\tvarM ) \times \tvarM, C_{\svarM}\rangle_{[1,2] \times B^3 \times S^2}.$$
Set $$C_{\svarM,\partial D(\tvarM)}=C_{\svarM} \cap \left( [1,2] \times\pi \times \partial D(\tvarM ) \times S^2 \right)$$
and $C_{\svarM,D(\tvarM)}=C_{\svarM} \cap \left( [1,2] \times\pi \times D(\tvarM ) \times S^2 \right)$. Without loss, assume that the above intersections are transverse.
$\delta(p_{SO(3)})$ is the degree of the projection of $C_{\svarM,\partial D(\tvarM)}$ on $S^2$ at $\tvarM$.
$$\partial C_{\svarM,\partial D(\tvarM)}= \left(2 \times \pi \times \partial D(\tvarM )\times_{\tau^{\prime}} \svarM\right) - \left(1 \times \pi \times \partial D(\tvarM )
\times \svarM\right).$$
Assume that $\svarM \in \partial D(\tvarM )$. Then $\pi \times \partial D(\tvarM )\times_{\tau^{\prime}} \svarM$ projects to $S^2$ as the double cover of $\partial D(\tvarM )$.
In particular the degree of the projection of $C_{\svarM,\partial D(\tvarM)}$ on $S^2$ is constant and equal to $\delta(p_{SO(3)})$ on the northern hemisphere, and, it is equal to $(\delta(p_{SO(3)})-2)$ on the southern hemisphere.

Since the degree of the projection on $S^2$ (via $\tau$) of
$$\partial C_{\svarM,D(\tvarM)}= \pm C_{\svarM} \cap \left((\partial [1,2]) \times\pi \times D(\tvarM ) \times S^2\right) \pm C_{\svarM,\partial D(\tvarM)}$$
vanishes, the degree of the projection of $C_{\svarM,\partial D(\tvarM)}$ on $S^2$ coincides up to sign with the degree of the projection of $\pi \times  D(\tvarM )\times_{\tau^{\prime}} \svarM$ that completely covers $S^2$ generically once.
Therefore, $\delta(p_{SO(3)})=1$. \eop

Similarly, define $p^{\prime}_1(g) = p_1(G \circ \tau)-p_1(\tau)$.
It is easy to see that this expression, that can be seen as an obstruction to extend a trivialisation on $[0,3]\times M$, does not depend on $\tau$ and that
$p^{\prime}_1$ is additive under the multiplication in $SO(3)$, so that $p^{\prime}_1(g)=\frac{k}{2}p^{\prime}_1(p_{SO(3)})$.
According to Proposition~1.8 in \cite{lesconst}, $p^{\prime}_1(p_{SO(3)}) =4$.
This ends the proof of Proposition~\ref{propvartau}.
\eop

\newpage
\section{Connected sum with a rational homology sphere}
\label{secconcas}
\setcounter{equation}{0}

This section is devoted to prove Proposition~\ref{propconcas} that
asserts that, for a rational homology sphere $N$, $$\CQ(\KK \subset M \sharp N)=\CQ(\KK \subset M)+6 \lambda(N).$$ The proof is based 
on the configuration space definition of the Casson-Walker invariant $\lambda$, that we summarize now.

\subsection{A configuration space definition of the Casson-Walker invariant}
\label{subdefcas}

For $r\in \RR$, let $B(r)$ denote the ball of radius $r$ in $\RR^3$
that is equipped with its standard parallelisation $\tau_s$.
A rational homology sphere $N$ may be written as $B_N \cup_{B(1) \setminus \mbox{\small Int}(B(1/2))} B^3$ where
$B_N$ is a {\em rational homology ball\/}, that is a connected compact (oriented) smooth $3$--manifold with boundary $S^2$ with the same rational homology as a point, $B^3$ is a $3$-ball, $B_N$ contains $(B(1) \setminus \mbox{Int}(B(1/2)))$ as a neighborhood of its boundary $\partial B_N=\partial B(1)$, and $B^3$ contains $(B(1) \setminus\mbox{Int}(B(1/2)))$ as a neighborhood of its boundary $\partial B^3=-\partial B(1/2)$.
Let $B(N)=B_N(3)$ be obtained from $B(3)$ by replacing the unit ball $B(1)$ of $\RR^3$ by $B_N$. Equip $B(N)$ with a trivialisation $\tau_N$ that coincides with $\tau_s$ outside $B_N$.

Let $W$ be a compact connected $4$-manifold with signature $0$ and with boundary 
$$\partial W =B_N(3) \cup_{\{1\} \times \partial B(3)} \left( -[0,1] \times \partial B(3) \right) \cup_{\{0\} \times \partial B(3)} (-B(3)).$$
Define $p_1(\tau_N) \in (H^4(W,\partial W;\pi_3(SU(4))) = \ZZ)$ as the obstruction to extend the trivialisation of $TW \otimes \CC$
induced by $\tau_s$ and $\tau_N$ on $\partial W$ to $W$. 
Again, we use the notation and conventions of \cite{milnorsta}, see also \cite[Section 1.5]{lesconst}.

Let $\RR^3(N)$ be obtained from $\RR^3$ by replacing its unit ball $B(1)$ by $B_N$. Let $C_2(\RR^3(N))$ be obtained from $\RR^3(N)^2$ by blowing-up the diagonal as in Subsection~\ref{subtcm}. Let $\projconf \colon C_2(\RR^3(N)) \rightarrow \RR^3(N)^2$ be the associated canonical projection and let $C_2(B(N))=\projconf^{-1}(B(N)^2)$.
Consider a smooth map 
$ \chi \colon \RR \rightarrow [0,1]$ that maps $]-\infty,-2]$ to $0$ and
$[-1,\infty[$ to $1$.
Define $$\begin{array}{llll}p_{B(3)} \colon &B(3)^2\setminus \mbox{diagonal}& \rightarrow &S^2\\
&(U,V) &\mapsto &\frac{\chi(\parallel V\parallel-\parallel U \parallel)V-\chi(\parallel U\parallel-\parallel V \parallel)U}{\parallel \chi(\parallel V\parallel-\parallel U \parallel)V-\chi(\parallel U\parallel-\parallel V \parallel)U \parallel} \end{array}
$$
This map extends to $C_2(B(3))$ to a map still denoted by $p_{B(3)}$, that reads as the projection to $S^2$ induced by $\tau_s$ (see Subsection~\ref{subtcm}) on the unit tangent bundle of $B(3)$.
A similar map $p_N$ can be defined on the boundary $\partial C_2(B(N))$: The map $p_N$ is the projection to $S^2$ induced by $\tau_N$ on the unit tangent bundle of $B(N)$, and the map $p_N$ is 
given by the above formula, where we set $\parallel U\parallel=1$ when $U \in B_N$, for the other points of the boundary that are pairs
$(U,V)$ of $\left(B(N)^2\setminus \mbox{diagonal}\right)$ where $U$ or $V$ belongs to $\partial B(3)$ (therefore a possible point of $B_N$ is replaced by $0 \in \RR^3$ in the formula).

The following theorem, that gives a configuration space definition for the Casson-Walker invariant, is due to Kuperberg and Thurston \cite{kt} for the case of integral homology spheres, though it is stated in other words. It has been generalised to rational homology spheres in \cite[Section 6]{sumgen}.

\begin{theorem}
\label{thmdefcasconf}
Let $\fvarM$, $\svarM$ and $\tvarM$ be three distinct points of $S^2$.
Under the above assumptions, for $\cvarM=\fvarM$, $\svarM$ or $\tvarM$, the $3$--cycle $p_N^{-1}(\cvarM)$ of $\partial C_2(B(N))$ bounds a rational chain $F_{N,\cvarM}$
in $C_2(B(N))$, and
$$\lambda(N)=\frac{\langle F_{N,\fvarM},F_{N,\svarM}, F_{N,\tvarM}\rangle_{C_2(B(N))}}{6} -\frac{p_1(\tau_N)}{24}.$$
\end{theorem}
\bp This essentially comes from \cite[Sections 1.1 to 1.5]{lesconst}
and \cite[Section 6.5]{sumgen} where is is shown that
$$\lambda(N)=\frac{\langle \Sigma_{\fvarM},\Sigma_{\svarM}, \Sigma_{\tvarM}\rangle_{C_2(N)}}{6} -\frac{p_1(\tau_N)}{24}$$
for chains $\Sigma_{\cvarM}$ of a configuration space $C_2(N)$ with prescribed boundaries $p_N(\tau_N)^{-1}(\cvarM)$
associated to a given map $p_N(\tau_N) \colon \partial C_2(N) \rightarrow S^2$.

The configuration space $C_2(N)$ contains $C_2(B(N))$ and it has the same homotopy type; the map $p_N(\tau_N)$ plays the same role as the above map $p_N$. That is essentially why
$$\langle F_{N,\fvarM},F_{N,\svarM}, F_{N,\tvarM}\rangle_{C_2(B(N))}
=\langle \Sigma_{\fvarM},\Sigma_{\svarM}, \Sigma_{\tvarM}\rangle_{C_2(N)}.$$
We nevertheless prove this with more details below.

The configuration space $C_2(N)$ is a compactification  of $$\mbox{Int}(C_2(N))=\left(\left(\RR^3 \setminus \mbox{Int}(B(1))\right) \cup_{\partial B(1)} B_N\right)^2 \setminus \mbox{diagonal}.$$ 
A point of $\partial C_2(N)$ may be viewed as a limit of a converging sequence $((U_n,V_n) \in \mbox{Int}(C_2(N)))_{n \in \NN}$.
The map $p_N(\tau_N)$
maps such a limit to $\lim_{n \rightarrow \infty} \frac{V_n-U_n}{\parallel V_n-U_n \parallel }$. The compactification is defined so that this limit makes sense.

In order to finish our proof, it suffices to extend the map $p_N$ to $C_2(N)\setminus \mbox{Int}(C_2(B(N)))$ so that $p_N=p_N(\tau_N)$ on
$\partial C_2(N)$. Indeed, this will allow us to set 
$$\Sigma_{\cvarM}=F_{N,\cvarM} \cup p_{N|C_2(N)\setminus \mbox{\small Int}(C_2(B(N)))}^{-1}(\cvarM)$$ and will make clear that the two algebraic intersections coincide.

Define a smooth map $\chi_2 \colon \RR^+ \rightarrow [0,1]$ such that
$\chi_2(x)=0$ if $x \leq 3$, and $\chi_2(x)=1$ if $x \geq 4$.
and $$\begin{array}{llll}\chi_3 \colon & \left(\RR^+\right)^2 &\rightarrow &[0,1]\\
       &(x,y)&\mapsto&(1-\chi_2(x))\chi(x-y) +\chi_2(x)
      \end{array}$$
so that $\chi_3(x,y) = 1$ if $x \geq 4$ or if $(x-y) \geq -1$.
Then define 
$$\begin{array}{llll}p_{N} \colon &\mbox{Int}(C_2(N))\setminus \mbox{Int}(C_2(B(N)))& \rightarrow &S^2\\
&(U,V) &\mapsto &\frac{\chi_3(\parallel V\parallel,\parallel U \parallel)V-\chi_3(\parallel U\parallel,\parallel V \parallel)U}{\parallel \chi_3(\parallel V\parallel,\parallel U \parallel)V-\chi_3(\parallel U\parallel,\parallel V \parallel)U \parallel} \end{array}
$$
where $\parallel U\parallel=1$ when $U \in B_N$.

Then $p_N$ coincides with the former $p_N$ on $\partial C_2(B(N))$
because $\chi_2=0$ there.
The pairs near the boundary of $C_2(N)$ are the pairs near the diagonal or the pairs of points where at least one point goes to $\infty$. For these pairs $\chi_3$ is $1$, unless one point is in $B_N(4)$ and the other one is near $\infty$, but in this case the limit does not see the point in $B_N(4)$, anyway.
\eop

\subsection{Proof of Proposition~\ref{propconcas}}
\label{proofconcas}

We keep all the notation from the previous subsection.

Consider the ball $B(3)$ of the previous subsection as a small ball embedded in $M$, outside a fixed neighborhood of $K$
that contains all the needed parallels of $K$, so that the trivialisation $\tau$ of $M$ coincides with the trivialisation $\tau_s$ on $B(3)$.
We shall perform the connected sum by letting the rational homology ball $B_N$ of the previous subsection replace $B(1)$.
Therefore we normalize our chains $F_{\cvarM}$ for $M$, over pairs of points that contain a point in $B(1)$.

Let $P \in \ZZ[t,t^{-1}]$ be such that $G_{\cvarM}=P F_{\cvarM}$ is an integral chain for $\cvarM=\fvarM$, $\svarM$ and $\tvarM$.
Let $m$ be the image of $0 \in B(3)$ under the implicit embedding of $B(3)$.
Without loss, we successively assume that 
\begin{itemize}
\item $G_{\cvarM}$ coincides with $P.p_{B(3)}^{-1}(\cvarM)$ on $\coprod_{k\in \ZZ} \theta^k(C_2(B(3))) \subset \TCM$,
\item
$G_{\cvarM}$ is transverse to the closures of ${m \times (\tilde{M} \setminus p_M^{-1}(B(3)))}$ and $(\tilde{M} \setminus p_M^{-1}(B(3))) \times m$,
\item 
the intersection with these pieces read $\pm m \times \gamma_1(m;\cvarM)$ and $\pm \gamma_2(m;\cvarM) \times m$ for
$1$-manifolds $\gamma_1(m;\cvarM)$ and $\gamma_2(m;\cvarM)$ of $\tilde{M}$ whose boundaries are supported in $p_M^{-1}(3\cvarM)$ and 
$p_M^{-1}(-3\cvarM)$, respectively,
\item $G_{\cvarM}$ intersects the closures of ${B(1) \times (\tilde{M} \setminus p_M^{-1}(B(3)))}$ and $(\tilde{M} \setminus p_M^{-1}(B(3))) \times B(1)$ as $ B(1) \times \gamma_1(m;\cvarM)$ and as $ \gamma_2(m;\cvarM) \times B(1)$, respectively,
\item the paths $\gamma_1(m;\fvarM)$, $\gamma_1(m;\svarM)$ and $\gamma_1(m;\tvarM)$ are pairwise disjoint,
\item the paths $\gamma_2(m;\fvarM)$, $\gamma_2(m;\svarM)$ and $\gamma_2(m;\tvarM)$ are pairwise disjoint.
\end{itemize}
These conditions can be successively achieved by small perturbations without losing the former ones.
Then let $F_{\sharp N,\cvarM}$ be the $4$-chain of $\tilde{C}_2(M\sharp N)$ that coincides
\begin{itemize}
\item with $F_{N,\cvarM} \subset
\left( C_2(B(N))=\{0\} \times C_2(B(N)) \right) \subset \left(\ZZ \times  C_2(B(N)) =\tilde{C}_2(B(N)=B_N(3)) \right)$ on $\tilde{C}_2(B(N))$,
\item with $\frac{1}{P}\left(B_N \times \gamma_1(m;\cvarM)\right)$ on the closure of $B_N \times (\tilde{M} \setminus p_M^{-1}(B(3)))$,
\item with $\frac{1}{P}\left(\gamma_2(m;\cvarM) \times B_N\right)$ on the closure of $(\tilde{M} \setminus p_M^{-1}(B(3)))\times B_N$,
\item with $F_{\cvarM}$, on the closure $E$ of the complement of the above subsets.
\end{itemize}

Then $\CQ(\KK \subset M)=\langle F_{\fvarM}, F_{\svarM},F_{\tvarM}\rangle_{e} -\frac{p_1(\tau)}{4}$ 
where the $F_{\cvarM}$ only intersect on $E$ under the above assumptions, and
$$\CQ(\KK \subset M \sharp N)=\langle F_{\sharp N,\fvarM}, F_{\sharp N,\svarM},F_{\sharp N,\tvarM}\rangle_{e} -\frac{p_1(\tau_{M \sharp N})}{4}$$
where $\tau_{M \sharp N}$ is the trivialisation induced by $\tau$ and $\tau_N$ that satisfies $p_1(\tau_{M \sharp N})=p_1(\tau)+p_1(\tau_N)$ and
$$\langle F_{\sharp N,\fvarM}, F_{\sharp N,\svarM},F_{\sharp N,\tvarM}\rangle_{e}=\langle F_{\fvarM}, F_{\svarM},F_{\tvarM}\rangle_{e,E}+
\langle F_{N,\fvarM}, F_{N,\svarM},F_{N,\tvarM} \rangle.$$

This shows that $\CQ(\KK \subset M \sharp N)=\CQ(\KK \subset M)+6 \lambda(N)$.
\eop

\newpage 
\section{First variations}
\setcounter{equation}{0}
\label{secvar}

\subsection{Preliminary computations}

To compute these triple intersections or their variations, we summarize some arguments that we shall frequently use.
Recall that a neighborhood of $\partial \TCM$ in $\TCM$ is diffeomorphic to $\ZZ \times ST(M) \times [0,1]$ where $[0,1]$ is the inward normal, and $ST(M)=_{\tau} M \times S^2$. (All the diffeomorphisms preserve the orientation, and the order of appearance of coordinates induces the orientation.) Near its boundary, $F_{\fvarM}$ can also be assumed to be diffeomorphic to 
$\partial F_{\fvarM}  \times [0,1]$, that is $s_{\tau}(M;\fvarM) \times [0,1] -I_{\Delta}(t)ST(K_{\fvarM}) \times [0,1]$ where the oriented normal of $s_{\tau}(M;\fvarM)$ in $\partial \TCM $ coincides with the oriented normal of  $s_{\tau}(M;\fvarM) \times [0,1]$ in $\TCM$ and is (the oriented tangent space of) $ST(\ast)$ and the oriented normal of $ST(K_{\fvarM})$ is the surface $S$. In particular, the oriented intersection of $s_{\tau}(M;\fvarM) \times [0,1]$ and $ST(K_{\svarM}) \times [0,1]$ is $s_{\tau}(K_{\svarM};\fvarM)\times [0,1]$, that is cooriented by $ST(S)$.

The following lemmas easily follow from these considerations.
\begin{lemma}
\label{leminsts}
$$\langle ST(S),F_{\svarM},F_{\tvarM} \rangle_e=-\ID(\svar) - \ID(\tvar).$$
$$\langle F_{\fvarM},ST(S),F_{\tvarM} \rangle_e=-\ID(\fvar) - \ID(\tvar)(=\ID(\svar\tvar) - \ID(\tvar)).$$
$$\langle F_{\fvarM},F_{\svarM},ST(S) \rangle_e=-\ID(\fvar)-\ID(\svar)(=\ID(\svar\tvar)-\ID(\svar)). $$
\end{lemma}
\eop
\begin{lemma}
\label{lemintca}
Let $C_{\fvarM}$ be a $4$-dimensional chain in $\partial \TCM$ whose boundary does not meet $$\ZZ \times \left(s_{\tau}(K_{\svarM};\tvarM) \cup s_{\tau}(K_{\tvarM};\svarM) \cup s_{TK}(K) \cup s_{-TK}(K)\right)$$
$$\langle C_{\fvarM},F_{\svarM},F_{\tvarM} \rangle_e=-\ID(\svar)\langle C_{\fvarM},s_{\tau}(K_{\svarM};\tvarM)\rangle_{e,\ZZ \times ST(M)} - \ID(\tvar)\langle C_{\fvarM},s_{\tau}(K_{\tvarM};\svarM)\rangle_{e,\ZZ \times ST(M)}$$
where the variable $t$ in the equivariant intersection numbers of the right-hand side should be replaced by $\fvar$.
$$\langle C_{\fvarM},A(K) \rangle_e= \langle C_{\fvarM},s_{TK}(K) \rangle_{e,\ZZ \times ST(M)} -t\langle C_{\fvarM},s_{-TK}(K) \rangle_{e,\ZZ \times ST(M)}$$
\end{lemma}
\bp
Recall that $\partial A(K)=s_{TK}(K)-t^{-1}s_{-TK}(K)$ from Lemma~\ref{lemak}.
\eop

In the above statement, the roles of $\fvarM$, $\svarM$ and $\tvarM$ can be permuted to give similar expressions for 
$$\langle F_{\fvarM},C_{\svarM},F_{\tvarM} \rangle_e=-\ID(\fvar)\langle C_{\svarM},s_{\tau}(K_{\fvarM};\tvarM)\rangle_{e,\partial \TCM} - \ID(\tvar)\langle C_{\svarM},s_{\tau}(K_{\tvarM};\fvarM)\rangle_{e,\partial \TCM}$$ and
$$\langle F_{\fvarM},F_{\svarM},C_{\tvarM} \rangle_e=-\ID(\fvar)\langle C_{\tvarM},s_{\tau}(K_{\fvarM};\svarM)\rangle_{e,\partial \TCM} - \ID(\svar)\langle C_{\tvarM},s_{\tau}(K_{\svarM};\fvarM)\rangle_{e,\partial \TCM}$$
since the triple intersection of codimension $2$ chains does not depend on the order other than for assignments of the variables $\fvar$, $\svar$ and $\tvar$.

\begin{lemma}
\label{lemvarsurf}
Let $\Sigma_{\fvarM}$ be a rational chain in $M$ whose boundary $(K^{\prime}_{\fvarM}-K_{\fvarM})$ is
disjoint from $(K_{\svarM},K_{\tvarM},K)$.
Let
 $$F^{\prime}_{\fvarM}=F_{\fvarM} -\ID(t) ST(\Sigma_{\fvarM})+\ID(t)\langle \Sigma_{\fvarM},K\rangle_M ST(S).$$
Then 
$$\partial F^{\prime}_{\fvarM} = s_{\tau}(M;\fvarM)-\ID K^{\prime}_{\fvarM} \times S^2\;\;\mbox{and} \;\;\langle F^{\prime}_{\fvarM}, A(K)\rangle_e=0.$$
$$\langle F^{\prime}_{\fvarM},F_{\svarM},F_{\tvarM} \rangle_e -\langle F_{\fvarM},F_{\svarM},F_{\tvarM} \rangle_e=$$
$$ \ID(\fvar)\langle\Sigma_{\fvarM},K\rangle_M(-\ID(\svar) - \ID(\tvar))+\ID(\fvar)(\ID(\svar)\langle\Sigma_{\fvarM},K_{\svarM}\rangle_M + \ID(\tvar)\langle\Sigma_{\fvarM},K_{\tvarM}\rangle_M).$$
\end{lemma}
\bp
Recall that $\langle ST(S), A(K)\rangle_e=1-t$.
\eop

\subsection{Changing the knot trivialisation}
\label{subknotpar}
Consider the map $$\begin{array}{llll} g \colon & S^1&\rightarrow &SO(3)\\
&\exp(i \alpha)& \mapsto &\rho(\qvarM, \alpha + \pi)\end{array}$$
where $\rho(\qvarM, \alpha)$ is the rotation with oriented axis $\qvarM$ and with angle $\alpha$. The homotopy class of $g$ generates $\pi_1(SO(3))$. Also consider the map $f$ from $M$ to $S^1$ whose restriction on $K$ has degree one, such that $f(M\setminus (S\times ]-1,1[))=-1$ and $f(\sigma \in S,\beta \in]-1,1[)= \exp(i \pi \beta)$.
This induces $$\begin{array}{llll}h\colon &M \times\RR^3 &\rightarrow& M \times\RR^3\\
&(m,v)& \mapsto &(m,(g\circ f (m))^{-1}(v)).\end{array}$$
Let $\tau^{\prime}=h\circ \tau \colon TM \rightarrow M \times \RR^3$.
The parallel $K^{\prime}_{\parallel}$ induced by $\tau^{\prime}$ is obtained from the parallel induced by $K$ by adding a positive meridian of $K$.
Note that $\tau$ and $\tau^{\prime}$ have the same Pontrjagin class since the change of parallelisations is associated with a map from $S^1$ to $SO(3)$ and $\pi_1(SU(4))$ is trivial.

In this subsection, we shall prove the following proposition.
\label{varpar}
\begin{proposition}
\label{propvarpar}
Under the above hypotheses, $$\CQ(K,K^{\prime}_{ \parallel})-\CQ(K,K_{\parallel})
=\CQ(K,K^{\prime}_{ \parallel},\tau^{\prime})-\CQ(K,K_{\parallel},\tau)
=\sum_{\mathfrak{S}_3(\fvar,\svar,\tvar)}\left(-\frac{\fvar \Delta^{\prime}(\fvar)}{2\Delta(\fvar)}\right)\ID(\svar)$$
\end{proposition}

\begin{remark}
The triple $(S^1\times S^2,S^1 \times\{\qvarM\},S^1 \times\{\qvarM^{\prime}\})$ is equivalent to the triple $(S^1\times S^2,S^1 \times\{\qvarM\},\{(\exp(i\beta)\in S^1,R_{\beta,\qvarM}(\qvarM^{\prime} ))\})$.

This is consistent with Proposition~\ref{propvarpar} that ensures that $\CQ$ does not vary when $\Delta=1$.
\end{remark}
$$\CQ(K,K_{\parallel},\tau) =\CQ(K_{\fvarM},K_{\svarM},K_{\tvarM},K,K_{\parallel},\tau).$$
To prove the proposition, we compute $$\delta_1=\CQ(K_{\fvarM},K_{\svarM},K_{\tvarM},K,K^{\prime}_{ \parallel},\tau^{\prime}) - \CQ(K,K_{\parallel},\tau)$$
in Lemma~\ref{lemdeltaone} and
$$\delta_2=\CQ(K,K^{\prime}_{\parallel},\tau^{\prime}) - \CQ(K_{\fvarM},K_{\svarM},K_{\tvarM},K,K^{\prime}_{ \parallel},\tau^{\prime})$$ in Lemma~\ref{lemdeltatwo}.

\begin{lemma}
\label{lemdeltaone}
Let $\delta_1=\CQ(K_{\fvarM},K_{\svarM},K_{\tvarM},K,K^{\prime}_{ \parallel},\tau^{\prime}) - \CQ(K,K_{\parallel},\tau)$
$$\delta_1=\sum_{\mathfrak{S}_3(\fvar,\svar,\tvar)}\frac{1+\fvar}{2(1-\fvar)}I_{\Delta}(\svar).$$
\end{lemma}
\bp We first change $F_{\fvarM}$ to $F_{\fvarM} + C_{\fvarM}$ where $\partial C_{\fvarM}=s_{\tau^{\prime}}(M;\fvarM) -s_{\tau}(M;\fvarM)$. To construct such a $C_{\fvarM}$,
define $$\begin{array}{llll}g_{\fvarM}  \colon &[0,1] \times S^1&\rightarrow &S^2\\
&(t,z)& \mapsto &g_{\fvarM}(t,z)\end{array}$$
where $g_{\fvarM}(0,z)=g_{\fvarM}(t,-1)=\fvarM$ and $g_{\fvarM}(1,z)=g(z)(\fvarM )$ and the image of $g_{\fvarM}$ is the disk $D_{\fvarM}(\qvarM)$ of $S^2$ bounded by the circle with center $\qvarM$ through $\fvarM$, $g_{\fvarM}$ has degree one at $\qvarM$.
Then the cobordism $$\begin{array}{llll}C_{\fvarM}\colon&[0,1] \times M &\rightarrow &M \times S^2\\
&(t,m) & \mapsto & (m,g_{\fvarM}(t,f(m)))\end{array}$$
in $ST(M)=_{\tau}M \times S^2$ satisfies
$$\partial C_{\fvarM}=s_{\tau^{\prime}}(M;\fvarM) -s_{\tau}(M;\fvarM),$$ and, according to Lemma~\ref{lemak},
$$\langle C_{\fvarM}, A(K) \rangle_e=\langle C_{\fvarM}, s_{\tau}(K;\qvarM) \rangle_{e,\partial \TCM} -t\langle C_{\fvarM}, s_{\tau}(K;-\qvarM) \rangle_{e,\partial \TCM}=\langle C_{\fvarM}, s_{\tau}(K;\qvarM) \rangle_{e,\partial \TCM}$$ since the image of $g_{\fvarM}$ does not meet $-\qvarM$.
Assume without loss that $K$ meets $S\times[-1,1]$ as $\{\sigma\} \times[-1,1]$. Then
$\langle C_{\fvarM}, s_{\tau}(K;\qvarM) \rangle_{e,\partial \TCM} $ is the degree at $\qvarM$ of the restriction of $C_{\fvarM}$ to $[0,1] \times \{\sigma\} \times[-1,1]$ that is $1$.

Let 
$$F^{\prime}_{\fvarM}= F_{\fvarM} + C_{\fvarM} -\frac{1}{1-t} ST(S).$$
Then $\langle F^{\prime}_{\fvarM}, A(K) \rangle=0$.

Assume that $\fvarM$, $\svarM$ and $\tvarM$ lie on different circles around $\qvarM$ so that $D_{\fvarM}(\qvarM) \subsetneq D_{\svarM}(\qvarM) \subsetneq D_{\tvarM}(\qvarM)$ and define $F^{\prime}_{\svarM}$ and $F^{\prime}_{\tvarM}$, similarly.
Then $$\begin{array}{ll}\delta_1&=\langle F^{\prime}_{\fvarM}, F^{\prime}_{\svarM},F^{\prime}_{\tvarM} \rangle_e -\langle F_{\fvarM}, F_{\svarM},F_{\tvarM} \rangle_e\\&=\langle F_{\fvarM}, F_{\svarM},C_{\tvarM} \rangle_e + \langle F_{\fvarM}, C_{\svarM},F^{\prime}_{\tvarM} \rangle_e +\langle C_{\fvarM}, F^{\prime}_{\svarM},F^{\prime}_{\tvarM} \rangle_e +\delta^{\prime}_1\end{array}$$
where $\delta^{\prime}_1 = \sum_{\circlearrowleft} \frac{1}{1-\fvar}(I_{\Delta}(\svar) + I_{\Delta}(\tvar))$, according to Lemma~\ref{leminsts}.
The symbol $\sum_{\circlearrowleft}$ stands for the sum of the three terms obtained from the written one by permuting $\fvarM$, $\svarM$ and $\tvarM$, cyclically.

Now, according to Lemma~\ref{lemintca}, $$\langle C_{\fvarM}, F^{\prime}_{\svarM},F^{\prime}_{\tvarM}\rangle_e=
-\ID(\svar)\langle C_{\fvarM},s_{\tau^{\prime}}(K_{\svarM};\tvarM)\rangle_{e,\ZZ \times ST(M)} - \ID(\tvar)\langle C_{\fvarM},s_{\tau^{\prime}}(K_{\tvarM};\svarM)\rangle_{e,\ZZ \times ST(M)}=0$$
since the image of $g_{\fvarM}$ does not meet the circles centered at $\qvarM$ through $\svarM$ and $\tvarM$.
$$\begin{array}{ll}\langle F_{\fvarM}, F_{\svarM},C_{\tvarM} \rangle_e&=-\ID(\fvar)\langle C_{\tvarM},s_{\tau}(K_{\fvarM};\svarM)\rangle_{e,\ZZ \times ST(M)} - \ID(\svar)\langle C_{\tvarM},s_{\tau}(K_{\svarM};\fvarM)\rangle_{e,\ZZ \times ST(M)}\\&=-\ID(\fvar)- \ID(\svar)\end{array}$$
since $\langle C_{\tvarM},s_{\tau}(K_{\fvarM};\svarM)\rangle_{e,\ZZ \times ST(M)}$ and  $\langle C_{\tvarM},s_{\tau}(K_{\svarM};\fvarM)\rangle_{e,\ZZ \times ST(M)}$ are as before the degrees of $g_{\tvarM}$, at $\svarM$ and $\fvarM$, that are one.
Similarly,
$$\langle F_{\fvarM}, C_{\svarM},F^{\prime}_{\tvarM} \rangle_e=-\ID(\tvar)\langle C_{\svarM},s_{\tau}(K_{\tvarM};\fvarM)\rangle_{e,\ZZ \times ST(M)}=-\ID(\tvar).$$
$$\delta_1=\sum_{\circlearrowleft}\left( \frac{1}{1-\fvar} -\frac{1}{2}\right)(I_{\Delta}(\svar) + I_{\Delta}(\tvar)).$$
\eop

\begin{lemma}
\label{lemdeltatwo}
Let $\delta_2=\CQ(K,K^{\prime}_{\parallel},\tau^{\prime}) - \CQ(K_{\fvarM},K_{\svarM},K_{\tvarM},K,K^{\prime}_{ \parallel},\tau^{\prime})$
$$\delta_2=-\sum_{\mathfrak{S}_3(\fvar,\svar,\tvar)}\frac{1}{2}I_{\Delta}(\fvar)I_{\Delta}(\svar).$$
\end{lemma}
\bp
Let $D^2 \times K$ be a tubular neighborhood of $K$ trivialised by $\tau_{|K}$ where $D^2$ is the unit disk of $\CC$. Assume $K_{\fvarM}=1/4 \times K$,
$K_{\svarM}=1/2 \times K$ and $K_{\tvarM}=3/4 \times K$.
Then let $K^{\prime}_{\fvarM} \subset 1/4S^1 \times K$ be a parallel of $K$ such that
$(K^{\prime}_{\fvarM}-K_{\fvarM})$ is homologous to a positive meridian of $K$ in $D^2 \times K$, and let $\Sigma_{\fvarM}$ be an annulus transverse to $K$ in $1/4D^2 \times K$ whose boundary is $(K^{\prime}_{\fvarM}-K_{\fvarM})$. Then 
$\langle \Sigma_{\fvarM},K\rangle_M =1$.
Similarly define $(\Sigma_{\svarM},K^{\prime}_{\svarM})$ and $(\Sigma_{\tvarM},K^{\prime}_{\tvarM})$ by replacing $1/4$ by $1/2$ and $3/4$, respectively.
$$\CQ(K,K^{\prime}_{ \parallel},\tau^{\prime})=\CQ(K^{\prime}_{\fvarM},K^{\prime}_{\svarM},K^{\prime}_{\tvarM},K,K^{\prime}_{\parallel},\tau^{\prime})$$
Use Lemma~\ref{lemvarsurf} and define 
$$F^{\prime \prime}_{\fvarM}=F^{\prime}_{\fvarM} -\ID(t) ST(\Sigma_{\fvarM})+\ID(t)\langle \Sigma_{\fvarM},K\rangle_M ST(S).$$
and its twin brothers $F^{\prime \prime}_{\svarM}$, $F^{\prime \prime}_{\tvarM}$.
Then $\delta_2=\sum_{\circlearrowleft}I_{\Delta}(\fvar)(-I_{\Delta}(\svar) - I_{\Delta}(\tvar)) +\delta^{\prime}_2$
where 
$$\begin{array}{lll}\delta^{\prime}_2&=& \ID(\fvar)(\ID(\svar)\langle\Sigma_{\fvarM},K_{\svarM}\rangle_M + \ID(\tvar)\langle\Sigma_{\fvarM},K_{\tvarM}\rangle_M)
\\&&+\ID(\svar)(\ID(\fvar)\langle\Sigma_{\svarM},K^{\prime}_{\fvarM}\rangle_M + \ID(\tvar)\langle\Sigma_{\svarM},K_{\tvarM}\rangle_M)
\\&&+\ID(\tvar)(\ID(\fvar)\langle\Sigma_{\tvarM},K^{\prime}_{\fvarM}\rangle_M + \ID(\svar)\langle\Sigma_{\tvarM},K^{\prime}_{\svarM}\rangle_M)
\\&=& \ID(\svar)\ID(\fvar)+\ID(\tvar)(\ID(\fvar)+ \ID(\svar)).\end{array}$$
\eop

\begin{remark}
 Proposition~\ref{propvarpar} could also be seen as a consequence of the Dehn surgery formula (Theorem~\ref{thmDehn} that will be proved independently later) by seeing the variation of $\CQ$ under the addition of a meridian to the parallel of $K$ as the result of a $(-1)$--surgery along the meridian of $K$. Then the Dehn surgery formula could be applied with the Seifert surface $S\setminus \mbox{Int}(D^2)$. 
I checked that this gives the same result when the genus of $S$ is one, but the above proof of Proposition~\ref{propvarpar} looks simpler even for the genus one case.
\end{remark}

\newpage
\section{Variation of $\CQ(\KK)$ under a general two--dimensional cobordism}
\setcounter{equation}{0}
\label{secvarcob}

This section is devoted to the proof of Theorem~\ref{thmfrakcha}.

\subsection{Beginning the proof of Theorem~\ref{thmfrakcha}}

In this subsection and the next one, we prove the first part of this theorem that is the following proposition, thanks to Proposition~\ref{propsym3}.

\begin{proposition}
\label{propvargen}
Let $\KK=(K,K_{\parallel})$ and $\KK^{\prime}=(K^{\prime},K^{\prime}_{\parallel})$ be two framed knots that represent the preferred generator of $H_1(M)/\mbox{Torsion}$. Then there exists an antisymmetric polynomial ${\cal V}(\KK,\KK^{\prime})$ in $\QQ[t,t^{-1}]$ such that 
$$\CQ(\KK^{\prime}) - \CQ(\KK)=\sum_{\mathfrak{S}_3(\fvar,\svar,\tvar)}\frac{{\cal V}(\KK,\KK^{\prime})(\fvar)}{\delta(M)(\fvar)}I_{\Delta}(\svar).$$
\end{proposition}

Let $K$ and $K^{\prime}$ be two knots that represent the preferred generator of $H_1(M)/\mbox{Torsion}$, and let $\tau$ be a trivialisation of $TM$ that maps the unit tangent vectors of $K$ and $K^{\prime}$ that induce the orientation of $K$ and $K^{\prime}$
to $\qvarM \in S^2$. Let $K_{\parallel}$ and $K^{\prime}_{\parallel}$
be the parallels of $K$ and $K^{\prime}$ induced by $\tau$. (The parallels can be chosen arbitrarily, thanks to Proposition~\ref{propvarpar}.)
Let $\KK=(K,K_{\parallel})$ and $\KK^{\prime}=(K^{\prime},K^{\prime}_{\parallel})$.

This subsection is devoted to the computation of 
$$\CQ(\KK^{\prime}) - \CQ(\KK)=\CQ(K^{\prime},K^{\prime}_{ \parallel},\tau)-\CQ(K,K_{\parallel},\tau).$$

There exists a rational $2$--chain $\Bor$ such that $$\partial \Bor =K^{\prime}-K.$$
Let $N(K)$ and $N(K^{\prime})$ be tubular neighborhoods of $K$ and $K^{\prime}$.
Then $$[\partial (\Bor \cap \left(M \setminus \mbox{Int}(N(K) \sqcup N(K^{\prime}))\right)] = [K^{\prime}_{\parallel} + r^{\prime}m(K^{\prime})] - [K_{\parallel} + rm(K)]$$
where $m(K)$ and $m(K^{\prime})$ denote meridians of $K$ and $K^{\prime}$, respectively and the brackets stand for homology classes in $\partial(N(K) \sqcup N(K^{\prime}))$.

Let $C(\Bor;\tau)$ be the following $2$--cycle of $\TCM$, $$C(\Bor;\tau)=A(K^{\prime}) - A(K) +t^{-1}s_{\tau}(\Bor;-\qvarM) -s_{\tau}(\Bor;\qvarM)$$ where $A(K)$ is defined in Lemma~\ref{lemak}.
Define $Q(\Bor;\tau) \in \QQ(t)$ so that $$[C(\Bor;\tau)]=(1-t^{-1})Q(\Bor;\tau)(t^{-1})[ST(\ast)]$$
in $H_2(C_2(M);\QQ(t))$.

\begin{proposition}
\label{propvarone} Under the assumptions above,
$$\CQ(\KK^{\prime}) - \CQ(\KK)=\sum_{\mathfrak{S}_3(\fvar,\svar,\tvar)}\left(Q(\Bor;\tau)(\fvar)+\frac{r^{\prime}-r}{2}I_{\Delta}(\fvar)\right)I_{\Delta}(\svar).$$
\end{proposition}
\bp
Observe that $\langle \Bor, K_{\fvarM} \rangle=\langle \Bor, K_{\svarM} \rangle=\langle \Bor, K_{\tvarM} \rangle =r$ and that $\langle \Bor, K^{\prime}_{\fvarM} \rangle=\langle \Bor, K^{\prime}_{\svarM} \rangle=\langle \Bor, K^{\prime}_{\tvarM} \rangle =-r^{\prime}$. Recall $\partial A(K)=s_{TK}(K)-t^{-1}s_{-TK}(K)$ and 
$\langle ST(S), A(K)\rangle_e=1-t$.
Let $\Sigma_{\fvarM}$ be a $2$--chain obtained from $\Bor$ by a natural small isotopy such that $\partial \Sigma_{\fvarM} =K_{\fvarM}^{\prime}-K_{\fvarM}$. Then $\langle \Sigma_{\fvarM}, K \rangle=r$ and $\langle \Sigma_{\fvarM}, K^{\prime} \rangle=-r^{\prime}$.

Let $F^1_{\fvarM}=F_{\fvarM}-I_{\Delta}(t)ST(\Sigma_{\fvarM})$ where $\langle A(K), F_{\fvarM} \rangle_e= 0$.
Then $$\begin{array}{lll}\langle A(K^{\prime}), F^1_{\fvarM} \rangle&=&\langle A(K^{\prime}), F_{\fvarM} \rangle_e + I_{\Delta}(t)\langle A(K^{\prime}), ST(\Sigma_{\fvarM}) \rangle_e \\
&=&\langle A(K^{\prime}), F_{\fvarM} \rangle_e +I_{\Delta}(t)(1-t^{-1})\langle \Sigma_{\fvarM}, K^{\prime} \rangle
\\&=&\langle C(\Bor;\tau)+ A(K) -t^{-1}s_{\tau}(\Bor;-\qvarM) +s_{\tau}(\Bor;\qvarM), F_{\fvarM} \rangle_e -r^{\prime}I_{\Delta}(t)(1-t^{-1})
\\&=&(1-t^{-1})\left(Q(\Bor;\tau)(t^{-1})-r^{\prime}I_{\Delta}(t)
+I_{\Delta}(t) \langle \Bor, K_{\fvarM} \rangle\right)
\\&=&(1-t^{-1})\left(Q(\Bor;\tau)(t^{-1})-(r-r^{\prime})I_{\Delta}(t^{-1})
\right).
\end{array}$$

Then $$F^{\prime}_{\fvarM}= F^1_{\fvarM} +\left((r-r^{\prime})I_{\Delta}(t)-Q(\Bor;\tau)(t)\right)ST(S) $$
satisfies
$\langle A(K^{\prime}),F^{\prime}_{\fvarM}\rangle =0$.
Set $Q(t)=Q(\Bor;\tau)(t)+(r^{\prime}-r)I_{\Delta}(t)$.
Define $\Sigma_{\svarM}$, $\Sigma_{\tvarM}$,
$F^{1}_{\svarM}=F_{\svarM}-I_{\Delta}(t)ST(\Sigma_{\svarM})$
and $F^{1}_{\tvarM}=F_{\tvarM}-I_{\Delta}(t)ST(\Sigma_{\tvarM})$, similarly. Then
$$\begin{array}{lll}\CQ(\KK^{\prime}) - \CQ(\KK)&=&\langle F^{\prime}_{\fvarM},F^{1}_{\svarM}-Q(t)ST(S),F^{1}_{\tvarM}-Q(t)ST(S) \rangle_e -\langle F_{\fvarM},F_{\svarM},F_{\tvarM} \rangle_e\\
&=& \sum_{\circlearrowleft}Q(\fvar)(I_{\Delta}(\svar)+I_{\Delta}(\tvar)) +\delta_3=\sum_{\mathfrak{S}_3(\fvar,\svar,\tvar)}Q(\fvar)I_{\Delta}(\svar)+\delta_3
  \end{array}$$
according to Lemma~\ref{leminsts}, where $$\begin{array}{lll}\delta_3&=&-I_{\Delta}(\fvar)\langle ST(\Sigma_{\fvarM}),F^{1}_{\svarM},F^{1}_{\tvarM} \rangle_e
-I_{\Delta}(\svar)\langle F_{\fvarM},ST(\Sigma_{\svarM}),F^{1}_{\tvarM} \rangle_e
-I_{\Delta}(\tvar)\langle F_{\fvarM},F_{\svarM},ST(\Sigma_{\tvarM}) \rangle_e
\\&=&I_{\Delta}(\fvar)\left(I_{\Delta}(\svar)\langle \Sigma_{\fvarM}, K^{\prime}_{\svarM} \rangle +I_{\Delta}(\tvar)\langle \Sigma_{\fvarM}, K^{\prime}_{\tvarM} \rangle\right)\\
&&+I_{\Delta}(\svar)\left(I_{\Delta}(\fvar)\langle \Sigma_{\svarM}, K_{\fvarM} \rangle +I_{\Delta}(\tvar)\langle \Sigma_{\svarM}, K^{\prime}_{\tvarM} \rangle\right)
\\
&&+I_{\Delta}(\tvar)\left(I_{\Delta}(\fvar)\langle \Sigma_{\tvarM}, K_{\fvarM} \rangle +I_{\Delta}(\svar)\langle \Sigma_{\tvarM}, K_{\svarM} \rangle\right)
\\&=& - r^{\prime}I_{\Delta}(\fvar)(I_{\Delta}(\svar)+I_{\Delta}(\tvar))
 + r I_{\Delta}(\svar)I_{\Delta}(\fvar)- r^{\prime}I_{\Delta}(\svar)I_{\Delta}(\tvar)
+r I_{\Delta}(\tvar)(I_{\Delta}(\fvar)+I_{\Delta}(\svar))
\\&=&\frac{r-r^{\prime}}{2}\sum_{\mathfrak{S}_3(\fvar,\svar,\tvar)}I_{\Delta}(\fvar)I_{\Delta}(\svar).
\end{array}$$
\eop

\begin{remark}
When $\Bor$ is changed to $\Bor + \alpha S$, then
$r^{\prime}=r^{\prime}(\Bor;\tau)$ is changed to $(r^{\prime}(\Bor + \alpha S ;\tau)=r^{\prime} -\alpha)$, $r$ is changed to $(r+\alpha)$ and 
$Q(\Bor + \alpha S ;\tau)=Q(\Bor;\tau) + \alpha I_{\Delta}(t)$, according to Theorem~\ref{thmstauS},
so that
$Q(\Bor;\tau)(t)+\frac{r^{\prime}-r}{2}I_{\Delta}(t)$ is invariant
as it must be.
\end{remark}

\begin{lemma}
\label{lemdegstwoso}
Let $\qvarM \in S^2$.
 Let $\Phi$ be a map from the unit disk $D^2$ of $\CC$ to $SO(3)$ such that $\Phi(\exp(i\beta))$ is the rotation $R_{2\beta,\qvarM}$ with axis directed by $\qvarM$ and with angle $2\beta$.
Then the map $$\begin{array}{llll}\Phi_{\qvarM}=\Phi(\cdot)(\qvarM) \colon &D^2  &\rightarrow & S^2\\
  &z&\mapsto &\Phi(z)(\qvarM)
 \end{array}$$
sends $\partial D^2$ to $\qvarM$ and the degree of the induced map
from $D^2/\partial D^2$ to $S^2$ is $(-1)$.
\end{lemma}
\bp
First note that the above degree does not depend on $\Phi$ on the interior of $D^2$, since $\pi_2(SO(3))=0$.
See the restriction of $\Phi_{\qvarM}$ to $\partial D^2$ as the path composition of the maps $\left(\beta \in [0,2\pi] \mapsto  R_{\beta,\qvarM} \right)$ and the inverse of $\left(\beta \in [0,2\pi] \mapsto R_{\beta,-\qvarM} \right)$ (that is twice the first map).
Consider an arc $\alpha$ of a great circle of $S^2$ from $-\qvarM$ to $\qvarM$, then $\Phi$ can be seen as the map from $\alpha \times [0,2\pi]$ to $SO(3)$ that maps $(\cvarM,\beta)$ to $R_{\beta,\cvarM}$ so that
the only preimage of $-\qvarM$ under $\Phi_{\qvarM}$ reads $(\cvarM_0,\pi)$ where $\cvarM_0 \perp \qvarM$ and the local degree is easily seen to be $(-1)$.
\eop

\begin{remark}
\label{rkhomso3}
Consider our framed knot $(K,K_{\parallel})$ in $M$ that generates $H_1(M;\ZZ)/\mbox{Torsion}$. Let $K^{\prime\prime}_{ \parallel}$ be the parallel of $K$ obtained from $K_{\parallel}$ by adding two positive meridians.
According to Proposition~\ref{propvarpar}, $$\CQ(K,K^{\prime\prime}_{\parallel})-\CQ(K,K_{ \parallel})
=\sum_{\mathfrak{S}_3(\fvar,\svar,\tvar)}\left(-\frac{\fvar \Delta^{\prime}(\fvar)}{\Delta(\fvar)}\right)\ID(\svar)$$

We can see this fact as a consequence of Proposition~\ref{propvarone} as follows.
Assume that $\tau$ has the following natural properties on $K \times D^2$ (trivialised with respect to $\tau$), it maps $TK \times \{z\}$ to $\RR \qvarM$ for $z \in D^2$, and it does not depend on $k \in K$ on $\{k\} \times TD^2$. Perform the following radial
homotopy of $\tau$ that changes $\tau$ into a trivialisation $\tau^{\prime}$ that coincides with $\tau$ outside $K \times D^2$, and that induces the same parallelisation as $ (K_{\parallel} +2m(K))$.
$$\tau \circ \tau^{\prime-1}(k,z)(v)=\nu(k,|z|)(v)$$ where $\nu(k,|z|) \in SO(3)$, $\nu(k,1) =\mbox{Id}$, $\nu(\exp(i\beta),0) =R_{2\beta,\qvarM}$.
Set $K=K\times\{0\}$, $K^{\prime}=K\times\{1\}$ and $B=-K\times [0,1]$.
Then $r(\tau^{\prime})=r-2$ and 
$s_{\tau^{\prime}}(\Bor;\qvarM)$ and $s_{\tau}(\Bor;\qvarM)$ coincide outside $K \times D^2$ and on $K$. The difference $[s_{\tau^{\prime}}(\Bor;\qvarM)-s_{\tau}(\Bor;\qvarM)]$ reads
$\degr  [ST(\ast)]$ where $\degr $ is the degree at $(-\qvarM)$ of $\tau \circ \tau^{\prime-1}(k,z)(\qvarM)$ restricted to $\Bor$. Then $\degr $ is the degree at $(-\qvarM)$ of the map
$$\begin{array}{lll}[0,1] \times K  &\rightarrow & S^2\\
  (t,k)&\mapsto &\nu(k,t)(\qvarM)
 \end{array}$$
that is $1$ according to Lemma~\ref{lemdegstwoso}.
Then $$[s_{\tau^{\prime}}(\Bor;\qvarM) - s_{\tau}(\Bor;\qvarM)]=[ST(\ast)].$$

Since $s_{\tau}(\Bor;-\qvarM)$ and $s_{\tau^{\prime}}(\Bor;-\qvarM)$ are obtained from $s_{\tau}(\Bor;\qvarM)$ and $s_{\tau^{\prime}}(\Bor;\qvarM)$
by the involution $\iota$, and since $[\iota(ST(\ast))]=-[ST(\ast)]$, 
we find that 
$[C(\Bor;\tau^{\prime})]-[C(\Bor;\tau)]= -(1+t^{-1})[ST(\ast)]$
so that
$Q(\Bor;\tau^{\prime})-Q(\Bor;\tau)=-\frac{1+t}{1-t}$.
In particular, 
$$\left(\CQ(\KK^{\prime},\tau^{\prime}) - \CQ(K,K^{\prime\prime}_{\parallel},\tau^{\prime})\right) -\left(\CQ(\KK^{\prime},\tau) - \CQ(\KK,\tau)\right)=\sum_{\mathfrak{S}_3(\fvar,\svar,\tvar)}\left(-\frac{1+\fvar}{1-\fvar}+I_{\Delta}(\fvar)\right)I_{\Delta}(\svar).$$
This is consistent with the expression of $\left(\CQ(\KK,\tau)-\CQ(K,K^{\prime\prime}_{\parallel},\tau^{\prime})\right)$ coming from Proposition~\ref{propvarpar}.
\end{remark}

The computation of $\left(\CQ(\KK^{\prime}) - \CQ(\KK)\right)$ is now reduced to the computation of $Q(\Bor;\tau)$ that only depends on the parallelisations of $K$ and $K^{\prime}$ and on the
homology class of $\Bor$.

\subsection{Computation of $Q(\Bor;\tau)$}

Without loss, we now assume that $\Bor$ induces the framing of $\KK$ i.e. that $r=0$. (It suffices to change $\Bor$ to $\Bor-rS$.)

See $K^{\prime}$ as a band sum of $K$ and of a rationally null-homologous knot. More precisely, write $K$ as the union of two oriented intervals $I=[\alpha,\beta]$ and $I_2$ glued along their boundaries 
$$K= I \cup_{\partial I_2} I_2 \;\;\;\mbox{and} \;\;\;K^{\prime}= I^{\prime} \cup_{\partial I_2} I_2$$
Then $$\partial \Bor =I^{\prime} \cup_{\partial I} (-I)$$
where the gluing is not smooth at $\partial I$, it has two cusps there.
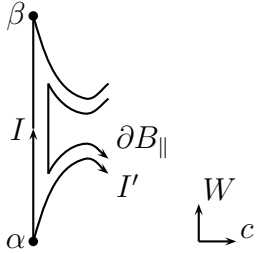
\begin{center}
\begin{pspicture}[shift=-0.1](-1,-1.6)(4,1.6)
\psline{*->}(0,-1.5)(0,0)
\psline{-*}(0,0)(0,1.5)
\rput[r](-.1,0){$I$}
\rput[r](-.1,-1.5){$\alpha$}
\rput[r](-.1,1.5){$\beta$}
\psecurve{->}(0,-2)(0,-1.5)(.7,-.4)(1,-.6)(1.5,-1.1)
\psecurve{-}(0,2)(0,1.5)(.7,.4)(1,.6)(1.5,1.1)
\rput[tl](1.1,-.6){$I^{\prime}$}
\psline(.2,-.6)(.2,.6)
\psecurve{->}(.2,-1.5)(.2,-.6)(.7,-.2)(1,-.4)(1.5,-.9)
\psecurve{-}(.2,1.5)(.2,.6)(.7,.2)(1,.4)(1.5,.9)
\rput[bl](1.1,-.4){$\partial \Bor_{\parallel}$}
\psline{->}(2.2,-1.5)(2.7,-1.5)
\psline{->}(2.2,-1.5)(2.2,-1)
\rput[bl](2.25,-.95){$\qvarM$}
\rput[bl](2.75,-1.45){$c$}
\end{pspicture}
\end{center}

Near $K$, $\Bor$ reads $[0,1[ \times I $ and its unit tangent space is mapped to a circle of $S^2$ that contains $\qvarM$ and $(-\qvarM)$, and that is made of two half-great circles that we shall denote by  $[-\qvarM,\qvarM]_{c}$ and $[-\qvarM,\qvarM]_{-c}$ when they are oriented from $-\qvarM$ to $\qvarM$. The half-circle $[-\qvarM,\qvarM]_{c}$ contains the direction $c$ of the inward normal $[0,1[$ of $\Bor$, while the half-circle $[-\qvarM,\qvarM]_{-c}$ contains the direction of the outward normal.

When $a$ and $b$ are real numbers such that $a<b$,
$[a,b]\times_{\leq}[a,b]=\{(t,u)\in[a,b]^2;t\leq u \}$ and $[a,b]\times_{\geq}[a,b]=\{(t,u)\in[a,b]^2;t\geq u \}$.
We shall use the notation $I \tilde{\times}_{\leq} I$ for the 
lift in $\TCM$ of the closure of $\{(t,u)\in[\alpha,\beta]^2;t < u \}$
in $C_2(M)$
that contains points of the blow-up of the preferred lift of the diagonal. The chains $(I^{\prime}\cup_{\beta}(-I))\tilde{\times}_{\leq}(I^{\prime}\cup_{\beta}(-I))$ and $(I^{\prime}\cup_{\beta}(-I))\tilde{\times}_{\geq}(I^{\prime}\cup_{\beta}(-I))$ should be understood similarly, they contain $\pm s_{\tau}(I^{\prime};\qvarM)$ and $\pm s_{\tau}(I^{\prime};-\qvarM)$,
respectively. Furthermore, $(I^{\prime}\cup_{\beta}(-I))\tilde{\times}_{\leq}(I^{\prime}\cup_{\beta}(-I))$ contains the (lift of the) closure of the set of points $(t,u)$ when $t$ approaches
$\alpha$ in $I^{\prime}$ and $u$ approaches $\alpha$ in $I$. Hence, it contains $\pm s_{\tau}(\alpha;[-\qvarM,\qvarM]_{-c})$. It similarly contains $\pm s_{\tau}(\beta;[-\qvarM,\qvarM]_{-c})$.

We want to get rid of $I_2$ by letting $\partial \Bor$ bound.
In general, $\Bor$ does not lift in $\tilde{M}$. However, $\delta(t_M)\partial \Bor$ bounds a rational chain denoted by $\delta(t_M)\Borp$ in $\tilde{M}$ that reads $\delta(t) [0,1[ \times I$ near the lifts of $K$. (Again, it suffices to add a combination of $\alpha_kt^k S$ to achieve this.) Let $\tilde{\alpha}$ be the preimage of $\alpha$ in $\tilde{M}$ that reads $\{0\} \times \alpha$ above. Then $\alpha \times \Borp$ denotes the closure in $\TCM$ of the projection in $\widetilde{M^2}$ of $\alpha \times (\Borp \setminus\tilde{\alpha}) $.
The chain $\Borp \times \alpha$ is also the most natural chain in $\TCM$ that fits with the notation, and that has the symmetric heavy definition.

\begin{lemma}
\label{lempremc2}
 $$[C(\Bor;\tau)]= t^{-1}[C_1(\Borp;\tau)]-[C_0(\Borp;\tau)]$$
in $H_2(C_2(M);\QQ(t))$, where
$$C_0(\Borp;\tau)=s_{\tau}(\Bor;\qvarM)-s_{\tau}(I;[-\qvarM,\qvarM]_{-c})-\alpha\times\Borp
-\Borp \times \alpha -(I^{\prime}\cup_{\beta}(-I))\tilde{\times}_{\leq}(I^{\prime}\cup_{\beta}(-I)).$$ and
$$C_1(\Borp;\tau)=s_{\tau}(\Bor;-\qvarM)+s_{\tau}(I;[-\qvarM,\qvarM]_{c})-\alpha\times\Borp
-\Borp \times \alpha +(I^{\prime}\cup_{\beta}(-I))\tilde{\times}_{\geq}(I^{\prime}\cup_{\beta}(-I)).$$
\end{lemma}
\bp
$$A(K)=[\alpha,\beta] \times I_2 + t^{-1}\left(I_2\times[\alpha,\beta]\right) + \left(I \tilde{\times}_{\leq} I\right)
+t^{-1}\left(I \tilde{\times}_{\geq} I\right)+ \left(I_2 \tilde{\times}_{\leq} I_2\right)
+t^{-1}\left(I_2 \tilde{\times}_{\geq} I_2\right)$$
where $s_{\tau}(I;\qvarM) \subset\partial( I \tilde{\times}_{\leq} I)$, $s_{\tau}(I;-\qvarM) \subset -\partial (I \tilde{\times}_{\geq} I)$, and some attention should also be paid to points $(\alpha,\alpha)$ and $(\beta,\beta)$. For example, $(\alpha,\alpha)$ goes to
$t^{-1}s_{\tau}(\alpha;-\qvarM)$ in $[\alpha,\beta] \times I_2$ and to
$t^{-1}s_{\tau}(\alpha;\qvarM)$ in $t^{-1}I_2\times[\alpha,\beta]$.

$$\begin{array}{ll}C(\Bor;\tau)&=\partial \Bor\times I_2 + t^{-1}I_2\times\partial \Bor\\&+ I^{\prime} \tilde{\times}_{\leq} I^{\prime}
+t^{-1}I^{\prime} \tilde{\times}_{\geq} I^{\prime}- I \tilde{\times}_{\leq} I
-t^{-1}I \tilde{\times}_{\geq} I \\&+t^{-1}s_{\tau}(\Bor;-\qvarM) -s_{\tau}(\Bor;\qvarM).
\end{array}$$

Consider the closure $\left(\Borp \times (I \cup_{\beta} I_2)\right)_{\TCM}$ in $\TCM$ of $\left(\Borp \times (I \cup_{\beta} I_2)\right) \setminus \mbox{diag}(I^2)$.

A point in $(\Bor \times I) \setminus \mbox{diag}(I^2)$ near the diagonal reads $(u +\varepsilon v, u)$, where $u\in I$, $\varepsilon \in [0,1[$, $v \in [-\qvarM,\qvarM]_{c}$, where the orientation is induced by $(\varepsilon,v,u)$, or $(-(\varepsilon , (u,-v)))$ since
$(v \mapsto(-v))$ preserves the orientation of $S^1$.
Then $\left(\Borp \times (I \cup_{\beta} I_2)\right)_{\TCM}$ is a $3$--chain whose boundary is
$$\partial \left(\Borp \times (I \cup_{\beta} I_2)\right)_{\TCM}=
 \partial \Bor\times(I \cup_{\beta} I_2) + (t^{-1}-1)(\Borp \times \alpha) -s_{\tau}(I;[-\qvarM,\qvarM]_{-c}).$$
(Again, $ \partial \Bor\times(I \cup_{\beta} I_2)$  abusively stands for the preferred lift in $\TCM$ of the closure in $C_2(M)$ of 
$ \left(\partial \Bor\times(I \cup_{\beta} I_2) \setminus \mbox{diag}(I^2)\right)$.)
Similarly, the boundary of the $3$--chain
$t^{-1}\left((I \cup_{\beta} I_2)\times \Borp\right)_{\TCM} $ is
$$\partial \left(t^{-1}\left((I \cup_{\beta} I_2)\times \Bor\right)_{\TCM}\right)
=-t^{-1}(I \cup_{\beta} I_2)\times\partial \Bor
+(1-t^{-1})(\alpha\times\Borp  )+t^{-1}s_{\tau}(I;[-\qvarM,\qvarM]_{c}),$$
and $C(\Bor;\tau)$ is homologous to
$$\begin{array}{l}-\partial \Bor\times I
+(1-t^{-1})(\Borp \times \alpha)
+s_{\tau}(I;[-\qvarM,\qvarM]_{-c})\\
- t^{-1}I\times\partial \Bor+ (1-t^{-1})(\alpha\times\Borp  ) +t^{-1}s_{\tau}(I;[-\qvarM,\qvarM]_{c})\\
+I^{\prime} \tilde{\times}_{\leq} I^{\prime}
+t^{-1}I^{\prime} \tilde{\times}_{\geq} I^{\prime}- I \tilde{\times}_{\leq} I
-t^{-1}I \tilde{\times}_{\geq} I+t^{-1}s_{\tau}(\Bor;-\qvarM) -s_{\tau}(\Bor;\qvarM)\end{array}$$
and to
$$\begin{array}{l}I^{\prime}\times (-I)
+(1-t^{-1})(\Borp \times \alpha)
+s_{\tau}(I;[-\qvarM,\qvarM]_{-c})\\
+ t^{-1}(-I)\times I^{\prime}+ (1-t^{-1})(\alpha\times\Borp  ) +t^{-1}s_{\tau}(I;[-\qvarM,\qvarM]_{c})\\
+I^{\prime} \tilde{\times}_{\leq} I^{\prime}
+t^{-1}I^{\prime} \tilde{\times}_{\geq} I^{\prime}+ I \tilde{\times}_{\geq} I
+t^{-1}I \tilde{\times}_{\leq} I+t^{-1}s_{\tau}(\Bor;-\qvarM) -s_{\tau}(\Bor;\qvarM).\end{array}$$
It is now easy to conclude. \eop

We shall assume that $\delta(t)\Borp$ reads $\delta(t) I^{\prime} \times [0,1[ - R^{\prime} M(I^{\prime})$ near the lifts of $K^{\prime}$, where $M(I^{\prime})$ is a small meridian disk of $I^{\prime}$, $R^{\prime} \in \QQ[t,t^{-1}]$ satisfies $R^{\prime}(1)=r^{\prime} \in \QQ$ and $I^{\prime} \times [0,1[$ induces the same parallelisation of $K^{\prime}$ as $\tau$.
Let $\partial \Bor_{\parallel}$ be the parallel of $\partial \Bor$ induced by the
parallels of $K$ and $K^{\prime}$, obtained by pushing $\partial \Bor$ inside $\Bor$ near $I$ and induced by the parallel of $K^{\prime}$ along $I^{\prime}$.
$$lk_e(\partial \Bor, \partial \Bor_{\parallel})=- \frac{R^{\prime}}{\delta}=-\frac{\overline{R^{\prime}}}{\delta}.$$

\begin{lemma}
\label{lemc2}
Set $$C_2(\Borp;\tau)=C_0(\Borp;\tau) - \frac{R^{\prime}}{2\delta} ST(\ast).$$ 
Then $\delta C_2(\Borp;\tau)$ is a rational cycle, and
$$[C_2(\Borp;\tau)]=[C_1(\Borp;\tau)]+ \frac{R^{\prime}}{2\delta} [ST(\ast)].$$
In particular,
$$[C(\Bor;\tau)]=(t^{-1}-1)[C_2(\Borp;\tau)]-(t^{-1} +1)\frac{R^{\prime}}{2\delta} [ST(\ast)].$$
\end{lemma}
\bp
We only (need to) prove that the two expressions of $[C_2(\Borp;\tau)]$ coincide.
$$\begin{array}{ll}
\partial (\partial \Bor \times (\Borp \setminus \mbox{Int}(-\frac{R^{\prime}}{\delta}M(I^{\prime})))_{\TCM}=&-s_{\tau}(I^{\prime};[-\qvarM,\qvarM]_{-c}) +s_{\tau}((-I);[-\qvarM,\qvarM]_{c})\\
 &-(I^{\prime}\cup_{\beta}(-I))\tilde{\times}_{\leq}(I^{\prime}\cup_{\beta}(-I))\\
&-(I^{\prime}\cup_{\beta}(-I))\tilde{\times}_{\geq}(I^{\prime}\cup_{\beta}(-I))
-\partial \Bor \times \frac{R^{\prime}}{\delta}m(I^{\prime})\end{array}$$

$$\begin{array}{ll}[C_0(\Borp;\tau)-C_1(\Borp;\tau)-\frac{R^{\prime}}{\delta}ST(\ast)]&= [s_{\tau}(\partial \Bor;[-\qvarM,\qvarM]_{-c})+s_{\tau}(\Bor;\qvarM)-s_{\tau}(\Bor;-\qvarM)]\\&= [\partial(s_{\tau}( \Bor;[-\qvarM,\qvarM]_{-c}))]=0.\end{array}$$
\eop

\begin{lemma}
\label{lemc3}
Define $Q_2(\Borp;\tau)(t)$ so that $[C_2(\Borp;\tau)]=Q_2(\Borp;\tau)(t)[ST(\ast)]$.
Then $$Q(\Bor;\tau)(t)=Q_2(\Borp;\tau)(t) -\frac{1+t}{1-t}\frac{R^{\prime}(t)}{2\delta}.$$
$$Q(\Bor;\tau)(t)=-Q(\Bor;\tau)(t^{-1}).$$
\end{lemma}
\bp
The involution that exchanges the two factors in $\tilde{M}^2$ induces an involution $\iota$ of $\TCM$ such that $\iota_{\ast}(P(t)[ST(\ast)])=-P(t^{-1})[ST(\ast)]$ for any $P\in \Lambda$.
Since $\iota(C_0(\Borp;\tau))=C_1(\Borp;\tau)$,
$$\iota_{\ast}([C_2(\Borp;\tau)])=[C_2(\Borp;\tau)]$$
and $Q_2(\Borp;\tau)(t)=-Q_2(\Borp;\tau)(t^{-1}).$
Now, according to Lemma~\ref{lemc2}, and to the definition of $Q(\Bor;\tau)$ before Proposition~\ref{propvarone},
$$Q(\Bor;\tau)(t^{-1})=-Q_2(\Borp;\tau)(t) -\frac{1+t^{-1}}{1-t^{-1}}\frac{R^{\prime}(t)}{2\delta}.$$
This proves the formula that implies the antisymmetry of $Q$.
\eop

\begin{lemma}
\label{lemQBrat}
Any cycle $\tvarM$ of $H_2(\TCM)$ reads $\frac{Q_{\tvarM}}{\delta(t-1)}[ST(\ast)]$
in $H_2(C_2(M);\QQ(t))$ for some $Q_{\tvarM}$ in $\QQ[t,t^{-1}]$.
\end{lemma}
\bp Since $\delta(t-1)\tvarM$ vanishes in $H_2(\TCMD)$ according to Proposition~\ref{proptilMtwocomp}, it bounds a rational $3$-chain, and $\pm Q_{\tvarM}$ is the equivariant intersection of this $3$-chain with the preferred lift of the diagonal of $M$.
\eop

\begin{lemma} 
\label{lemc0q0}
Under the assumptions above, where $r=0$,
set $$[C_0(\Borp;\tau)]=Q_0(\Borp;\tau)[ST(\ast)].$$
Then Proposition~\ref{propvargen} is satisfied with
$${\cal V}(\KK, \KK^{\prime})= \delta Q(\Bor;\tau) + \delta\frac{r^{\prime}}{2} \ID= \delta Q_0(\Borp;\tau) -\frac{R^{\prime}}{2} +\frac{1+t}{t-1}\frac{R^{\prime}}{2} + \delta\frac{r^{\prime}}{2} \ID.$$
\end{lemma}
\bp The first formula for ${\cal V}(\KK, \KK^{\prime})$ comes from Proposition~\ref{propvarone}, and the second one easily follows from Lemmas~\ref{lemc2} and \ref{lemc3}. It is enough to
see that ${\cal V}(\KK, \KK^{\prime}) \in \QQ[t,t^{-1}]=\Lambda$. Recall that $\delta=\delta(M)$ is the annihilator of $H_1(\tilde{M})$.
Since it has the same roots as $\Delta(M)$, $\delta \frac{\Delta^{\prime}}{\Delta} \in \Lambda$.
Since $(r^{\prime}\delta -R^{\prime})(1)=0$, $(t-1)$ divides $(r^{\prime}\delta -R^{\prime})$ and $\frac{1+t}{2(1-t)}(r^{\prime}\delta -R^{\prime}) \in \Lambda$.
Therefore, it is enough to prove that $\delta Q_0(\Borp;\tau) \in \Lambda$.
Thanks to Lemma~\ref{lemQBrat}, $(1-t)^2\delta(t)Q(\Bor;\tau)(t^{-1}) \in \Lambda$. Therefore, $(1-t)^2\delta Q_0(\Borp;\tau)\in \Lambda$.
According to Proposition~\ref{proptilMtwocomp}, since $\delta C_0(\Borp;\tau)$ is a rational cycle, it is 
rationally homologous in $\TCMD$ to $u [S\times \ast]$ + $v [\ast \times S]$
up to some elements of $\delta$--torsion, where $u$ and $v$ are rational numbers that can be computed as the algebraic intersection with the preimage of $M \times K_{\parallel}$ (up to sign) or $K_{\parallel} \times M$, respectively. Since our assumptions imply that $\Borp$ does not meet the preimage of $K_{\parallel}$, and since  $\Bor$ does not meet $K_{\parallel}$ either, $u$ and $v$ are
zero, and $\delta^2C_0(\Borp;\tau)$ is null-homologous in 
$H_2(\TCMD)$. Therefore $\delta^2Q_0(\Borp;\tau) \in \Lambda$, and since $\delta$ and $(1-t)^2$ are coprime, $\delta Q_0(\Borp;\tau) \in \Lambda$. \eop 

\subsection{An independent lemma}

\begin{lemma} 
\label{lemhomdiagSbry} Let $\Sigma$ be a compact oriented surface with one boundary component $J(S^1)$ equipped with a basepoint $\ast=J(1)$.
Let $\Sigma$ and $\Sigma^+$ be two copies of $\Sigma$, let $(z_i)_{i=1, \dots 2g}$ and $(z^{\ast}_i)_{i=1, \dots, 2g}$ be two dual bases
of $H_1(\Sigma;\ZZ)$ such that 
$$\langle z_i, z^{\ast}_j\rangle=\delta_{ij}.$$
Set $$J\times_{\ast,\leq}J^+=\{(J(\exp(2i\pi t)),J(\exp(2i\pi u))); (t,u)\in [0,1]^2,t\leq u\}$$ and
$$J\times_{\ast,\geq}J^+=\{(J(\exp(2i\pi t)),J(\exp(2i\pi u))); (t,u)\in [0,1]^2,t\geq u\}.$$
Let $\mbox{diag}(\Sigma \times \Sigma^+)=\{(x,x); x \in \Sigma\}$.
Then the chains $$C_{\ast,\leq}(\Sigma,\Sigma^+)=\mbox{diag}(\Sigma \times \Sigma^+)-\ast \times \Sigma^+ -\Sigma \times \ast^+ -J\times_{\ast,\leq}J^+$$ and
$$C_{\ast,\geq}(\Sigma,\Sigma^+)=\mbox{diag}(\Sigma \times \Sigma^+)-\ast \times \Sigma^+ -\Sigma \times \ast^+ +J\times_{\ast,\geq}J^+$$
are cycles and
we have the following equality in $H_2(\Sigma \times \Sigma^+)$
$$[C_{\ast,\leq}(\Sigma,\Sigma^+)]=[C_{\ast,\geq}(\Sigma,\Sigma^+)] =\sum_{i=1}^{2g} [z_i \times z^{\ast+}_i].$$
\end{lemma}
\bp Since $$\partial (J\times_{\ast,\leq}J^+) =\mbox{diag}(J \times J^+) -\ast \times J^+ -J \times \ast^+$$
 $C_{\ast,\leq}(\Sigma,\Sigma^+)$ is a cycle. 
Consider the closed surface $S$ obtained from $\Sigma$ by gluing a disk $D$ along $J$. According to Proposition~\ref{prophomdiagS}, in $H_2(S\times S^+)$,
$$[\mbox{diag}(S \times S^+)]=[\ast \times S^+] + [S \times \ast^+]+ \sum_{i=1}^{2g} [z_i \times z^{\ast +}_i].$$
This implies that $$[C_{\ast,\leq}(\Sigma,\Sigma^+)-C_{\ast,\leq}(-D,(-D)^+)]=\sum_{i=1}^{2g} [z_i \times z^{\ast +}_i]$$ in $H_2(S \times S^+)$.
Since the cycle $C_{\ast,\leq}(-D,(-D)^+)$ lives in $D\times D^+$, it is null-homologous there, and
since $H_2(\Sigma \times \Sigma^+)$ naturally injects into 
$H_2(S \times S^+)$, we can conclude that $[C_{\ast,\leq}(\Sigma,\Sigma^+)] =\sum_{i=1}^{2g} [z_i \times z^{\ast +}_i]$. The proof for $C_{\ast,\geq}(\Sigma,\Sigma^+)$ is the same.
\eop

\subsection{Examples of variations}

In this subsection, we shall keep the notation of the first two subsections of this section and assume that $r=r^{\prime}=0$.
We shall compute $Q(\Bor)=Q(\Bor;\tau)$ in some examples where $\Bor$ lifts in $\tilde{M}$, so that $\Bor=\Borp$ and $R^{\prime}=0$.

\begin{lemma}
\label{lemkeyvarsur}
If $\Bor$ lifts in $\tilde{M}$, (if $R^{\prime}=r^{\prime}=r=0$,) and if  $(z_i)_{i=1, \dots 2g}$ and $(z^{\ast}_i)_{i=1, \dots, 2g}$ are two dual bases
of $H_1(\Bor;\ZZ)$ such that 
$\langle z_i, z^{\ast}_j\rangle=\delta_{ij},$ then $$[C_0(\Bor)]= [C_0(\Bor;\tau)]=g[ST(\ast)]+\sum_{i=1}^{2g} [z_i \times z^{\ast +}_i]$$ and 
$$Q(\Bor)(t)=\sum_{i=1}^{2g}lk_e(z_i,z^{\ast +}_i) +g.$$
\end{lemma}
\bp The expression of $Q(\Bor)$ follows from the expression of $[C_0(\Bor)]$, thanks to Lemma~\ref{lemc0q0}. Thus, we are left with the computation of $[C_0(\Bor)]= [C_0(\Bor;\tau)]$ that was defined in Lemma~\ref{lempremc2}. Our computation will rely on Lemma~\ref{lemhomdiagSbry}.

Set $d=\qvarM$.
Along $I$, $\tau$ maps the inward normal of $\Bor$ to $c \in S^2$ and the positive normal of $\Bor$ to $e$ so that $(c,d,e)$ is an oriented orthonormal basis for $\RR^3$. Let $\langle c,d,e\rangle$ denote the intersection of the sphere 
with the convex hull of $0$, $3c$, $3d$ and $3e$ in $\RR^3$ with the orientation of $S^2$
while $\langle d,c,e\rangle=-\langle c,d,e\rangle$.
$[-\qvarM,\qvarM]_{c}=\langle -d,c\rangle + \langle c,d\rangle$.
Consider $\Bor$ as $\Bor \times 0$ in $\Bor \times [0,1] \subset M$
and consider the closure $C_{\ast,\leq}(\Bor,\Bor \times [0,1])$ in $\TCM$
of $\cup_{t\in ]0,1]}C_{\ast,\leq}(\Bor,\Bor \times t)$, where the basepoint $\ast$ is $\alpha$, with the notation of Lemma~\ref{lemhomdiagSbry}. This is a $3$--chain whose boundary contains $(-s_+(\Bor))$ in the closure of $\cup_{t\in ]0,1]}\mbox{diag}(\Bor \times \Bor \times t)$.
Let us consider the ways of approaching $ST(\alpha)$ and determine $ST(\alpha) \cap \partial C_{\ast,\leq}(\Bor,\Bor \times [0,1])$.
\begin{itemize}
\item In the closure of $\alpha \times (\Bor \times ]0,1])$,
we only find $s_{\tau}(\alpha;\langle d,e\rangle)$.
\item  In the closure of $\Bor \times(\alpha\times]0,1])$, we only have $s_{\tau}(\alpha;\langle -d,e\rangle) $. 
\item When a pair of points of $I$, (resp. of $I^{\prime}$) in $(I^{\prime}\cup_{\beta}(-I))\tilde{\times}_{\leq}((I^{\prime}\cup_{\beta}(-I)) \times ]0,1]$ approaches $ST(\alpha)$, it approaches $s_{\tau}(\alpha;\langle \pm d,e\rangle) $, again.
\item When the first point is in $I^{\prime}$ and the second one is in $(-I)$, the vector between them is in $\langle d,-c,e\rangle \cup \langle -c,-d,e\rangle$
where $I^{\prime}$ gives the inward normal, $(-I)$ moves along $[d,-d]_{-c}$, and $[0,1]$ gives the $e$ direction. In particular, the boundary of $-I^{\prime} \times (-I)\times [0,1]$ is oriented like $(\langle d,-c,e\rangle \cup \langle -c,-d,e\rangle)$, that is like the sphere.
\end{itemize}

The only part that will matter to us is this last $2$--dimensional part.

We similarly determine the $2$--dimensional part of $ST(\beta) \cap \partial C_{\ast,\leq}(\Bor,\Bor \times [0,1])$ that comes from pairs of points
where the first point is in $I^{\prime}$ and the second one is in $(-I)$, the vector between them is again in $\langle d,-c,e\rangle \cup \langle -c,-d,e\rangle$. Now,
$I^{\prime}$ gives the outward normal, $(-I)$ moves along $[d,-d]_{-c}$, and $[0,1]$ gives the $e$ direction. In particular the boundary of $-I^{\prime} \times (-I)\times [0,1]$ is oriented like $(-(\langle d,-c,e\rangle \cup \langle -c,-d,e\rangle))$ unlike the sphere.

We now compute
$$C_{0,+}=C_{\ast,\leq}(\Bor,\Bor \times 1)-\partial C_{\ast,\leq}(\Bor,\Bor \times [0,1])$$ forgetting the $1$--dimensional parts,

$$\begin{array}{lll}C_{0,+}&=&s_+(\Bor)-\alpha \times (\Bor \setminus \alpha) - (\Bor \setminus \alpha)\times\alpha \\ &&-(I^{\prime}\cup_{\beta}(-I))\tilde{\times}_{\leq}(I^{\prime}\cup_{\beta}(-I)) 
+ s_{\tau}(I^{\prime};\langle d,e\rangle)+ s_{\tau}(-I;\langle -d,e\rangle))\\&&- s_{\tau}(\alpha;\langle d,-c,e\rangle \cup \langle -c,-d,e\rangle)+s_{\tau}(\beta;\langle d,-c,e\rangle \cup \langle -c,-d,e\rangle).
\end{array}$$
where $d=W$ and $e$ denotes the positive normal to $\Bor$ along $I^{\prime}$, and
$$[s_+(\Bor) +s_{\tau}(\partial \Bor;\langle d,e\rangle)-s_{\tau}(\Bor;\qvarM)]=-g[ST(\ast)]$$
since this is an obstruction to homotope the positive normal to a point in $S^2$ via $\tau$ that would be zero if $\Bor$ were a disk, and that is similar to the obstruction computed in Proposition~\ref{propsplusstau}. 

Then $$\begin{array}{lll}[C_0(\Bor)]&=&[s_+(\Bor) +s_{\tau}(\partial \Bor;\langle d,e\rangle)+gST(\ast)-s_{\tau}(I;[-\qvarM,\qvarM]_{-c})\\&&-\alpha\times\Bor
-\Bor \times \alpha-(I^{\prime}\cup_{\beta}(-I))\tilde{\times}_{\leq}(I^{\prime}\cup_{\beta}(-I))]\\
&=&[C_{0,+}+g ST(\ast)+s_{\tau}(\partial \Bor;\langle d,e\rangle)-s_{\tau}(I;[-\qvarM,\qvarM]_{-c})\\&&- s_{\tau}(I^{\prime};\langle d,e\rangle)- s_{\tau}(-I;\langle -d,e\rangle)\\
&&+ s_{\tau}(\alpha;\langle d,-c,e\rangle \cup \langle -c,-d,e\rangle)-s_{\tau}(\beta;\langle d,-c,e\rangle \cup \langle -c,-d,e\rangle)].
\end{array}$$
where $$\partial s_{\tau}((-I);\langle d,-c,e\rangle \cup \langle -c,-d,e\rangle)$$
$$=s_{\tau}(\alpha;\langle d,-c,e\rangle \cup \langle -c,-d,e\rangle)-s_{\tau}(\beta;\langle d,-c,e\rangle \cup \langle -c,-d,e\rangle)
+s_{\tau}((-I);[-\qvarM,\qvarM]_{-c} -[-\qvarM,\qvarM]_{e}).$$
Since $[C_{0,+}]=[C_{\ast,\leq}(\Bor,\Bor \times 1)]=
\sum_{i=1}^{2g} [z_i \times z^{\ast +}_i]$, according to Lemma~\ref{lemhomdiagSbry},
we find
$$[C_0(\Bor)]=g[ST(\ast)]+\sum_{i=1}^{2g} [z_i \times z^{\ast +}_i].$$
\eop

\noindent{\sc Proof of Proposition~\ref{propcorcalvarb}:} It is a direct consequence of Proposition~\ref{propvarone} and Lemma~\ref{lemkeyvarsur}, using the same arguments as in Remark~\ref{remlog}.
\eop

Note that Proposition~\ref{propcorcalvarb} implies that if the $a_i$ and the $b_i$ bound in the complement of $S$, then $Q(\Bor)(t)=0$.

\begin{example}
According to Lemma~\ref{lemblanchnondeg}, there exist a two-component link $(a,b)$ in $\tilde{M} \setminus p_M^{-1}(K)$, and a rational number $q\neq 0$ such that $lk_e(a,b)=\frac{qt^k}{\delta(M)}$.
Plumb two annuli around $a$ and $b$ so that $lk_e(a,b^+)=\frac{qt^k}{\delta(M)}$. Let $K^{\prime}$ be a band sum of $K$ with the boundary of the genus one surface obtained by this plumbing. Then $K^{\prime}$ is homologous to $K$ and $Q(t)=q\frac{t^k-t^{-k}}{\delta(M)}$.
\end{example}

This example concludes the proof of Theorem~\ref{thmfrakcha}.
\eop

\begin{remark}
Consider a cobordism $\Bor$ between $K$ and $1 \times K$ constructed from $[0,1] \times K +S $ by surgery around the intersection segment.  $$\partial \Bor=1 \times K-m(1 \times K) -K -m(K).$$
Then $\CQ(K,K-m(1 \times K))-\CQ(K,K+m(1 \times K))$
is computed both in Proposition~\ref{propvarpar}, and in Proposition~\ref{propcorcalvarb}.
Thanks to Remark~\ref{remlog}, these two computations lead to the same result and show that
$Q(\Bor)(t)= \frac{t\Delta^{\prime}(t)}{\Delta(t)}$ in this case.
\end{remark}

\newpage
\section{The Lagrangian-preserving surgery formula}
\setcounter{equation}{0}
\label{secsur}

\subsection{The statement of the Lagrangian-preserving surgery formula}
\label{substateLP}

A {\em genus $g$ $\QQ$--handlebody}\/ is an (oriented, compact) 3--manifold $A$
with 
the same homology with rational coefficients as the standard (solid) handlebody
$H_g$ of the following figure.
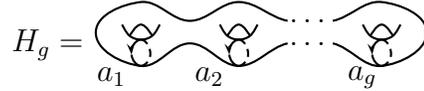
\begin{center}
$$H_g = \begin{pspicture}[shift=-0.75](0,-.5)(4.5,.95) 
\psset{xunit=.5cm,yunit=.5cm}
\psecurve{-}(5.7,1.3)(5.2,1.3)(3.9,1.8)(2.6,1.3)(1.3,1.8)(.1,1)(1.3,.1)(2.6,.7)
(3.9,.1)(5.2,.7)(5.7,.7) 
\pscurve{-}(.8,1.2)(1,.9)(1.3,.8)(1.6,.9)(1.8,1.2) 
\pscurve{-}(1,.9)(1.3,1.2)(1.6,.9) 
\rput[r](.9,-.2){$a_1$} 
\psecurve{->}(1.6,.4)(1.3,.8)(1.05,.4)(1.3,.1) 
\psecurve{-}(1.3,.8)(1.05,.4)(1.3,.1)(1.6,.4) 
\psecurve[linestyle=dashed,dash=3pt 2pt](1,.4)(1.3,.8)(1.55,.4)(1.3,.1)(1,.4) 
\pscurve{-}(3.4,1.2)(3.6,.9)(3.9,.8)(4.2,.9)(4.4,1.2) 
\pscurve{-}(3.6,.9)(3.9,1.2)(4.2,.9) 
\rput[r](3.5,-.2){$a_2$} 
\psecurve{->}(4.2,.4)(3.9,.8)(3.65,.4)(3.9,.1) 
\psecurve{-}(3.9,.8)(3.65,.4)(3.9,.1)(4.2,.4) 
\psecurve[linestyle=dashed,dash=3pt
2pt](3.6,.4)(3.9,.8)(4.15,.4)(3.9,.1)(3.6,.4) 
\rput(5.8,1.3){\dots} 
\rput(5.8,.7){\dots} 
\psecurve{-}(5.8,1.3)(6.4,1.3)(7.7,1.8)(8.9,1)(7.7,.1)(6.4,.7)(5.8,.7)
\pscurve{-}(7.2,1.2)(7.4,.9)(7.7,.8)(8,.9)(8.2,1.2) 
\pscurve{-}(7.4,.9)(7.7,1.2)(8,.9) 
\rput[l](6.8,-.2){$a_g$} 
\psecurve{->}(8,.4)(7.7,.8)(7.45,.4)(7.7,.1) 
\psecurve{-}(7.7,.8)(7.45,.4)(7.7,.1)(8,.4) 
\psecurve[linestyle=dashed,dash=3pt
2pt](7.4,.4)(7.7,.8)(7.95,.4)(7.7,.1)(7.4,.4) 
\end{pspicture}$$ 
\label{fig1}
\end{center}
The boundary of such a $\QQ$--handlebody $A$ (or {\em rational homology handlebody\/}) is homeomorphic to the 
boundary $(\partial H_g =\Sigma_g)$ of $H_g$.

For a $\QQ$--handlebody $A$,
the kernel ${\cal L}_A$ of the map induced by the inclusion
$$ H_1(\partial A;\QQ) \longrightarrow H_1( A;\QQ)$$ 
is a Lagrangian of $(H_1(\partial A;\QQ),\langle,\rangle_{\partial A})$, that is called the {\em {Lagrangian}}\/ of $A$.
A {\em Lagrangian-preserving surgery\/} or {\em LP--surgery} $(A,A^{\prime})$ is the replacement of a 
$\QQ$--handlebody $A$ embedded in a 3--manifold by
another such $A^{\prime}$ with identical (identified via a homeomorphism) boundary and Lagrangian. The manifold obtained by such 
an LP surgery $(A,A^{\prime})$ from a manifold $M$ (with $A \subset M$) is denoted by $M(A^{\prime}/A)$.

There is a canonical isomorphism 
$$\partial_{MV} \colon H_2(A \cup_{\partial A}
-A^{\prime};\QQ) \rightarrow {\cal L}_A$$
that maps the class of a closed surface in the closed 3--manifold $(A \cup_{\partial A}
-A^{\prime})$ to the boundary of its intersection with $A$.
This isomorphism carries the algebraic triple intersection of surfaces to a trilinear antisymmetric form $\CI_{AA^{\prime}}$ on $\CL_A$.
$$\CI_{AA^{\prime}}(a_{i},a_{j},a_{k})=\CI_{AA^{\prime}}(a_{i}\wedge a_{j}\wedge a_{k})=\langle \partial_{MV}^{-1}(a_i), \partial_{MV}^{-1}(a_j), \partial_{MV}^{-1}(a_k)\rangle_{A \cup
-A^{\prime}}$$

\begin{theorem}
\label{thmLP}
Let $A$ and $B$ be two disjoint rational homology handlebodies of $M$ disjoint from $K$. Assume that $H_1(A)$ and $H_1(B)$ are sent to $0$ in $H_1(M;\QQ)$ by the maps induced by the inclusions.
Let $(a_i,z_i)_{i=1, \dots, g_A}$  be a basis of $H_1(\partial A)$  such that $(a_i)_{i=1, \dots, g_A}$ is a basis of $\CL_A$
and $\langle a_i,z_j \rangle_{\partial A}=\delta_{ij}$.
Let $(b_i,y_i)_{i=1, \dots, g_B}$  be a basis of $H_1(\partial B) $  such that $(b_i)_{i=1, \dots, g_B}$ is a basis of $\CL_B$
and $\langle b_i,y_j \rangle_{\partial B}=\delta_{ij}$.
Let $A^{\prime}$ and $B^{\prime}$ be two rational homology handlebodies
whose boundaries are identified to $\partial A$ and to $\partial B$, respectively, so that 
$\CL_{A^{\prime}}=\CL_A$ and $\CL_{B^{\prime}}=\CL_B$. Set
$M_A=M(A^{\prime}/A)$, $M_B=M(B^{\prime}/B)$, $M_{AB}=M(A^{\prime}/A,B^{\prime}/B)$ and
let $$\CS=
\CQ(\KK \subset M_{AB})+\CQ(\KK \subset M) -\CQ(\KK \subset M_A) -\CQ(\KK \subset M_B).$$
Then $$ \CS= -
\sum_{(i,j,k) \in  \{1, \dots, g_A\}^3,(\ell,m,n) \in  \{1, \dots, g_B\}^3} lk(i,j,k,\ell,m,n) \CI(i,j,k,\ell,m,n)$$
where
$$ \begin{array}{ll}lk(i,j,k,\ell,m,n)=&
lk_e(z_i, y_{\ell})(\fvar)lk_e(z_j, y_{m})(\svar)lk_e(z_k, y_{n})(\tvar)\\ &+lk_e(z_i, y_{\ell})(\fvar^{-1})lk_e(z_j, y_{m})(\svar^{-1})lk_e(z_k, y_{n})(\tvar^{-1})\\
&+\sum_{\mathfrak{S}_3(\fvar,\svar,\tvar)}lk(z_i, y_{\ell})lk_e(z_j, z_k)(\fvar)lk_e(y_{m}, y_{n})(\svar)
,\end{array} $$
and
 $$\CI(i,j,k,\ell,m,n)=\CI_{AA^{\prime}}(a_i\wedge a_j \wedge a_k)\CI_{BB^{\prime}}(b_{\ell}\wedge b_m \wedge b_n).$$
Here, when the $z_i$ (resp. the $y_{\ell}$) are arguments of equivariant linking numbers, they are lifts of the $z_i$ (resp. of the $y_{\ell}$) sitting all in the same lift of $A$ (resp. of $B$) in $\tilde{M}$.
\end{theorem}

The proof of this theorem will be given in Section~\ref{secproofLP}.

\subsection{Relations with the clasper theory and with the Kontsevich integral}

For $P$, $Q$ and $R$ in $\QQ[t^{\pm 1},\frac{1}{\delta(t)}]$, and for $\lambda \in \QQ$, see the two {\em beaded graphs} below as the following  elements of
$R_{\delta}$ $$\begin{pspicture}[shift=-0.4](0,-.1)(1.2,1.3) 
\psccurve(0.1,.4)(.2,.9)(1,.9)(1.1,.4)(1,-.1)(.2,-.1)
\psline{*-*}(0.1,.4)(1.1,.4) 
\rput[b](.6,1.15){$P$}
\rput[b](.6,.45){$Q$}
\rput[b](.6,-.1){$R$} \end{pspicture}=
\sum_{\mathfrak{S}_3(\fvar,\svar,\tvar)} \left(P(\fvar)Q(\svar)R(\tvar)+P(\fvar^{-1})Q(\svar^{-1})R(\tvar^{-1})\right)
$$

$$\begin{array}{ll}\begin{pspicture}[shift=-.4](-.2,-.1)(1.7,1.1) 
\pscircle(0.1,0.4){.3}
\psline{*-*}(0.4,.4)(1,.4) 
\pscircle(1.3,0.4){.3}
\rput(.1,.9){$P$}
\rput(1.3,.9){$Q$}
\rput[b](.7,.5){$\lambda$} \end{pspicture}&=
\lambda\sum_{\mathfrak{S}_3(\fvar,\svar,\tvar)} \left(P(\fvar)Q(\svar^{-1})+P(\fvar^{-1})Q(\svar)-P(\fvar)Q(\svar)-P(\fvar^{-1})Q(\svar^{-1})\right)
\\&=\lambda\sum_{\mathfrak{S}_3(\fvar,\svar,\tvar)}\left(P(\fvar)-P(\fvar^{-1})\right)\left(Q(\svar^{-1})-Q(\svar)\right)
\end{array}$$
so that the following IHX or Jacobi relation is satisfied.
$$\begin{pspicture}[shift=-.4](-.2,-.1)(1.7,1.1) 
\pscircle(0.1,0.4){.3}
\psline{*-*}(0.4,.4)(1,.4) 
\pscircle(1.3,0.4){.3}
\rput(.1,.9){$P$}
\rput(1.3,.9){$Q$}
\rput[b](.7,.5){$\lambda$} \end{pspicture}=\begin{pspicture}[shift=-0.4](0,-.1)(1.2,1.3) 
\psccurve(0.1,.4)(.2,.9)(1,.9)(1.1,.4)(1,-.1)(.2,-.1)
\psline{*-*}(0.1,.4)(1.1,.4) 
\rput[b](.6,1.15){$\overline{P}$}
\rput[b](.6,.5){$\lambda$}
\rput[b](.6,-.15){$Q$} \end{pspicture}  -
\begin{pspicture}[shift=-0.4](0,-.1)(1.2,1.3) 
\psccurve(0.1,.4)(.2,.9)(1,.9)(1.1,.4)(1,-.1)(.2,-.1)
\psline{*-*}(0.1,.4)(1.1,.4) 
\rput[b](.6,1.15){$P$}
\rput[b](.6,.5){$\lambda$}
\rput[b](.6,-.15){$Q$} \end{pspicture}$$

Via the identification
above, the invariants $\delta(M)(\fvar)\delta(M)(\svar)\delta(M)(\tvar)\CQ(M)$ and $$\delta(M)(\fvar)\delta(M)(\svar)\delta(M)(\tvar)\hat{\CQ}(\hat{K})$$ take their values in the space $\CA_2(\QQ[t,t^{-1}])$ of $2$-loop trivalent graphs with beads described in \cite[Definition 1.6]{GR}. (Here, the edges are oriented so that their polynomial labels are on their left-hand sides.) The lack of denominators follows from Proposition~\ref{propdenwithoutz} that will be shown in an independent section.
This space $\CA_2(\QQ[t,t^{-1}])$ maps to unitrivalent graphs via the map Hair described in \cite[Section 1.4]{GR}. The following result is a direct consequence of Theorem~\ref{thmLP} and of the results in \cite{GR}.
\begin{theorem}
\label{thmKont}
For knots with trivial Alexander polynomials in integral homology spheres, the map
 $\mbox{Hair} \circ \hat{\CQ}$ coincides with the $2$-loop primitive part of the Kontsevich integral.
\end{theorem}
Before saying more precisely how this theorem follows from Theorem~\ref{thmLP} and \cite{GR},
let us relate Theorem~\ref{thmLP} with clasper calculus and illustrate it by some examples.

A special case of a Lagrangian-preserving surgery is the case where $A$ is the regular neighborhood of a $Y$-graph, that is a genuine genus $3$ handlebody and $A^{\prime}$ is obtained from $A$ by a borromean surgery. See \cite[Section 2.2]{al}. Then the three curves $z_1$, $z_2$ and $z_3$ can be chosen as
longitudes of the leaves, the three curves $a_1$, $a_2$ and $a_3$ are the corresponding meridians and 
$$\CI_{AA^{\prime}}(a_i\wedge a_j \wedge a_k)=\pm 1$$
See \cite[Lemma 4.2]{al}.
The fact that $H_1(A)$ is sent to $0$ in $H_1(M;\QQ)$ is equivalent to the fact that the leaves are null-homologous. 
Assume that $(A,A^{\prime})$ and $(B,B^{\prime})$ are both such borromean surgeries and that $y_1$, $y_2$ and $y_3$ are the longitudes of the leaves of the underlying $Y$--graph for $B$.

\begin{lemma}
\label{lemnullclasp}
Under the hypotheses above, the expression $\CS$ of Theorem~\ref{thmLP} reads
$$\CS=\sum_{\sigma \in \mathfrak{S}_3}\mbox{sign}(\sigma)\begin{pspicture}[shift=-0.6](-.6,-.3)(2.6,1.6)
\psline{*-*}(-.4,.4)(2.4,.4) 
\psframe[linearc=.4](-.4,-.3)(2.4,1.1)
\rput[b](1,1.15){$lk_e(z_3, y_{\sigma(1)})$}
\rput[b](1,.45){$lk_e(z_2, y_{\sigma(2)})$}
\rput[b](1,-.25){$lk_e(z_1, y_{\sigma(3)})$} 
\end{pspicture}
+ \sum_{(i,m) \in \{1,2,3\}^2}\begin{pspicture}[shift=-.4](-.7,0)(5.4,.8) 
\pscircle(0.2,0.4){.4}
\psline{*-*}(0.6,.4)(4.2,.4) 
\pscircle(4.6,0.4){.4}
\rput[b](.2,.85){$lk_e(y_{m+1},y_{m+2})$}
\rput[b](4.6,.85){$lk_e(z_{i+1},z_{i+2})$}
\rput[b](2.4,.45){$lk(y_m,z_i)$} \end{pspicture}$$
where the indices are integers mod $3$ in the latter expression.
\end{lemma}
\eop

\noindent{\sc Proof of Theorem~\ref{thmKont}: }
It follows from Theorem~\ref{thmLP} that $\CQ$ and $\hat{\CQ}$ are of null-type $2$ in the sense of \cite[Definition 1.1]{GR} since the considered LP-surgeries do not change the equivariant linking numbers.
According to \cite[Theorem 2]{GR}, a null-type $2$ invariant of knots with trivial Alexander polynomial in homology spheres, that is valued in a rational vector space, is determined by its value on the trivial knot of $S^3$ and by its values on the degree $2$ null-claspers (that are special cases of 
the surgeries considered in Theorem~\ref{thmLP}). Then since for such a clasper $G$, $\CS$ is the contraction of $G$ in the sense of \cite[Theorem 4]{GR}, the value of $\mbox{Hair} \circ \hat{\CQ}(G)$ coincides with the primitive $2$-loop part of the Kontsevich integral of $G$ according to \cite[Theorem 4]{GR}, and we are done since $\hat{\CQ}$ vanishes for the trivial knot, and since the Kontsevich integral of the trivial knot is made of wheels \cite{BLT}.
\eop

Let us see some more examples of applications.
Applying Theorem~\ref{thmLP} to the case of a clasper in the complement of a (gray in the picture) Seifert surface of $\hat{K}$ like 
\begin{pspicture}[shift=-0.4](-.9,0)(2.3,.8) 
\pscircle(0.7,0.4){.2} 
\psline{*-*}(-.2,.4)(.5,.4) 
\psline{*-*}(0.9,.4)(1.6,.4)
\pspolygon*[linecolor=lightgray](-.65,.7)(-.65,0)(-.35,0)(-.35,.7)
\psecurve{->}(-.3,.5)(-.7,.5)(-.8,.4)(-.4,.2)(-.2,.4)(-.3,.5)(-.7,.5)
\rput[bl](-.25,.5){$m$}
\psline(-.65,.7)(-.65,0)
\psline(-.35,.7)(-.35,0)
\pspolygon*[linecolor=lightgray](1.75,.7)(1.75,0)(2.05,0)(2.05,.7)
\psecurve{->}(2.1,.5)(1.7,.5)(1.6,.4)(2,.2)(2.2,.4)(2.1,.5)(1.7,.5)
\rput[bl](2.15,.5){$n$}
\psline(2.05,.7)(2.05,0)
\psline(1.75,.7)(1.75,0)
\end{pspicture}
like in Proposition 4.17 in \cite{oht2} yields
$$\CS=
\hat{\CQ}(\hat{K} \subset (M_{\KK})_{AB})-\hat{\CQ}(\hat{K} \subset M_{\KK})=\begin{pspicture}[shift=-0.4](0,-.1)(1.2,1.3) 
\psccurve(0.1,.4)(.2,.9)(1,.9)(1.1,.4)(1,-.1)(.2,-.1)
\psline{*-*}(0.1,.4)(1.1,.4) 
\rput[b](.6,1.15){$lk_e(m,n)$}
\rput[b](.6,.5){$1$}
\rput[b](.6,-.1){$1$} \end{pspicture}$$

The following lemma gives an expression for $lk_e(m,n)$ that finishes to show that the two-loop polynomial described in \cite{oht2} behaves like $12\Delta(\fvar)\Delta(\svar)\Delta(\tvar)\hat{\CQ}$ in this case.

\begin{lemma}
Let $m$ and $n$ be two meridian curves of $1$-handles of a Seifert surface $\Sigma$ of $\hat{K}$ as above.
Consider a basis $(z_i)_{i\in\{1,\dots,2g_{\Sigma}\}}$ of $H_1(\Sigma;\ZZ)$, the matrix $C=[(c_{ij})_{(i,j)\in\{1,\dots,2g_{\Sigma}\}^2}]$
where $$c_{ij}=t_M^{1/2}lk(z_i^+,z_j)-t_M^{-1/2}lk(z_i^-,z_j)$$ and its inverse $D=C^{-1}=[(d_{ij})_{(i,j)\in\{1,\dots,2g_{\Sigma}\}^2}]$.
Then $$lk_e(m,n)=-(t_M^{1/2}-t_M^{-1/2})\sum_{(i,j)\in\{1,\dots,2g_{\Sigma}\}^2}lk(m,z_i)d_{ij}lk(n,z_j)$$
where $m$ and $n$ stand for lifts of $m$ and $n$ in the same lift of $M_{\KK}\setminus \Sigma$, when the former ones are arguments of $lk_e(m,n)$.
\end{lemma}
\bp
There exist $x$ and $y$ in $H_1(\Sigma)$ such that $m$ and $n$ are respectively homologous to $(x^+ - x^-)$ and $(y^+-y^-)$ in $M_{\KK}\setminus \Sigma$.
Then for any $z \in H_1(\Sigma)$, $lk(m,z)=\langle x,z\rangle$.
Without loss, assume that $(x^+ - x^-)$ sits in the boundary of a regular neighborhood of $\Sigma$
and that $(y^+ - y^-)$ sits in the boundary of a regular neighborhood of $\Sigma$ inside the former one.
Then with similar notation as in Section~\ref{subsecderal}, 
$$lk_e(m,n)= lk_e(\hat{x}^{+}-\hat{x}^{-},\hat{y}^{+}-\hat{y}^{ -})=(t_M-1)lk_e(\hat{x}^{+}-\hat{x}^{-},\tilde{y})$$
where $$\hat{x}=\sum_{i=1}^{2g}\langle z_i, x\rangle \hat{z}_i^{\ast}$$
and like in Proposition~\ref{proplogderlk}
$$t_M^{-1/2}(\tilde{z}^+_j-\tilde{z}^-_j)=-
\sum_{i=1}^{2g}c_{ji}(\hat{z}_i^{\ast +}-\hat{z}_i^{\ast -}).$$
Then $$lk_e(m,n)=(t_M-1)\sum_{i=1}^{2g}\langle z_i, x\rangle\left( \sum_{j=1}^{2g}d_{ij}lk_e(-t_M^{-1/2}(\tilde{z}^+_j-\tilde{z}^-_j),\tilde{y})\right)$$
where $$lk_e(\tilde{z}^+_j-\tilde{z}^-_j, \tilde{y})=\langle z_j,y\rangle.$$
\eop

Similarly, the two-loop polynomial described in \cite{oht2} behaves like $12\Delta(\fvar)\Delta(\svar)\Delta(\tvar)\hat{\CQ}$ in the case of \cite[Proposition 4.18]{oht2}.

\subsection{Proving the Dehn surgery formula from the LP surgery formula}
\label{subDehnproof}
In this subsection, we prove Theorem~\ref{thmDehn} from Theorem~\ref{thmLP}.
Since the proof is completely similar to the proof given \cite[Section 9]{surfor}, we shall be very sketchy and refer to \cite[Section 9]{surfor} for details.

Let $F=(\Sigma \setminus \mbox{Int}(D^2))$ be obtained from the Seifert surface $\Sigma$ of $J$ by removing an open disk with boundary $c$,
let $F\times[-1,2]$ be embedded in a neighborhood of this surface and let
$U$ be a trivial knot in $M$ located in the exterior of $F\times[-1,2]$, that is a meridian of $c$.
Let $U_{\parallel}$ denote the parallel of $U$ in the exterior of $F\times[-1,2]$ that does not link $U$, $J_{\parallel}$ denotes the parallel of $J$ that does not link $J$.
In \cite[Section 9]{surfor},
we define two $LP$-surgeries $(A,A^{\prime})=(A_F=F\times[-1,0],A^{\prime}_F)$ and $(B,B^{\prime})=(B_F=F\times[1,2],B^{\prime}_F)$ such that
\begin{itemize}
\item performing both surgeries does not change $(M,U,U_{\parallel})$,
\item performing one of the surgeries changes $(M,U,U_{\parallel})$ to $(M,\pm J,\pm J_{\parallel})$.
\end{itemize}

This implies that $$M(J;p/q)=M(U;p/q)(A^{\prime}/A)=M(U;p/q)(B^{\prime}/B)$$
and that $$M \sharp S^3(U;p/q)=M(U;p/q)=M(U;p/q)(A^{\prime}/A,B^{\prime}/B).$$

Then, according to Proposition~\ref{propconcas},
$$\CQ(\KK \subset M)=\CQ(\KK \subset M \sharp S^3(U;p/q)) -6\lambda(S^3(U;p/q))$$
and Theorem~\ref{thmLP} allows us to write
$$2\CQ(\KK \subset M(J;p/q)) - 2\CQ(\KK \subset M)=-\CS +12\lambda(S^3(U;p/q))$$
where
$$(-\CS)=\sum_{(i,j,k) \in  \{1, \dots, g_{A}\}^3,(\ell,m,n) \in  \{1, \dots, g_{B}\}^3} lk(i,j,k,\ell,m,n) \CI(i,j,k,\ell,m,n)$$
with notation consistent with the statement of Theorem~\ref{thmLP}.
Here, we use the basis $$\left(c=c \times \{0\},(c^-_j=c_j \times \{-1\} ,d_j=d_j \times \{0\})_{j=1,\dots,g}\right)$$ that plays the role of the basis $z_i$ for $H_1(A)$, and we use the basis $$\left(c^+=c \times \{1\},(c^+_j=c_j \times \{1\} ,d^{++}_j=d_j \times \{2\})_{j=1,\dots,g}\right)$$ that plays the role of the basis $y_i$ for $H_1(B_F)$. Then according to \cite[Lemma 9.3]{surfor} and to the notation of \cite[Section 3]{surfor}, the only triples $(z_i,z_j,z_k)$ for which 
$\CI_{AA^{\prime}}(a_i,a_j,a_k)\neq 0$ are the triples that are obtained from a triple of type $(c,c^-_r,d_r)$ by permutation, and for $(z_i,z_j,z_k)=(c,c^-_r,d_r)$,
$\CI_{AA^{\prime}}(a_i,a_j,a_k)=-1,$ while the only triples $(y_{\ell},y_m,y_n)$ for which 
$\CI_{BB^{\prime}}(b_{\ell},b_m,b_n)\neq 0$ are the triples that are obtained from a triple of type $(c^+,c^+_s,d^{++}_s)$ by permutation, and for $(y_{\ell},y_m,y_n)=(c^+,c^+_s,d^{++}_s)$,
$\CI_{BB^{\prime}}(b_{\ell},b_m,b_n)=1.$
All the lifts in $\tilde{M}$ of the mentioned curves except $c$ and $c^+$ bound in the complement
of the preimage of a ball that contains $U$, $c$ and $c^+$. Therefore the equivariant linking number of a curve among the $c^{(+)}_j$, $d^{(+)}_j$ with another mentioned curve is not affected by the surgery on $U$, and these curves do not link $c$ or $c^+$. Furthermore, $$lk_e((c,c^+)\subset M(U;p/q) )=-q/p$$ Therefore, with the above bases, if $(i,j,k,\ell,m,n)$ contributes, then either $(z_i,y_{\ell})=(c,c^+)$,
or $(z_j,y_{m})=(c,c^+)$, or $(z_k,y_{n})=(c,c^+)$,
and $$(-\CS)=\frac{q}{p}
\sum_{\mathfrak{S}_3(\fvar,\svar,\tvar)}
\sum_{r=1}^g\sum_{s=1}^g \left(\alpha_{rs}(\fvar,\svar) + \alpha_{rs}(\fvar^{-1},\svar^{-1})+\beta_{rs}(\fvar,\svar)\right)$$
where
$\alpha_{rs}(\fvar,\svar)=lk_e(c_r,c^+_s )(\fvar)lk_e(d_r,d^+_s)(\svar) -lk_e(c_r,d^+_s )(\fvar)lk_e(d_r,c^+_s)(\svar)$\\
and $\beta_{rs}(\fvar,\svar)=\left(lk_e(c_r,d^+_r )(\fvar)-lk_e(d^+_r,c_r)(\fvar)\right)\left(lk_e(c_s,d^+_s )(\svar)-lk_e(d^+_s,c_s)(\svar)\right)$.
\eop

An example of application of Theorem~\ref{thmDehn} is the following one.

\begin{example}
By Blanchfield duality (Lemma~\ref{lemblanchnondeg}), for any polynomial $P$ of $\QQ[t_M^{\pm 1}]$, there exists a link $(d_1,d_2)$ such that 
$lk_e(d_1,d_2) = q \frac{P}{\delta(M)}$ for some $q \in \QQ \setminus \{0\}$.
Let $(c_1,c_2)$ be a Hopf link that does not link $(d_1,d_2)$. Frame $c_1$ and $c_2$ trivially so that $lk_e(c_1 ,c_2)=1$, $lk_e(c_1 ,c_1^+)=0$, $lk_e(c_2 ,c_2^+)=0$.
Construct a surface $\Sigma$ by plumbing bands around $c_1$ and $d_1$ and bands around $c_2$ and $d_2$ and by connecting them. This can be achieved so that
$$\begin{array}{ll}\lambda_e^{\prime}(\partial \Sigma)=&\frac{1}{3}( lk_e(d_1 ,d_2)(\fvar^{-1}) + lk_e(d_1 ,d_2)(\svar^{-1}) + lk_e(d_1 ,d_2)(\tvar^{-1}) )\\
 &+\frac{1}{3}( lk_e(d_1 ,d_2)(\fvar) + lk_e(d_1 ,d_2)(\svar) + lk_e(d_1 ,d_2)(\tvar) ) . \end{array}
$$
\end{example}

The next section contains independent preliminaries for the proof of Theorem~\ref{thmLP} that will be given in Section~\ref{secproofLP}.

\newpage
\section{Pseudo-trivialisations}
\setcounter{equation}{0}
\label{secpseudotriv}

\subsection{Introduction to pseudo-trivialisations}

Let $\tau$ be a trivialisation of $M$. Let $A$ be a rational homology handlebody embedded in $M$, and let $A^{\prime}$ be another rational homology handlebody such that $\partial A=\partial {A^{\prime}}$ and $\CL_A=\CL_{A^{\prime}}$. The restriction of $\tau$ on $M \setminus \mbox{Int}(A)$ does not necessarily extend as an actual trivialisation to $A^{\prime}$. (See \cite[Section 4.2]{sumgen}.) It does for the LP surgeries that are involved in the proof of the Dehn surgery formula, and more generally for the Torelli sugeries, or when $A$ and $A^{\prime}$ are 
integral homology handlebodies.

To prove Theorem~\ref{thmLP} in its full generality, we shall make the chains $F_{\cvarM}$ coincide as much as possible. Therefore, like in \cite{sumgen}, we shall
introduce {\em pseudo-trivialisations\/} $\tilde{\tau}$ that will
\begin{itemize}
\item generalize trivialisations,
\item always extend to rational homology handlebodies,
\item induce genuine trivialisations $\tilde{\tau}_{\CC}$ of $TM \otimes_{\RR} \CC$ that have Pontrjagin classes $p_1(\tilde{\tau}_{\CC})$,
\item define $3$-dimensional pseudo-sections $s_{\tilde{\tau}}(M;\cvarM)$ of $ST(M)$, for $\cvarM \in S^2$, so that $$\langle s_{\tilde{\tau}}(M;\fvarM) \cap ST(\Sigma), s_{\tilde{\tau}}(M;\svarM)\rangle_{ST(M)} =0$$ for a two cycle $\Sigma$ when $\fvarM \neq \svarM$ ,
\item thus provide $3$-dimensional cycles 
$$\partial F_{\cvarM}(\tilde{\tau})=s_{\tilde{\tau}}(M;\cvarM) -\ID ST(M)_{|K_{\cvarM}}$$ for any $\cvarM \in S^2$,
and allow us to define $\CQ(\KK;\tilde{\tau})=\langle F_{\fvarM}(\tilde{\tau}),F_{\svarM}(\tilde{\tau}), F_{\tvarM}(\tilde{\tau})\rangle_e$ for generic $\fvarM$, $\svarM$, $\tvarM$ of $S^2$, as in Proposition~\ref{propuninvtripl},
\end{itemize}
so that for such pseudo-trivialisations $\tilde{\tau}$ $$\CQ(\KK)=\CQ(\KK,\tilde{\tau}) - \frac{p_1(\tilde{\tau}_{\CC})}{4}.$$

In this section, we define and prove all we need about these pseudo-trivialisations. This section can be avoided by the reader who is only interested by the proof of Theorem~\ref{thmLP} in the cases where the trivialisations extend to the replacing rational homology handlebodies, except for the independent lemma~\ref{lemtruetriv} that will be used in Section~\ref{secproofLP}.
This section will also be used in Section~\ref{secaug}.

\subsection{Definitions of pseudo-trivialisations and associated notions}
\begin{definition}
\label{defpseudotriv}
A {\em pseudo-trivialisation\/} $\tilde{\tau}=(N(c);\tau_e,\tau_b)$ of $(M,\KK)$ consists of 
\begin{itemize}
\item a link $c$ of $M$ equipped with a neighborhood $N(c)=[a,b] \times c \times [-1,1]$
that avoids $K$,
\item a trivialisation $\tau_e$ of 
$M$ outside $N(c)$ that sends the oriented tangent vectors of $K$ to $\RR^+\qvarM$ for a fixed $\qvarM \in S^2$,
\item a trivialisation $\tau_b$ of $T(N(c))$ such that
$$\tau_b=\left\{\begin{array}{ll} \tau_e & \mbox{over}\; \partial([a,b] \times c \times [-1,1])\setminus (\{a\} \times c \times [-1,1]) \\
\CT_c^{-1} \circ \tau_e  & \mbox{over}\; \{a\} \times c \times [-1,1].
\end{array} \right.$$
\end{itemize}
where 
$$\CT_c^{-1} \circ \tau_e \left(\tau_e^{-1}(t,\gamma \in c,u \in [-1,1];\cvarM\in S^2)\right)=(t,\gamma,u,R_{\alpha(-u)}(\cvarM))$$
where 
$R_{\alpha(-u)}=R_{\alpha(-u),(1,0,0)}$ denotes the rotation of $\RR^3$ 
with axis directed by $(1,0,0)$ and with angle $\alpha(-u)$, and $\alpha$ is a smooth map from $[-1,1]$ to $[0,2 \pi]$ 
that maps $[-1,-1+\varepsilon]$ to $0$,
that increases from $0$ to $2\pi$ on $[-1+\varepsilon, 1-\varepsilon]$, and such that $\alpha(-u)+\alpha(u)=2\pi$ for any $u \in [-1,1]$.
\end{definition}

\begin{lemma}
 Let $A$ be a $\QQ$-handlebody and let $\tau$ be a trivialisation of $TA$ defined on a collar $[-4,0]\times \partial A$ of $\partial A$.
Then there is a pseudo-trivialisation of $A$ that extends the restriction of $\tau$ to $[-1,0]\times \partial A$.
\end{lemma}
\bp
Indeed, there exists a trivialisation $\tau^{\prime}$ of $TA$ on $A$.
There exists a curve $c$ of $\{-2\} \times \partial A$ such that $\tau^{\prime} = \CT_c \circ \tau$ on $\{-2\} \times \partial A$ (after a homotopy of $\tau$ around $\{-2\} \times \partial A$).
Then equip $c$ with the neighborhood $[-2,-1]\times c\times[-1,1]$ and define $\tau_b$ as $\tau$ on $N(c)$, define $\tau_e$ as $\tau$ on $([-2,0]\times \partial A) \setminus \mbox{Int}(N(c))$ and as $\tau^{\prime}$ on $A \setminus (]-2,0]\times \partial A)$.
\eop

\begin{definition}
\label{defpseudotrivpone}[Trivialisation $\tilde{\tau}_{\CC}$ of $TM \otimes_{\RR} \CC$]
Let $F_U$ be a smooth map such that
$$\begin{array}{llll}F_U:&[a,b] \times [-1,1] &\longrightarrow &SU(3)\\
& (t,u) & \mapsto & \left\{\begin{array}{ll}\mbox{Identity}&
 \mbox{if}\; |u|>1-\varepsilon\\
R_{\alpha(u)} &  \mbox{if}\; t<a+\varepsilon\\
\mbox{Identity} &  \mbox{if}\; t>b-\varepsilon.\end{array}\right.\end{array}$$
$F_U$ extends to $[a,b] \times [-1,1]$  because $\pi_1(SU(3))$ is trivial.
Define the trivialisation $\tilde{\tau}_{\CC}$ of $TM \otimes_{\RR} \CC$ as follows.
\begin{itemize}
\item On $T(M\setminus N(c))$, $\tau_{\CC} =\tau_e \otimes 1_{\CC}$,
\item Over $[a,b] \times c \times [-1,1]$, 
$\tau_{\CC}\tau_b^{-1}(t,\gamma,u;\cvarM) =(t,\gamma,u;F_U(t,u)(\cvarM))$.
\end{itemize}
Since $\pi_2(SU(3))$ is trivial, the homotopy class of $\tau_{\CC}$
is well-defined. The definition of $p_1$ for genuine trivialisations in Section~\ref{submanpar} naturally extends to define $p_1(\tilde{\tau}_{\CC})$.
\end{definition}

\begin{definition}[Pseudo-sections $s_{\tilde{\tau}}(.;\cvarM)$]
\label{defpseudosec}
 Let $\varepsilon>0$ be a small positive number and let $F$ be a smooth map such that
$$\begin{array}{llll}F:&[a,b] \times [-1,1] &\longrightarrow &SO(3)\\
& (t,u) & \mapsto & \left\{\begin{array}{ll}\mbox{Identity}&
 \mbox{if}\; |u|>1-\varepsilon\\
R_{\alpha(u)} &  \mbox{if}\; t<a+\varepsilon\\
R_{-\alpha(u)} &  \mbox{if}\; t>b-\varepsilon\end{array}\right.\end{array}$$
where $\alpha$ has been defined in Definition~\ref{defpseudotriv}.
The map $F$ extends to $[a,b] \times [-1,1]$  because its restriction to the boundary
is trivial in $\pi_1(SO(3))$.

Let ${F}(c,\tau_b)$ be defined on $ST(N(c)) \stackrel{\tau_b}{=}   [a,b] \times   c \times [-1,1] \times S^2$ as follows
$$\begin{array}{llll} {F}(c,\tau_b): &  [a,b] \times   c \times [-1,1] \times S^2 & \longrightarrow & S^2\\
&(t,\gamma,u; \cvarM) & \mapsto &  F(t,u)(\cvarM).\end{array}$$

Let $\cvarM \in S^2$ and let $S^1(\cvarM)$ be the circle (or point) in $S^2$ that lies in a plane orthogonal to the axis generated by $(1,0,0)$ and that contains $\cvarM$.
There is a $2$-dimensional chain $C_2(\cvarM)$ in $[-1,1] \times S^1(\cvarM)$ whose
boundary is 
$\{\left(u,R_{-\alpha(u)}(\cvarM)\right), u \in[-1,1]\}
+\{\left(u,R_{\alpha(u)}(\cvarM)\right), u \in[-1,1]\} - 2[-1,1]\times \{\cvarM\}.$
Then $s_{\tilde{\tau}}(M;\cvarM)$ is the following $3$--cycle
$$s_{\tilde{\tau}}(M;\cvarM)=s_{\tau_e}(M\setminus \mbox{Int}(N(c));\cvarM) + 
\frac{s_{\CT_c \circ \tau_b}(N(c);\cvarM) + F(c,\tau_b)^{-1}(\cvarM) + \{b\} \times c \times C_2(\cvarM)}{2}.
$$
When $\Sigma$ is a $2$--chain that intersects $N(c)$ along sections $N_{\gamma}(c)=[a,b]\times\{\gamma\}\times [-1,1],$ $$s_{\tilde{\tau}}(\Sigma;\cvarM)=s_{\tilde{\tau}}(M;\cvarM)\cap ST(M)_{|\Sigma}$$ so that $$s_{\tilde{\tau}}(N_{\gamma}(c);\cvarM)= 
\frac{s_{\CT_c \circ \tau_b}(N_{\gamma}(c);\cvarM) + F(c,\tau_b)^{-1}(\cvarM)\cap ST(M)_{|N_{\gamma}(c)} - \{b\} \times \{\gamma\} \times C_2(\cvarM)}{2}.
$$
\end{definition}

\subsection{Two properties of pseudo-sections.}

The obvious property that {\em genuine sections \/$s_{\tau}(\Sigma;\fvarM)$ and $s_{\tau}(\Sigma;\svarM)$ corresponding to distinct \/$\fvarM$ and \/$\svarM$ of \/$S^2$ are disjoint\/} generalizes as follows for pseudo-sections.

\begin{lemma}
\label{lemintpseudo}
Let $\Sigma$ be a surface embedded in a \/$3$-manifold equipped with a pseudo-trivialisation $\tilde{\tau}=(N(c);\tau_e,\tau_b)$, such that $\Sigma$ only intersects $N(c)$ along sections 
$N_{\gamma_i}(c)=[a,b] \times \{\gamma_i\} \times [-1,1]$ in the interior of $\Sigma$.
Then if $\svarM \in S^2$ and if $\tvarM \in S^2 \setminus S^1(\svarM)$, $s_{\tilde{\tau}}(\Sigma;\svarM)$ and $s_{\tilde{\tau}}(\Sigma;\tvarM)$ do not intersect algebraically in $ST(\Sigma)$.
\end{lemma}
\bp Recall that $S^1(\cvarM)$ denotes the circle in $S^2$ that lies in a plane orthogonal to the axis generated by $e_1=(1,0,0)$ and that contains $\cvarM$.
Let us consider the contribution to $\langle s_{\tilde{\tau}}(\Sigma;\svarM),s_{\tilde{\tau}}(\Sigma;\tvarM)\rangle_{ST(\Sigma)}$ of an intersection point $\gamma$ of $c$ with $\Sigma$.
Such a contribution will read $$\pm \frac{1}{4} \left(\langle s_{\CT_c \circ \tau_b}(N_{\gamma}(c);\svarM), F(c,\tau_b)^{-1}(\tvarM) \rangle_{N_{\gamma}(c) \times S^2} +\langle F(c,\tau_b)^{-1}(\svarM),s_{\CT_c \circ \tau_b}(N_{\gamma}(c);\tvarM) \rangle_{N_{\gamma}(c) \times S^2}\right).$$
We show that both intersection numbers are $\pm 1$ and that their signs are opposite.
Consider an arc $\xi$ of great circle from $e_1$ to $-e_1$.
When $\xi$ stands for $[a,b]$, the map $F$ of Definition~\ref{defpseudosec} above can be seen as the map that maps $(V \in \xi,u \in [-1,1])$ to the rotation $R_{\alpha(u),V}$ with axis $V$ and angle $\alpha(u)$.
Then
$\langle F(c,\tau_b)^{-1}(\tvarM),s_{\CT_c \circ \tau_b}(N_{\gamma}(c);\svarM) \rangle$ is the degree of the map 
$$\begin{array}{llll}f_{\svarM}\colon &\xi \times [-1,1]&\rightarrow &S^2\\&(V,u)&\mapsto &R_{\alpha(u),V}\circ R_{\alpha(-u),e_1} (\svarM)\end{array}$$ at $\tvarM$ while the other intersection
is the degree of the map $f_{\tvarM}$ at $\svarM$.
Since the boundary of the image of $f_{\svarM}$ is $\pm 2 S^1(\svarM)$, the degree jumps by $\pm 2$ from one component of $S^2 \setminus S^1(\svarM)$ to the other one. On the other hand, the degree of $f_{\svarM}$ on the component of $e_1$ is independent of $\svarM\neq e_1$, and the degree on the component of $-e_1$ is independent of $\svarM \neq (-e_1)$. Therefore, using Lemma~\ref{lemdegstwoso}, we easily see that the degrees are $\pm 1$ and that they are opposite.
\eop

\begin{lemma}
\label{lemtruetriv} 
Let $e_1=(1,0,0) \in \RR^3$.
Let $\Sigma$ be a surface immersed in a \/$3$-manifold $M$ equipped with a trivialisation \/$\tau :TM \rightarrow M\times \RR^3$ such
that \/$\tau^{-1}(.;e_1)$ is a positive normal to $\Sigma$ along the boundary \/$\partial \Sigma$ of $\Sigma$.
Let \/$a^{(1)}$, \dots, $a^{(k)}$ denote the \/$k$ connected components of the boundary \/$\partial \Sigma$ of $\Sigma$. 
For \/$i=1,\dots,k$, the unit bundle of $T\Sigma_{|a^{(i)}}$
is an \/$S^1$-bundle over \/$a^{(i)}$ with a canonical trivialisation induced by \/$Ta^{(i)}$. For a trivialisation \/$\tau$ as above, let \/$d(\tau,a^{(i)})$ be the degree of the projection on the fiber \/$S^1$ of this bundle of the section \/$\tau^{-1}(a^{(i)}\times e_2)$.
Let \/$\mbox{diag}(n)(\Sigma) \subset ST(M)$ be the section of \/$ST(M)_{|\Sigma}$ in \/$ST(M)$ associated with the positive normal field \/$n$ to \/$\Sigma$, that
coincides with \/$\tau^{-1}(\partial \Sigma \times e_1)$ on $\partial \Sigma$.
Then $$2(\mbox{diag}(n)(\Sigma)-s_{\tau}(\Sigma; e_1)) - \left(\sum_{i=1}^k d(\tau,a^{(i)}) +\chi(\Sigma)\right)ST(\ast) $$ is a cycle that is null-homologous in \/$ST(\Sigma)$.
\end{lemma}
\bp
This comes from the fact that the relative Euler class of $s_{\tau}(\partial \Sigma; e_2)$ in $T\Sigma$ is $$\chi(\Sigma)+\sum_{i=1}^kd(\tau,a^{(i)}),$$ up to sign. Let us give some more details. Proposition~\ref{propsplusstau} gives the result when $\partial \Sigma =\emptyset$.
Since the generator of $\pi_1(SO(3))$ can be realized by rotations around the $e_1$-axis and since $\pi_2(SO(3))=0$, the homology class of
$\left(\mbox{diag}(n)(\Sigma) -\tau^{-1}(\Sigma\times e_1)\right)$ does not depend on the trivialisation $\tau$ of $(T\Sigma \oplus \RR)$ when the $d(\tau,a^{(i)})$ are fixed.
Then it is easy to find an embedding of
$\Sigma$ in $\RR^3$ for which $d(\tau,a^{(i)})=1$ for all $i$ with the standard trivialisation of $\RR^3$ (a standard one obtained from a closed surface by drilling small holes) and another one for which $d(\tau,a^{(i)})=-1$ for all $i$ by gluing the following pieces to the hole boundaries.

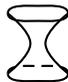
\begin{center}
\begin{pspicture}[shift=0.5](0,1)(0,1) 
\psecurve(.5,.2)(.1,.1)(.5,0)(.9,.1)(.5,.2)
\psecurve[linestyle=dashed](.5,0)(.1,.1)(.5,.2)(.9,.1)(.5,0)
\psccurve(.5,1)(.1,.9)(.5,.8)(.9,.9)
\pscurve(.1,.1)(.4,.5)(.1,.9)
\pscurve(.9,.1)(.6,.5)(.9,.9)
\end{pspicture}
\end{center}

Thus, the lemma is easy to prove when $\tau$ is an actual trivialisation and when $d(\tau,a^{(i)})$ does not depend on $i$ and is $\pm 1$, by computing degrees of Gauss maps.

Let $A$ be an annulus that is a regular neighborhood of a boundary component of $\Sigma$. Let $h_A \colon A \rightarrow SO(3)$ be a homotopy between the trivial loop of $SO(3)$ and a loop of rotations of $SO(3)$ around $e_1$ that has degree $2$.
Since the composition of $\tau$ by $h_A$ sends $\tau^{-1}(A\times e_1)$ to the sphere with degree $\pm 1$ (see Lemma~\ref{lemdegstwoso}), and because the previous argument fixes the correct sign, we easily conclude that the cycle of the statement is null-homologous when the $d(\tau,a^{(i)})$ are odd.

Then twisting the trivialisation around $e_1$ across paths in $\Sigma$ from one boundary component to another one allows us to prove that it is null-homologous for arbitrary $d(\tau,a^{(i)})$.
\eop

\begin{lemma}
\label{lempseudotriv}
Lemma~\ref{lemtruetriv} is true if $\tau$ is a pseudo-trivialisation that is a genuine trivialisation satisfying the assumptions of Lemma~\ref{lemtruetriv} around $\partial \Sigma$, such that 
$\Sigma$ meets $N(c)$ along sections $N_{\gamma}(c)$, where $s_{\tau}(\Sigma;e_1)$ is defined in Definition~\ref{defpseudosec}.
\end{lemma}
\bp
Since the formula of Lemma~\ref{lemtruetriv} behaves well under gluing (or cutting) surfaces along curves that satisfy the assumptions of Lemma~\ref{lemtruetriv} and since these assumptions are easily satisfied by isotopy (when the cutting process is concerned), we are left with the proof of the lemma for a meridian disk $\Sigma$ of
$c$ in Definition~\ref{defpseudotriv}. For such a disk $s_{\tau}(\Sigma; e_1)$
is the average of two genuine sections corresponding to trivialisations $\tau_1$ and $\tau_2$ so that
$$2(\mbox{diag}(n)(\Sigma)-s_{\tau}(\Sigma; e_1)) - \left(\frac{d(\tau_1,\partial \Sigma)+d(\tau_2,\partial \Sigma)}{2}  +\chi(\Sigma)\right)ST(\ast)$$
is a null-homologous cycle.
Since the exterior trivialisation $\tau_e$ is such that
$$d(\tau_e,-\partial \Sigma) =\frac{d(\tau_1,-\partial \Sigma)+d(\tau_2,-\partial \Sigma)}{2},$$ we are done.
\eop

\subsection{Making the definition of $\CQ$ more flexible with pseudo-sections}

\begin{lemma}
\label{lemhomhthreepseudo}
Let $\tilde{\tau}$ be a pseudo-trivialisation as in Definition~\ref{defpseudotriv}, and let $\tau$ be a genuine trivialisation of $M$ that coincides with $\tau_e$ on a tubular neighborhood of $K$.
For any $\cvarM \in S^2$, there exists a $4$-chain $C_4(\cvarM)$
of $[0,1] \times ST(E=M \setminus \mbox{Int}(N(K)))$ such that $\partial C_4(\cvarM) =c_3(\cvarM)$ where
$$c_3(\cvarM)= \{1\} \times s_{\tilde{\tau}}(ST(E);\cvarM)  -\{0\} \times s_{\tau}(ST(E);\cvarM)
+[0,1] \times \partial N(K) \times_{\tau} \{\cvarM\}.$$
\end{lemma}
\bp
Like in the proof of Lemma~\ref{lemhomhthree}, it suffices to prove that the algebraic intersection
$$\langle c_3(\cvarM),[0,1] \times (S \setminus \mbox{Int}(D^2))\times_{\tau} \qvarM \rangle_{[0,1] \times ST(E)}$$ is zero for some $\qvarM \in S^2\setminus S^1(\cvarM)$, where $S$ is assumed to meet $N(c)$ along sections $N_{\gamma}(c)$. Up to sign, this intersection is 
$$\langle s_{\tilde{\tau}}(S \setminus \mbox{Int}(D^2);\cvarM),(S \setminus \mbox{Int}(D^2))\times_{\tau} \qvarM \rangle_{ST(S \setminus \mbox{\small Int}(D^2))}.$$ Since
Lemma~\ref{lempseudotriv} implies that the cycle $\left(s_{\tilde{\tau}}(S \setminus \mbox{Int}(D^2);\cvarM)- s_{\tau}(S \setminus \mbox{Int}(D^2);\cvarM)\right)$ bounds in $ST(S \setminus \mbox{Int}(D^2))$, this algebraic intersection vanishes.
\eop

Lemma~\ref{lemhomhthreepseudo} allows us to construct a $4$--chain $F_{\cvarM}(\tilde{\tau})$
with the wanted boundary $$\partial F_{\cvarM}(\tilde{\tau})=s_{\tilde{\tau}}(M;\cvarM) -\ID ST(M)_{|K_{\cvarM}}$$
from our $4$-chain $F_{\cvarM}(\tau)$ previously associated to $\tau$,
(using modifications near $ST(E) \subset \partial C_2(M)$) such that $\langle F_{\cvarM}(\tilde{\tau}), A(K) \rangle_e=0$.

Note that as soon as the circles $S^1(\fvarM)$, $S^1(\svarM)$ and $S^1(\tvarM)$ are disjoint, the chains $F_{\fvarM}(\tilde{\tau})$, $F_{\svarM}(\tilde{\tau})$ and $F_{\tvarM}(\tilde{\tau})$ do not have triple intersection on the boundary.

As before, we define
$$\CQ(\KK,\tilde{\tau})=\langle F_{\fvarM}(\tilde{\tau}), F_{\svarM}(\tilde{\tau}), F_{\tvarM}(\tilde{\tau}) \rangle_e.$$

Because of the specific construction outlined above,
$$\CQ(\KK,\tilde{\tau})-\CQ(\KK,\tau)=\langle C_4(\fvarM), C_4(\svarM), C_4(\tvarM)\rangle_{[0,1] \times ST(E)}.$$

We shall prove that the following lemma follows from \cite{sumgen}.
\begin{lemma}
\label{lempsetriv}
$$\langle C_4(\fvarM), C_4(\svarM), C_4(\tvarM)\rangle_{[0,1] \times ST(E)}=\frac{p_1(\tilde{\tau}_{\CC})-p_1(\tau)}{4}.$$
\end{lemma}
Together with Proposition~\ref{propvartau}, this lemma obviously implies the following proposition.
\begin{proposition}
Under the assumptions above, 
$$\CQ(\KK)=\CQ(\KK,\tilde{\tau}) - \frac{p_1(\tilde{\tau}_{\CC})}{4}.$$
\end{proposition}

\begin{definition}
\label{defpseudotrivrat}
A {\em pseudo-trivialisation\/} $\tilde{\tau}=(N(c);\tau_e,\tau_b)$ {\em of a rational homology ball\/} $B(N)$ associated with a rational homology sphere $N$ (as in Subsection~\ref{subdefcas}) consists of a framed link $c$ of $B(N)$ equipped with a neighborhood $N(c)$ in the interior of $B_N$, a trivialisation $\tau_e$ of $\left(B(N) \setminus \mbox{Int}(N(c))\right)$ that is standard on $\left(B(N) \setminus \mbox{Int}(B_N)=B(3) \setminus \mbox{Int}(B(1))\right)$, and a trivialisation $\tau_b$ of $T(N(c))$ where $\tau_e$ and $\tau_b$ are related like in Definition~\ref{defpseudotriv}.

Then the {\em trivialisation\/} $\tilde{\tau}_{\CC}$ of $TB(N) \otimes_{\RR} \CC$
is defined like in Definition~\ref{defpseudotrivpone} and the definition
of the {\em Pontrjagin class\/} for genuine trivialisations standard on $B(N) \setminus \mbox{Int}(B_N)$, outlined in Subsection~\ref{subdefcas} and given in details in \cite[Section 1.5]{lesconst}, naturally extends to pseudo-trivialisations of rational homology balls.

Finally, define a {\em $3$-cycle\/} $\partial F_{N,\cvarM}(\tilde{\tau})$ as $p_N^{-1}(\cvarM)$ on $\partial C_2(B(N)) \setminus \mbox{Int}(ST(N(c)))$, with respect to the map $p_N$ of Subsection~\ref{subdefcas}, and as $s_{\tilde{\tau}}(N(c);\cvarM)$ on $ST(N(c))$, like in Definition~\ref{defpseudosec}.
\end{definition}

Recall that $S^1(\cvarM)$ denotes the circle through $\cvarM$ in $S^2$ that lies in a plane orthogonal to the axis generated by $(1,0,0)$.

\begin{proposition}
\label{propdefcaspseudo}
Let $\fvarM$, $\svarM$ and $\tvarM$ be three distinct points of $S^2$ such that the circles $S^1(\fvarM)$, $S^1(\svarM)$ and $S^1(\tvarM)$ are disjoint.
With the above notation, for $\cvarM=\fvarM$, $\svarM$ or $\tvarM$, $\partial F_{N,\cvarM}(\tilde{\tau})$ bounds a rational chain $F_{N,\cvarM}(\tilde{\tau})$
in $C_2(B(N))$, and
$$\lambda(N)=\frac{\langle F_{N,\fvarM}(\tilde{\tau}),F_{N,\svarM}(\tilde{\tau}), F_{N,\tvarM}(\tilde{\tau})\rangle_{C_2(B(N))}}{6} -\frac{p_1(\tilde{\tau})}{24}.$$
\end{proposition}
\bp
It is essentially a consequence of \cite[Proposition 4.8]{sumgen}. We nevertheless give a few details to see how this proposition applies.
According to \cite[Theorem 2.6]{sumgen}, $Z_1(N)=\frac{\lambda(N)}{2}[\theta]$.
According to \cite[Proposition 4.8]{sumgen}, 
$$Z_1(N)=Z_1(N;\omega(c,\tau_e,\tau_b))+\frac{p_1(\tilde{\tau})}{4}\xi_1=\frac{\int_{C_2(B(N))}\omega(c,\tau_e,\tau_b)^3}{12}[\theta]+\frac{p_1(\tilde{\tau})}{4}\xi_1$$
where $\xi_1=-\frac{1}{12}[\theta]$ according to \cite[Proposition 2.45]{lesconst}, and
$\omega(c,\tau_e,\tau_b)$ is defined as $p_N^{\ast}(\omega_{S^2})$ on $\partial C_2(B(N)) \setminus \mbox{Int}(ST(N(c)))$, and like in \cite[Notation 4.9]{sumgen} on $ST(N(c))$ (that is like in \cite[Notation 4.9]{sumgen}, except that $C_2(N)$ is replaced by $C_2(B(N))$ like in the proof of Theorem~\ref{thmdefcasconf}).
Therefore $\lambda(N)=\frac{\int_{C_2(B(N))}\omega(c,\tau_e,\tau_b)^3}{6}-\frac{p_1(\tilde{\tau})}{24}.$
Now, we can change $\omega(c,\tau_e,\tau_b)^3$ in the above integral to
$\omega_{\fvarM} \wedge \omega_{\svarM} \wedge \omega_{\tvarM}$ where the closed form $\omega_{\cvarM}$ is Poincar\'e dual to $F_{N,\cvarM}(\tilde{\tau})$ and supported in a small neighborhood of $F_{N,\cvarM}(\tilde{\tau})$, like in \cite[Lemma 6.15]{sumgen} for instance.
\eop

\noindent{\sc Proof of Lemma~\ref{lempsetriv}:}
The trivialisation $\tau$ of $M \setminus \mbox{Int}(N(K))$ extends to a rational homology sphere $M^{\prime}$ obtained by Dehn filling with respect to a parallel of $K$ (that differs from the given parallel of $K$ by adding a meridian). Call this trivialisation $\tau(M^{\prime})$ and assume without loss that, after removing from $M^{\prime}$ an open ball inside the new solid torus $T^{\prime}$ of $M^{\prime}$, the boundary of the obtained rational homology ball $B(M^{\prime})$ can be identified with the boundary of the ball of radius $3$ of $\RR^3$ so that $\tau(M^{\prime})$ coincides with the standard trivialisation of $\RR^3$ near the boundary.
Then $\left(p_1(\tilde{\tau}_{\CC})-p_1(\tau)\right)$ is the obstruction to extend
the $SU(4)$-trivialisation 
induced by 
\begin{itemize}
\item $\tilde{\tau}_{\CC}$ (and $\tau(M^{\prime})$ on $T^{\prime}$) on $\{1\} \times M$, $\{1\} \times M^{\prime}$ or $\{1\} \times B(M^{\prime})$,
\item $\tau$ (and $\tau(M^{\prime})$ on $T^{\prime}$) on $\{0\} \times M$, $\{0\} \times M^{\prime}$ or $\{0\} \times B(M^{\prime})$, 
\end{itemize}
to $[0,1] \times M$, $[0,1] \times M^{\prime}$ or $[0,1] \times B(M^{\prime})$, respectively.
Now, as before
$$
\langle F_{M^{\prime},\fvarM}(\tilde{\tau}),F_{M^{\prime},\svarM}(\tilde{\tau}), F_{M^{\prime},\tvarM}(\tilde{\tau})\rangle_{C_2(B(M^{\prime}))} -\langle F_{M^{\prime},\fvarM}({\tau}),F_{M^{\prime},\svarM}({\tau}), F_{M^{\prime},\tvarM}({\tau})\rangle_{C_2(B(M^{\prime}))}$$
$$=\langle C_4(\fvarM), C_4(\svarM), C_4(\tvarM)\rangle_{[0,1] \times ST(E)}.$$
Thus, Proposition~\ref{propdefcaspseudo} allows us to conclude.
\eop

\newpage
\section{Proof of the LP surgery formula}
\setcounter{equation}{0}
\label{secproofLP}

\subsection{Sketch of proof of Theorem~\ref{thmLP}}

The restriction of any trivialisation $\tau$ of $M$ to $M \setminus \mbox{Int}(A\cup B)$ extends as a pseudo-trivialisation to $A^{\prime} \cup B^{\prime}$. This induces a pseudo-trivialisation $\tau(M_{AB})$
on $M_{AB}$, and pseudo-trivialisations
$\tau(M_A)$ and $\tau(M_B)$ on $M_A$ and $M_B$, respectively,
that coincide with each other and with $\tau$ as much as possible.
The chains $F_{\cvarM}(M_S)$ associated with these manifolds $M_S$ will be associated with these pseudo-trivialisations so that
\begin{itemize}
\item $\partial F_{\cvarM}(M_S)=\partial F_{\cvarM}(\tau(M_S))$,
\item $\CQ(\KK \subset M_S)=\langle F_{\fvarM}(M_S),F_{\svarM}(M_S),F_{\tvarM}(M_S) \rangle_e -\frac{p_1(\tau(M_S)_{\CC})}{4}$
\item $p_1(\tau=\tau(M))-p_1(\tau(M_A)_{\CC})=p_1(\tau(M_B)_{\CC})-p_1(\tau(M_{AB})_{\CC})$
\end{itemize}
The last equality is proved like
\cite[Lemma 3.2]{sumgen} (where $Y$ is a $3$-manifold with fixed rational homology --instead of a $\QQ$-handlebody-- and fixed Lagrangian). It implies that the Pontrjagin classes can be forgotten
in the expression of the alternate sum 
$\CS$ that we are about to compute using the above expressions of the
$\CQ(\KK \subset M_S)$.

For $i=1,\dots,g_A$, let $\Sigma(a_i)$ and $\Sigma^{\prime}(a_i)$ be rational chains of $A$ and $A^{\prime}$, respectively, such that
$\partial(\Sigma(a_i))=\partial(\Sigma^{\prime}(a_i))=a_i$.
Then $\CI_{AA^{\prime}}(a_i\wedge a_j \wedge a_k)$ is the algebraic intersection $$\langle \Sigma(a_i) \cup (-\Sigma^{\prime}(a_i)),\Sigma(a_j) \cup (-\Sigma^{\prime}(a_j)),\Sigma(a_k) \cup (-\Sigma^{\prime}(a_k)) \rangle_{A\cup_{\partial A}-A^{\prime}}.$$
Similarly, for $i=1,\dots,g_B$, let $\Sigma(b_i)$ and $\Sigma^{\prime}(b_i)$ be rational chains of $B$ and $B^{\prime}$, respectively, such that
$\partial(\Sigma(b_i))=\partial(\Sigma^{\prime}(b_i))=b_i$.
For a curve $c$ among the $a_i$ or the $b_i$, let $\Sigma_{\fvarM}(c)$, $\Sigma_{\svarM}(c)$ and $\Sigma_{\tvarM}(c)$
denote three parallel copies of $\Sigma(c)$, and let $\Sigma^{\prime}_{\fvarM}(c)$, $\Sigma^{\prime}_{\svarM}(c)$ and $\Sigma^{\prime}_{\tvarM}(c)$ similarly
denote three parallel copies of $\Sigma^{\prime}(c)$ such that $\partial(\Sigma_{\cvarM}(c))=\partial(\Sigma^{\prime}_{\cvarM}(c))$ is parallel to $c$, for $\cvarM = \fvarM$, $\svarM$ or
$\tvarM$. Fix preferred lifts for $A$ and $B$ in $\tilde{M}$, the preferred lifts of the curves $c$ and the chains $\Sigma_{\cvarM}(c)$ in the preferred lifts of $A$ and $B$ are also denoted by 
$c$ or $\Sigma_{\cvarM}(c)$.

For $i=1, \dots, g_A$, $\delta(t_M)z_i$ bounds a rational chain in $\left(\tilde{M} \setminus p_M^{-1}(K)\right)$. (The possible intersections with $p_M^{-1}(K)$ can be removed with the help of $S$.) Therefore, it cobounds a rational cycle $\delta(t_M)\Sigma(\check{z}_i)$ in $\left(\tilde{M} \setminus p_M^{-1}(K \cup \mathring{A})\right)$ with a combination of $a_j$ that has its coefficients in $\Lambda_M$.
$$\partial \Sigma(\check{z}_i)= z_i - \sum_{j=1}^{g_A} lk_e(z_i,z^-_j)a_j=\check{z}_i.$$
We can furthermore assume that
$$\Sigma(\check{z}_i) \cap p_M^{-1}(B)=\sum_{j=1}^{g_B}lk_e(z_i,y_j)\Sigma(b_j).$$
Define $\Sigma^{\prime}(\check{z}_i) \subset \left(\widetilde{M_B} \setminus p_{M_B}^{-1}(K \cup \mathring{A})\right)$ from $\Sigma(\check{z}_i)$
by replacing the pieces $\Sigma(b_j)$ by pieces $\Sigma^{\prime}(b_j)$.

Similarly, construct $\Sigma(\check{y}_i)$ in $\left(\widetilde{M} \setminus p_M^{-1}(K \cup \mathring{B})\right)$, for $i=1, \dots, g_B$, such that
$$\partial \Sigma(\check{y}_i)= y_i - \sum_{j=1}^{g_B} lk_e(y_i,y^-_j)b_j=\check{y}_i$$ and
$$\Sigma(\check{y}_i) \cap p_M^{-1}(A)=\sum_{j=1}^{g_A}lk_e(y_i,z_j)\Sigma(a_j)$$
and define $\Sigma^{\prime}(\check{y}_i) \subset \left(\widetilde{M_A} \setminus p_{M_A}^{-1}(K \cup \mathring{B})\right)$ from $\Sigma(\check{y}_i)$
by replacing the pieces $\Sigma(a_j)$ by pieces $\Sigma^{\prime}(a_j)$.

When $A$ is a submanifold of a closed $3$-manifold $M$ with first Betti number one, $C_2(A)$ will denote the preimage of $A^2$ under the natural projection $\projconf \colon C_2(M) \rightarrow M^2$, and $\tilde{C}_2(A)$ will be the preimage of $A^2$ under the natural projection $\tilde{\projconf} \colon \TCM \rightarrow M^2$.

\begin{proposition}
\label{proplprecap}
We can assume that the fundamental chains $F_{\cvarM}(M_S)$ for the manifolds
$M$, $M_A=M(A^{\prime}/A)$, $M_B=M(B^{\prime}/B)$, $M_{AB}=M(A^{\prime}/A,B^{\prime}/B)$, and for $\cvarM = \fvarM$, $\svarM$ or
$\tvarM$, satisfy the following properties.

\begin{enumerate}
\item The chains $F_{\cvarM}(M)$ intersect
\begin{itemize}
 \item $A \times B$ as
$\sum_{(i,j) \in \{1, \dots, g_A\} \times \{1, \dots, g_B\}} \overline{lk_e(z_i,y_j)} \Sigma_{\cvarM}(a_i) \times \Sigma_{\cvarM}(b_j)$
and \item $B \times A$ as
$\sum_{(i,j) \in \{1, \dots, g_A\} \times \{1, \dots, g_B\}} lk_e(z_i,y_j) \Sigma_{\cvarM}(b_j) \times \Sigma_{\cvarM}(a_i)$.
\end{itemize}
and the chains 
$F_{\cvarM}(M_A)$, $F_{\cvarM}(M_{B})$ and 
$F_{\cvarM}(M_{AB})$ have similar intersections with $A^{(\prime)} \times B^{(\prime)}$ and $B^{(\prime)} \times A^{(\prime)}$ where appropriate primes are added.

\item Set $$ST(\Sigma_{AA^{\prime}}) = \sum_{(i,j,k)\in \{1,\dots g_A\}^3} \CI_{AA^{\prime}}(a_i\wedge a_j \wedge a_{k})lk_e(z_{k},z_j) ST(M)_{|p_M(\Sigma(\check{z}_{i}))} $$
and $$ST(\Sigma_{BB^{\prime}}) = \sum_{(\ell,m,n)\in \{1,\dots g_B\}^3} \CI_{BB^{\prime}}(b_{\ell}\wedge b_m \wedge b_n)lk_e(y_{n},y_m) ST(M)_{|p_M(\Sigma(\check{y}_{\ell}))}$$
Similarly define $$ST(\Sigma^{\prime}_{AA^{\prime}}) \subset \tilde{C}_2\left(M_B \setminus (\mathring{A} \cup K)\right)$$ and 
$$ST(\Sigma^{\prime}_{BB^{\prime}}) \subset \tilde{C}_2\left(M_A \setminus (\mathring{B} \cup K)\right)$$
by replacing $\Sigma$ by $\Sigma^{\prime}$ in the right-hand sides.
Let $ \ZZ \times ST(M)\times [0,1]=\partial \tilde{C}_2(M)\times [0,1]$ denote an equivariant neighborhood of $\partial \tilde{C}_2(M)$ in $\tilde{C}_2(M)$ and, for $\eta \in ]0,1]$, let $ST_{\eta}(\Sigma_{..})$ denote
$ST(\Sigma_{..}) \times \{\eta\}$ in this neighborhood.\\
Let $\eta \in ]0,1/6[$, let $r(\fvarM)=0$, $r(\svarM)=1$ and $r(\tvarM)=2$.

\begin{itemize}
\item $F_{\cvarM}(M_A)=F_{\cvarM}(M)+ST_{(3-r(\cvarM))\eta}(\Sigma_{AA^{\prime}})$ on $\tilde{C}_2(M \setminus \mathring{A})$,
\item $F_{\cvarM}(M_{AB})=F_{\cvarM}(M_{B})+ST_{(3-r(\cvarM))\eta}(\Sigma^{\prime}_{AA^{\prime}})$ on $\tilde{C}_2(M_{B} \setminus \mathring{A})$,
\item $F_{\cvarM}(M_{B})=F_{\cvarM}(M)+ST_{(6-r(\cvarM))\eta}(\Sigma_{BB^{\prime}})$ on $\tilde{C}_2(M \setminus \mathring{B})$,
\item $F_{\cvarM}(M_{AB})=F_{\cvarM}(M_A)+ST_{(6-r(\cvarM))\eta}(\Sigma^{\prime}_{BB^{\prime}})$ on $\tilde{C}_2(M_{A} \setminus \mathring{B})$.
\item For any permutation $\sigma$ of $\{\fvarM,\svarM,\tvarM\}$,
\begin{itemize} 
\item $p(ST_{(3-r(\sigma(\fvarM)))\eta}(\Sigma_{AA^{\prime}})) \cap
p(F_{\sigma(\svarM)}(M) )\cap p( F_{\sigma(\tvarM)}(M)) =\emptyset$
\item $p(ST_{(6-r(\sigma(\fvarM)))\eta}(\Sigma_{BB^{\prime}})) \cap
p(F_{\sigma(\svarM)}(M) )\cap p( F_{\sigma(\tvarM)}(M)) =\emptyset$
\item The intersection of $ST_{(3-r(\sigma(\fvarM)))\eta}(\Sigma^{\prime}_{AA^{\prime}})$, 
$F_{\sigma(\svarM)}(M_B) $ and $ F_{\sigma(\tvarM)}(M_B)$ is located in the interior of $p^{-1}(C_2(B^{\prime}))$ and their algebraic equivariant intersection vanishes.
\item The intersection of $ST_{(6-r(\sigma(\fvarM)))\eta}(\Sigma^{\prime}_{BB^{\prime}})$, 
$F_{\sigma(\svarM)}(M_A)$ and $ F_{\sigma(\tvarM)}(M_A)$ is located in the interior of $p^{-1}(C_2(A^{\prime}))$, their algebraic equivariant intersection is
$$\sum_{(i,j,k)\in \{1,\dots g_A\}^3,(\ell,m,n)\in \{1,\dots g_B\}^3} \CI(i,j,k,\ell,m,n)L_{\sigma}(i,j,k,\ell,m,n)$$
where $\CI(i,j,k,\ell,m,n)=\CI_{AA^{\prime}}(a_i\wedge a_j \wedge a_{k})\CI_{BB^{\prime}}(b_{\ell}\wedge b_m \wedge b_n)
$, $\sigma$ permutes $\fvar$, $\svar$ and $\tvar$ like $\fvarM$, $\svarM$ and $\tvarM$ and 
$$L_{\sigma}(i,j,k,\ell,m,n)=$$
$$-lk_e(y_{n},y_m)(\sigma(\fvar))\left(lk_e(z_{k},z_j)(\sigma(\svar))+lk_e(z_{k},z_j)(\sigma(\tvar))\right)lk(z_i,y_{\ell}).$$
\end{itemize}
\end{itemize}
\end{enumerate}
\end{proposition}
The proof of this key proposition will start in the next subsection and end at the end of the section.
Assuming it, let us conclude the proof of Theorem~\ref{thmLP}.

For a part $C$ of $C_2(M_S)$, let $$\CI_S(C)=\langle F_{\fvarM}(M_S),F_{\svarM}(M_S),F_{\tvarM}(M_S)\rangle_{e,p^{-1}(C)}$$
denote the restriction of the intersection to this part.
Then 
$$\CI_B(C_2(M_B \setminus \mathring{B}^{\prime}))=\CI(C_2(M \setminus \mathring{B}))$$ and the equivariant intersection points are actually the same in both sides
because the pieces $ST_{(6-r(\sigma(\fvarM)))\eta}(\Sigma_{BB^{\prime}})$ are pairwise disjoint and they do not meet $\left(F_{\sigma(\svarM)}(M) \cap  F_{\sigma(\tvarM)}(M)\right)$.
Similarly, $$\CI_{AB}(C_2(M_{AB} \setminus \mathring{B}^{\prime}))-\CI_{A}(C_2(M_A \setminus \mathring{B}))=(\CI_{AB}-\CI_A)(C_2(A^{\prime}))$$
$$ =-\sum_{(i,j,k) \in  \{1, \dots, g_A\}^3,\,(\ell,m,n) \in  \{1, \dots, g_B\}^3,\,\sigma \in \mathfrak{S}_3(\fvar,\svar,\tvar)} \CI(i,j,k,\ell,m,n)lk(z_i, y_{\ell})lk_e(z_j, z_k)(\fvar)lk_e(y_{m}, y_{n})(\svar)$$
and $$(\CI_{AB}-\CI_A-\CI_B+\CI)\left(C_2(M_S \setminus \mathring{B}^{(\prime)})\right)
=(\CI_{AB}-\CI_A)(C_2(A^{\prime})).$$
(The parentheses around the prime mean that there can be a prime.)

Consider the part $$C(B^{(\prime)},M_S \setminus \mathring{A}^{(\prime)})=C_2(M_S \setminus \mathring{A}^{(\prime)})\cap \overline{C_2(M_S)\setminus C_2(M_{S} \setminus {B}^{(\prime)})}$$ of $C_2(M_S)$
made of the configurations that project to $M_S^2$ as pairs of points such that
at least one of the points is in $B^{(\prime)}$ and both of them are in $M_S \setminus \mbox{Int}(A^{(\prime)})$.

Note that $\CI_A$ and $\CI$ coincide there with the same actual equivariant intersection like in the whole $C_2(M_S \setminus \mbox{Int}(A^{(\prime)}))$, while the equivariant intersection points contributing to $(\CI_{AB}-\CI_B)(C(B^{\prime},M_B \setminus \mbox{Int}(A)))$ project to the interior of $C_2(B^{\prime})$ and yield a trivial contribution. Thus, 
$$(\CI_{AB}-\CI_A-\CI_B+\CI)\left(C_2(M_S \setminus \mathring{A}^{(\prime)}) \cup C_2(M_S \setminus \mathring{B}^{(\prime)})\right)=(\CI_{AB}-\CI_A)(C_2(A^{\prime})).$$

Since the closure of the complement of $C_2(M_S \setminus \mbox{Int}(A^{(\prime)})) \cup C_2(M_S \setminus \mbox{Int}(B^{(\prime)}))$
in $C_2(M_S)$
is $\left(A^{(\prime)} \times B^{(\prime)} \right)\cup \left(B^{(\prime)} \times A^{(\prime)}\right)$, 
$$\left(\CS-(\CI_{AB}-\CI_A)(C_2(A^{\prime}))\right)$$
is the contribution of the intersection points that project to $A^{(\prime)} \times B^{(\prime)}$ or
$B^{(\prime)}\times A^{(\prime)}$.
Then the remaining contributions can be seen as contributions of intersections
in\\ $p^{-1}\left(\left(A\cup(-A^{\prime})\right) \times \left(B\cup(-B^{\prime})\right)\right)$ and intersections
in $p^{-1}\left(\left(B\cup(-B^{\prime})\right) \times \left(A\cup(-A^{\prime})\right)\right)$.
The intersections in $p^{-1}\left(\left(A\cup(-A^{\prime})\right) \times \left(B\cup(-B^{\prime})\right)\right)$ contribute as
$$\sum_{((i,\ell),(j,m),(k,n)) \in  (\{1, \dots, g_A\} \times \{1, \dots, g_B\})^3} lk_e(z_i, y_{\ell})(\fvar^{-1})lk_e(z_j, y_{m})(\svar^{-1})lk_e(z_k, y_{n})(\tvar^{-1}) \CI_{ijk\ell mn}$$

where
$\CI_{ijk\ell mn}$ is the algebraic triple intersection of
$(\Sigma(a_i) \cup -\Sigma^{\prime}(a_i)) \times (\Sigma(b_{\ell}) \cup -\Sigma^{\prime}(b_{\ell}))$,
$(\Sigma(a_j) \cup -\Sigma^{\prime}(a_j)) \times (\Sigma(b_{m}) \cup -\Sigma^{\prime}(b_{m}))$ and
$(\Sigma(a_k) \cup -\Sigma^{\prime}(a_k)) \times (\Sigma(b_{n}) \cup -\Sigma^{\prime}(b_{n}) )$
that is 
$$\CI_{ijk\ell mn}=-\CI_{AA^{\prime}}(a_i\wedge a_j \wedge a_k)\CI_{BB^{\prime}}(b_{\ell}\wedge b_m \wedge b_n).$$

The intersections in $p^{-1}\left(\left(B\cup(-B^{\prime})\right) \times \left(A\cup(-A^{\prime})\right)\right)$ give the conjugate expression. This gives the wanted expression for $\CS$.

The proof of Theorem~\ref{thmLP} is now reduced to the proof of Proposition~\ref{proplprecap} that will occupy the rest of this section.
\eop

\subsection{Normalizing $F_{\fvarM}$ with respect to one handlebody $A$}
\label{subproofbeg}

Let $A$ be a rational homology handlebody embedded in $M$ (outside a tubular neighborhood of $K$) whose $H_1$ goes to torsion in $H_1(M)$ so that $A$ lifts in $\tilde{M}$. Again, we shall consider a preferred lift of $A$ in $\tilde{M}$ and often identify $A$ with its preferred lift, and submanifolds of $A$ with submanifolds of the preferred lift, without mentioning it.

Let $(a_i,z_i)_{i=1, \dots g}$ be a basis of $H_1(\partial A) $ such that $a_i=\partial(\Sigma(a_i) \subset A)_{i=1,\dots, g}$ where $\Sigma(a_i)$ is a rational chain of $A$ and $\langle a_i,z_j \rangle=\delta_{ij}$.

Let $[-4,4] \times \partial A \subset M$ be a neighborhood of $\partial A= 0 \times\partial A$. ($[-4,0] \times \partial A \subset A$.)

For $s \in [-4,4]$, set
$$A_s=\left\{\begin{array}{ll} A \cup  (\partial A \times [0,s]) & \mbox{if}\;\; s \geq 0\\
A \setminus ( ]s,0]\times \partial A ) & \mbox{if} \;\;s \leq 0\end{array} \right.$$
and $(M \setminus A)_s = M \setminus \mathring{A}_s$.
$$ \partial A_s =\{s\}  \times\partial A=-\partial (M \setminus A)_s.$$

Assume that $\Sigma(a_i)$ intersects $ [-4,0]\times\partial A $ as
$[-4,0]\times a_i $ and construct a rational chain 
$\Sigma_4(a_i)=\Sigma(a_i)\cup [0,4] \times a_i $.
Set $\Sigma_s(a_i)=A_s \cap \Sigma_4(a_i)$.

Let $p_A$ be a point of $\partial A$ outside the $a_i$ and the $z_i$.
Let $(1-t_M)\gamma_A$ denote a path in $p_M^{-1}((M \setminus A)_0)$ that intersects $\left(p_M^{-1}(A_4)\subset \tilde{M}\right)$ as $(1-\theta_M)( [0,4] \times p_A)$ so that its boundary is $(\theta_M-1)p_A$. The path $(1-t_M)\gamma_A$ can be constructed as follows: First join $\partial A$ to a parallel of $K$ by a path, and 
let $(1-t_M)\gamma_A$ be the conjugate of this parallel of $K$ by this path. 

We shall assume that $\Sigma(\check{z}_i)$ (defined just before Proposition~\ref{proplprecap}) avoids the preimage of $\gamma_A$ and
(by removing the possible intersections with the help of $\partial A$)
that $\Sigma(\check{z}_i)$ intersects $p_M^{-1}([0,4] \times\partial A )$
as $ -[0,4]\times\check{z}_i $ and set 
$\Sigma_s(\check{z}_i)=\Sigma(\check{z}_i)\cap p_M^{-1}((M \setminus A)_s)$.

Set $\gamma_{A,s}=\gamma_A \cap p_M^{-1}((M \setminus A)_s)$.

Set $E=A \times (M\setminus A)_3$ and
$$\tilde{E}=p^{-1}(E)=A \times_{e} (M\setminus A)_3=A \times p_M^{-1}((M\setminus A)_3).$$

\begin{proposition}
\label{propnorone}
We can assume that $F_{\fvarM}$ intersects $\tilde{E}=A \times_{e} (M\setminus A)_3$  like 
$$A \times \gamma_{A,3} + \sum_{(i,j) \in \{1, \dots, g\}^2} \overline{lk_e(z_i,a_j^+)}\Sigma(a_i)\times \Sigma_3(\check{z_j})$$ 
where $\overline{lk_e(z_i,a_j^+)}=lk_e(z_i,a_j^+)(t^{-1})$.
\end{proposition}
\bp
We shall check that the class of $F_{\fvarM}$ in 
$$H_4(C_2(M),C_2(M) \setminus \mbox{Int}(E);\QQ(t))=H_4(E,\partial E;\QQ(t))$$
coincides with the class $F_1$ of the chain of the statement by proving first that $H_4(E,\partial E;\QQ(t))$ is the dual of $H_2( E;\QQ(t))$ -as Poincar\'e duality prescribes- and next that $F_{\fvarM}$ and $F_1$ intersect a basis of $H_2(E;\QQ(t))$ in the same way. This is enough to produce a $5$--chain $W$ in $\tilde{E}$, or rather in the homotopically equivalent
$A_1 \times_e (M\setminus A)_2$, such that $(\partial W +  F_{\fvarM} - F_1) \subset \partial \left(A_1 \times_e (M\setminus A)_2\right)$
and to conclude by adding $\partial W$ to $F_{\fvarM}$.
Details follow.
\eop

\begin{lemma}
 $$H_{\ast}(M \setminus A; \QQ(t_M))=H_{\ast+1}(A,\partial A;\QQ)\otimes_{\QQ} \QQ(t_M)$$
$$H_{2}(M \setminus A; \QQ(t_M))=\QQ(t_M)[\partial A].$$
$$H_{1}(M \setminus A; \QQ(t_M))=\oplus_{i=1}^g \QQ(t_M)[a_i].$$
$$H_{0}(M \setminus A; \QQ(t_M))=0.$$
\end{lemma}
\bp Since $H_{\ast}(M ; \QQ(t_M))=0$ by Lemma~\ref{lemhomtilM},
$$H_{\ast}(M \setminus A; \QQ(t_M))=H_{\ast+1}(M,M \setminus A; \QQ(t_M)).$$
\eop

\begin{lemma}
$$H_{q}(E; \QQ(t))=
\left(H_{\ast}(A;\QQ)\otimes  \left(\QQ[\partial A_3] \oplus\oplus_{i=1}^g \QQ[a_i\subset \partial A_3]\right)\right)_q \otimes_{\QQ} \QQ(t)$$
 $$H_{q}(\partial A \times (M \setminus A)_3; \QQ(t))=
\left(H_{\ast}(\partial A;\QQ)\otimes  \left(\QQ[\partial A_3] \oplus\oplus_{i=1}^g \QQ[a_i\subset \partial A_3]\right)\right)_q \otimes_{\QQ} \QQ(t)$$
\end{lemma}
\bp $\tilde{E}=p^{-1}(E)=A \times p_M^{-1}(M\setminus \mbox{Int}(A_3))$.
\eop

\begin{lemma}
 $$H_{4}(E,\partial E; \QQ(t))=\oplus_{(i,j)\in\{1,2,\dots,g\}^2}\QQ(t) [\Sigma(a_i)\times\Sigma_3(\check{z}_j)]\oplus \QQ(t)[A \times \gamma_{A,3}].$$
\end{lemma}
\bp
Since $H_4(E; \QQ(t))=0$, $H_{4}(E,\partial E; \QQ(t))$ is isomorphic to the kernel of the natural map
$$H_3(\partial E; \QQ(t)) \rightarrow H_3(E; \QQ(t)).$$
Let $A_{ij}=\partial \left(\Sigma(a_i)\times\Sigma_3(\check{z}_j)\right)$ and $A_{00}=\partial (A \times \gamma_{A,3})$.
It suffices to prove that the $A_{ij}$ and $A_{00}$ form a basis of this kernel.

Let us compute $H_3(\partial E=A\times \partial A_3 \cup_{\partial A\times \partial A_3} \partial A\times  (M\setminus A)_3; \QQ(t))$ with the help of the Mayer-Vietoris sequence.
The cokernel of the Mayer-Vietoris map 
$$H_3(\partial A\times \partial A_3; \QQ(t))\rightarrow H_3(A\times \partial A_3; \QQ(t)) \oplus H_3(\partial A \times  (M\setminus A)_3; \QQ(t))$$
is $$ \oplus_{i=1}^g \QQ(t)[z_i\times \partial A_3].$$
The kernel of the Mayer-Vietoris map 
$$H_2(\partial A\times \partial A_3; \QQ(t))\rightarrow H_2(A\times \partial A_3; \QQ(t)) \oplus H_2(\partial A \times  (M\setminus A)_3; \QQ(t))$$
is freely generated over $\QQ(t)$ by the $[a_i\times \check{z}_j]$, for  $(i,j)\in\{1,2,\dots,g\}^2$ and $[\partial A \times p_{A_3}]$ that are the images of the $\pm A_{ij}$ and $\pm A_{00}$
via the Mayer-Vietoris boundary map.
Therefore $$H_3(\partial E; \QQ(t))=\oplus_{i=1}^g \QQ(t)[z_i\times \partial A_3] \oplus \QQ(t)A_{00} \oplus\oplus_{(i,j)\in \{1,2,\dots,g\}^2} \QQ(t)A_{ij}.$$
On the other hand, we know $$H_3(E;\QQ(t))=\oplus_{i=1}^g \QQ(t)[z_i\times \partial A_3].$$
This easily leads to the result.
\eop

\noindent{\sc End of proof of Proposition~\ref{propnorone}:}
Enlarge $\tilde{E}$ to the homeomorphic $$\tilde{E}_1 = A_1 \times_e (M\setminus A)_2=p^{-1}(E_1)$$
where $E_1=A_1 \times (M\setminus A)_2$.
The basis of $H_4(E_1,\partial E_1;\QQ(t))$ is dual to the basis of $H_2(E_1;\QQ(t))$
made of the $z_i \times (\{3\} \times a_j )$ and $p_{A} \times \partial A_3$, with respect to the equivariant algebraic intersection. 
According to Proposition~\ref{propdeflkeq} and to Theorem~\ref{thmstauM},
$\langle z_i \times (\{3\} \times a_j ),F_{\fvarM}\rangle_e=lk_e(z_i,\{3\} \times a_j)$ and $\langle p_{A} \times \partial A_3,F_{\fvarM}\rangle_e=1$.
Therefore the initial chain $F_{\fvarM}$ and the wanted one $F_1$ (naturally extended to $\tilde{E}_1$) intersect
this basis in the same way algebraically, and their classes in $H_4(E_1,\partial E_1;\QQ(t))$ coincide. Then there exists a $5$--chain
$W$ of $\tilde{E}_1$ such that $(\partial W +  F_{\fvarM} - F_1) \subset \partial \tilde{E}_1$.
Replace $F_{\fvarM}$ by $\partial W +  F_{\fvarM}$.
\eop

Set $$(M\setminus A)_3\times_{e}A= p^{-1}\left((M\setminus A)_3\times A\right)=-\iota\left(A \times_{e} (M\setminus A)_3\right).$$

\begin{proposition}
\label{propnoronesym}
We can furthermore assume that 
$$\begin{array}{ll}F_{\fvarM} \cap\left((M\setminus A)_3\times_{e}A\right)&=\iota\left(F_{\fvarM} \cap\left(A \times_{e} (M\setminus A)_3\right)\right)\\&=- \gamma_{A,3} \times A+ \sum_{(i,j) \in \{1, \dots, g\}^2} lk_e(z_i,a_j^+)\Sigma_3(\check{z_j})\times \Sigma(a_i).\end{array}$$ 

\end{proposition}
\bp Since $\left((M\setminus A)_3\times_{e}A\right)$ and $\left(A \times_{e} (M\setminus A)_3\right)$ are disjoint, this statement does not interfer with the statement of Proposition~\ref{propnorone}. Its proof is obtained from the other one by symmetry. \eop

\begin{lemma}
\label{lemintFSig}
 We can furthermore assume that 
$\langle F_{\fvarM}, \Sigma_3(a_i)\times p_A \rangle_e=0$ and that
$\langle F_{\fvarM}, p_A\times \Sigma_3(a_i) \rangle_e=0$ for any $i$.
\end{lemma}
\bp
Indeed, adding to $F_{\fvarM}$ multiples of $\partial (\Sigma_2(\check{z}_i) \times A_1)$
allows us to fix
the first intersections arbitrarily and independently. Similarly, adding to $F_{\fvarM}$ multiples of $\partial (A_1 \times \Sigma_2(\check{z}_i))$ allows us to fix
the second intersections without changing the first ones, and these additions do not change $F_{\fvarM}$ on $\left((M\setminus A)_3\times_{e}A\right) \cup \left(A\times_{e}(M\setminus A)_3\right)$.
\eop

\subsection{Normalization with respect to a finite set of handlebodies}

Assume that a set of $n$ disjoint rational homology handlebodies $A^{(i)}$ in $M$ is given, outside a tubular neighborhood of $K$.

In this article, $n=2$ is enough and $A^{(2)}=B$, but it is not harder
to study the general case that will be useful in a later study of 
the higher loop degree case.
 For each $i$, let $p_{A^{(i)}} $ be a point of $\partial A^{(i)}$ outside the $(a_j^{(i)})$ and the $(z_j^{(i)})$ that play similar roles as before for each $A^{(i)}$. 
Construct disjoint paths $(1-t_M) \gamma_{A^{(i)}}$ in $p_M^{-1}(M \setminus (\coprod_{j=1}^{n}\mbox{Int}(A^{(j)})))$ as before. The path
$(1-t_M) \gamma_{A^{(i)}}$ will intersect $p_M^{-1}(A^{(i)}_{4})$ as $ (1-\theta_M)[0,4] \times p_{A^{(i)}}$ so that its boundary is $(\theta_M-1)p_{A^{(i)}}$.
All the $\gamma_{A^{(i)}}$ are supposed to be disjoint from each other and from a tubular neighborhood of $K$.

Then construct the $\Sigma(\check{z}_j^{(i)})$ as before Proposition~\ref{proplprecap} so that they are transverse and disjoint from the $(1-t_M) \gamma_{A^{(k)}}$. (Again, their possible intersections with the paths
and $K$ can be removed with the help of the $\partial A^{(k)}$ and $S$, respectively.)
Also assume that the $\Sigma(\check{z}_k^{(i)})$ for $A^{(i)}$ intersects $A_4^{(j)}$ as copies of $\Sigma_4(a_{\ell}^{(j)})$.

\begin{proposition}
\label{propnoronesim}
 The normalizations of the previous subsection can be achieved simultaneously. In other words, we can assume that
\begin{itemize}
\item For any $j=1,\dots, n$,
$F_{\fvarM}$ intersects $A^{(j)} \times_{e} (M\setminus A^{(j)})_3$  like 
$$A^{(j)} \times \gamma_{A^{(j)},3} + \sum_{(i,k) \in \{1, \dots, g(A^{(j)})\}^2} \overline{lk_e(z_i^{(j)},a^{(j)+}_k)}\Sigma(a_i^{(j)})\times \Sigma_3(\check{z}_k^{(j)}).$$
\item  For any $j=1,\dots, n$,
$F_{\fvarM}$ intersects $(M\setminus A^{(j)})_3 \times_{e} A^{(j)}$ like 
$$-\gamma_{A^{(j)},3} \times A^{(j)} + \sum_{(i,k) \in \{1, \dots, g(A^{(j)})\}^2} lk_e(z_i^{(j)},a^{(j)+}_k)\Sigma_3(\check{z}^{(j)}_k)\times \Sigma(a_i^{(j)}).$$
\item For any $j=1,\dots, n$, for any $i =1,\dots, g(A^{(j)})$,
$\langle F_{\fvarM}, \Sigma_3(a^{(j)}_i)\times p_{A^{(j)}} \rangle_e=0$ and
$$\langle F_{\fvarM}, p_{A^{(j)}}\times \Sigma_3(a^{(j)}_i) \rangle_e=0.$$
\end{itemize}
As usual, $A^{(j)} \times p_{A^{(j)},3}$, $\Sigma(a_i^{(j)})\times \check{z}_k^{(j)}$, $p_{A^{(j)},3} \times A^{(j)}$, $\check{z}_k^{(j)}\times \Sigma(a_i^{(j)})$ are in the preferred lift of $C_2(A^{(j)}_4)$, the one that contains $ST(A^{(j)}_4)$.

\end{proposition}
\bp
Let $F_{k-1}$ be a $4$--chain that satisfies the
 normalization conditions -with respect to $(A^{(j)}_1,(M\setminus A^{(j)})_2)$ instead of $(A^{(j)},(M\setminus A^{(j)})_3)$ for the first two ones without loss-  for $j<k$ and let us try to modify it into a $4$--chain $F_{k}$ that also satisfies them with respect to $A^{(k)}$. To do that, 
we would like to add some $\partial W$ to $F_{k-1}$ where $W$ is a $5$--chain of
$A^{(k)}_1 \times_e (M\setminus A^{(k)})_2$.
Unfortunately, we need to control $\partial W$ on $\partial \left(A^{(k)}_1 \times_e (M\setminus A^{(k)})_2\right)$ since this boundary contains $\partial A^{(k)}_1 \times_e A^{(j)}$ for $j<k$ where we cannot afford to lose our former modifications.

We know that $F_{k-1}$ is already as wanted on $A^{(k)}_1 \times_e \coprod_{j=1}^{k-1}A^{(j)}_1$. 
Set $$D=A^{(k)}_1 \times \left((M\setminus A^{(k)})_2 \setminus (\coprod_{j=1}^{k-1}\mbox{Int}(A^{(j)}_1))\right).$$
$$E_1=A^{(k)}_1 \times (M\setminus A^{(k)})_2.$$
We shall consider the difference of  $F_{k-1}$ with the wanted chain $F_{k,D}$
in $H_4(D,D\cap \partial E_1; \QQ(t))$
where $$D\cap \partial E_1=\partial A^{(k)}_1 \times\left( (M\setminus A^{(k)})_2 \setminus (\coprod_{j=1}^{k-1}\mbox{Int}(A^{(j)}_1))\right) \cup A^{(k)}_1 \times \partial A^{(k)}_2$$
$$=\partial E_1 \setminus \left( \partial A^{(k)}_1 \times(\coprod_{j=1}^{k-1}\mbox{Int}(A^{(j)}_1)) \right).$$

We shall prove the following lemma at the end of this subsection.
\begin{lemma}
\label{lemhomchiant}
 $H_4(D,D\cap \partial E_1; \QQ(t))$ is generated by chains of the following forms for $j<k$, \\
$\Sigma_1(a_i^{(k)}) \times \partial A_1^{(j)}$, $A^{(k)}_1 \times \left(\{1\} \times (a_i^{(j)})\right)$, $A^{(k)}_1 \times \gamma_{A^{(k)},2}$,
and combinations of $\Sigma_1(a_i^{(k)}) \times \Sigma_2(\check{z}_{\ell}^{(k)})$ that are in $D$.
\end{lemma}

Assume this lemma for a moment and let us conclude the proof.

The chains of the statement will be detected by the following dual chains $\check{z}_i^{(k)} \times \gamma_{A^{(j)}}$, $p_{A^{(k)}} \times \Sigma(\check{z}_{i}^{(j)})$, $p_{A^{(k)}} \times \partial A^{(k)}_3$ and $z_i^{(k)} \times a_{\ell,3}^{(k)}$.
The chains $F_{k,D}$ and $F_{k-1}$ intersect the cycles of the last two kinds as prescribed by the linking form condition that is consistent with the
prescription of the statement.
Let us show that their algebraic intersections with the chains of the  first two kinds also coincide.
Since we know $F_{k,D}$  on $A^{(k)}_1  \times_e (M\setminus A^{(k)})_2$, we 
see that $\gamma_{A^{(j)}}$ is not in the interaction locus of $A^{(k)}$, and $\langle F_{k,D},\check{z}_i^{(k)} \times \gamma_{A^{(j)}} \rangle_e=0$. Similarly, $p_{A^{(k)}}$ interacts only with $\gamma_{A^{(k)}}$ that does not meet $\Sigma(\check{z}_{i}^{(j)})$, and $\langle F_{k,D},p_{A^{(k)}} \times \Sigma(\check{z}_{i}^{(j)}) \rangle_e=0$.
On the other hand,
$$\langle F_{k-1},\partial\left(\Sigma(\check{z}_i^{(k)}) \times \gamma_{A^{(j)}} \right) =\check{z}_i^{(k)} \times \gamma_{A^{(j)}} - \Sigma(\check{z}_i^{(k)}) \times p_{A^{(j)}}\rangle_e=0.$$
Therefore, since 
$$\langle F_{k-1},\Sigma(\check{z}_i^{(k)})\times p_{A^{(j)}}\rangle_e=0$$
(because $\langle F_{k-1},\Sigma_3(a_{\ell}^{(j)}) \times p_{A^{(j)}}\rangle_e=0$ for all $\ell$),
$\langle F_{k-1},\check{z}_i^{(k)} \times \gamma_{A^{(j)}} \rangle_e=0$.

Similarly, 
$$\langle F_{k-1},\partial\left(\gamma_{A^{(k)}} \times \Sigma(\check{z}_{i}^{(j)}) \right) =-p_{A^{(k)}} \times \Sigma(\check{z}_{i}^{(j)}) - \gamma_{A^{(k)}} \times \check{z}_{i}^{(j)}\rangle_e=0,$$
and the prescriptions for $F_{k-1}$ impose that
$\langle F_{k-1},\gamma_{A^{(k)}} \times \check{z}_{i}^{(j)}\rangle_e=0$.
Therefore, $\langle F_{k-1},p_{A^{(k)}} \times \Sigma(\check{z}_{i}^{(j)})\rangle_e=0$.

In particular the class of $(F_{k,D}-F_{k-1})$ vanishes in $H_4(D,D\cap \partial E_1; \QQ(t))$. Therefore, there exists a $5$--chain $W$ in 
$p^{-1}(D)$ such that $\partial W +F_{k-1}-F_{k,D} \subset p^{-1}(D\cap \partial E_1)$, and
$F_k=\partial W + F_{k-1}$ satisfies the normalizations conditions on $p^{-1}(E_1)$. It can similarly be assumed to satisfy them on the symmetric part $\iota(p^{-1}(E_1))$. Now, we independently deal with the additional conditions $\langle F_k, \Sigma_3(a^{(k)}_i)\times p_{A^{(k)}} \rangle_e=0$ and 
$\langle F_k,p_{A^{(k)}} \times \Sigma_3(a^{(k)}_i) \rangle_e=0$ for any $i$, as follows.

We shall add a combination $C$ of $\partial (\Sigma_2(\check{z}_i^{(k)}) \times A^{(k)}_1)$ to $F_k$ in order to satisfy the first condition, for example. But this combination must not intersect some $A^{(j)} \times_e (M \setminus A^{(j)})_3$, for $j<k$. Under our assumptions, this is equivalent to say that the algebraic equivariant intersection
$\langle \check{z}^{(j)}_r \times \gamma_{A^{(k)}},C \rangle_e$ must vanish.

Therefore, in order to prove that we may furthermore assume that $\langle F_k, \Sigma_3(a^{(k)}_i)\times p_{A^{(k)}} \rangle_e=0$, it suffices to prove the following lemma.

\begin{lemma}
 The combination $C$ of $\partial (\Sigma_2(\check{z}_i^{(k)}) \times A^{(k)}_1)$ such that
$\langle F_k +C , \Sigma_3(a^{(k)}_i)\times p_{A^{(k)}} \rangle_e=0$ satisfies $$\langle \check{z}^{(j)}_r \times \gamma_{A^{(k)}},C \rangle_e=0$$ 
for any $j<k$ and for any $r=1, \dots, g^{(j)}$.
\end{lemma}
\bp
Define a linear form $f$ on $H_1(M\setminus \mathring{A}^{(k)};\QQ(t_M))$
by its values on the generators $a_i^{(k)}$
$$f(a_i^{(k)})=\langle  \Sigma_3(a^{(k)}_i)\times p_{A^{(k)}},F_k \rangle_e.$$
Then $$f(a_i^{(k)})=\langle  a^{(k)}_{i,3} \times \gamma_{A^{(k)}},F_k \rangle_e.$$
Now, $$\langle \check{z}^{(j)}_r \times \gamma_{A^{(k)}},F_k \rangle_e
=\langle \left(\Sigma(\check{z}^{(j)}_r) \cap (M\setminus A^{(k)})_3\right) \times p_{A^{(k)}},F_k \rangle_e + f(z^{(j)}_r)$$
where $\langle \check{z}^{(j)}_r \times \gamma_{A^{(k)}},F_k \rangle_e=0$ because of the form of $F_k$ on $A^{(j)} \times_e (M\setminus A^{(j)})_3 $ and $$\langle\left( \Sigma(\check{z}^{(j)}_r )\cap (M\setminus A^{(k)})_3\right) \times p_{A^{(k)}},F_k \rangle_e=0$$ because of the form of $F_k$ on $(M\setminus A^{(k)})_3 \times_e A^{(k)} $.
Therefore, $f(z^{(j)}_r)=0$.

Our combination $C$ is a boundary that is defined so that
$$f(a_i^{(k)})=-\langle  \Sigma_3(a^{(k)}_i)\times p_{A^{(k)}},C \rangle_e=-\langle  a^{(k)}_{i,3} \times \gamma_{A^{(k)}},C \rangle_e.$$
Thus,
$$\langle \check{z}^{(j)}_r \times \gamma_{A^{(k)}},C \rangle_e=\langle \left(\Sigma(\check{z}^{(j)}_r) \cap (M\setminus A^{(k)})_3\right) \times p_{A^{(k)}},C \rangle_e-f(z^{(j)}_r).$$
Since
$\partial (\Sigma_2(\check{z}_i^{(k)}) \times A^{(k)}_1)$ does not meet
$\left(\Sigma(\check{z}^{(j)}_r) \cap (M\setminus A^{(k)})_3 \right)\times p_{A^{(k)}}$, $\langle \check{z}^{(j)}_r \times \gamma_{A^{(k)}},C \rangle_e=0$, and we are done.

\eop

Similarly, we may furthermore assume that $\langle F_k,p_{A^{(k)}} \times \Sigma_3(a^{(k)}_i) \rangle_e=0$ for any $i$.
This concludes the proof of Proposition~\ref{propnoronesim} up to the proof of Lemma~\ref{lemhomchiant} that we give now.
\eop

\noindent{\sc Proof of Lemma~\ref{lemhomchiant}:}
$D=A^{(k)}_1 \times \left((M\setminus A^{(k)})_2 \setminus (\coprod_{j=1}^{k-1}\mbox{Int}(A^{(j)}_1))\right)$.
$$E_1=A^{(k)}_1 \times (M\setminus A^{(k)})_2.$$

Since $H_4(D; \QQ(t))=0$, $H_4(D,D\cap \partial E_1; \QQ(t))$
is the kernel of the natural map
$H_3(D\cap \partial E_1; \QQ(t))\rightarrow H_3(D; \QQ(t))$
where
$$H_3(D; \QQ(t))=H_1(A^{(k)}_1;\QQ) \otimes_{\QQ}\left( \oplus_{j=1}^{k-1}\QQ(t)[\partial(A^{(j)}_1)] \oplus \QQ(t)[\partial(A^{(k)}_2)]\right).$$
We compute $H_3(D\cap \partial E_1; \QQ(t))$ using the Mayer-Vietoris
sequence associated with the decomposition
$$D\cap \partial E_1=D_1 \cup_{\partial A^{(k)}_1 \times \partial A^{(k)}_2} D_2$$
where
$$D_1=\partial A^{(k)}_1 \times\left( (M\setminus A^{(k)})_2 \setminus (\coprod_{j=1}^{k-1}\mbox{Int}(A^{(j)}_1))\right) \;\; \mbox{and}\;\;
D_2=A^{(k)}_1 \times \partial A^{(k)}_2.$$
The cokernel of the map from $H_3(D_1\cap D_2=\partial A^{(k)}_1 \times \partial A^{(k)}_2; \QQ(t))$ to $H_3(D_1; \QQ(t)) \oplus H_3(D_2; \QQ(t))$  is
$$\left(H_1(A^{(k)}_1;\QQ) \otimes \QQ(t)[\partial A^{(k)}_2]\right)
\oplus
\left(H_1(\partial A^{(k)}_1;\QQ) \otimes \oplus_{j=1}^{k-1}\QQ(t)[\partial A^{(j)}_1]\right)
\oplus \bigoplus_{j=1}^{k-1} \bigoplus_{i=1}^{g(A^{(j)})} \QQ(t)[\partial A^{(k)}_1 \times a_i^{(j)}]
.$$
The kernel of the map from $H_2(D_1\cap D_2
; \QQ(t))$ to $H_2(D_1; \QQ(t)) \oplus H_2(D_2; \QQ(t))$  is 
$$(\oplus_i \QQ a_i^{(k)}) \otimes \mbox{Ker}\left(H_1(\partial A^{(k)}_2; \QQ(t_M))
\rightarrow
H_1\left((M\setminus A^{(k)})_2 \setminus (\coprod_{j=1}^{k-1}\mbox{Int}(A^{(j)}_1));\QQ(t_M)\right)\right)$$ $$
\oplus
\QQ(t) [\partial A^{(k)}_1 \times p_{A^{(k)}_2}].
$$
Therefore $H_3(D\cap \partial E_1; \QQ(t))$ is generated by 
classes of the following form that vanish in $H_3(D;\QQ(t))$
\begin{itemize}
\item $\partial \left(A^{(k)}_1 \times \gamma_{A^{(k)},2}\right)$,
\item combinations  $\sum_{(i,j)}\alpha_{ijk}\partial (\Sigma_1(a_i^{(k)}) \times \Sigma_2(\check{z}_j^{(k)}))$ such that $\sum_{(i,j)}\alpha_{ijk}\Sigma_1(a_i^{(k)}) \times \Sigma_2(\check{z}_j^{(k)}) \subset D$
\item $(\{1\} \times a_i^{(k)}) \times \partial A^{(j)}_1$, and 
\item $\partial A^{(k)}_1 \times (\{1\} \times a_i^{(j)})$
\end{itemize}
and classes of $H_1(A^{(k)}_1;\QQ) \otimes_{\QQ} \left( \oplus_{j=1}^{k-1}\QQ(t)[\partial A^{(j)}] \oplus \QQ(t)[\partial A^{(k)}_2]\right)$ that survive in $H_3(D;\QQ(t))$. Only the first four ones contribute to $H_4(D,D\cap \partial E_1; \QQ(t))$ as in the statement.
\eop

\subsection{Almost concluding the proof of Proposition~\ref{proplprecap}}
\label{subconclLP}

With the notation of the previous subsection, for any $i=1,\dots, n$, let $A^{(i)\prime}$ be another rational homology handlebody with the same 
boundary and the same Lagrangian as $A^{(i)}$, let $\Sigma_4^{\prime}(a^{(i)}_j)$ be a rational $4$-chain with the same boundary $4 \times a^{(i)}_j$ as $\Sigma_4(a^{(i)}_j)$ and that intersects $ [-4,4] \times \partial A^{(i)}$ like $\Sigma_4(a^{(i)}_j)$. 

Without loss, in addition to the requirements of Proposition~\ref{propnoronesim}, assume that
\begin{itemize}
\item For any $j=1,\dots, n$, for any $t\in[-4,0]$,
$F_{\fvarM}$ intersects $A^{(j)}_t \times_{e} (M\setminus A^{(j)})_{t+3}$  like\\ 
$A^{(j)}_t \times \gamma_{A^{(j)},t+3} + \sum_{(i,k) \in \{1, \dots, g(A^{(j)})\}^2} \overline{lk_e(z_i^{(j)},a^{(j)+}_k)}\Sigma_t(a_i^{(j)})\times \Sigma_{t+3}(\check{z}_k^{(j)}).$
\item  For any $j=1,\dots, n$, for any $t\in[-4,0]$,
$F_{\fvarM}$ intersects $(M\setminus A^{(j)})_{t+3} \times_{e} A^{(j)}_t$ like\\ 
$-\gamma_{A^{(j)},t+3} \times A^{(j)}_t + \sum_{(i,k) \in \{1, \dots, g(A^{(j)})\}^2} lk_e(z_i^{(j)},a^{(j)+}_k)\Sigma_{t+3}(\check{z}^{(j)}_k)\times \Sigma_t(a_i^{(j)}).$
\item $F_{\fvarM}=\partial F_{\fvarM} \times [0,1]$ in an equivariant neighborhood $\partial \TCM \times [0,1]$ of $\partial \TCM$ in $\TCM$.
\end{itemize}

Let $M_i=M(A^{(i)\prime}/A^{(i)})$.
Define $\tilde{F}(M_i)$ on $(\tilde{C}_2(M_i) \setminus \mbox{Int}(\tilde{C}_2(A^{(i)\prime}_{-1})))$ so that
\begin{itemize}
 \item $\tilde{F}(M_i)= F_{\fvarM}$ on 
$\tilde{C}_2((M \setminus A^{(i)\prime})_{-4})$
\item
For any $t\in[-4,0]$, $\tilde{F}(M_i)$ intersects $A^{(i)\prime}_t \times_{e} (M\setminus A^{(i)})_{t+3}$  like 
$$\sum_{(j,k) \in \{1, \dots, g(A^{(i)})\}^2} \overline{lk_e(z_j^{(i)},a^{(i)+}_k)}\Sigma^{\prime}_t(a^{(i)}_j)\times \Sigma_{t+3}(\check{z}^{(i)}_k) + A^{(i)\prime}_t \times \gamma_{A^{(i)},t+3}.$$
\item  
For any $t\in[-4,0]$, $\tilde{F}(M_i)$ intersects $(M\setminus A^{(i)})_{t+3} \times_{e} A^{(i)\prime}_t$ like 
$$\sum_{(j,k) \in \{1, \dots, g(A^{(i)})\}^2} lk_e(z_j^{(i)},a^{(i)+}_k)\Sigma_{t+3}(\check{z}^{(i)}_k)\times \Sigma_t^{\prime}(a_j^{(i)}) - \gamma_{A^{(i)},t+3} \times A^{(i)\prime}_t.$$
\item $\tilde{F}(M_i)=\partial F_{\fvarM}(M_i)$ on $ST(M_i)$.
\end{itemize}

In particular, $\tilde{F}(M_i)$ is well-defined outside the interior of
$\tilde{C}_2(A^{(i)\prime}_{-1})$ and we wish to extend it there.
Assume without loss that $\tilde{F}(M_i)$ is transverse to $\partial \tilde{C}_2(A^{(i)\prime})$, and
consider the class of the $3$--cycle $$\tilde{F}_3= \partial \tilde{C}_2(A^{(i)\prime}) \cap_e \tilde{F}(M_i).$$

The following lemma will be proved in Subsection~\ref{subsecprephalt}.
\begin{lemma}
\label{lemprephalt}
Under the hypotheses above $$\tilde{F}_3-\sum_{(j,k,\ell)\in \{1,\dots g(A^{(i)})\}^3} \CI_{A^{(i)}A^{(i)\prime}}(a_j^{(i)}\wedge a_k^{(i)} \wedge a_{\ell}^{(i)})lk_e(z^{(i)}_{\ell},z^{(i)}_k) ST(p_M(\check{z}_{j}^{(i)}))$$ 
bounds a $4$--chain $\tilde{F}_{\fvarM}(A^{(i)\prime})$ in $\tilde{C}_2(A^{(i)\prime})$.
\end{lemma}
Let $\eta <\frac{1}{3n}$. Define $ST([-4+3i\eta,0]\times \partial \Sigma_{A^{(i)}A^{(i)\prime}})$ in $\tilde{C}_2([-4,0]\times \partial A^{(i)\prime})$ as 
$$
 \sum_{(j,k,\ell)\in \{1,\dots g(A^{(i)})\}^3} \CI_{A^{(i)}A^{(i)\prime}}(a^{(i)}_j\wedge a^{(i)}_k \wedge a^{(i)}_{\ell})lk_e(z^{(i)}_{\ell},z^{(i)}_k) ST(M)_{|[-4+3i\eta,0]\times p_M(\check{z}_{j}^{(i)})},$$
and set $$ST(\Sigma_{A^{(i)}A^{(i)\prime},-4+3i\eta})=ST(\Sigma_{A^{(i)}A^{(i)\prime}}) \cup -ST([-4+3i\eta,0]\times \partial \Sigma_{A^{(i)}A^{(i)\prime}}).$$

Assuming Lemma~\ref{lemprephalt}, we shall rather assume without loss that
$$\partial \tilde{F}_{\fvarM}(A^{(i)\prime})=\tilde{F}_3 - \partial ST_1( \Sigma_{A^{(i)}A^{(i)\prime},-4+3i\eta}),$$
with notation consistent with the statement of Proposition~\ref{proplprecap},
and that $\tilde{F}_{\fvarM}(A^{(i)\prime})$ coincides with $\tilde{F}(M_i)$ on $\tilde{C}_2(A^{(i)\prime}) \setminus \tilde{C}_2(A^{(i)\prime}_{-1})$.
Set
$$F_{\fvarM}(A^{(i)\prime})=\tilde{F}_{\fvarM}(A^{(i)\prime})
+\left(ST_{3i\eta}( \Sigma_{A^{(i)}A^{(i)\prime},-4+3i\eta}) \cap \tilde{C}_2(A^{(i)\prime})\right)
-\left(ST(\partial \Sigma_{A^{(i)}A^{(i)\prime},-4+3i\eta})\times[3i\eta,1]\right)$$
(where $ST(\partial \Sigma_{A^{(i)}A^{(i)\prime},-4+3i\eta})\times[3i\eta,1] \subset \partial \TCM \times [0,1]$).
Then 
$$F_{\fvarM}(M_i)=\tilde{F}(M_i) \cap \overline{\tilde{C}_2(M)\setminus \tilde{C}_2(A^{i})}
+F_{\fvarM}(A^{(i)\prime}) + \left(ST_{3i\eta}( \Sigma_{A^{(i)}A^{(i)\prime},-4+3i\eta}) \cap \overline{\tilde{C}_2(M) \setminus \tilde{C}_2(A^{(i)\prime})}\right)$$
is a $4$-chain of $\tilde{C}_2(M_i)$ such that $\partial F_{\fvarM}(M_i)$
is the wanted fixed boundary.
Without loss, we shall furthermore assume that
$$\begin{array}{ll}\partial F_{\fvarM}(M_i) \cap \left(\partial \tilde{C}_2(M_i)\times [0,1]\right)
=&\partial F_{\fvarM}(M_i) \times [0,1]\\& +ST_{3i\eta}( \Sigma_{A^{(i)}A^{(i)\prime},-4+3i\eta})
-ST(\partial \Sigma_{A^{(i)}A^{(i)\prime},-4+3i\eta})\times[3i\eta,1].\end{array}$$
Now, for any subset $I$ of $\{1,2,\dots,n\}$, set $M_{I}=M((A^{(j)\prime}/A^{(j)})_{j\in I})$ and consider the chain $ST(\Sigma_{A^{(i)}A^{(i)\prime}}^{I})$ in $\tilde{C}_2(M_{I})$ obtained from $ST(\Sigma_{A^{(i)}A^{(i)\prime}})$ by replacing the pieces $\Sigma(a_k^{(j)})$ by $\Sigma^{\prime}(a_k^{(j)})$, for $j \in I$.

Define $\tilde{F}_{\fvarM}(M_I)$ on $\tilde{C}_2(M_I) \setminus \coprod_{i \in I}\mathring{\tilde{C}}_2(A^{(i)\prime})$ as 
\begin{itemize}
 \item $F_{\fvarM}$ on $\tilde{C}_2\left(M \setminus \coprod_{i \in I}\mathring{A}^{(i)}_{-4}\right)$,
\item $\sum_{(j,k) \in \{1, \dots, g(A^{(i)})\}^2} \overline{lk_e(z_j^{(i)},a^{(i)+}_k)}\Sigma_t^{\prime}(a^{(i)}_j)\times \Sigma^{I}_{t+4}(\check{z}^{(i)}_k) + A^{(i)\prime}_t \times \gamma_{A^{(i)},t+4}$ on $A^{(i)\prime}_{t} \times_e (M_I \setminus A^{(i)\prime})_{t+4}$, for $i\in I$ and $t\in [-4,0]$,
\item $\sum_{(j,k) \in \{1, \dots, g(A^{(i)})\}^2}
 lk_e(z_j^{(i)},a^{(i)+}_k)\Sigma^{I}_{t+4}(\check{z}^{(i)}_k)\times \Sigma_t^{\prime}(a_j^{(i)}) - \gamma_{A^{(i)},t+4} \times A_t^{(i)\prime}$ on $(M_I \setminus A^{(i)\prime})_{t+4} \times_e A^{(i)\prime}_t$, for $i\in I$ and $t\in [-4,0]$.
\end{itemize}

$$\begin{array}{ll}F_{\fvarM}(M_I)=&\tilde{F}_{\fvarM}(M_I) \cap \left(\tilde{C}_2(M_I) \setminus \coprod_{i \in I}\mathring{\tilde{C}}_2(A^{(i)\prime}) \right)\\&
+\sum_{i\in I}ST_{3i\eta}(\Sigma_{A^{(i)}A^{(i)\prime},-4+3i\eta}^I)
\cap \left(\tilde{C}_2(M_I) \setminus \mathring{\tilde{C}}_2(A^{(i)\prime}) \right)
+\sum_{i\in I}F_{\fvarM}(A^{(i)\prime}).\end{array}$$

We construct $F_{\svarM}$ and $F_{\tvarM}$, similarly, taking parallel
copies of the pieces $\Sigma$ in the description of $\tilde{F}(M_I)$, and changing $3i\eta$ to $(3i-1)\eta$ and $(3i-2)\eta$, respectively.
$F_{\cvarM}(M_I)$ reads $\partial F_{\cvarM}(M_I)\times [0,1]$ in
$\partial \tilde{C}_2(M_I) \times [0,1]$, with the following two exceptions
\begin{itemize}
 \item the pairwise disjoint $ST_{(3i-r(\cvarM))\eta}(\Sigma_{A^{(i)}A^{(i)\prime}}^I)$ parts
\item the additional parts
$ST_{(3i-r(\cvarM))\eta}( \Sigma_{A^{(i)}A^{(i)\prime},-4+(3i-r(\cvarM))\eta}) \cap \tilde{C}_2(A^{(i)\prime})$ and
$$-ST(\partial \Sigma_{A^{(i)}A^{(i)\prime},-4+(3i-r(\cvarM))\eta})\times[(3i-r(\cvarM))\eta,1]$$ in the $p^{-1}(ST(A^{(i)\prime})\times [0,1])$.
\end{itemize}

Then since the $ST(\Sigma_{A^{(i)}A^{(i)\prime}})$ avoid $\tilde{C}_2(N(K))$, the $ST_{.}(\Sigma_{A^{(i)}A^{(i)\prime}})$ avoid the loci of double intersections of two $\left(\tilde{F}_{\cvarM}(M_I) \cap \left(\tilde{C}_2(M_I) \setminus \coprod_{i \in I}\mathring{\tilde{C}}_2(A^{(i)\prime}) \right)\right)$ parts, and they can only meet double intersections of $F_{\cvarM}(A^{(j)\prime})$, with $j\neq i$.

Now, forget about the superscript $(i)$ and let us come back to our case, where $n=2$, and to the notation $A=A^{(1)}$, $B=A^{(2)}$.
The $ST_{(3-r(\cvarM))\eta}(\Sigma_{AA^{\prime}})$-pieces do not meet the loci of double intersections of two $F_{\cvarM}(B^{\prime})$ because
the $F_{\cvarM}(B^{\prime})$ read as disjoint products on $ST(B^{\prime}) \times [0,3.5\eta]$, if we are dealing with actual trivialisations.
For pseudo-trivialisations, they read as non disjoint products on $ST(B^{\prime})\times [0,3.5\eta]$ whose double intersections will yield algebraically cancelling intersections with the $ST_{(3-r(\cvarM))\eta}(\Sigma_{AA^{\prime}})$, thanks to Lemma~\ref{lemintpseudo}.

Then we are left with the computation of the algebraic intersection of $ST_{(6-r(\sigma(\fvarM)))\eta}(\Sigma^{\prime}_{BB^{\prime}})$,
$F_{\sigma(\svarM)}(M_A)$ and $ F_{\sigma(\tvarM)}(M_A)$, in $p^{-1}(ST(A^{\prime})\times [0,1])$.
We shall assume that $\sigma$ is the trivial permutation and it will be clear
that the other permutations yield a similar result.
The possible intersections caused by pseudo-trivialisations algebraically cancel as above and we shall forget about them.
Then the only remaining intersection
will involve $ST_{(6-r(\fvarM))\eta}(\Sigma^{\prime}_{BB^{\prime}})$,
$-ST(\partial \Sigma_{AA^{\prime},-4+(3-r(\svarM))\eta})\times[(3-r(\svarM))\eta,1]$
and 
$s_{\tau_{A^{\prime}}}(A^{\prime};\tvarM)\times[0,1]$,
up to a permutation of $\svarM$ and $\tvarM$.
Since $ST_{6\eta}(\Sigma^{\prime}_{BB^{\prime}})$ intersects $p^{-1}(ST(A^{\prime})\times [0,1])$
as
$$\sum_{(\ell,m,n)\in \{1,\dots g_B\}^3, i\in \{1,\dots g_A\} } \CI_{BB^{\prime}}(b_{\ell}\wedge b_m \wedge b_n)lk_e(y_{n},y_m)lk({y}_{\ell},z_i) ST_{6\eta}(\Sigma^{\prime}(a_i)) $$
and since
$$\langle ST_{6\eta}(\Sigma^{\prime}(a_i)),-ST(z_{j,-4+2\eta})\times[2\eta,1],s_{\tau_{A^{\prime}}}(A^{\prime};\tvarM)\times[0,1]\rangle_e=-\delta_{ij},$$
we are done.

This ends the proof of Proposition~\ref{proplprecap}, where $ST_{(3i-r(\cvarM))\eta}( \Sigma_{A^{(i)}A^{(i)\prime}})$ should rather be defined as $ST_{(3i-r(\cvarM))\eta}( \Sigma_{A^{(i)}A^{(i)\prime},-4+(3i-r(\cvarM))\eta}) \cap \left(\tilde{C}_2(M_I) \setminus \mathring{\tilde{C}}_2(A^{(i)\prime}) \right)$, up to the proof of Lemma~\ref{lemprephalt}.
\eop

\subsection{Proof of Lemma~\ref{lemprephalt}}
\label{subsecprephalt}

To prove this lemma, we shall first refine and restate in other words some results of \cite[Section 5.3]{sumgen}.

Again, consider a rational homology handlebody $A$ equipped with a collar
$[-4,0] \times \partial A$.
For $s \in [-4,0]$, 
$A_s= A \setminus  ( ]s,0] \times \partial A )$, $ \partial A_s =\{s\} \times \partial A$.
Let $(a_i,z_i)_{i=1, \dots g_A}$ be a basis of $H_1(\partial A) $ such that $a_i=\partial(\Sigma(a_i) \subset A)_{i=1,\dots, g}$ where $\Sigma(a_i)$ is a rational chain of $A$ and $\langle a_i,z_j \rangle=\delta_{ij}$.

Consider a curve $a$ representing an element of $\CL_A$ of order $k$ in $H_1(A;\ZZ)$, $k \in \NN \setminus \{0\}$.
Let $\Sigma=k\Sigma(a)$ be a surface of $A$ immersed in $A$ bounded by $ka$ that
intersects $ [-1,0] \times \partial A $ as $k$ copies of $[-1,0] \times a$, and that intersects $\mbox{Int}(A_{-1})$ as an embedded surface.

\begin{lemma}
\label{homctimec}
Let $(c_i)_{i=1, \dots 2g}$ and $(c^{\ast}_i)_{i=1, \dots, 2g}$ be two dual bases of $H_1(\Sigma;\ZZ)/H_1( \partial \Sigma;\ZZ)$, $\langle c_i, c^{\ast}_j\rangle=\delta_{ij}$.

Then $\sum_{i=1}^{2g} c_i \times c^{\ast}_i$ is homologous to $\sum_{(j,{\ell}) \in \{1, \dots g_A\}^2}\langle \Sigma,\Sigma(a_j),\Sigma(a_{\ell}) \rangle z_j \times z_{\ell}$ in $A^2$.

Furthermore, $\sum_{i=1}^{2g} c_i \times c^{\ast +}_i$ is homologous to $\sum_{(j,{\ell}) \in \{1, \dots g_A\}^2}\langle \Sigma,\Sigma(a_j),\Sigma(a_{\ell}) \rangle z_j \times z_{\ell} -gST(\ast)$ in $C_2(A)$, where the $z_j$ are pairwise disjoint representatives of the $[z_j]$ on $\partial A$.
\end{lemma}
\bp For $(j,{\ell}) \in \{1, \dots g_A\}^2$, set $S_{j\ell}=\Sigma(a_j)\cap\Sigma(a_{\ell})$, $S_{\Sigma j}=\Sigma\cap\Sigma(a_j)$ and $S_{\Sigma \ell}=\Sigma\cap\Sigma(a_{\ell})$.
Then in $H_1(A)$,
$c_i=\sum_{j=1}^{g_A}\langle c_i, \Sigma(a_j) \rangle_A z_j=\sum_{j=1}^{g_A}\langle c_i , S_{\Sigma j} \rangle_{\Sigma} z_j$ and 
similarly,
$$c^{\ast}_i=\sum_{\ell=1}^{g_A}\langle c^{\ast}_i, S_{\Sigma \ell} \rangle_{\Sigma} z_{\ell}.$$
Thus, in $H_2(A^2)$,
$$c_i \times c^{\ast}_i
=\sum_{(j,\ell) \in \{1, \dots g_A\}^2}\langle c_i, S_{\Sigma j} \rangle_{\Sigma}\langle c^{\ast}_i, S_{\Sigma \ell} \rangle_{\Sigma} z_j \times z_{\ell} $$
On the other hand in $H_1(\Sigma)/H_1(\partial \Sigma)$,
$$S_{\Sigma j}=\sum_{i=1}^{2g}\langle c_i, S_{\Sigma j} \rangle_{\Sigma}c^{\ast}_i$$
and, 
$$S_{\Sigma \ell}=-\sum_{i=1}^{2g}\langle c^{\ast}_i, S_{\Sigma \ell}\rangle_{\Sigma}c_i.$$
Then $$\langle \Sigma,\Sigma(a_j),\Sigma(a_{\ell}) \rangle=\langle S_{\Sigma j},S_{\Sigma \ell} \rangle_{\Sigma}=\sum_{i=1}^{2g}\langle c_i, S_{\Sigma j} \rangle_{\Sigma}\langle c^{\ast}_i,  S_{\Sigma \ell} \rangle_{\Sigma}$$
and the first assertion is proved.
Let us now prove that
$\alpha=\sum_{i=1}^{2g} c_i \times c^{\ast+}_i$ is homologous to 
$$\beta=\sum_{(j,{\ell}) \in \{1, \dots g_A\}^2}\langle \Sigma(a),\Sigma(a_j),\Sigma(a_{\ell}) \rangle z_j \times z_{\ell}-g ST(\ast)$$ in $C_2(A)$. 
First note that the homology class of $\alpha$ in $C_2(A)$ is independent of the dual bases $(c_i)$ and $(c^{\ast}_i)$. Indeed, since both $a \times \sigma^+$ and $\sigma \times a^+$ are null-homologous in $C_2(A)$, the class of $\alpha$ in $C_2(A)$ only depends on the class of $\sum_{i=1}^{2g} c_i \times c^{\ast+}_i$ in $H_1(\Sigma)/H_1(\partial \Sigma) \otimes H_1(\Sigma^+)/H_1(\partial \Sigma^+)$ that is determined by the following property: For any two closed curves $e$ and $f$ of $\Sigma$, $\langle e \times f^+,\sum_{i=1}^{2g} c_i \times c^{\ast+}_i \rangle_{\Sigma \times \Sigma^+}=-\langle e,f \rangle_{\Sigma}$.

In particular, $[\alpha]=[\sum_{i=1}^{2g}  c^{\ast}_i \times (-c^{+}_i)]$.
The previous computation tells us that the difference $[\beta-\alpha]$ of the two classes
is a rational multiple of $[ST(\ast)]$.

Then
$$\begin{array}{lll}[\beta-\alpha]&=&+\frac{1}{2} \sum_{(j,{\ell}) \in \{1, \dots g_A\}^2}\langle \Sigma(a),\Sigma(a_j),\Sigma(a_{\ell})\rangle(lk_e(z_j,z_{\ell})+lk_e(z_{\ell},z_j))[ST(\ast)]-g[ST(\ast)]
\\&&-\frac{1}{4}\sum_{i=1}^{2g}\left(lk_e(c_i,c^{\ast+}_i)-lk_e(c^{\ast-}_i \times c_i))
+\overline{lk_e(c_i,c^{\ast+}_i)-lk_e(c^{\ast-}_i \times c_i)}\right)[ST(\ast)]
\\&=&-g[ST(\ast)]-\frac{1}{4} (-2g-2g)[ST(\ast)]=0.\end{array}$$
\eop

We now define a cycle $F^2(a)$ of $\partial C_2(A)$ that is associated to $\Sigma$.
Let $(a \times [-1,1])$ be a tubular neighborhood of $a$ in $\partial A$.
Let $p(a) \in a$ and see $a$ as the image of a map $a\colon [0,1] \rightarrow a$ such that $a(0)=a(1)=p(a)$.
For $s \in [-2,0]$, $\Sigma_{s}=\Sigma \cap A_s$
Let
$\Sigma^+=\Sigma_{-1} \cup k\{(t-1, a(\alpha),t);(t,\alpha)\in[0,1]^2\}$
so that $\partial \Sigma^+= k (a \times \{1\})$.

Let
$p(a)^+=(p(a),1)=(0,p(a),1) \in a \times [-1,1] \subset (\partial A =\{0\} \times \partial A)$.

Let $T(a)=\{((a(v),0),(a(w),+1));(v,w)\in [0,1]^2,v\geq w \}.$

Let $A(a)$ be the closure of $\{((a(v),0),(a(v),t));(t,v)\in ]0,1] \times [0,1]\}$.

Let $\mbox{diag}(n)\Sigma$ be the positive normal section of $ST(A)_{|\Sigma}$, and let $e(\Sigma(a)=\frac{\Sigma}{k})= \frac{g+k-1}{k}$ where $g$ is the genus of $\Sigma$.
$$e(\Sigma(a))=\frac{-\chi(\Sigma)}{2k} + \frac{1}{2}.$$

\begin{lemma}
\label{lemf(a)}
With the notation above
$$\begin{array}{ll}F^2(a)=& A(a) + T(a) - p(a) \times \frac{1}{k}\Sigma^+ -\frac{1}{k}\Sigma \times p(a)^+ \\&+\frac{1}{k}\mbox{diag}(n)\Sigma + e(\Sigma(a))[ST(\ast)]
\\&-\sum_{(j,{\ell}) \in \{1, \dots g_A\}^2}\langle \Sigma(a),\Sigma(a_j),\Sigma(a_{\ell}) \rangle z_j \times z_{\ell}\end{array}$$
is null-homologous in $C_2(A)$.
\end{lemma}
\bp
For $k=1$, (when we are dealing with integral homology handlebodies, for example) it is a direct consequence of Lemma~\ref{lemhomdiagSbry} and Lemma~\ref{homctimec} above.
Note that $T(a)=a\times_{p(a),\geq}(a\times\{1\})$
with the notation of Lemma~\ref{lemhomdiagSbry}.

Let us now focus on the case $k>1$.
(This is similar to \cite[Lemma 5.3]{sumgen}. 
However, we present a simpler independent proof below.)
Without loss, assume that 
$$\Sigma \cap ([-2,-1] \times \partial A)=\{(t-2,a(\alpha),\frac{(j-1)(1-t)}{k});(t,\alpha)\in [0,1]^2;j\in{1,2,\dots,k}\}$$
and change the definition of $\Sigma^+$ for the proof so that
$$\Sigma^+ \cap ([-1,0] \times \partial A)=k([-1,0] \times a \times \{1\})$$
and
$$\Sigma^+ \cap ([-2,-1]\times\partial A)=\{(t-2,a(\alpha),\frac{(j-\frac12)(1-t)}{k}+t);(t,\alpha)\in [0,1]^2;j\in{1,2,\dots,k}\}$$
as in the picture below that presents $\Sigma \cap ([-2,-1]\times p(a) \times[-1,1]) $ as the thick lines and $\Sigma^+ \cap ([-2,-1] \times p(a) \times[-1,1])$ as the thin lines when $k=3$.

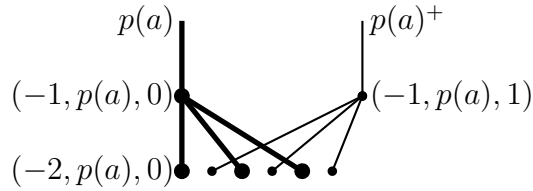
\begin{center}
\begin{pspicture}[shift=-0.4](0,-.7)(2.4,2.9) 
\psline[linewidth=2pt]{*-}(0,0)(0,2) 
\psline[linewidth=2pt]{*-*}(0.8,0)(0,1)
\psline[linewidth=2pt]{*-}(1.6,0)(0,1)
\psline{*-}(0.4,0)(2.4,1)(2.4,2)
\psline{*-*}(1.2,0)(2.4,1)
\psline{*-*}(2,0)(2.4,1)
\rput[r](-.1,2){$p(a)$}
\rput[r](-.1,1){$(-1,p(a),0)$}
\rput[r](-.1,0){$(-2,p(a),0)$}
\rput[l](2.5,1){$(-1,p(a),1)$} 
\rput[l](2.5,2){$p(a)^+$}
\end{pspicture}
\end{center}

Let $\Sigma_{-2}=\Sigma \cap A_{-2}$, $\partial \Sigma_{-2} =\cup_{j=1}^k \left(\{-2\}\times a \times\{\frac{j-1}{k}\}\right)$, and let $\Sigma_{-2}^+$ be a parallel copy of $\Sigma$ on its positive side with boundary $\partial(-\Sigma^+ \cap ([-2,-1]\times \partial A))$.

Glue abstract disks $D_j$ with respective boundaries $\{-2\}\times (-a) \times\{\frac{j-1}{k}\}$ on 
$\partial \Sigma_{-2}$
(resp. $D_j^+$ with boundaries $\{-2\}\times (-a) \times\{\frac{j-\frac12}{k}\}$ on 
$\partial \Sigma_{-2}^+$), and let $S$ (resp. $S^+$) be the obtained closed surface.
For $j=1,\dots, k$, set $p_j=(-2,p(a),\frac{j-1}{k}) \in \partial A_{-2}$ and $p^+_j=(-2,p(a),\frac{j-\frac12}{k})$.
Then it follows from Proposition~\ref{prophomdiagS} that
$$C(S)=\mbox{diag}(S\times S^+)- p_1 \times S^+ -S \times p^+_k - \sum_{i=1}^{2g} c_i \times c^{\ast +}_i$$
is null-homologous in $H_2(S\times S^+)$.
Let $[p_1,p_j]$ (resp. $[p_j^+,p_k^+]$) denote a path in $\left(\Sigma_{-2} \setminus \cup_{i=1}^{2g}c_i\right)$ from $p_1$ to $p_j$ (resp. in $\left(\Sigma^+_{-2} \setminus \cup_{i=1}^{2g}c^+_i\right)$ from $p_j^+$ to $p_k^+$).
Adding the null-homologous cycles
$$\partial (-[p_1,p_j] \times D_j^+)=p_1\times D_j^+ -p_j\times D_j^+
+[p_1,p_j]\times \partial D_j^+,$$ 
$$\partial (D_j \times [p_j^+,p_k^+])=D_j \times p^+_k -D_j \times p^+_j +\partial D_j \times [p_j^+,p_k^+],$$
for $j=1,\dots,k$, and the null-homologous cycles of Lemma~\ref{lemhomdiagSbry}
$$(-C_{\ast,\leq}(D_j,D_j^+)) $$ to $C(S)$ transforms it to the still null-homologous cycle
$$\begin{array}{ll}C(\Sigma_{-2})=&\mbox{diag}(\Sigma_{-2}\times \Sigma_{-2}^+)- p_1 \times \Sigma_{-2}^+ -\Sigma_{-2} \times p^+_k - \sum_{i=1}^{2g} c_i \times c^{\ast +}_i\\& +\sum_{j=1}^k(\partial D_j \times [p_j^+,p_k^+]+[p_1,p_j]\times \partial D_j^+ +\partial D_j \times_{p(a),\leq} \partial D_j^+).\end{array}$$
This cycle can be naturally continuously extended from the level $\{-2\} \times \Sigma$ to the level $\{0\} \times \Sigma$ to become naturally homologous to $kF^2(a)$ that will be therefore null-homologous
provided that 
$$\sum_{j=1}^k\left(lk(-a, \lim_{s\rightarrow -1}[\{s\} \times p_j^+,\{s\} \times p_k^+]) + lk(\lim_{s\rightarrow -1} [\{s\} \times p_1,\{s\} \times p_j],-a^+)\right)=k-1.$$
To conclude, it suffices to prove this equality. When $s$ approaches $(-1)$,
$[\{s\} \times p_j^+,\{s\} \times p_k^+]$ becomes a loop on $\Sigma^+_{-1}=\Sigma^+ \cap A_{-1}$, its linking number with $(-a)$ is its intersection with $\frac{-1}{k}\Sigma$ that only occurs where $\Sigma$ and $\Sigma^+$ intersect, in $([-2,-1] \times \partial A)$. The intersection can be seen in the picture below that $\lim_{s\rightarrow -1}[\{s\} \times p_j^+,\{s\} \times p_k^+]$ intersects as the oriented paths (for $j=1$ and $k=3$).

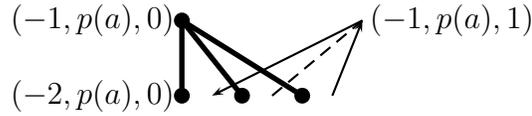
\begin{center}
\begin{pspicture}[shift=-0.4](0,-.2)(2.4,1.2) 
\psline[linewidth=2pt]{*-*}(0,0)(0,1) 
\psline[linewidth=2pt]{*-}(0.8,0)(0,1)
\psline[linewidth=2pt]{*-}(1.6,0)(0,1)
\psline{->}(2.4,1)(0.4,0)
\psline[linestyle=dashed](1.2,0)(2.4,1)
\psline{->}(2,0)(2.4,1)
\rput[r](-.1,1){$(-1,p(a),0)$}
\rput[r](-.1,0){$(-2,p(a),0)$}
\rput[l](2.5,1){$(-1,p(a),1)$} \end{pspicture}
\end{center}

In particular, since the positive normal to $\Sigma$ goes from left to right,
we see that $$lk(-a, \lim_{s\rightarrow -1}[\{s\} \times p_j^+,\{s\} \times p_k^+])=\frac{k-j}{k}.$$
Similarly, $lk(\lim_{s\rightarrow -1} [\{s\} \times p_1,\{s\} \times p_j],-a^+)=\frac{j-1}{k}$ and we are done.
\eop

\begin{lemma}
If $A$ is a rational homology handlebody such that $H_1(A)=\oplus_{j=1}^{g(A)}[z_j]$, then
$$H_3(C_2(A);\QQ(t))=\oplus_{j=1}^{g(A)} \QQ(t)[ST(z_j)].$$
\end{lemma}
\bp
$C_2(A)$ and $C_2(\mathring{A})$ have the same homotopy type, that is the homotopy type of $\mathring{A}^{2} \setminus \mbox{diag}$.
$H_3(\mathring{A}^2)=0$, 
$H_4(\mathring{A}^{2},\mathring{A}^{2} \setminus \mbox{diag})=\oplus_{j=1}^{g(A)} \QQ(t)[z_j \times B^3]$.
\eop

\noindent{\sc Proof of Lemma~\ref{lemprephalt}:}
We drop the useless superscripts $(i)$ in this proof.
Since $$\langle ST(p_M(\check{z}_{k})),F^2(a_j)\rangle_{e,\partial \tilde{C}_2(A)}=\delta_{jk},$$
according to the above two lemmas, it suffices to prove that
$$\langle\tilde{F}_3-\sum_{(i,k,\ell)\in \{1,\dots g(A)\}^3} \CI_{AA^{\prime}}(a_i\wedge a_k \wedge a_{\ell})lk_e(z_{\ell},z_k) ST(p_M(\check{z}_{i})),F_{A^{\prime}}^2(a_j) \rangle_{e,\partial \tilde{C}_2(A^{\prime})}=0$$
for all $j$. Fix $j$ and set $a=a_j$.
Since $F^2(a)$ vanishes according to Lemma~\ref{lemf(a)},
$$\langle F_{\fvarM}, F^2(a) \rangle_{e,\TCM} =0.$$

$$\langle F_{\fvarM},\mbox{diag}(n)\Sigma(a) + e(\Sigma(a))[ST(\ast)]\rangle_{e, \TCM}
=\langle F_{\fvarM},\frac{2k\mbox{diag}(n)\Sigma(a)}{2k} -\frac{\chi(k\Sigma(a))}{2k}[ST(\ast)]\rangle_{e, \TCM} + \frac{1}{2}$$
Assume without loss that $\tau^{-1}(.,e_1)$ is a positive normal to $k\Sigma(a)$ along $a$.
Then, according to Lemma~\ref{lemtruetriv},
$$\begin{array}{ll}\langle F_{\fvarM},\mbox{diag}(n)\Sigma(a) + e(\Sigma(a))[ST(\ast)]\rangle_{e, \TCM}
&=\langle F_{\fvarM},\frac{2s_{\tau}(k\Sigma(a);e_1)}{2k} \rangle_{e, \TCM} + \frac{d(\tau,a)+1}{2}
\\&=\frac{d(\tau,a)+1}{2}.\end{array}$$
Similarly, according to Lemma~\ref{lempseudotriv},
$$\langle \tilde{F}_3,\mbox{diag}(n)\Sigma^{\prime}(a) + e(\Sigma^{\prime}(a))[ST(\ast)]\rangle_{e,\partial \tilde{C}_2(A^{\prime})}
=\langle \tilde{F}_3,\frac{2s_{\tau(A^{\prime})}(k^{\prime}\Sigma^{\prime}(a);e_1)}{2k^{\prime}} \rangle_{e,\partial \tilde{C}_2(A^{\prime})} + \frac{d(\tau,a)+1}{2}.$$
Since this is equal to 
$\frac{d(\tau,a)+1}{2}$
according to Lemma~\ref{lemintpseudo},
$$\langle F_{\fvarM},\mbox{diag}(n)\Sigma(a) + e(\Sigma(a))[ST(\ast)]\rangle_{e, \TCM}
=\langle \tilde{F}_3,\mbox{diag}(n)\Sigma^{\prime}(a) + e(\Sigma^{\prime}(a))[ST(\ast)]\rangle_{e,\partial \tilde{C}_2(A^{\prime})}.$$

Since $A(a) + T(a) \subset C_2(\partial A)$, $\langle\tilde{F}_3,A(a) + T(a)\rangle_{e,\partial \tilde{C}_2(A^{\prime})} =\langle F_{\fvarM},A(a) + T(a)\rangle_{e, \TCM}$.
Similarly, setting $\Sigma_{[-4,0]}=\Sigma_0 \setminus \mbox{Int}(\Sigma_{-4})$,
and assuming here $\Sigma_{[-4,0]}=\Sigma^{\prime}_{[-4,0]}=[-4,0] \times a$ without loss,
$$\langle\tilde{F}_3, - p(a) \times \Sigma_{[-4,0]}^{\prime+}(a) -\Sigma_{[-4,0]}^{\prime}(a) \times p(a)^+\rangle_{e,\partial \tilde{C}_2(A^{\prime})}=\langle F_{\fvarM}, - p(a) \times \Sigma_{[-4,0]}^+(a) -\Sigma_{[-4,0]}(a) \times p(a)^+\rangle_{e, \TCM}$$
while the normalization conditions imply that
$$\langle\tilde{F}_3, - p(a) \times \Sigma_{-4}^{\prime+}(a) -\Sigma_{-4}^{\prime}(a) \times p(a)^+\rangle_{e,\partial \tilde{C}_2(A^{\prime})}=\langle F_{\fvarM}, - p(a) \times \Sigma_{-4}^+(a) -\Sigma_{-4}(a) \times p(a)^+\rangle_{e, \TCM}.$$

Therefore,
$$\langle \tilde{F}_3,F_{A^{\prime}}^2(a_j)\rangle_{e,\partial \tilde{C}_2(A^{\prime})}-\langle F_{\fvarM}, F^2(a_j) \rangle_{e, \TCM}=\sum_{(k,{\ell}) \in \{1, \dots g_A\}^2}\CI_{AA^{\prime}}(a_j\wedge a_k \wedge a_{\ell})lk_e(z_{\ell},z_k),$$
and we are done.
\eop

\newpage 
\section{Interactions with $K$}
\setcounter{equation}{0}
\label{secspecK}

\subsection{Introduction}
\label{subintroint}
In this section, we shall prove 
\begin{proposition}
\label{propdenwithoutz}
 $$\delta(M)(\fvar)\delta(M)(\svar)\delta(M)(\tvar)\CQ(\KK) \in \frac{\QQ[\fvar^{\pm 1},\svar^{\pm 1}, \tvar^{\pm 1}]}{(\fvar \svar \tvar=1)}.$$
\end{proposition}

In order to do this, we shall fix the chains $F_{\fvarM}$, $F_{\svarM}$ and $F_{\tvarM}$
outside $\tilde{C}_2(M\setminus T(K))$, or equivalently we shall
fix the interactions of $T(K)$ via $F$.
More precisely, (in Proposition~\ref{propnorFK}) we shall define chains $\Phi_{\fvarM}$, $\Phi_{\svarM}$ and $\Phi_{\tvarM}$ that are similar to $F_{\fvarM}$, $F_{\svarM}$ and $F_{\tvarM}$ and that can replace them without loss (by Proposition~\ref{propnorFKb}) such that 
$$\CQ(\KK,\tau)=\langle F_{\fvarM} , F_{\svarM}, F_{\tvarM} \rangle_{e,\TCM}=\langle \Phi_{\fvarM} , \Phi_{\svarM}, \Phi_{\tvarM} \rangle_{e,\TCM}.$$

We shall fix the chains $\Phi_{\cvarM}$ for $\cvarM=\fvarM$, $\svarM$ and $\tvarM$, so that $\Phi_{\fvarM}$, $\Phi_{\svarM}$ and $\Phi_{\tvarM}$ intersect only on $\tilde{C}_2(M\setminus T(K))$, and  $$C_{\cvarM}=\delta(M)(t)(\Phi_{\cvarM}\cap\tilde{C}_2(M\setminus T(K))$$ is rational for all $\cvarM$, so that
$$\delta(M)(\fvar)\delta(M)(\svar)\delta(M)(\tvar)\CQ(\KK,\tau)=\langle C_{\fvarM}, C_{\svarM} , C_{\tvarM} \rangle_{e,\tilde{C}_2(M\setminus T(K))} \in \QQ[\fvar^{\pm 1}, \svar^{\pm 1}, \tvar^{\pm 1}]$$
 where $C_{\cvarM}$ will be defined by boundary conditions (see Proposition~\ref{propdefbord}) (like in the work \cite{Ju} of Julien March\'e, so that
we get a definition of $\CQ(\KK)$ in the spirit of this work). The precise definition is stated in Subsection~\ref{subdefbord}.

\subsection{An alternative definition of $\CQ(\KK)$ with boundary conditions}
\label{subdefbord}

Let us first introduce some notation.
Write the sphere $S^2$ as the quotient of $[0,8] \times S^1$ where $\{0\}\times S^1$ is identified to a single point (the North Pole of $S^2$) and $\{8\}\times S^1$ is identified to another single point (the South Pole of $S^2$).
When $\alpha \subset [0,8]$, $D^2_{\alpha}$ denotes the image of $\alpha \times S^1$ via the quotient map $q$. For example, $D^2_{[1,8]}$ is a disk. Embed $D^2_{[1,8]} \times S^1$ as a tubular
neighborhood of $K$, so that $K=\{\ast_{\qvarM}\}\times S^1$ for some $\ast_{\qvarM} \in \partial D^2_{[0,6]}$, and $K_{\parallel}=\{q_{\fvarM}\}\times S^1$ for some $q_{\fvarM} \in \partial D^2_{[0,2]}$, and let $M_{[0,1]}=M\setminus (D^2_{]1,8]} \times S^1)$.
More generally, let 
$$\begin{array}{llll}
   r \colon & M & \rightarrow & [1,8]\\
  & x \in M_{[0,1]} & \mapsto & 1\\
& (q(t,z_{\qvarM}),z) \in D^2_{[1,8]} \times S^1&\mapsto &t.\\
\end{array}
$$
When $\alpha \subset [0,8]$, $M_{\alpha}=r^{-1}(\alpha)$
and $\tilde{M}_{\alpha}=p_M^{-1}(M_{\alpha})$.

Consider a map $f\colon M \rightarrow S^1$ that coincides with the projection
onto $S^1$ on $D^2_{[1,8]} \times S^1$, and a lift of this map $\tilde{f}\colon \tilde{M} \rightarrow \RR$.
Embed $\tilde{M}_{]1,6]}=p_M^{-1}(D^2_{]1,6]} \times S^1)$ in $\RR^3$, seen as $\CC \times \RR$, 
as $\{z \in \CC; 1< |z| \leq 6\} \times \RR$ naturally so that the projection on $\RR$ is $\tilde{f}$. Here $\CC$ is thought of as horizontal and $\RR$ is vertical.
This embedding induces a trivialisation $\tau$ on $TM_{]1,6]}$ that we extend on $TM_{[0,6]}$.
(The only obstruction to do so would be the obstruction to extend it to $(S\setminus D^2_{]1,8]})$ that vanishes.) This trivialisation respects the product structure with
$\RR$ on $\tilde{M}_{]1,6]}$, we also extend it on $\tilde{M}_{]1,8]}$
so that it still respects the product structure with
$\RR$, there.

\noindent{ \em Construction of a map $\pi \colon p^{-1}\left((M_{[0,5]}^2 \setminus M_{[0,3[}^2)\setminus \mbox{\rm diag}({M}_{[3,5]}^2)  \right)\rightarrow S^2$}\\
Let $$\begin{array}{llll}\chi \colon &[-4,4]& \rightarrow& [0,1]\\
&t \in [-0.5,4] &\mapsto &1\\
&t \in [-4,-1] &\mapsto &0
\end{array}$$
be a smooth map. Recall that $\tilde{M}_{]1,5]}$ is embedded in $\RR^3$ that is seen as $\CC \times \RR$.
When $(u,v) \in \left(\tilde{M}_{[0,5]}^2 \setminus \tilde{M}_{[0,3[}^2\right) \setminus \mbox{diag}(\tilde{M}_{[0,5]}^2)$, set
$$U(u,v)=(1-\chi(r(u) -r(v)))(0,\tilde{f}(u)) + \chi(r(u)-r(v))u$$
$$V(u,v)=(1-\chi(r(v) -r(u)))(0,\tilde{f}(v)) + \chi(r(v)-r(u))v$$
so that $(U(u,v),V(u,v)) \in (\RR^3)^2 \setminus \mbox{diag}$.

Define $$\begin{array}{llll}\pi \colon &p^{-1}\left((M_{[0,5]}^2 \setminus M_{[0,3[}^2)\setminus \mbox{diag}({M}_{[3,5]}^2) \right)&\rightarrow& S^2\\
 &(u,v)&\mapsto & \frac{V(u,v)-U(u,v)}{\parallel V(u,v)-U(u,v)\parallel}.
         \end{array}$$
The map $\pi$ extends naturally to $\tilde{C}_2(M_{[0,5]})\setminus \tilde{C}_2(M_{[0,3[})$.

\begin{proposition}
 \label{propdefbord}
Let $q_{\fvarM}$, $q_{\svarM}$ and $q_{\tvarM}$ be three distinct points
on $\partial D^2_{[0,2]}$, and let $\fvarM$, $\svarM$ and $\tvarM$ be three distinct vectors
of $S^2$ whose vertical coordinate is in $]0,\frac{1}{50}[$.
For $\cvarM=\fvarM, \svarM$ or $\tvarM$, let $K_{\cvarM}=q_{\cvarM}\times S^1 $,
then
there exist an element $J_{\Delta}$ of $\frac{1}{\delta(M)}\QQ[t,t^{-1}]$, and a $4$-dimensional rational chain $C_{\cvarM}$ of $\tilde{C}_2(M_{[0,3]})$ whose boundary is
$$\delta(M)\left(\pi_{|\partial \tilde{C}_2(M_{[0,3]}) \setminus \partial \tilde{C}_2(M_{[0,3[})}^{-1}(\cvarM) \cup s_{\tau}(M_{[0,3]};\cvarM) \cup (-J_{\Delta}) K_{\cvarM} \times_{\tau}S^2\right),$$
and that is transverse
to $\partial \tilde{C}_2(M_{[0,3]})$.
\end{proposition}
Note that this proposition is true when $M=S^1 \times S^2$, with $J_{\Delta}=0$, because in this case, $\pi$ extends to $\tilde{C}_2(M_{[0,3]})$ and $C_{\cvarM}=\pi_{| \tilde{C}_2(M_{[0,3]})}^{-1}(\cvarM)$ fulfills the conditions. It will be proved in general in Subsection~\ref{subsecproofC}.

\begin{proposition}
 \label{propdefbord2}
With the notation of Proposition~\ref{propdefbord},
$$\delta(M)(\fvar)\delta(M)(\svar)\delta(M)(\tvar)\CQ(\KK,\tau)=\langle C_{\fvarM}, C_{\svarM} , C_{\tvarM} \rangle_{e,\tilde{C}_2(M_{[0,3]})}$$
\end{proposition}

\begin{proposition}
 \label{propdefbord3}
With the notation of Proposition~\ref{propdefbord},
$$J_{\Delta}=\frac{t\Delta^{\prime}(t)}{\Delta(t)}.$$
\end{proposition}

These propositions, that will be proved later in this section, imply Proposition~\ref{propdenwithoutz}. They give an alternative definition of $\CQ(\KK)$ together with Proposition~\ref{propvartau}.

\subsection{Fixing the interactions with $K$}

Let $\ast_{\fvarM}$, $\ast_{\svarM}$, $\ast_{\tvarM}$ and $\ast_{\qvarM}$ be $4$ distinct points
in $D^2_{[6,8]}$, $\ast_{\qvarM} \in \partial D^2_{[6,8]}$, let $p_{\fvarM}$, $p_{\svarM}$ and $p_{\tvarM}$ be three distinct points
on $\partial D^2_{[0,5]}$ close to each other and such that the three distances between two of them are pairwise distinct, and 
let $h$
be a positive number smaller than $1/4$. Consider the three points
$\fvarM=\frac{(p_{\fvarM},h)}{\parallel (p_{\fvarM},h) \parallel}$, $\svarM=\frac{(p_{\svarM},h)}{\parallel (p_{\svarM},h) \parallel}$, $\tvarM=\frac{(p_{\tvarM},h)}{\parallel (p_{\tvarM},h) \parallel}$ of $S^2$.

Let $P \colon \tilde{M}^2 \rightarrow \TCMD$ be the canonical quotient map.
Define $C^5_{[4,8]}$ as the closure in $\tilde{C}_2(M_{[4,8]})$ of 
$P\left(\{((z\in D^2_{[4,8]},0);(z_2\in D^2_{[4,8]},t \in \RR\setminus \ZZ))\}\right)$
where $M_{[4,8]} = D^2_{[4,8]} \times S^1$.
Then
$$\tilde{C}_2(M_{[4,8]})=S^1 \times C^5_{[4,8]} $$
where $$\left(\exp(2i\pi u),P\left((z,0);(z_2,t\notin \ZZ) \right)\right)=P\left((z,u);(z_2,t+u) \right) $$ for $u \in[0,1]$.

The proof of the following proposition will be given in Subsection~\ref{subpfnorFK} assuming Proposition~\ref{propdefbord}.

\begin{proposition}
\label{propnorFK}
For $\cvarM=\fvarM$, $\svarM$ or $\tvarM$, 

\begin{itemize}
\item The map $\pi$ is regular from $\tilde{C}_2(M_{[0,5]})\setminus \tilde{C}_2(M_{[0,3[})$ to $S^2$, it factors through $C^5_{[4,5]}$ on $(\tilde{C}_2(M_{[4,5]})=S^1 \times C^5_{[4,5]})$, and
$\pi^{-1}(\cvarM )$ is a $4$--submanifold of $\tilde{C}_2(M_{[0,5]})\setminus \tilde{C}_2(M_{[0,3[})$.
\item Let $G(\cvarM )$ be the closure in $\TCM$ of $$P\left(\{(m;(p_{\cvarM},\tilde{f}(m) + h +u));m \in \tilde{f}^{-1}(]0,1]) \subset \tilde{M},u\in [0,1] \setminus \{1-h\} \}\right).$$
\item Let $G_{\iota}(\cvarM )$ be the closure in $\TCM$ of $$P\left(\{((-p_{\cvarM},\tilde{f}(m) - h+u);m);u\in [0,1] \setminus \{h\},m \in \tilde{f}^{-1}(]0,1]) \subset \tilde{M} \}\right)$$ oriented as $]0,1[ \times \tilde{M}$.
\end{itemize}
There exists a rational $3$--chain $E(\cvarM )$ of $C^5_{[4,8]}$ 
such that the boundary of
$$\Phi_{\cvarM} = \frac{1}{t-1}(tG(\cvarM )+G_{\iota}(\cvarM )) +\pi^{-1}(\cvarM ) + \frac{1}{\delta}C_{\cvarM} + S^1 \times E(\cvarM )$$
is
$$s_{\tau}(M;\cvarM) - J_{\Delta} ST(q_{\cvarM} \times S^1)
-\frac{1}{1-t}ST(\{p_{\cvarM}\} \times S^1) -\frac{1}{1-t} ST(\{-p_{\cvarM}\} \times S^1)+ST(\{\ast_{\cvarM}\} \times S^1)$$
where the $4$--chain $S^1 \times E(\cvarM )$ is in $\left(\tilde{C}_2(M_{[4,8]})=S^1 \times C^5_{[4,8]}\right)$, and $C_{\cvarM}$ is the rational chain of Proposition~\ref{propdefbord}.
\end{proposition}

\begin{proposition}
\label{propnorFKb}
\begin{enumerate}
\item The chains $\Phi_{\fvarM}$, $\Phi_{\svarM}$ and $\Phi_{\tvarM}$
of Proposition~\ref{propnorFK} have no algebraic triple intersection outside $\tilde{C}_2(M_{[0,3]})$.
\item Set $K=\{\ast_{\qvarM}\}  \times S^1$, and recall Notation $A(K)$ from Lemma~\ref{lemak} then for $\cvarM=\fvarM$, $\svarM$ or $\tvarM$,
$$\langle \Phi_{\cvarM},A(K) \rangle_{e,\TCM}=0.$$
\item $$\CQ(\KK,\tau)=\langle \Phi_{\fvarM} , \Phi_{\svarM}, \Phi_{\tvarM} \rangle_{e,\TCM}=\frac{1}{\delta(M)(\fvar)\delta(M)(\svar)\delta(M)(\tvar)}\langle C_{\fvarM} , C_{\svarM}, C_{\tvarM} \rangle_{e,\TCM}.$$

\end{enumerate}
\end{proposition}

\begin{proposition}
 \label{propGrat}
For the chains $F_{\cvarM}$ (for $\cvarM=\fvarM$, $\svarM$ or $\tvarM$) of Proposition~\ref{propexistF} and Proposition~\ref{propuninvtripl},
$\delta(M)(t-1)F_{\cvarM}$ can be assumed to be rational.
\end{proposition}
\noindent {\sc Proof of Propositions~\ref{propdefbord2}, \ref{propnorFKb}, \ref{propGrat} and \ref{propdefbord3} assuming Propositions~\ref{propdefbord} and \ref{propnorFK}: }
Proposition~\ref{propdefbord2} is contained in Proposition~\ref{propnorFKb}.
Observe that the boundaries of the $\Phi$ have no triple intersection on $\ZZ \times ST(M)$.
The supports of the chains $C_{\fvarM}$, $C_{\svarM}$ and $C_{\tvarM}$ on $\tilde{C}_2(M_{[0,3]}) \setminus \tilde{C}_2(M_{[0,3[})$ are in $\pi^{-1}(\fvarM )$, $\pi^{-1}(\svarM )$ and  $\pi^{-1}(\tvarM )$, respectively, so that there is no triple intersection there.
Clearly, the $C$ pieces cannot intersect the other ones, and two pieces
of the same kind among $\pi^{-1}$, $G$, $G_{\iota}$ cannot intersect each other so that the only triple intersections that could occur outside $\tilde{C}_2(M_{[0,3]})$ should
involve the three kinds or they should live in $\tilde{C}_2(M_{[4,8]})$.
Finally, they should live in $\tilde{C}_2(M_{[4,8]})$ where all the chains factor through $S^1$ as codimension $2$-chains in a $5$--dimensional manifold so that the algebraic intersection between the chains vanishes there.

Similarly, since $A(K)$ factors through $S^1$, the algebraic intersections $\langle \Phi_{\cvarM}, A(K) \rangle_{e,\TCM}$ vanish, too, and the first two assertions of Proposition~\ref{propnorFKb} are proved.

To construct chains $F_{\fvarM}$, $F_{\svarM}$ and $F_{\tvarM}$ such that $$\CQ(\KK,\tau)=\langle F_{\fvarM} , F_{\svarM}, F_{\tvarM} \rangle_{e,\TCM},$$ from the chains $\Phi_{\cvarM}$, use a product neighborhood $[-2,0] \times \partial \TCM $ of $\partial \TCM=\{0\} \times  \partial \TCM$, shrink the part of $\Phi$ inside $[-2,0] \times \partial \TCM$ into $[-2,-1] \times \partial \TCM$ using $(t,u)\mapsto ((t+2)/2-2,u)$. Consider three disjoint connected graphs $\Gamma_{\fvarM}$, $\Gamma_{\svarM}$ and $\Gamma_{\tvarM}$ embedded in $\left(D^2_{[2,8]} \setminus \{\ast_{\qvarM}\}\right)$ such that, for $\cvarM=\fvarM$, $\svarM$ and $\tvarM$, $\Gamma_{\cvarM}$ contains $p_{\cvarM}$, $(-p_{\cvarM})$, $q_{\cvarM}$ and $\ast_{\cvarM}$.
Then construct $F_{\cvarM}$ by completing the shrinked $\Phi_{\cvarM}$ by a cobordism 
between its boundary $\{-1\} \times \partial \Phi_{\cvarM}$ and 
$$\partial F_{\cvarM} = s_{\tau}(M;\cvarM) - \left(J_{\Delta} +\frac{1+t}{1-t} \right) ST(\{p_{\cvarM}\} \times S^1)$$
supported in $[-1,0] \times \left(s_{\tau}(M;\cvarM) \cup ST(\Gamma_{\cvarM} \times S^1)\right)$. It is clear that $$\langle F_{\cvarM},A(K) \rangle_{e,\TCM}=\langle \Phi_{\cvarM},A(K) \rangle_{e,\TCM}=0.$$ This proves Proposition~\ref{propGrat}.
According to Proposition~\ref{prophomTCM} and Theorem~\ref{thmstauM},
$J_{\Delta}$ must be equal to $\frac{t\Delta^{\prime}(t)}{\Delta(t)}$.
This proves Proposition~\ref{propdefbord3}.

Now since
$$\langle \Phi_{\fvarM} , \Phi_{\svarM}, \Phi_{\tvarM} \rangle_{e,\TCM}=\langle F_{\fvarM} , F_{\svarM}, F_{\tvarM} \rangle_{e,\TCM},$$ the last assertion of Proposition~\ref{propnorFKb} is also proved according to Proposition~\ref{propuninvtripl}.
\eop

Thus, we are left with the proofs of Proposition~\ref{propnorFK} and  Proposition~\ref{propdefbord}.

\subsection{Proof of Proposition~\ref{propnorFK} assuming Proposition~\ref{propdefbord}}
\label{subpfnorFK}

For $\alpha \subset [0,8]$,
let $M_{\alpha}\times_{\tilde{f},h} (\{p_{\cvarM}\} \times \RR)$ be $$P\left(\{(m;(p_{\cvarM},\tilde{f}(m) + h));m \in \tilde{f}^{-1}(]0,1]) \cap \tilde{M}_{\alpha} \}\right)$$
oriented by $\tilde{M}_{\alpha}$. Similarly,
let $(\{-p_{\cvarM}\} \times \RR)\times_{\tilde{f},h} M_{\alpha}$ be $$P\left(\{((-p_{\cvarM},\tilde{f}(m) - h);m);m \in \tilde{f}^{-1}(]0,1]) \cap \tilde{M}_{\alpha}\}\right)$$
oriented by $\tilde{M}_{\alpha}$.

$$\parallel(p_{\cvarM},h) \parallel \cvarM = (p_{\cvarM},h) = (5,z(\cvarM),h).$$

The boundary of $G(\cvarM )$ is
$$\partial G(\cvarM )= (1-t^{-1})\left(M \times_{\tilde{f},h} (\{p_{\cvarM}\} \times \RR)\right) \cup t^{-1} ST(\{p_{\cvarM}\} \times \RR).$$
(About the sign in front of $t^{-1} ST(\{p_{\cvarM}\} \times \RR)$,
$M$ reads $(\tilde{f}^{-1}(s), s)$ where $s$ is the parameter for $\RR$ so that we have the orientation $(-(u,\tilde{f}^{-1}(s),s))$ for $G(\cvarM )$. When $u>1-h$, $u$ is an inward normal, $\tilde{f}^{-1}(s)$ has the orientation opposite to that of a surrounding sphere, but it is first instead of second.)

The chain $G_{\iota}(\cvarM )$ is the closure in $\tilde{C}_2(M)$ of $$P\left(\{((-p_{\cvarM},\tilde{f}(m) - h+u);m);u\in [0,1]\setminus\{h\},m \in \tilde{f}^{-1}(]0,1]) \subset \tilde{M} \}\right)$$ oriented as $[0,1] \times \tilde{M}$.
$$\partial G_{\iota}(\cvarM )= (t-1)(\{-p_{\cvarM})\} \times \RR)\times_{\tilde{f},h}M)   \cup  ST(\{-p_{\cvarM}\} \times S^1).$$

Assuming Proposition~\ref{propdefbord}, the boundary of the $4$-chain $\left(\pi^{-1}(\cvarM ) \cup \frac{1}{\delta}C_{\cvarM}\right)$ is 
$$\pi_{|\partial \tilde{C}_2(M_{[0,5]}) \setminus \partial \tilde{C}_2(M_{[0,5[})}^{-1}(\cvarM ) \cup s_{\tau}(M_{[0,5]};\cvarM) \cup (-J_{\Delta})K_{\cvarM} \times_{\tau}S^2.$$
Set
$$\partial_{[0,5]\leftrightarrow 5}=\pi_{|\partial \tilde{C}_2(M_{[0,5]}) \setminus \partial \tilde{C}_2(M_{[0,5[})}^{-1}(\cvarM ).$$

$$\partial_{[0,5]\leftrightarrow 5}=\partial_{[0,4],5} \cup \partial_{5,[0,4]} \cup
\partial_{[4,5]\leftrightarrow 5}$$
where $$\partial_{[4,5]\leftrightarrow 5}= \pi_{|\partial C_2(M_{[4,5]}) \setminus \partial C_2(M_{[4,5[})}^{-1}(\cvarM )$$
$$\partial_{[0,4],5}=-(M_{[0,4]}\times_{\tilde{f},h} (\{p_{\cvarM}\} \times \RR)$$
(since $M^2$ is oriented as $(M,\mbox{outward normal to} \;M_{[0,5]}, \partial M_{[0,5]})$, and since $\pi$ maps $\partial M_{[0,5]}$ to $S^2$ in an orientation-preserving way, we get the minus sign above)
and 
$$\partial_{5,[0,4]}=-((\{-p_{\cvarM}\} \times \RR)\times_{\tilde{f},h} M_{[0,4]}.$$

Now, the boundary of the chain
$$\frac{1}{t-1}(tG(\cvarM )+G_{\iota}(\cvarM )) +\pi^{-1}(\cvarM ) + \frac{1}{\delta}C_{\cvarM}$$
is 
$\partial_r +\partial_e $
where $$\partial_r=s_{\tau}(M_{[0,5]};\cvarM) \cup (-J_{\Delta}) ST(K_{\cvarM}) +\frac{1}{t-1}ST(\{p_{\cvarM}\} \times S^1) +\frac{1}{t-1} ST(\{-p_{\cvarM}\} \times S^1)$$ is part of $\partial \Phi_{\cvarM}$ and where the remaining part
is 
$$\partial_e=\left(M_{[4,8]} \times_{\tilde{f},h} (\{p_{\cvarM}\} \times \RR)\right) +\left(\{-p_{\cvarM}\} \times \RR)\times_{\tilde{f},h}M_{[4,8]}\right) + \partial_{[4,5]\leftrightarrow 5} .$$

Recall that $$\tilde{C}_2(M_{[4,8]})=S^1 \times C^5_{[4,8]}$$
where
the orientation of $C^5_{[4,8]}$ is induced by this product structure.
Note that $\partial_e =S^1 \times \partial^2_e$, and that 
$s_{\tau}(M_{[5,8]};\cvarM)=S^1 \times s_{\tau}(M_{[5,8]}\cap f^{-1}(1);\cvarM)$.
Proposition~\ref{propnorFK} is now the consequence of the following lemma.

\begin{lemma}
There exists a rational $3$--dimensional chain $E(\cvarM )$ in $C^5_{[4,8]}$ such that
$$\partial E(\cvarM ) =\partial^2_e - s_{\tau}(M_{[5,8]}\cap f^{-1}(1);\cvarM) -ST(\ast_{\cvarM}).$$
\end{lemma}
\bp
Let $s_+(M_{[5,8]}\cap f^{-1}(1))$ be a section that is in the hemisphere of the $\RR$
direction and that coincides with $s_{\tau}(M_{[5,8]}\cap f^{-1}(1);\cvarM)$
on $\partial M_{[5,8]}\cap f^{-1}(1)$.
Then $(\partial^2_e-s_+(M_{[5,8]}\cap f^{-1}(1)))$ is a $2$--cycle in the closure of
$P\left(\{((z\in D^2_{[4,8]},0);(z_2\in D^2_{[4,8]},t \in ]0,1/2[))\}\right)$
that is homotopy equivalent to $D^2 \times D^2 \times [1/2]$. Therefore this $2$-cycle bounds a rational $3$--chain in $C^5_{[4,8]}$.
Now, $s_+(M_{[5,8]}\cap f^{-1}(1))- s_{\tau}(M_{[5,8]}\cap f^{-1}(1);\cvarM)
=ST(\ast_{\cvarM})$ since the difference of the degrees of the Gauss maps of the sections of $T(S^2 \subset \RR^3)$ given by the outward normal section, on one hand, and a trivialisation, on the other hand, is one.
\eop

\subsection{Homology of $\tilde{C}_2(M_{[0,3]})$}

Let $E=M_{[0,3]}$, $E$ is the exterior of $K$, $E$ is a rational homology torus, and we are going to compute the homology of $\tilde{C}_2(E)$
as in Section~\ref{sechomtcm}.
First note that
$H_2(E;\ZZ)=H_3(E;\ZZ)=0$, $H_0(E;\ZZ)=\ZZ$ and $H_1(E;\ZZ)/\mbox{Torsion}=\ZZ[K_{\fvarM}]$.

Let $\Sigma =S\cap E$, $\Sigma$ is obtained from $S$ by removing an open disk.
The pair $(M \setminus S, E \setminus \Sigma)$ has the same homology as the pair $(D^2 \times \RR, S^1\times \RR)$ by excision. Therefore, Lemma~\ref{lemhomMsetminusS}
implies the following similar lemma:

\begin{lemma}
\label{lemhomEsetminusSigma}
$H_i(E \setminus \Sigma;\ZZ)=\{0\}$ for any $i \geq 2$,
$H_0(E \setminus \Sigma;\ZZ)=\ZZ[\ast^+]$.\\
Let $(z_i)_{i=1, \dots 2g}$ and $(z^{\ast}_i)_{i=1, \dots, 2g}$ be two dual bases
of $H_1(\Sigma;\ZZ)$ such that 
$\langle z_i, z^{\ast}_j\rangle=\delta_{ij}.$
Then $$H_1(E \setminus \Sigma;\QQ) =\bigoplus_{i=1}^{2g} \QQ[z_i^{+} -z_i^{-}]$$
and for any $v\in H_1(E \setminus \Sigma;\QQ)$,
$v=\sum_{i=1}^{2g}lk(v,z^{\ast}_i)(z_i^{+} -z_i^{-}).$
\end{lemma}

The bases in the statements of Lemmas~\ref{lemhomMsetminusS} and \ref{lemhomEsetminusSigma} can and will be assumed to be the same.

Let $\tilde{E}$ be the infinite cyclic covering of $E$.

\begin{lemma}
\label{lemhomtilE}
$$H_{0}(\tilde{E}) = \frac{\Lambda_M}{(t_M-1)}, H_{2}(\tilde{E}) = \{0\}, H_{3}(\tilde{E}) = \{0\}$$
and
$$H_{1}(\tilde{E}) = \oplus_{i=1}^k\frac{\Lambda_M}{(\delta_i(M))}$$
where $\prod_{i=1}^k\delta_i(M)=\Delta(M)$ is the Alexander polynomial of $M$ and $E$.
\end{lemma}
\bp Compare the homology of $\tilde{E}$ with
the homology of $\tilde{M}$ given by Lemma~\ref{lemhomtilM} using the
long exact sequence associated with the pair $(\tilde{M},\tilde{E})$ whose homology is the homology of $(D^2\times \RR, S^1\times \RR)$.
Again, $H_i(D^2\times \RR, S^1\times \RR) =0$ when $i\neq 2$ and $H_2(D^2\times \RR, S^1\times \RR) =\QQ [D^2]$ where $[D^2]$ is the image
of the generator $[S]$ of $H_2(\tilde{M})$ under the composition of the natural map to $H_2(\tilde{M},\tilde{E})$ with the excision map.
\eop

Now, the homology of $\widetilde{E^2}$ can be computed as in Subsection~\ref{subhomTCMDmore} to find the following proposition.
Assume without loss that the $\Sigma_i$, the $c_i$ and the $C(\Sigma_i \times \Sigma_j)$ of Proposition~\ref{proptilMtwocomp} live in $\tilde{E}$ or in $\widetilde{E^2}$.

\begin{proposition}
\label{proptilEtwocomp}
The rational homology of $\widetilde{E^2}$ reads as follows\\
$\begin{array}{lll}
H_0(\widetilde{E^2})&=&\frac{\Lambda}{(t-1)}[\ast \times \ast]\\
H_1(\widetilde{E^2})&=&\frac{\Lambda}{(t-1)}[\mbox{\rm diag}(K^2)]\\
H_2(\widetilde{E^2})&=& \oplus_{(i,j) \in \{1,2,\dots,k\}^2}\frac{\Lambda}{\left(\delta_{\mbox{\tiny \rm min}(i,j)}\right)}[c_i \times c_j]\\
H_3(\widetilde{E^2})&=& \oplus_{(i,j) \in \{1,2,\dots,k\}^2}\frac{\Lambda}{\left(\delta_{\mbox{\tiny \rm min}(i,j)}\right)}[C(\Sigma_i \times \Sigma_j)]\\
H_4(\widetilde{E^2})&=&H_5(\widetilde{E^2})=0
\end{array}$\\
where the $\Sigma_i$, the $c_i$ and the $C(\Sigma_i \times \Sigma_j)$ are the same as in Proposition~\ref{proptilMtwocomp}.
\end{proposition}

\subsection{Proof of Proposition~\ref{propdefbord}}
\label{subsecproofC}

We must find
a $4$-dimensional rational chain $C_{\cvarM }$ of $\tilde{C}_2(E)$ whose boundary is
$$\delta(M)\left(\pi_{|\partial \tilde{C}_2(E) \setminus \partial \tilde{C}_2(M_{[0,3[})}^{-1}(\cvarM ) \cup s_{\tau}(E;\cvarM) \cup \lambda K_{\cvarM } \times_{\tau}S^2\right),$$
for some $\lambda$ in $\frac{1}{\delta(M)}\QQ[t,t^{-1}]$, and that is transverse
to $\partial \tilde{C}_2(E)$.

Set $$A=\delta(M)\left(\pi_{|\partial \tilde{C}_2(M_{[0,3]}) \setminus \partial \tilde{C}_2(M_{[0,3[})}^{-1}(\cvarM ) \cup s_{\tau}(E;\cvarM)\right).$$
$A$ is a $3$-cycle of $\partial \tilde{C}_2(E)$ that can be first pushed
inside $\tilde{C}_2(E)$ using a product cobordism in the neighborhood of the boundary towards a $3$-cycle $B$ that is parallel to $A$ but that does not meet the boundary anymore.
According to Proposition~\ref{proptilEtwocomp}, $\delta(M)B$ bounds a $4$-chain in $\widetilde{E^2}$. This chain can be assumed to be transverse to the preimage of the diagonal of $E^2$ under the covering map.
This chain can be modified by the addition of $4$-cycles of the form $\partial (B^3 \times G)$ where $G$ is a 
$2$--cobordism in the preimage of the diagonal of $E^2$, and $B^3$ is a normal section of the diagonal of $E^2$. Thus, our chain can be furthermore assumed to intersect this preimage along the preimage of the diagonal of $K_{\cvarM}^2$. This shows that there exists a $4$-dimensional rational chain $C_{\cvarM }$ of $\tilde{C}_2(E)$ whose boundary is
$$\delta(M)\left(\pi_{|\partial \tilde{C}_2(E) \setminus \partial \tilde{C}_2(M_{[0,3[})}^{-1}(\cvarM ) \cup s_{\tau}(E;\cvarM) \cup \lambda K_{\cvarM } \times_{\tau}S^2\right),$$ and that is transverse
to $\partial \tilde{C}_2(E)$, for some $\lambda$ of $\frac{1}{\delta(M)}\QQ[t,t^{-1}]$.
\eop

\newpage 

\section{On the augmentation of $\CQ$.}
\setcounter{equation}{0}
\label{secaug}

\subsection{The result and the sketch of the proof}

This section is devoted to the proof of the following theorem.

\begin{theorem}
\label{thmaugcas}
 Let $M_{\KK}$ denote the rational homology sphere obtained from $M$ by surgery along
$\KK$. Then $$\CQ(\KK \subset M)(1,1,1)=6\lambda(M_{\KK}).$$
\end{theorem}

The proof of this proposition heavily relies on Propositions~\ref{propdefbord} and \ref{propdefbord2}. 
Using the same notation as in these propositions, and denoting the meridian disk of the torus reglued during the surgery on $(M,\KK)$ by $D_{[7,8]}^2(\hat{K})$,  write 
$$M_{\KK} = M_{[0,7]} \cup_{\partial M_{[0,7]} \sim (-\partial (S^1 \times D_{[7,8]}^2(\hat{K}) ))} S^1 \times D_{[7,8]}^2(\hat{K})$$
where $$p_M(7z,u) \in \partial M_{[0,7]} = \{p_M(7z,u); z \in S^1, u \in [0,1]\}$$ 
(with the notation of Subsection~\ref{subdefbord})
is identified to 
$(z, \exp(2i\pi u) \in \partial D_{[7,8]}^2(\hat{K})) \in S^1 \times D_{[7,8]}^2(\hat{K})$.

Referring to the configuration space construction of the Walker invariant in Subsection~\ref{subdefcas}, assume that $B(M_{\KK})=B_{M_{\KK}}(3)$ is embedded in $M_{\KK}$ and that the ball 
$\left(M_{\KK} \setminus B(M_{\KK})\right)$ is a small ball inside $M_{]6,7[}$.
Recall $E=M_{[0,3]}$ and set $B(M_{\KK})_{[3,8]}=B(M_{\KK}) \setminus M_{[0,3[}$.
Note that $B(M_{\KK})_{[3,8]}$ is a solid torus minus a small open ball that is independent of $(M,\KK)$.

For $\cvarM=\fvarM, \svarM$ or $\tvarM$ that will be assumed to be almost horizontal, let $p(C_{\cvarM})$ be the image of the chain $C_{\cvarM}$ of Proposition~\ref{propdefbord} in 
$C_2(E)$ by the covering map. Since the covering map sends $t$ to $1$,
$$\partial p(C_{\cvarM}) = \pi_{|\partial \tilde{C}_2(M_{[0,3]}) \setminus \partial \tilde{C}_2(M_{[0,3[})}^{-1}(\cvarM) \cup s_{\tau}(M_{[0,3]};\cvarM).$$
We are going to construct $4$--chains $D_{\cvarM}$ of $C_2(B(M_{\KK}))$
\begin{itemize}
 \item with support outside
the interiors of $C_2(E=M_{[0,3]})$ and $C_2(B(M_{\KK})_{[3,8]})$
\item such that $D_{\fvarM}$, $D_{\svarM}$ and $D_{\tvarM}$ have no triple intersection,
\item and such that there exist transverse $4$--chains $E_{\cvarM}$ of $C_2(B(M_{\KK})_{[3,8]})$ independent of $M_{\KK}$

such that
$$6\lambda(M_{\KK})=\langle p(C_{\fvarM}) + D_{\fvarM} + E_{\fvarM},p(C_{\svarM}) + D_{\svarM} + E_{\svarM} ,p(C_{\tvarM}) + D_{\tvarM} + E_{\tvarM} \rangle + C(p_1)-\frac{p_1(\tau)}{4}$$
where $C(p_1)$ is a rational number that does not depend on $(M,\KK)$.\end{itemize}
Therefore, we shall have

$$6\lambda(M_{\KK})=\langle p(C_{\fvarM}),p(C_{\svarM}),p(C_{\tvarM}) \rangle_{C_2(E)} + \langle E_{\fvarM}, E_{\svarM}, E_{\tvarM} \rangle + C(p_1)-\frac{p_1(\tau)}{4}$$
where $$\CQ(\KK \subset M)(1,1,1)=\langle p(C_{\fvarM}),p(C_{\svarM}),p(C_{\tvarM}) \rangle_{C_2(E)}-\frac{p_1(\tau)}{4}$$
so that $$\CQ(\KK \subset M)(1,1,1)-6\lambda(M_{\KK})$$ will be independent of $(M,\KK)$.
The theorem will follow since $\CQ(S^1 \times S^2,S^1 \times u)=0=\lambda(S^3)$.

\subsection{The construction of $D_{\cvarM}$}
Set $B(M_{\KK})_{[5,8]} = B(M_{\KK})\setminus M_{[0,5[}$.
The chain $D_{\cvarM}$ will read
 $$p\left(\pi^{-1}_{|\tilde{C}_2(M_{[0,5]})\setminus \left( \tilde{C}_2(M_{[0,3[}) \cup \tilde{C}_2(M_{]3,5]})\right)}(\cvarM)\right) + d_{\cvarM} + d_{\cvarM}^{\prime}$$
where $d_{\cvarM}$ will be supported in $M_{[0,3]} \times B(M_{\KK})_{[5,8]}$
and $d_{\cvarM}^{\prime}$ that will be close to $\iota(d_{\cvarM})$ will be supported in $B(M_{\KK})_{[5,8]} \times M_{[0,3]}$.
The chain $d_{\cvarM}$
is made of the following $3$ pieces:
\begin{itemize}
\item
$-\bigcup_{u \in [0,1]}f_{|E}^{-1}(\exp(2i\pi (u-h))) \times p_M\left(\{p_{\cvarM}\} \times [0,u]\right),$
\item
$M_{[0,3]} \times [p_M(p_{\cvarM},0),3{\cvarM}]$\\
where $[p_M(p_{\cvarM},0),3{\cvarM}]$ is a path in $B(M_{\KK})_{[5,8]}$
that goes from $(p_M(p_{\cvarM},0) \in \partial M_{[0,5]}) $ to $(3\cvarM \in \partial B(M_{\KK}))$, and the paths $[p_M(p_{\fvarM},0),3{\fvarM}]$, $[p_M(p_{\svarM},0),3{\svarM}]$ and $[p_M(p_{\tvarM},0),3{\tvarM}]$ are pairwise disjoint,
\item
$f_{|E}^{-1}(\exp(2i\pi (-h))) \times \left(p_{\cvarM} \times D_{[5,8]}^2(\hat{K})\right)$\\
where $\left(p_{\cvarM} \times D_{[5,8]}^2(\hat{K})\right)$ is a meridian disk of $\hat{K}$ in $B(M_{\KK})_{[5,8]}$ with boundary $p_M(p_{\cvarM} \times [0,1])$.
\end{itemize}

Split

$$\partial p\left(\pi^{-1}_{|\tilde{C}_2(M_{[0,5]})\setminus \left( \tilde{C}_2(M_{[0,3[}) \cup \tilde{C}_2(M_{]3,5]})\right)}(\cvarM)\right)$$
as the following sum $$\partial_{[0,3],5} +  \partial_{5,[0,3]} + \partial_{[3,5]\leftrightarrow 3} + \partial_{[0,3]\leftrightarrow 3}$$
where the subscripts reflect the values of the map $r$ of the beginning of Subsection~\ref{subdefbord} of the points of an ordered pair of $C_2(M_{[0,5]})$.

$$\partial_{[0,3],5}=-\bigcup_{u \in [0,1]}f_{|E}^{-1}(\exp(2i\pi (u-h))) \times p_M\left(\{p_{\cvarM}\} \times \{u\}\right).$$
$$\partial_{5,[0,3]}=-\bigcup_{u \in [0,1]}p_M\left(\{-p_{\cvarM}\} \times \{u\}\right) \times f_{|E}^{-1}(\exp(2i\pi (u+h))).$$

\begin{lemma}
 $$\begin{array}{ll}\partial d_{\cvarM}=& -M_{[0,3]} \times \{3{\cvarM}\} - \partial_{[0,3],5}\\
&
+\bigcup_{u \in [0,1]}\partial f_{|E}^{-1}(\exp(2i\pi (u-h))) \times p_M\left(\{p_{\cvarM}\} \times [0,u]\right) \\
&+\partial M_{[0,3]} \times [p_M(p_{\cvarM},0),3{\cvarM}]\\
&+\partial f_{|E}^{-1}(\exp(2i\pi (-h))) \times \left(p_{\cvarM} \times D_{[5,8]}^2(\hat{K})\right).
\end{array}$$
\end{lemma}
\bp Note that the part $-f_{|E}^{-1}(\exp(2i\pi (1-h))) \times p_M\left(\{p_{\cvarM}\} \times [0,1]\right)$ cancels with $$f_{|E}^{-1}(\exp(2i\pi (-h))) \times \left(p_{\cvarM} \times \partial D_{[5,8]}^2(\hat{K})\right).$$
\eop

The chain $d_{\cvarM}^{\prime}$ is made of the following $3$ pieces:
\begin{itemize}
\item
$-\bigcup_{u \in [0,1]} p_M\left(\{-p_{\cvarM}\} \times [0,u]\right)\times f_{|E}^{-1}(\exp(2i\pi (u+h))) $,
\item
$[-3\cvarM,p_M(-p_{\cvarM},0)] \times M_{[0,3]}$
where $[-3\cvarM,p_M(-p_{\cvarM},0)]$ is a path in $B(M_{\KK})_{[5,8]}$
that goes from $(-3\cvarM)$ in the boundary of $B(M_{\KK})$ to $p_M(-p_{\cvarM},0)$,  and the paths $[-3\fvarM,p_M(-p_{\fvarM},0)]$, $[-3\svarM,p_M(-p_{\svarM},0)]$ and $[-3\tvarM,p_M(-p_{\tvarM},0)]$ are pairwise disjoint,
\item
$\left(\{-p_{\cvarM}\} \times D_{[5,8]}^2(\hat{K})\right) \times f_{|E}^{-1}(\exp(2i\pi h))$ where $\left(\{-p_{\cvarM}\} \times D_{[5,8]}^2(\hat{K})\right)$ is a meridian disk of $\hat{K}$ in $B(M_{\KK})_{[5,8]}$ with boundary $p_M(\{-p_{\cvarM}\} \times [0,1])$.
\end{itemize}

\begin{lemma}
 $$\begin{array}{ll}\partial d_{\cvarM}^{\prime}= &-\{-3{\cvarM}\} \times M_{[0,3]} - \partial_{5,[0,3]}\\
&-\bigcup_{u \in [0,1]}p_M\left(\{-p_{\cvarM}\}\times [0,u]\right) \times \partial f_{|E}^{-1}(\exp(2i\pi (u+h)))  \\
&-[-3\cvarM,p_M(-p_{\cvarM},0)] \times \partial M_{[0,3]}\\
&+\left(\{-p_{\cvarM}\} \times D_{[5,8]}^2(\hat{K})\right) \times \partial f_{|E}^{-1}(\exp(2i\pi (h))).
\end{array}$$
\end{lemma}
\bp The proof is similar, of course.
\eop

\begin{lemma}
 $$\partial (p(C_{\cvarM}) + D_{\cvarM})=\partial_r +\partial_e$$
where
$\partial_r=-M_{[0,3]} \times \{3{\cvarM}\}-\{-3{\cvarM}\} \times M_{[0,3]}+s_{\tau}(M_{[0,3]};\cvarM)$ and\\
$\begin{array}{ll}\partial_e=& \partial_{[3,5]\leftrightarrow 3}\\
&-\bigcup_{u \in [0,1]}p_M\left(\{-p_{\cvarM}\}\times [0,u]\right) \times \partial f_{|E}^{-1}(\exp(2i\pi (u+h))) \\
&+\bigcup_{u \in [0,1]}\partial f_{|E}^{-1}(\exp(2i\pi (u-h))) \times p_M\left(\{p_{\cvarM}\} \times [0,u]\right) \\
&-[-3\cvarM,p_M(-p_{\cvarM},0)] \times \partial M_{[0,3]}+\partial M_{[0,3]} \times [p_M(p_{\cvarM},0),3{\cvarM}]\\
&+\left(\{-p_{\cvarM}\} \times D_{[5,8]}^2(\hat{K})\right) \times \partial f_{|E}^{-1}(\exp(2i\pi (h)))+\partial f_{|E}^{-1}(\exp(2i\pi (-h))) \times \left(p_{\cvarM} \times D_{[5,8]}^2(\hat{K})\right).
\end{array}$
\end{lemma}

Note that the chains $D_{\cvarM}$ have no triple intersection. Indeed $D_{\cvarM}$ meets $M_{[0,3]} \times \partial M_{[0,5]}$ in 
$M_{[0,3]} \times p_M(p_{\cvarM} \times \RR)$ and the $p_M(p_{\cvarM} \times \RR)$ are pairwise disjoint so that two $D_{\cvarM}$ cannot intersect in $M_{[0,3]} \times \partial M_{[0,5]}$. Similarly, two $D_{\cvarM}$ cannot intersect in $\partial M_{[0,5]} \times M_{[0,3]}$. Then triple intersections of $D_{\fvarM}$, $D_{\svarM}$ and $D_{\tvarM}$ must be triple intersections of $d_{\fvarM}$, $d_{\svarM}$ and $d_{\tvarM}$ or triple intersections of $d_{\fvarM}^{\prime}$, $d_{\svarM}^{\prime}$ and $d_{\tvarM}^{\prime}$ that cannot occur under the further natural assumption that the meridian disks $((\pm p_{\cvarM}) \times D^2_{[5,8]}(\hat{K}))$ are disjoint.

\subsection{The existence of $E_{\cvarM}$}
The trivialisation $\tau$ of $M_{[0,7]}$ of Subsection~\ref{subdefbord} does not extend to $B(M_{\KK})$ as a genuine trivialisation but it extends as a pseudo-trivialisation $\tilde{\tau}$ of Definition~\ref{defpseudotrivrat} where $c$ is a parallel of $\hat{K}$, that is a meridian of $K$. Set $N=M_{\KK}$.
According to Proposition~\ref{propdefcaspseudo},
$$6\lambda(N)= \langle F_{N,\fvarM}(\tilde{\tau}),F_{N,\svarM}(\tilde{\tau}), F_{N,\tvarM}(\tilde{\tau})\rangle_{C_2(B(N))} -\frac{p_1(\tilde{\tau})}{4}$$
where the boundary of $F_{N,\cvarM}(\tilde{\tau})$ is defined in Definition~\ref{defpseudotrivrat}, for $\cvarM=\fvarM$, $\svarM$ or $\tvarM$.

Recall that $(B(M_{\KK})_{[3,8]}= B(M_{\KK})\setminus M_{[0,3[})$ is a solid torus minus a small open ball that is independent of $(M,\KK)$.

Note that $\partial_r$ is part of the wanted boundary $\partial F_{N,\cvarM}(\tilde{\tau})$, and that the cycle
 $$e_{\cvarM}=\left(\partial F_{N,\cvarM}(\tilde{\tau})-\partial_e - \partial_r \right)$$ is therefore supported in $\partial C_2(B(M_{\KK})_{[3,8]})$.

Thus, $e_{\cvarM}$ is a $3$--cycle of $ \partial C_2(B(M_{\KK})_{[3,8]})$ that is independent of $(M,\KK)$. 
Let $h_{\fvarM}$, $h_{\svarM}$ and $h_{\tvarM}$ be three distinct small positive numbers.
Let $\hat{K}_{\cvarM}$ denote a horizontal parallel of $\hat{K}$ that reads $p_M(3S^1 \times \{h_{\cvarM}\})$ in $\partial M_{[0,3]}$.
We are going to prove the following lemma.
\begin{lemma}
\label{lemmu}
 There exists $\mu \in \QQ$ such that $$e_{\cvarM}-\mu ST(\hat{K}_{\cvarM})$$
bounds a $4$--chain $E_{\cvarM}$ in $C_2(B(M_{\KK})_{[3,8]})$ transverse to $ \partial C_2(B(M_{\KK})_{[3,8]})$.
\end{lemma}

Let us show that Lemma~\ref{lemmu} implies Theorem~\ref{thmaugcas}.
Let $\Sigma_{\fvarM}$, $\Sigma_{\svarM}$ and $\Sigma_{\tvarM}$ be three disjoint
Seifert surfaces in $M_{[0,3]}$ for $\hat{K}_{\fvarM}$, $\hat{K}_{\svarM}$ and $\hat{K}_{\tvarM}$, respectively. Let $N=M_{\KK}$.

Then, we can use the chains
$$F_{N,\cvarM}(\tilde{\tau})=p(C_{\cvarM})+D_{\cvarM}+E_{\cvarM}+\mu ST(\Sigma_{\cvarM})$$
to compute $(6\lambda(N)-\CQ(\KK \subset M)(1,1,1))$.
Since the $ST(\Sigma_{\cvarM})$ are pairwise disjoint, and since the other pieces of $F_{N,\cvarM}(\tilde{\tau})$ only meet $ST(M_{[0,3]})$ on $s_{\tau}(M_{[0,3]};\cvarM)$, no triple intersection will involve the $ST(\Sigma_{\cvarM})$ that can be forgotten.
If $p(C_{\fvarM})$ met $D_{\svarM}$ or $E_{\svarM}$ outside $ST(\Sigma_{\svarM})$, it would be in $p(\pi^{-1}(\fvarM)\cap\pi^{-1}(\svarM))$, so they do not meet.
If $D_{\fvarM}$ met $E_{\svarM}$ outside $ST(\Sigma_{\svarM})$, it would be in $\partial E_{\svarM}$ and actually in $\partial D_{\svarM}$. Therefore, since the $D_{\cvarM}$ have no triple intersection,
$$6\lambda(N)=\langle p(C_{\fvarM}),p(C_{\svarM}),p(C_{\tvarM}) \rangle_{C_2(E)} + \langle E_{\fvarM}, E_{\svarM}, E_{\tvarM} \rangle_{C_2(B(M_{\KK})_{[3,8]})} -\frac{p_1(\tilde{\tau})}{4}.$$

\begin{lemma}
The difference $(p_1(\tau)-p_1(\tilde{\tau}))$ is independent of $(M,\KK)$.
\end{lemma}
\bp
Let $W$ be  $4$--manifold with boundary $B_N(3) \cup (-[0,1]\times \partial B(3)) \cup (-B(3))$ with signature $0$. $p_1(\tilde{\tau})$ is the obstruction to extend the trivialisation induced by $\tilde{\tau}_{\CC}$ on $\partial W$ to $W$.
Glue a $3$-ball $B$ along $\partial B(3)$, so that $M_{\KK}=B_N(3) \cup B$, extend $\tilde{\tau}$ as $\tau_B$ on $B$. Consider $W_B=W\cup_{[0,1]\times \partial B(3)}[0,1]\times B$. Then $p_1(\tilde{\tau})$ is the obstruction to extend the trivialisation induced by $\tilde{\tau}_{\CC}$ and $\tau_B$ on $\partial W_B$ to $W_B$, where
the signature of $W_B$ is zero. Glue a ball $B^4$ to the part
$\{0\} \times (B(3) \cup B)$ of $\partial W_B$ where $\tau_{B(3)}$ has been extended by $\tau_B$.
The obstruction to extend this trivialisation to $B^4$ will be denoted by
$p_1(\tau_{B(3)} \cup \tau_B)$.
Then $N=M_{\KK}=\partial (W_B \cup B^4)$ and the obstruction to extend $\tilde{\tau}$ to $W_B \cup B^4$ is $p_1(\tilde{\tau}) + p_1(\tau_{B(3)} \cup \tau_B)$.
Now, attaching a $2$--handle along $\hat{K} \times D^2 \subset \partial (W_B \cup B^4)$ equipped with its canonical trivialisation produces a $4$-manifold $W_M$ with signature $0$ and with boundary $M$ and the obstruction to extend the trivialisation induced by $\tau$ to $W_M$ is
$p_1(\tau)$.
Let $p_1(H)$ be the obstruction to extend the trivialisation induced by $\tau$ and $\tilde{\tau}_{\CC}$ on the $2$--handle, it is independent of $M$ and $\KK$. Then $p_1(\tau)=p_1(H)+p_1(\tilde{\tau}) + p_1(\tau_{B(3)} \cup \tau_B)$.
\eop

Since $$\CQ(\KK \subset M)(1,1,1)=\langle p(C_{\fvarM}),p(C_{\svarM}),p(C_{\tvarM}) \rangle_{C_2(E)}-\frac{p_1({\tau})}{4}$$
$$6\lambda(M_{\KK})-\CQ(\KK \subset M)(1,1,1)= \frac{p_1(H) + p_1(\tau_{B(3)} \cup \tau_B)}{4}+\langle E_{\fvarM}, E_{\svarM}, E_{\tvarM} \rangle_{C_2(B(M_{\KK})_{[3,8]})} $$ is independent of $(M,\KK)$.
Thus, we are left with the proof of Lemma~\ref{lemmu} to prove Theorem~\ref{thmaugcas}.

We shall use the following easy lemma.

\begin{lemma}
$H_3(C_2(B(M_{\KK})_{[3,8]});\QQ) = \QQ[ \partial B(M_{\KK}) \times \hat{K}] \oplus \QQ[ \hat{K} \times \partial B(M_{\KK})] \oplus \QQ[ST(\hat{K})]$. A $3$-cycle of $C_2(B(M_{\KK})_{[3,8]})$
reads $\mu ST(\hat{K})$
if and only if its 
algebraic intersections with the two chains
$[P,3\qvarM] \times D_{[3,8]}^2(\hat{K})_P$ and
$D_{[3,8]}^2(\hat{K})_P \times [P,3\qvarM]$ vanish,
where $D_{[3,8]}^2(\hat{K})_P$ is a meridian disk of $\hat{K}$ in $B(M_{\KK})_{[3,8]}$,
$3\qvarM \in \partial B(M_{\KK})$, $P \in \partial M_{[0,3]}$, and $[P,3\qvarM]$ is a path from $P$ to $3\qvarM$ that does not meet $D_{[3,8]}^2(\hat{K})_P$.
\end{lemma}
\eop

\noindent{\sc Proof of Lemma~\ref{lemmu}:}

It is enough to show that $D_{[3,8]}^2(\hat{K})_P$, $[P,3\qvarM]$, 
$\{\pm p_{\cvarM}\} \times D_{[5,8]}^2(\hat{K})$,
$[-3\cvarM,p_M(-p_{\cvarM},0)]$ and $ [p_M(p_{\cvarM},0),3{\cvarM}]$
can be fixed so that neither $\partial_e$ nor $\partial F_{N,\qvarM}(\tilde{\tau})$ intersects $$H=[P,3\qvarM] \times D_{[3,8]}^2(\hat{K})$$ or
 $\iota(H)=D_{[3,8]}^2(\hat{K})_P \times [P,3\qvarM]$.
To do that, fix $\{\pm p_{\cvarM}\} \times D_{[5,8]}^2(\hat{K})$
so that it intersects $M_{[5,6]}$ as $\{(q(t,\pm z(\cvarM)),z);t \in [5,6], z \in S^1\}$, with respect to  the notation of Subsection~\ref{subdefbord}, where $(p_{\cvarM},h_{\cvarM}) = (5,z(\cvarM),h_{\cvarM})$.
Then complete these two vertical annuli with two parallel meridian disks of the tubular neighborhood $N_{[6,8]}(\hat{K})$ of $\hat{K}$ in $M_{\KK}$.
Here is a picture of the projection on $\CC$ of $M_{[0,7]}\setminus M_{[0,3[}$
and the corresponding picture of $S^1 \times D_{[3,8]}^2(\hat{K})$. Note that, for $z\in S^1$, the preimage of $\RR^+ \{z\}$ in $(M_{[0,7]}\setminus M_{[0,3[})$
corresponds to $\{z\} \times D_{[3,8]}^2(\hat{K})$ in $(M_{\KK}\setminus M_{[0,3[})$.
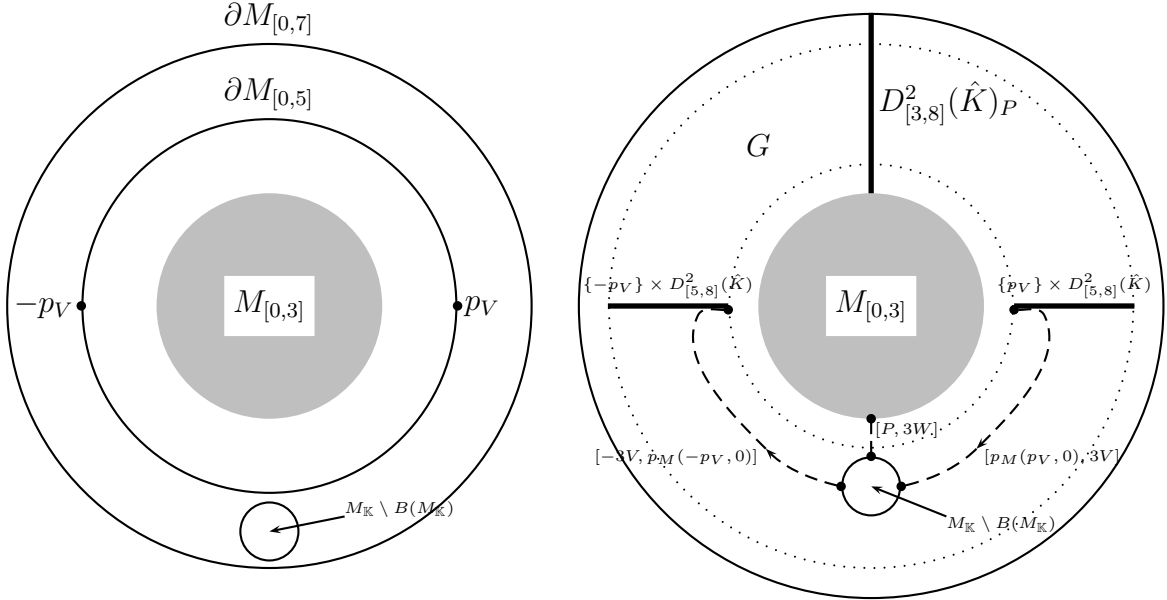
$$
\begin{pspicture}[shift=-.4](-4,-4)(12,4) 
\pscircle*[linecolor=lightgray](0,0){1.5}
\psframe*[linecolor=white](-.6,-.4)(.6,.4)
\rput(0,0){$M_{[0,3]}$}
\pscircle(0,0){2.5}
\rput[b](0,2.6){$\partial M_{[0,5]}$}
\pscircle(0,0){3.5}
\rput[b](0,3.6){$\partial M_{[0,7]}$}
\psdots(-2.5,0)(2.5,0)
\rput[l](2.6,0){$p_V$} \rput[r](-2.6,0){$-p_V$} 
\pscircle(0,-3){.4}
\psline{->}(1,-2.8)(0,-3)
\rput[bl](1,-2.8){\tiny $M_{\KK}\setminus B(M_{\KK})$}
\pscircle*[linecolor=lightgray](8,0){1.5}
\psframe*[linecolor=white](7.4,-.4)(8.6,.4)
\rput(8,0){$M_{[0,3]}$}
\psline[linewidth=2pt](8,1.5)(8,3.9)
\rput[l](8.1,2.7){$D_{[3,8]}^2(\hat{K})_P$}
\psline[linewidth=2pt](9.9,0)(11.5,0)
\rput[b](10.7,.1){\tiny $\{p_{\cvarM}\} \times D_{[5,8]}^2(\hat{K})$}
\psline[linewidth=2pt](6.1,0)(4.5,0)
\rput[b](5.3,.1){\tiny $\{-p_{\cvarM}\} \times D_{[5,8]}^2(\hat{K})$}
\rput[b](6.5,2){$G$}
\pscircle(8,0){3.9}
\pscircle[linestyle=dotted](8,0){1.9}
\pscircle[linestyle=dotted](8,0){3.5}
\psecurve[linestyle=dashed]{*->}(8,-.2)(9.9,-.05)(10.3,-.12)(9.4,-1.9)(8.4,-2.4)
\psecurve[linestyle=dashed]{-*}(10.3,-.1)(9.4,-1.9)(8.4,-2.4)(7.6,-2.4)
\rput[lt](9.5,-1.9){\tiny $[p_M(p_{\cvarM},0),3{\cvarM}]$}
\psecurve[linestyle=dashed]{*->}(8.4,-2.4)(7.6,-2.4)(6.6,-1.9)(5.7,-.1)
\psecurve[linestyle=dashed]{-*}(7.6,-2.4)(6.6,-1.9)(5.7,-.12)(6.1,-.05)(8,-.2)
\rput[rt](6.5,-1.9){\tiny $[-3\cvarM,p_M(-p_{\cvarM},0)]$}
\pscircle(8,-2.4){.4}
\psline[linestyle=dashed]{*-*}(8,-1.5)(8,-2)
\rput[l](8.05,-1.65){\tiny $[P,3\qvarM]$}
\psline{->}(9,-2.8)(8,-2.4)
\rput[tl](9,-2.8){\tiny $M_{\KK}\setminus B(M_{\KK})$}
\end{pspicture}
$$

The complement of the two meridian disks $\{\pm p_{\cvarM}\} \times D_{[5,8]}^2(\hat{K})$ in the solid torus $N_{[5,8]}(\hat{K})$ has two connected components. Let $G$ be the connected component that does not contain $M_{\KK} \setminus B(M_{\KK})$.

Let $\frac53P$ be a point of $\partial D^2_{[0,5]}$ such that $(\frac53P,1) \notin \overline{G}$.
Fix $D_{[3,8]}^2(\hat{K})_P$ so that it intersects $M_{[3,6]}$ as $\{(q(t,-\frac{P}{|P|}),z);t \in [3,6], z \in S^1\}$, and complete it 
into a meridian disk in $G$.
Fix the path $[P,3\qvarM]$ so that it intersects $M_{[3,6]}$ as $\{(q(t,\frac{P}{|P|}),1);t \in [3,6]\}$ and so that it does not intersect $\overline{G}$. 
Thus, $\pi$
maps $C_2(M_{[0,5]}) \cap \left(H \cup
 \iota(H)\right)$ to $(\lambda P,h)$ for two real numbers $\lambda$ and $h$, so that
$\partial_{[3,5]\leftrightarrow 3}$ does not meet $(H \cup
 \iota(H))$, because $(\lambda P,h)\neq \pm \cvarM$.
The path $[P,3\qvarM]$ is furthermore assumed to be transverse to
$\partial B(M_{\KK})$ with $\qvarM \neq \pm \cvarM$ so that $\partial F_{N,\cvarM}(\tilde{\tau})$  does not meet $(H \cup
 \iota(H))$.
It is clear that $(H \cup
 \iota(H))$ does not meet 
$$\begin{array}{l}
+\bigcup_{u \in [0,1]}p_M\left(\{-p_{\cvarM}\}\times [0,u]\right) \times \partial f_{|E}^{-1}(\exp(2i\pi (u+h))) \\
+\bigcup_{u \in [0,1]}\partial f_{|E}^{-1}(\exp(2i\pi (u-h))) \times p_M\left(\{p_{\cvarM}\} \times [0,u]\right) \\
+\left(\{-p_{\cvarM}\} \times D_{[5,8]}^2(\hat{K})\right) \times \partial f_{|E}^{-1}(\exp(2i\pi (h)))+\partial f_{|E}^{-1}(\exp(2i\pi (-h))) \times \left(p_{\cvarM} \times D_{[5,8]}^2(\hat{K})\right)
\end{array}$$
either, and since we can fix the paths
$[-3\cvarM,p_M(-p_{\cvarM},0)]$ and $[p_M(p_{\cvarM},0),3\cvarM]$ so that they meet neither $D_{[3,8]}^2(\hat{K})_P$ nor $[P,3\qvarM]$, we
conclude that $\partial_e$ does not meet $(H \cup
 \iota(H))$ at all. Then the class of $e_{\cvarM}$ reads $\mu ST(\hat{K}_{\cvarM})$ in $H_3(C_2(B(M_{\KK})_{[3,8]}))$ and Lemma~\ref{lemmu} and therefore Theorem~\ref{thmaugcas} are proved.
\eop

\newpage

\section{On the target vector space for $\overline{\CQ}(M)$}
\setcounter{equation}{0}
\label{sectarg}

\subsection{Notation}

Consider the polynomial ring 
$$\Lambda_f=\QQ[\fvar^{\pm 1},\svar^{\pm 1}]=\frac{\QQ[\fvar^{\pm 1},\svar^{\pm 1},\tvar^{\pm 1}]}{\fvar\svar\tvar=1}$$
and its subring
$$\QQ_{s}[\fvar^{\pm 1},\svar^{\pm 1}]= \{P \in \Lambda_f; P(\fvar,\svar)=P(\fvar, \tvar=(\fvar\svar)^{-1})=P(\svar,\fvar)=P(\fvar^{-1},\svar^{-1})\}.$$

Then $$\CQ_2(\KK)=\delta(\fvar)\delta(\svar)\delta(\tvar) \CQ(\KK)$$ takes its values in $\QQ_{s}[\fvar^{\pm 1},\svar^{\pm 1}]$, according to Proposition~\ref{propdenwithoutz} and to Subsection~\ref{subsecsym}.
Define
 $$P_k=P_k(M)=P_k(\Delta(M),\delta(M))=P_k(\Delta,\delta)=\delta(\fvar)\delta(\svar)\delta(\tvar)\sum_{\mathfrak{S}_3(\fvar,\svar,\tvar)}\frac{(\fvar^k-\fvar^{-k})}{\delta(\fvar)}\ID(\tvar).$$
\begin{lemma}
$P_k \in \QQ_{s}[\fvar^{\pm 1},\svar^{\pm 1}]$.
$$P_k=\sum_{\mathfrak{S}_3} (\fvar^k-\fvar^{-k})\delta(\svar)\delta(\tvar)\frac{\tvar\Delta^{\prime}(\tvar)}{\Delta(\tvar)}
+\sum_{\circlearrowleft}
\frac{\left((\fvar\tvar)^k-(\fvar\tvar)^{-k}\right)\delta(\fvar)-\left(\fvar^k-\fvar^{-k}\right)\delta(\fvar\tvar)}{\tvar - 1}     \delta(\tvar)(\tvar + 1)$$
\end{lemma}
\bp
The symmetry properties of $P_k$ are easy to check, and it is enough to check that $P_k$ can be expressed as in the statement to prove that
$P_k$  is actually a polynomial.
$$P_k-\sum_{\mathfrak{S}_3} (\fvar^k-\fvar^{-k})\delta(\svar)\delta(\tvar)\frac{\tvar\Delta^{\prime}(\tvar)}{\Delta(\tvar)}$$
$$=-\sum_{\mathfrak{S}_3}(\fvar^k-\fvar^{-k})\delta(\svar)\delta(\tvar)\frac{\tvar + 1}{\tvar - 1}$$
$$=-\sum_{\mathfrak{S}_3}(\fvar^k-\fvar^{-k})\delta(\fvar\tvar)     \delta(\tvar)\frac{\tvar + 1}{\tvar - 1}$$
$$=\sum_{\circlearrowleft}
\frac{\left((\fvar\tvar)^k-(\fvar\tvar)^{-k}\right)\delta(\fvar)-\left(\fvar^k-\fvar^{-k}\right)\delta(\fvar\tvar)}{\tvar - 1}     \delta(\tvar)(\tvar + 1).$$
\eop

Let $\Lambda(M)=\Lambda(\Delta,\delta)$ be the quotient of 
$\QQ_{s}[\fvar^{\pm 1},\svar^{\pm 1}]$ by the rational vector space generated by the $P_k(M)$ for $k \in \NN$, $k \geq 1$.
According to Theorem~\ref{thmfrakcha}, the class $\overline{\CQ}_2(M)$ of $\CQ_2(\KK)$ in $\Lambda(M)$ is independent of $\KK$. Therefore, it is worth studying this quotient.

I conjecture that $\overline{\CQ}_2(M)$ is equivalent to an invariant previously defined by Ohtsuki in \cite{ohtb} when $\delta=\Delta$, and that it is a refinement of this Ohtsuki invariant when $\delta \neq \Delta$.
Ohtsuki has proved more results on $\Lambda(\Delta,\Delta)$ in \cite{ohtb}.

\subsection{Detection of constants}

In this subsection we prove the following proposition.

\begin{proposition}
\label{propinj}
If $\Delta(M)$ has no multiple roots, then $\Delta(M)=\delta(M)$ and
$\delta(\fvar)\delta(\svar)\delta(\tvar)$ does not vanish in $\Lambda(M)$. In particular, $\overline{\CQ}_2(M \sharp N) \neq \overline{\CQ}_2(M)$ for any rational homology sphere $N$ with non trivial Walker invariant.
\end{proposition}

We have the following tautological lemma.

\begin{lemma}
In general, $\delta(\fvar)\delta(\svar)\delta(\tvar)$ vanishes in $\Lambda(M)$ if and only if the following equation is satisfied
\begin{equation}
\label{eqabs}
\delta(\fvar)\delta(\svar)\delta(\tvar)= \sum_{i=1}^s\beta_iP_i
\end{equation} where $s \in \NN \setminus \{0\}$, $\beta_i \in \QQ$, for $i=1, \dots, s$, and $\beta_s \neq 0$.
\end{lemma}
\eop

Set
$$p_k=p_k(\Delta,\delta)=P_k(\fvar, \fvar^{-1},1).$$

Of course, if Equation~\ref{eqabs} is satisfied, then  
\begin{equation}
\label{eqabsev}
\delta^2(\fvar) = \sum_{i=1}^s\beta_ip_i(\fvar).
\end{equation}

\begin{lemma}
\label{lempk}
$$p_k(\Delta,\delta)=2(\fvar^k-\fvar^{-k}) \delta(\fvar) \left(  - \frac{\fvar \Delta^{\prime}(\fvar)}{\Delta(\fvar)} - \frac{\fvar \delta^{\prime}(\fvar)}{\delta(\fvar)}
+ \frac{\fvar +1}{\fvar-1}\right)
+2k(\fvar^k+\fvar^{-k})\delta(\fvar).$$
$$p_k(\Delta,\Delta)=-4\fvar\Delta^{\prime}(\fvar)(\fvar^k-\fvar^{-k}) +2\Delta(\fvar)
\left( k(\fvar^k+\fvar^{-k})+(\fvar +1)\frac{\fvar^k-\fvar^{-k}}{\fvar-1}\right)$$
$$=-4\fvar\Delta^{\prime}(\fvar)(\fvar^k-\fvar^{-k}) +2\Delta(\fvar)
\left( (k+1) (\fvar^k +\fvar^{-k}) + 2 (\fvar^{k-1} + \fvar^{k-2} + \dots + \fvar^{2-k} +\fvar^{1-k}) \right).$$
\end{lemma}
\bp \\
$\begin{array}{lll}p_k&=&-2(\fvar^k-\fvar^{-k})\delta(\fvar)\frac{\fvar\Delta^{\prime}(\fvar)}{\Delta(\fvar)}\\
&&+\mbox{lim}_{\tvar \rightarrow 1}\frac{\left((\fvar\tvar)^k-(\fvar\tvar)^{-k}\right)\delta(\fvar)-\left(\fvar^k-\fvar^{-k}\right)\delta(\fvar\tvar)}{\tvar - 1}     \delta(\tvar)(\tvar + 1)\\
&&+\frac{\left(\fvar^k-\fvar^{-k}\right)}{\fvar - 1}     \delta(\fvar)(\fvar + 1)
-\frac{\left(\fvar^k-\fvar^{-k}\right)}{\svar - 1}     \delta(\svar)(\svar + 1)
\end{array}$\\
where $\svar=\fvar^{-1}$.
Of course, the mentioned limit is nothing but the derivative at $\tvar=1$ that is $$2\left(k\left(\fvar^k+\fvar^{-k}\right)\delta(\fvar)-\left(\fvar^k-\fvar^{-k}\right)\fvar\delta^{\prime}(\fvar)\right).$$
The second equality is an easy consequence of the first one.
\eop

If $\Delta=1$, then $p_k =2\left( (k+1) (\fvar^k +\fvar^{-k}) + 2 (\fvar^{k-1} + \fvar^{k-2} + \dots + \fvar^{2-k} +\fvar^{1-k})  \right)$.
Therefore, the degree of $p_k$ is $k$ for any $k$ and Equation~\ref{eqabs} cannot hold with some $\beta_s \neq 0$.
Therefore, Proposition~\ref{propinj} is true in this case.

\begin{lemma}
\label{lemdetcons}
Assume that $$\Delta=\delta=\sum_{i=1}^r \alpha_i (t^i + t^{-i}) + 1-2\sum_{i=1}^r \alpha_i,$$ with $\alpha_r \neq 0$ and $r>0$.
If Equation~\ref{eqabs} is satisfied, then $s=r$, $r>1$ and $\alpha_r=(2-2r)\beta_r$.
\end{lemma}
\bp
The leading term in
$$p_k=-4\fvar\Delta^{\prime}(\fvar)(\fvar^k-\fvar^{-k}) +2\Delta(\fvar)
\left( (k+1) (\fvar^k +\fvar^{-k}) + 2 (\fvar^{k-1} + \fvar^{k-2} + \dots + \fvar^{2-k} +\fvar^{1-k}) \right).$$
is $$(2(k+1)\alpha_r - 4r\alpha_r)\fvar^{k+r}$$ unless $k=2r-1$, and in this latter case the degree of $p_k$ is lower than $(k+r)$.

The degree of the left-hand side in Equation~\ref{eqabsev} is $2r$ while the degree of the right-hand side is $s+r$ if $s\neq 2r-1$ and is lower than $s+r$ if $s=2r-1$.
Therefore, if Equation~\ref{eqabsev} holds with $s\neq r$, then
$s=2r-1$ and $r<s$, and therefore $r>1$.

In particular, if $r=1$, Equation~\ref{eqabsev} could only hold with
$r=s=1$ but in this case $s=2r-1$, and the degree of $p_1$ is smaller than $2$.
Thus, if Equation~\ref{eqabsev} is satisfied, then $r>1$.

Then Lemma~\ref{lemdetcons} will be proved as soon as the following lemma is proved.

\begin{sublemma}
 If $\Delta=\delta=\sum_{i=1}^r \alpha_i (t^i + t^{-i}) + 1-2\sum_{i=1}^r \alpha_i,$ with $\alpha_r \neq 0$ then 
Equation~\ref{eqabs} cannot hold if $s=2r-1>r>1$.
\end{sublemma}
\noindent{\sc Proof of the sublemma:}
We compare the terms of lexicographically higher degrees in both sides
of the equation:
$$\Delta(\fvar)\Delta(\svar)\Delta(\tvar)(\fvar-1)(\svar-1)(\tvar-1) = (\fvar-1)(\svar-1)(\tvar-1)\sum_{i=1}^s\beta_iP_i$$
in $\QQ[\fvar^{\pm 1},\svar^{\pm 1}]$.
On the left-hand side we have
$\alpha_r^3 \fvar^{2r} \svar^{2r} (-\fvar\svar)$.
On the right-hand side, since $s>r>1$, the term of highest degree in
$$(\fvar-1)(\svar-1)(\tvar-1)P_s$$
$$=\sum_{\mathfrak{S}_3} (\fvar-1)(\fvar^s-\fvar^{-s}) (\svar-1) \delta(\svar)  \left((\tvar-1)  \tvar \delta^{\prime}(\tvar) - \delta(\tvar)(\tvar +1)\right)$$
is $2(\fvar\svar)^{s+1}(r-1)\alpha_r^2(\fvar\svar)^r$.
Then since $s>r$ the degree of the right-hand side is bigger than the degree of the left-hand side.
\eop

The fact than $\alpha_r=(2-2r)\beta_r$ can be deduced by identifying the leading terms of both sides in Equation~\ref{eqabsev} or in Equation~\ref{eqabs}.

\eop

\noindent{\sc End of proof of Proposition~\ref{propinj}:}
If $\Delta$ has no multiple root then $\Delta$ is coprime with $\Delta^{\prime}$ and therefore to $\fvar\Delta^{\prime}$.

Then if Equation~\ref{eqabsev} is satisfied, according to Lemma~\ref{lempk}, $\Delta$ must divide  $\sum_{i=1}^s\beta_i(\fvar^i-\fvar^{-i})$ and since $\Delta$ is symmetric, the ratio is of degree at least $1$.
Then $s \geq r+1$.
This concludes the proof of Proposition~\ref{propinj}.
\eop

\end{document}